\documentclass[11pt]{article}

\usepackage[margin=1in]{geometry}
\usepackage{amsthm,amsmath,amsfonts,amssymb}
\usepackage[authoryear,round]{natbib}
\usepackage[colorlinks,citecolor=blue,urlcolor=blue,linkcolor=blue]{hyperref}
\usepackage{graphicx}
\usepackage{tikz}
\usetikzlibrary{patterns}
\usepackage{mathtools}
\usepackage{xcolor}
\usepackage{float}
\usepackage{placeins}

\numberwithin{equation}{section}
\theoremstyle{plain}
\newtheorem{theorem}{Theorem}
\newtheorem{proposition}{Proposition}
\newtheorem{lemma}{Lemma}
\newtheorem{corollary}{Corollary}
\theoremstyle{definition}
\newtheorem{defn}{Definition}
\newtheorem{example}{Example}
\newtheorem{remark}{Remark}
\newcommand{\load}{s_+}

\newcommand{\votemapwidth}{0.75}
\makeatletter
\newenvironment{figH}[1][!bp]{\edef\@tempa{\noexpand\figure[#1]}\@tempa}{\endfigure}
\makeatother
\newcommand{\figposrisk}{!t}
\newcommand{\figposconv}{!t}
\newcommand{\figposmap}{!tp}

\title{Reweighting for Transfer Classification:\\ Exact Risk Conversion and Minimax Adaptation}
\author{Aoran Chen\\
  Department of Statistics, The Pennsylvania State University\\
  \texttt{azc6609@psu.edu}
  \and
  Yang Feng\thanks{Corresponding author.}\\
  Department of Biostatistics, New York University\\
  \texttt{yf31@nyu.edu}}
\date{}

\begin{document}
\maketitle

\begin{abstract}
A source can make a training mixture easier to learn while making its accuracy less valuable for the target. We study how the training problem a source creates and the finite-sample error of learning it jointly determine target-risk guarantees under reweighting. For binary classification, for each pair of laws and each positive target weight, we derive the sharpest target excess hinge-risk bound implied by mixture excess alone, uniformly over clipped predictors. This exact conversion curve is determined by density-weighted local votes. A two-sample oracle inequality for weighted support vector machines supplies the mixture excess, and the composition becomes explicit for Gaussian kernels. On a tilted-boundary class, we establish the all-procedure minimax classification rate \([\,n+\load^2/(m^{-1}+\tau^2)\,]^{-1}\) in the transferable strength \(\load\) and the tilt \(\tau\), separating source sampling uncertainty (\(m^{-1}\)) from relation uncertainty (\(\tau^2\)). A target-checked Gaussian-SVM selector adapts to unknown \(\load\) and \(\tau\) up to a logarithmic factor for fixed dimension, and retains the target-only rate up to the same factor against arbitrary independent sources under the stated target model. On a separate smooth opposition class, the minimax risk within the family of exact Gaussian SVMs fitted to the original labels with nonnegative weights is exponentially larger than the unrestricted minimax risk whenever the target sample is at most half the source sample, despite arbitrary data-dependent weighting and tuning. The value of a source thus depends on the training problem it creates, on what accuracy there means for the target, and on whether the permitted use of the source can implement the needed correction.
\end{abstract}



\section{Introduction}

In transfer classification, a learner is trained partly on source
data and judged on the target
\citep{pan2009survey,zhuang2021survey,blitzer2007learning}.  Consider two source distributions
that share the Bayes classifier of the target: at every feature value, each would
predict the same label as the target does.  Borrowing from either
source might be expected to help equally.  Not always.  Pooling
combines label evidence, and the evidence a source
contributes at a feature value carries more than its predicted
label: it carries how confident that label is and how often the
source is observed there.
Two sources with identical predicted
labels can therefore strengthen or weaken the pooled objective's
margin at the same location.  Moreover, even when the pooled objective
and the target agree on the optimal prediction everywhere, they
can assign different costs to the same mistake.  A fitted
classifier always makes some mistakes, so what the source is worth
depends on this cost.

Now give one of these sources confident labels.  The training
mixture can become easier to learn: its margins are stronger, so a
given accuracy is reached with fewer samples.  The target does not
grade the mixture.  It grades where the learner's remaining
mistakes fall and how much each one costs under the target law.  We study
these effects through reweighting, with a support vector machine
fitted to the weighted union of independent target and source
samples.
Two quantities then have to be known together: how accurately the
mixture can be learned from the samples at hand, and what
that accuracy is worth under the target distribution.

Let \(P_T\) and \(P_S\) be the target and source laws, let \(n,m\)
be their sample sizes, and let \(\alpha\in(0,1]\) be the total
target weight, with \(u=1-\alpha\).  The population training
distribution is \(P_\alpha=\alpha P_T+uP_S\).  Write
\(m_D(x)=2P_D(Y=1\mid X=x)-1\) for the posterior margin of domain
\(D\in\{T,S,\alpha\}\), and let \(p_D\) be its covariate density
with respect to a common dominating measure.  The signed votes
combine as
\begin{equation}
  p_\alpha(x)\,m_\alpha(x)
  =\alpha\,p_T(x)\,m_T(x)+u\,p_S(x)\,m_S(x).
  \label{eq:intro_vote}
\end{equation}
Density determines how often a domain votes near \(x\), and the
margin gives the strength and sign of its vote.  A single global
weight can
preserve the target label in one region and reverse it in another;
this local reversal is the mechanism behind negative transfer
\citep{rosenstein2005transfer,wang2019characterizing}.

Write \(\mathcal E_\alpha\) and \(\mathcal E_T\) for hinge excess
risks relative to the respective Bayes risks, and suppose a fitted
score has mixture excess at most \(t\).
The largest target excess compatible with this budget, over all
clipped scores, is a function \(\Psi_\alpha(t)\) of \(t\ge0\):
every measurable \(\widehat g:\mathcal X\to[-1,1]\) satisfies
\begin{equation}
  \mathcal E_\alpha(\widehat g)\le t
  \quad\Longrightarrow\quad
  \mathcal E_T(\widehat g)\le\Psi_\alpha(t).
  \label{eq:intro_value}
\end{equation}
The samples determine how small \(t\) can be: with \(n\) target
and \(m\) source observations a high-probability bound is
\(t_{\alpha,n,m}\).  The target guarantee is their composition
\(\Psi_\alpha(t_{\alpha,n,m})\).  A source with confident labels
can lower \(t_{\alpha,n,m}\) by making small mixture margins
rarer, so the mixture is easier to learn.  The same source can raise \(\Psi_\alpha\), because the
mixture assigns little weight to locations where its votes nearly
cancel while the target vote there remains large, so each unit of
mixture error permits more target error.
The source and
its weight act on both, so neither quantity alone decides whether
borrowing helped.

The laws behind these quantities are unknown, so the worth of a
source must be judged from data.  One might hope that enough
source data settles the matter: learn the source law exactly, then
decide.  Knowing the source law exactly still leaves open how that
law relates to the target.  On a structured class we separate the
uncertainty that additional source observations remove from the
uncertainty they cannot touch.  Target data can nevertheless
support a useful choice without first computing the class-informed
optimal weight, the weight that would be optimal if the relation
were known.  Once weighting works as well as it can, we ask
whether any choice of nonnegative weights can extract the value of
a source whose labels become informative only when reversed.

\subsection{Main results}
\label{subsec:intro_overview}

Theorems~\ref{thm:affine_frontier} and~\ref{thm:weighted_crerm}
supply the two quantities: the first the conversion, the second
the budget.  Theorem~\ref{thm:affine_frontier} computes
\(\Psi_\alpha\) exactly for each fixed pair of laws and weight,
together with the full family of affine bounds
\(\mathcal E_T(g)\le c\,\mathcal E_\alpha(g)+\mathcal B_\alpha(c)\),
where \(c\ge0\) is the slope applied to mixture excess and
\(\mathcal B_\alpha(c)\) is the smallest intercept valid for every
clipped score \(g\).  At zero mixture excess,
\(\Psi_\alpha(0)=H_{\rm flip}(u)\), the flip profile: on the set
where the source vote opposes the target's, the mass under
\(P_{X,T}\), weighted by \(2|m_T|\), of locations where source
weight \(u\) has tied or reversed the target label
(Section~\ref{sec:problem_setup}).
This is the target excess that remains possible when the mixture is
learned perfectly, and it is the limit of the intercepts \(\mathcal B_\alpha(c)\) as \(c\to\infty\), which no finite slope need attain (Section~\ref{sec:adaptive}).  We keep the whole family
because the two domains' votes can nearly cancel: there the
mixture is almost indifferent while the target is not, so the
local multiplier is unbounded and no single multiplier bounds
target error.  A finite slope with a positive intercept still
does, and this is the form that enters the fitted-SVM analysis of
Section~\ref{sec:adaptive}.

Theorem~\ref{thm:weighted_crerm} bounds the mixture excess of the
weighted SVM over a general reproducing kernel Hilbert space through an effective sample size that combines the two samples (the Supplement keeps their contributions separate).  Source sample size and source weight are
different interventions: more source
observations make the source part of the empirical risk more
precise, whereas a larger source weight makes that part count for
more in the objective being minimized.  With high probability the
clipped fitted score \(\bar f_D\) satisfies
\begin{equation}
  \mathcal E_T(\bar f_D)
  \le\Psi_\alpha\bigl(t_{\alpha,n,m}\bigr)
  =\inf_{c\ge0}\bigl\{c\,t_{\alpha,n,m}+\mathcal B_\alpha(c)\bigr\},
  \label{eq:intro_composition}
\end{equation}
where \(t_{\alpha,n,m}\) is a two-sample oracle bound on mixture
excess.  At the same confidence level, borrowing improves on the
target-only guarantee \(t_{1,n,m}\), the same bound at \(\alpha=1\),
when \(\Psi_\alpha(t_{\alpha,n,m})<t_{1,n,m}\), in the range where
\(\Psi_\alpha\) has not saturated (it is constant once the budget
exceeds the mixture's total vote capacity);
Corollary~\ref{thm:main_target_bound} makes both ingredients
explicit for Gaussian kernels.  Example~\ref{ex:constant_Q} shows
what the composition does: two mixtures with the same Bayes rule
as the target, at the same budget, receive sharply different
target guarantees.

We then ask what finite data must reveal about the relation to
guide borrowing.
On a tilted-boundary class, in which each posterior margin is a
linear form in \(x\) scaled by a radial profile, the source signal
splits
into a component along the target direction, whose coefficient,
when positive, is the transferable strength \(\load\) (\(\load=0\) covers
every nonpositive aligned strength), and a remainder, whose size is
at most the tilt \(\tau\).  Write
\(\mathcal P_1(\load,\tau)\) for the resulting quality class and
\(\mathcal E_T^{01}\) for target classification excess.
Theorem~\ref{thm:target_safe_adaptive} establishes, for fixed
dimension and for \(\load\) and \(\tau\) in the ranges stated there, the
all-procedure minimax rate
\begin{equation}
  \inf_{\widetilde h}\,
  \sup_{P\in\mathcal P_1(\load,\tau)}
  \mathbb E_P\,\mathcal E_T^{01}(\widetilde h)
  \;\asymp\;
  \Bigl[\,n+\frac{m\load^2}{1+m\tau^2}\Bigr]^{-1},
  \label{eq:intro_adaptive}
\end{equation}
up to constants, with the infimum over all randomized measurable
classifiers.  The source contribution can be written as
\(\load^2/(m^{-1}+\tau^2)\).  As \(m\) grows without bound at fixed \(\tau>0\), the source
law becomes known exactly, yet the contribution stops at
\(\load^2/\tau^2\).  What remains is ambiguity: within the class,
distinct target directions are compatible with the same source
law, and source observations alone cannot separate them.

A single target-checked selector approaches this benchmark without
knowing \(\load\) or \(\tau\) and without reconstructing the
class-informed optimal weight.  It uses an estimated strength to
propose a source-assisted fit, then independent target data to
check its decisions.  For fixed dimension, at fixed kernel and
regularization, its risk is within a logarithmic factor \(\ell_N=\log(e(n+m))\) of the benchmark, with
no logarithmic loss on the clean class \(\tau=0\).  Under the same
target model, it also retains the target-only rate up to that
factor against arbitrary independent sources.
The fallback is
paid for only when it is used: its target-only risk is multiplied
by the rejection probability, which the proof controls on an
informative source, so the fit keeps its faster rate.

Having reached the benchmark, we ask a further question: is
sufficiently careful weighting always enough?  Theorem~\ref{thm:opposition_separation}
gives a negative answer on a separate smooth opposition class, in
which the source labels are systematically reversed relative to
the target.  Reversing the source labels gives an accurate target
rule.  Keeping them and fitting an exact Gaussian SVM with
nonnegative observation weights leaves an event on which the
weights are irrelevant: when every observed label, in both samples, is the label opposite to the target's Bayes label, any such fit returns that label, while reversing
the source labels returns the target's.  Whenever the target sample is at most half the source sample \((n\le m/2)\) this event
is rare but not rare enough.
  The minimax risk within the
exact-SVM family is exponentially larger than the unrestricted
minimax risk.  The gap survives arbitrary data-dependent weights,
retained observations and tuning, even when the opposition
relation is known; it lies in what the final fit can return.

A source is therefore worth what its training problem lets the
learner achieve, what that achievement means under the target law,
and what the permitted way of using the source can implement.

\subsection{Related work}
\label{subsec:intro_related}

The closest work compares excess risks across distributions.
\citet{hanneke2019value} do so through transfer exponents, and
\citet{hanneke2024aggregation} through a weak modulus of transfer,
the largest target excess risk of a hypothesis whose source excess
risk is within a budget; discrepancy-based bounds
\citep{ben2010theory,mansour2009domain} use a global distance
between the two laws.  In each case the comparison is a fixed
property of the source--target pair.  Ours starts from the
training mixture that a chosen weight creates, so it changes with
the weight; under the hinge loss it is an exact identity at each
feature value, and the curve and its affine frontier are computed
from the local vote \eqref{eq:intro_vote}.

Minimax rates for structured transfer classes are the second point
of contact.  On a sign-linked posterior-drift class,
\citet{cai2021transfer} obtain a rate in which the source acts as
an effective sample size set by a relative signal exponent,
attained by weighted \(k\)-NN with known parameters and
adaptively up to logarithmic factors; there the two Bayes
decisions agree.  \citet{kpotufe2021marginal} and
\citet{reeve2021adaptive} treat covariate shift and
covariate-dependent relations.  \citet{fan2025robust} establish minimax bounds
for logistic transfer that balance source estimation error and a
squared directional tilt against target-only learning, under a
fixed positive signal lower bound.  Our source-quality rate has
the same form up to constants, the source boundary may rotate away
from the target's, and a single globally weighted Gaussian SVM
attains the rate uniformly over unknown strength and tilt,
including vanishing positive strength.

\citet{hanneke2023limits} identify gaps between oracle and
adaptive model-selection rates under transfer; our limitation
comes from an exact-SVM output restriction and appears even when
the opposition class is known.  Procedures that adjust the source
rule itself, through a link-scale correction
\citep{maity2024linear} or low-dimensional transformations of a
source-trained decision rule \citep{wang2025decision}, are
complementary: our results quantify what retaining the original
labels with nonnegative weights can achieve.  Transfer through
contrasts between model parameters, in high-dimensional regression
and generalized linear models
\citep{bastani2021predicting,li2022transfer,tian2023transfer},
covariate-shift regression
\citep{shimodaira2000improving,ma2023optimally}, and multi-task
and representation learning
\citep{duan2023adaptive,tian2025representations,tian2024gmm},
concerns parametric models and leaves the local vote implicit.
The SVM analysis
builds on the fast Gaussian-kernel rates of
\citet{steinwart2007fast} and localized empirical-process theory
\citep{bartlett2005local,koltchinskii2006local,steinwart2008support},
keeping the two sampling contributions separate.

Section~\ref{sec:problem_setup} describes how reweighting combines
local decision evidence; Section~\ref{sec:mixture_to_target}
evaluates its target consequences; Section~\ref{sec:adaptive}
identifies what target data must learn to guide borrowing and asks
whether better weighting is always enough;
Section~\ref{sec:numerical} varies source sample size, source-label
strength and source-boundary orientation numerically; and Section~\ref{sec:discussion} returns to what
these distinctions imply for source use.  Proofs are in the
Supplementary Material.

\section{Setup and the Local Vote}
\label{sec:problem_setup}

Our goal is to classify observations from a target distribution.
We also have labeled observations from a source distribution,
which we can use in training.  Combining the two samples may add
information, and it also changes the distribution on which we
train.  We first examine the classification rule of that
distribution, and then ask how its errors compare with errors on
the target.

\subsection{Weighted SVMs and Target Excess Risk}

Let \(P_T\) and \(P_S\) be the target and source distributions on
\(\mathcal X\times\mathcal Y\), with \(\mathcal X\subset\mathbb R^d\)
and \(\mathcal Y=\{-1,+1\}\), and let
\(D_T^n=\{(X_i^T,Y_i^T)\}_{i=1}^n\sim P_T^{\otimes n}\) and
\(D_S^m=\{(X_j^S,Y_j^S)\}_{j=1}^m\sim P_S^{\otimes m}\) be
independent samples, and write \(D:=(D_T^n,D_S^m)\) for the combined sample (as the subscript of a law, margin, vote or risk, \(D\) instead indexes a domain, \(T\) or \(S\), or the mixture \(\alpha\)).  A score \(f\) assigns the label
\(\operatorname{sign}f(x)\) to \(x\), and we evaluate it by its
target classification risk
\(\mathcal R_{L_{01},P_T}(f):=P_T\{\operatorname{sign}f(X)\ne Y\}\).
No relation between \(P_T\) and \(P_S\) is assumed: arbitrary
covariate shift and posterior drift are allowed
\citep{pan2009survey,cai2021transfer}, and margin, noise, and
smoothness conditions enter only later as explicit inputs to
individual statements.

We fit a support vector machine with the hinge loss
\(L(y,t)=(1-yt)_+\) \citep{vapnik1995support,zhang2004statistical}.
For an RKHS \(H\) \citep{scholkopf2002learning}, regularization
\(\lambda>0\), and target weight \(\alpha\in[0,1]\), the weighted
transfer SVM \citep{blitzer2007learning} is any solution of
\begin{equation}
  f_{D,\lambda,\alpha}
  \in\arg\min_{f\in H}
  \left\{
  \lambda\|f\|_H^2
  +\frac{\alpha}{n}\sum_{i=1}^n L(Y_i^T,f(X_i^T))
  +\frac{1-\alpha}{m}\sum_{j=1}^m L(Y_j^S,f(X_j^S))
  \right\}.
  \label{eq:transfer_svm_empirical}
\end{equation}
Write \(u:=1-\alpha\) for the total source weight.  The losses are
averaged within each sample before they are combined, so \(m\)
controls how many source observations estimate the source term,
while \(u\) controls how much that term counts in the objective;
\(\alpha=1/2\), \emph{equal-domain weighting}, gives the two
samples equal total weight, and \(\alpha=n/(n+m)\),
\emph{proportional pooling}, gives every observation equal
weight.\footnote{Weighting the source empirical risk by \(q\ge0\)
against an unweighted target risk is equivalent to
\eqref{eq:transfer_svm_empirical} with \(\alpha=(1+q)^{-1}\) and
\(\lambda=\lambda_q/(1+q)\); this \(q=u/\alpha\) is the
source-to-target odds used in Section~\ref{sec:adaptive}.}  We
write \(f_D\) for \(f_{D,\lambda,\alpha}\) when \(\lambda\) and
\(\alpha\) (and, for Gaussian kernels, the bandwidth \(\gamma\))
are fixed by context, and display a parameter in the subscript
only where it varies.  Its population version replaces the two averages by the hinge
risk \(\mathcal R_{L,P}(f):=\mathbb E_PL(Y,f(X))\) under the mixture
\[
  P_\alpha:=\alpha P_T+uP_S .
\]
Changing \(u\) therefore changes the weights in the fitted
objective and can change the classification rule that its
population version favors.  To see what changes, we set sampling
and fitting aside and ask how \(P_\alpha\) classifies a point.

\subsection{Local Votes and the Mixture Label}
\label{subsec:mixture_geometry}
\label{subsec:local_identities}

Fix a feature value \(x\).  The Bayes rule of a distribution
predicts the label with the larger probability near \(x\).  Let
\(p_D(x)\) be the covariate density of domain \(D\in\{T,S\}\) with
respect to a common dominating measure \(\mu\) (one may take
\(\mu=P_{X,T}+P_{X,S}\)), and let
\(\eta_D(x):=P_D(Y=1\mid X=x)\).  The joint density of \((X,Y)\)
under domain \(D\) is \(p_D(x)\eta_D(x)\) at \((x,+1)\) and
\(p_D(x)\{1-\eta_D(x)\}\) at \((x,-1)\).  The mixture adds the
corresponding quantities of the two domains: its density at
\((x,+1)\) is \(\alpha p_T\eta_T+up_S\eta_S\), and at \((x,-1)\) it
is \(\alpha p_T(1-\eta_T)+up_S(1-\eta_S)\).  Their difference
determines the preferred label.  Write \(m_D(x):=2\eta_D(x)-1\)
for the \emph{posterior margin} of domain \(D\); the difference is
\(\alpha p_Tm_T+up_Sm_S\).  The mixture has density
\(p_\alpha:=\alpha p_T+up_S\), conditional label probability
\(\eta_\alpha\), and posterior margin \(m_\alpha:=2\eta_\alpha-1\),
and because mixing the \emph{joint} distributions mixes both label
densities, the same difference equals \(p_\alpha m_\alpha\):
\begin{equation}
  p_\alpha(x)m_\alpha(x)
  =
  \alpha p_T(x)m_T(x)+u p_S(x)m_S(x)
  \qquad \mu\text{-a.e.}
  \label{eq:margin_identity}
\end{equation}
We call \(s_D:=p_Dm_D\) the \emph{local vote} of domain \(D\), so
that \(s_\alpha:=p_\alpha m_\alpha=\alpha s_T+us_S\).  Its sign is
the label the domain prefers at \(x\), and its magnitude is the
difference between the two label densities; the corresponding
Bayes representative is \(f_D^*:=\operatorname{sign}m_D\), with
\(\operatorname{sign}(0)=1\).  With equal covariate densities the
mixture margin is the plain average \(\alpha m_T+um_S\); otherwise
the density factors set the local balance.  A vote takes the same value for every version of \(m_D\), because
\(p_D\) vanishes off the support of \(P_{X,D}\).

Now change \(u\) with the two distributions fixed.  Suppose at one
location the target vote is \(0.2\) and the source vote is
\(-0.3\).  Then
\[
  s_\alpha(x)=(1-u)(0.2)+u(-0.3)=0.2-0.5u,
\]
so the mixture favors \(+1\) for \(u<0.4\), has no preference at
\(u=0.4\), and favors \(-1\) for \(u>0.4\).  When the two domains
favor the same label,
increasing the source weight strengthens the mixture vote if the
source vote is the larger of the two in magnitude and weakens it
if the source vote is the smaller; the preferred label does not
change.  When
the domains favor different labels, increasing the source weight
always reduces the vote in the target's direction: the
target-directed vote becomes weaker and, once the source receives
enough weight, reverses.  Both the change in the preferred label
and the change in the strength of the vote matter when we compare
the mixture with the target.

\subsection{Flip Thresholds and Decision Regions}
\label{subsec:local_global_dictionary}

The example has a threshold at \(u=0.4\).  Another location can
have a different threshold, because the two distributions can
have different local votes there.

Let \(\mathcal D_T:=\{x:|s_T(x)|>0\}\) be the target-active set,
and let \(\Omega_-:=\{x:s_T(x)s_S(x)<0\}\) and
\(\Omega_+:=\{x:s_T(x)s_S(x)\ge0\}\) be the strict-disagreement
and agreement regions; neither depends on \(\alpha\), and
\(\Omega_-\subset\mathcal D_T\) because \(s_Ts_S<0\) forces
\(s_T\ne0\), so \(\Omega_-\cap\mathcal D_T\) below is \(\Omega_-\)
itself.  For
\(0\le u<1\), the mixture can strictly reverse the target label
only on \(\Omega_-\cap\mathcal D_T\).  At such a
point, write
\[
  r(x):=\frac{|s_S(x)|}{|s_T(x)|}
\]
for the \emph{local vote ratio}, the ratio of the two vote
magnitudes before the weights are applied.  Because the votes have opposite signs, the mixture vote
relative to the target vote is
\(s_\alpha(x)/s_T(x)=(1-u)-u\,r(x)\): the target label is retained
while this quantity is positive, the votes tie when it is zero,
and the label reverses when it is negative.  Solving
\((1-u)-u\,r(x)=0\) gives the \emph{critical source weight}
\begin{equation}
  u_c(x):=\frac{1}{1+r(x)}.
  \label{eq:pointwise_flip_threshold}
\end{equation}
For the location above, \(r(x)=0.3/0.2=1.5\) and \(u_c(x)=0.4\).

Why does \(r(x)\) differ across locations?  Substituting
\(s_D=p_Dm_D\) into the comparison gives, on
\(\Omega_-\cap\mathcal D_T\) and for \(0\le u<1\),
\begin{equation}
  \frac{u\,|s_S(x)|}{(1-u)\,|s_T(x)|}
  =
  \underbrace{\frac{u}{1-u}}_{\substack{\text{source-to-target}\\\text{odds}}}\;
  \underbrace{\frac{p_S(x)}{p_T(x)}}_{\text{density ratio}}\;
  \underbrace{\frac{|m_S(x)|}{|m_T(x)|}}_{\text{margin ratio}},
  \label{eq:three_factor_vote_ratio}
\end{equation}
and the mixture reverses the target label exactly when this
product exceeds one.  The first factor is the same at every
location; the other two vary with \(x\).  With \(u=1/3\), a source
density three times the target density and opposing margins of
equal size give \(\tfrac12\times3\times1=1.5\): the source decides
the local label although it receives the smaller total weight.  A
location with a larger opposing source vote relative to the
target vote has a smaller \(u_c(x)\) and reverses at a lower
source weight.  At the location just described, \(r(x)=3\) and
\(u_c(x)=1/4\): reducing the source weight below \(1/4\) would
preserve the target label there.  Because \(u\) is global, the
same reduction also lowers the source's contribution at every
other location, including locations where the source and target
agree.  Changing \(u\) leaves \(p_T\), \(p_S\), \(m_T\), and
\(m_S\) as they are; it changes only their relative influence in
the mixture.

At a fixed \(u\), the target-active locations that have tied or
reversed form
\begin{equation}
  \mathcal C(u)
  :=\{x\in\Omega_-\cap\mathcal D_T:u_c(x)\le u\},
  \label{eq:charged_as_threshold_sublevel}
\end{equation}
and these sets are nested as \(u\) grows.  On
\(\Omega_+\cap\mathcal D_T\) the target label is retained for
\(u<1\); on \(\Omega_-\cap\mathcal D_T\) it is retained below
\(u_c(x)\), tied at \(u_c(x)\), and reversed above it.  We call
these three parts of the disagreement region \emph{surviving
opposition}, \emph{ties}, and \emph{strict flips}, written
\(\Omega_-^{\rm no\text{-}flip}(\alpha)\),
\(\Omega_-^{\rm tie}(\alpha)\), and \(\Omega_-^{\rm flip}(\alpha)\),
so that
\(\mathcal C(u)=(\Omega_-^{\rm tie}(\alpha)\cup\Omega_-^{\rm flip}(\alpha))\cap\mathcal D_T\);
the three regions carry the target weight \(\alpha\) as their
index, while \(\mathcal C(u)\) and the flip profile
\(H_{\rm flip}(u)\) of Section~\ref{subsec:flip_profile} carry the
source weight \(u=1-\alpha\).
At the source-only endpoint \(u=1\), target-active locations with
a zero source vote become ties as well (Supplement,
Remark~\ref{cmp:source-only-endpoint}).  A single weight can thus
change the classification rule in some regions and leave it
unchanged in others.  When the population vote is continuous, the
boundary between regions with opposite preferred labels lies in
its zero set.

Raising the source weight changes the mixture vote continuously,
and the preferred label changes only where the vote crosses zero.
At a fixed location in \(\Omega_-\cap\mathcal D_T\), the vote in the
target's direction weakens as \(u\) rises, the label stays until
\(u\) reaches \(u_c(x)\), the weighted votes cancel there, and beyond it
the mixture favors the opposite label with a vote that grows
again.  Figure~\ref{fig:vote}(a) shows the same process across
locations.  Each location sits at its own \(r\) on the horizontal
axis, and raising \(u\) moves it up its vertical line until it
crosses the threshold curve \(u_c(r)\).  At a fixed weight
\(u_0\), the horizontal line through \(u_0\) separates the
locations that have reached their thresholds, \(r\ge u_0^{-1}-1\),
from those that have not; raising the weight moves this line up
and its crossing with the curve to the left, so more locations
enter \(\mathcal C(u)\).  A location below the curve keeps its
label while its target-directed vote has already weakened.
Section~\ref{subsec:local_cost} measures that change.

\begin{figH}[!t]
\centering
\begin{tikzpicture}[font=\footnotesize,>=stealth,line cap=round,line join=round]
\definecolor{ink}{HTML}{111111}
\definecolor{bblue}{HTML}{1D4E89}
\definecolor{bred}{HTML}{D62828}
\definecolor{byellow}{HTML}{F4C20D}
\begin{scope}[x=1.55cm,y=3.3cm]
  \def\uz{0.4}
  \draw[->,ink] (0,0) -- (3.6,0) node[anchor=south east,yshift=3pt,xshift=-2pt] {local vote ratio $r$};
  \draw[->,ink] (0,0) -- (0,1.12) node[anchor=south west,xshift=-5pt] {source weight $u$};
  \draw[bblue,very thick] plot[domain=0:3.4,samples=80] (\x,{1/(1+\x)});
  \node[bblue,anchor=south west] at (0.42,0.80) {$u_c(r)$: tie};
  \node[ink,anchor=west] at (1.75,0.78) {flip};
  \node[ink,anchor=west,align=left] at (0.3,0.24) {surviving\\opposition};
  \draw[ink,dashed] (0,\uz) -- (1.5,\uz);
  \draw[bred,very thick] (1.5,\uz) -- (3.4,\uz);
  \fill[ink] (1.5,\uz) circle (1.4pt);
  \node[bred,anchor=south] at (2.55,\uz) {$\mathcal C(u_0)$};
  \draw[ink,dashed] (1.5,\uz) -- (1.5,0);
  \draw[ink] (-0.06,1) -- (0,1) node[anchor=east,xshift=-2pt] {1};
  \draw[ink] (-0.06,\uz) -- (0,\uz) node[anchor=east,xshift=-2pt] {$u_0$};
  \node[anchor=east,xshift=-2pt] at (0,0) {0};
  \draw[ink] (1.5,-0.02) -- (1.5,0) node[anchor=north,yshift=-2pt] {$u_0^{-1}-1$};
  \node[anchor=north,yshift=-2pt] at (0,0) {0};
  \node[anchor=north west] at (0.08,1.1) {\textbf{(a)}};
\end{scope}
\begin{scope}[xshift=7.45cm,x=5.2cm,y=3.3cm]
  \def\uz{0.4}
  \draw[->,ink] (0,0) -- (1.08,0) node[anchor=south east,yshift=3pt,xshift=-2pt] {source weight $u$};
  \draw[->,ink] (0,0) -- (0,1.12) node[anchor=south west,xshift=-5pt] {target-margin mass};
  \draw[ink,dotted,thick] (0,0.8) -- (1,0.8);
  \draw[bblue,very thick] (0,0) -- (0.15,0) -- (0.35,0.45) -- (0.5,0.45) -- (0.62,0.8) -- (1,0.8);
  \node[bblue,anchor=north west] at (0.66,0.76) {$H_{\rm flip}(u)$};
  \draw[ink,dashed] (\uz,0) -- (\uz,0.45) -- (0,0.45);
  \fill[ink] (\uz,0.45) circle (1.4pt);
  \draw[ink] (-0.012,0.45) -- (0,0.45) node[anchor=east,xshift=-2pt] {$H_{\rm flip}(u_0)$};
  \draw[ink] (-0.012,0.8) -- (0,0.8) node[anchor=east,xshift=-2pt] {$2\nu_T(\Omega_-)$};
  \node[anchor=east,xshift=-2pt] at (0,0) {0};
  \draw[ink] (\uz,-0.02) -- (\uz,0) node[anchor=north,yshift=-2pt] {$u_0$};
  \node[anchor=north,yshift=-2pt] at (0,0) {0};
  \node[anchor=north west] at (0.025,1.1) {\textbf{(b)}};
\end{scope}
\end{tikzpicture}

\caption{Each opposed location flips at its own source weight,
and the flip profile is the largest target hinge excess the tied
and flipped locations can incur at zero mixture excess.  (a) A location with
local vote ratio \(r\) flips once the source weight passes
\(u_c(r)=(1+r)^{-1}\) (blue curve: ties on the curve, flips above
it, surviving opposition below it).  Raising \(u\) moves a location
up its vertical line; below the curve it keeps its label while its
target-directed vote weakens.  At a fixed weight \(u_0\), the
horizontal line picks out the locations with \(r\ge u_0^{-1}-1\),
the tie-or-flip set \(\mathcal C(u_0)\) of
\eqref{eq:charged_as_threshold_sublevel} (red segment); raising
\(u_0\) moves the line up and the segment to the left.  (b) The
flip profile \(H_{\rm flip}(u)\) of \eqref{eq:flip_profile} is the
mass \(2\nu_T\) of \(\mathcal C(u)\) under the target-margin
measure \(\nu_T\) of \eqref{eq:margin_measure} (blue).  It never
decreases as \(u\) grows, it is flat where the locations reaching
their thresholds carry zero target-margin mass, it is defined for
\(u<1\), and it stops at the mass \(2\nu_T(\Omega_-)\) of all
opposed locations (dotted).  Population quantities; the profile is
a schematic shape.}
\label{fig:vote}
\end{figH}

\subsection{Local Error Costs under a Retained Label}
\label{subsec:local_cost}

Two distributions that favor the same label at \(x\) can penalize
a wrong prediction there by different amounts.  If
\(P(Y=1\mid X=x)=0.9\), predicting \(-1\) instead of \(+1\) raises
the conditional classification error from \(0.1\) to \(0.9\), an
increase of \(0.8=|2(0.9)-1|\); a distribution with
\(P(Y=1\mid X=x)=0.6\) favors the same label and loses only
\(0.2\).

For a distribution \(P\) on \(\mathcal X\times\mathcal Y\), write
\(\eta_P\), \(m_P:=2\eta_P-1\), and \(f_P^*:=\operatorname{sign}m_P\)
as in Section~\ref{subsec:mixture_geometry}, and
\(\mathcal R_{L,P}\) and \(\mathcal R_{L_{01},P}\) for its hinge and
classification risks.  Clipping a score to
\(\bar t=\max\{-1,\min\{t,1\}\}\) satisfies
\((1-y\bar t)_+\le(1-yt)_+\) for both labels, so clipping cannot
increase hinge risk and preserves the predicted sign; the
regularization term of \eqref{eq:transfer_svm_empirical} uses the
RKHS score itself, since a clipped score need not belong to
\(H\).  We therefore compare risks for measurable clipped scores
\(g:\mathcal X\to[-1,1]\), without loss for the empirical SVM
after clipping.  Their conditional hinge risk is affine,
\(\mathbb E_P[L(Y,g(x))\mid X=x]=1-m_P(x)g(x)\), which gives the
excess-risk identity and the calibration inequality used
throughout:
\begin{equation}
  \mathcal R_{L,P}(g)-\mathcal R^*_{L,P}
  =
  \int |m_P(x)|\,|g(x)-f_P^*(x)|\,dP_X(x),
  \label{eq:hinge_pointwise_excess}
\end{equation}
\begin{equation}
  \mathcal R_{L_{01},P}(g)
  -\mathcal R^*_{L_{01},P}
  =
  \int |m_P(x)|
  \mathbf 1\{\operatorname{sign}g(x)\ne f_P^*(x)\}\,dP_X(x)
  \le
  \mathcal R_{L,P}(g)-\mathcal R^*_{L,P},
  \label{eq:hinge_classification_calibration}
\end{equation}
the inequality because a sign error forces \(|g-f_P^*|\ge1\)
\citep{zhang2004behavior,bartlett2006convexity}.  The margin sets
the cost of choosing the less likely label, and the covariate
distribution sets how often that location occurs; every target
hinge-excess bound below therefore also controls the
classification excess of the induced classifier.

Now compare these costs under the mixture and under the target.
At target-active points, define the \emph{normalized mixture vote}
\begin{equation}
  Q_\alpha(x)
  :=
  \frac{p_\alpha(x)m_\alpha(x)}{p_T(x)m_T(x)}
  =\alpha+u\,\frac{s_S(x)}{s_T(x)},
  \qquad x\in\mathcal D_T .
  \label{eq:Q_def}
\end{equation}
Its sign says whether the mixture retains the target label, and
where it is positive its size says how the two distributions
penalize the same score error: \(Q_\alpha>1\) means a larger
mixture excess-risk density than the target's, \(0<Q_\alpha<1\) a
smaller one, and \(Q_\alpha=1\) the same; \(Q_\alpha=0\) is a
vanishing mixture vote, and \(Q_\alpha<0\) a disagreement of the
preferred labels.  On \(\Omega_-\cap\mathcal D_T\),
\(Q_\alpha=(1-u)-u\,r(x)\).  The ratio \(Q_\alpha\) compares
density-weighted votes; Section~\ref{sec:mixture_to_target} also
uses the mixture margin \(m_\alpha=s_\alpha/p_\alpha\), the
posterior margin of \(P_\alpha\), which divides the mixture vote
by the mixture density.

Define the pointwise excess hinge-risk densities
\(e_D(g;x):=p_D(x)|m_D(x)|\,|g(x)-f_D^*(x)|\), \(D\in\{T,\alpha\}\),
whose integrals against \(\mu\) are the hinge excess risks by
\eqref{eq:hinge_pointwise_excess}.  On the \emph{stable set}
\(\mathcal D_\alpha^+:=\mathcal D_T\cap\{Q_\alpha>0\}\) the preferred
labels agree, \(f_\alpha^*=f_T^*\), and \(|s_\alpha|=Q_\alpha|s_T|\),
so that, with the \emph{local multiplier}
\(A_\alpha(x):=Q_\alpha(x)^{-1}\), every measurable
\(g:\mathcal X\to[-1,1]\) satisfies, \(\mu\)-a.e.\ on
\(\mathcal D_T\),
\begin{align}
  e_T(g;x)
  &=A_\alpha(x)\,e_\alpha(g;x),
  &&Q_\alpha(x)>0,
  \label{eq:pointwise_stable_conversion}\\
  e_T(g;x)
  &\le 2\,p_T(x)|m_T(x)|,
  &&Q_\alpha(x)\le0,
  \label{eq:pointwise_charged_conversion}
\end{align}
the second line because \(|g-f_T^*|\le2\).  If \(Q_\alpha(x)=0.1\),
the same score error has one tenth as much excess hinge cost under
the mixture as under the target.  As \(u\) approaches \(u_c(x)\)
from below, this mixture cost approaches zero while the target
cost of that fixed error stays the same.  At a tie or a strict
flip the score \(g(x)=-f_T^*(x)\) incurs zero mixture excess and
attains the upper bound in the second line, so a score that is
optimal for the mixture can be wrong for the target.  In terms of \(Q_\alpha\), the decision regions on
\(\mathcal D_T\) are, for \(0<\alpha<1\),
\begin{equation}
  \begin{array}{ccll}
  Q_\alpha\ge\alpha
  &:&
  \text{agreement or zero source vote}
  &(\Omega_+),\\[1mm]
  0<Q_\alpha<\alpha
  &:&
  \text{surviving opposition}
  &(\Omega_-^{\rm no\text{-}flip}(\alpha)),\\[1mm]
  Q_\alpha\le0
  &:&
  \text{tie or flip}
  &(\Omega_-^{\rm tie}(\alpha)\cup\Omega_-^{\rm flip}(\alpha)),
  \end{array}
  \label{eq:vote_cases}
\end{equation}
and the multiplier is at most \(1/\alpha\) under agreement and
exceeds \(1/\alpha\) under surviving opposition;
Corollary~\ref{cor:regime_flip} and Section~\ref{sec:numerical}
refer to the regimes by these names.

\subsection{Tie-or-Flip Sets and the Flip Profile}
\label{subsec:flip_profile}

The set \(\mathcal C(u)\) and the bound
\eqref{eq:pointwise_charged_conversion} say, location by location,
where the label can change and how much target error is possible
there.  To add these contributions over locations, define the
\emph{margin-weighted design measure}
\begin{equation}
  \nu_D(E):=\int_E|s_D|\,d\mu=\int_E|m_D(x)|\,dP_{X,D}(x),
  \qquad D\in\{T,S,\alpha\},
  \label{eq:margin_measure}
\end{equation}
and the \emph{total target margin}
\(J_T:=\nu_T(\mathcal X)=\mathbb E_{P_{X,T}}|m_T(X)|\).  The factor \(|m_T(x)|\) is the extra
conditional classification error from choosing the less likely
target label, and \(P_{X,T}\) is how often such locations occur
under the target.  For \(0\le u<1\), the largest target hinge
excess that the tie-or-flip set can incur without any mixture
excess is the \emph{flip profile}
\begin{equation}
  H_{\rm flip}(u)
  :=
  2\,\nu_T\{Q_\alpha\le0\}
  =2\int_{\Omega_-\cap\mathcal D_T}
  |m_T(x)|\mathbf 1\{u_c(x)\le u\}\,dP_{X,T}(x),
  \qquad 0\le u<1,
  \label{eq:flip_profile}
\end{equation}
that is, \(2\nu_T(\mathcal C(u))\), the target-margin-weighted
(unnormalized) distribution function of the threshold \(u_c(X)\)
(Figure~\ref{fig:vote}(b)).  The factor two is the largest
difference between a clipped score and the target-optimal score,
and \(\mathcal C(u)\) grows with \(u\), so the profile cannot
decrease.  Where it is flat, the locations reaching their
thresholds carry zero target-margin mass, while the votes of
locations below their thresholds still weaken;
\eqref{eq:pointwise_stable_conversion} records that change through
\(Q_\alpha\).  At a strict flip a mixture-optimal score takes the
target-opposed sign; at an exact tie the mixture is indifferent
among clipped scores (\(1-m_\alpha g=1\) when \(m_\alpha=0\)), so a
particular score need not incur the maximal target cost.  The
exact conversion of Section~\ref{subsec:risk_comparison} takes the
profile as an input.  The
other locations also matter: there the preferred label is
unchanged, and \(Q_\alpha(x)\) sets how much mixture excess a given
target error produces.  Section~\ref{sec:mixture_to_target}
combines the two and determines how much target error can remain
when a fitted score has mixture excess at most \(t\).

\textbf{Notation.}  A subscript \(D\in\{T,S,\alpha\}\) attaches an
object to the target, the source, or the training mixture.  Roman
capitals \(P\) denote laws and calligraphic \(\mathcal P\) classes of
laws.  Lowercase \(f\), \(g\), and \(h\) denote an RKHS score, a
clipped score, and a hard classifier; a hat marks a quantity
computed from the data and a bar the clipping of a score to
\([-1,1]\).

\section{From Mixture Learning to Target Guarantees}
\label{sec:mixture_to_target}
\label{sec:help}
\label{sec:comparison}

A classifier fitted on the mixture is evaluated on the target.
Section~\ref{sec:problem_setup} showed why these are different
classification problems: their preferred labels can differ, and
where they agree, a prediction error can have a different cost.
The classifier is also learned from finite samples.  Two questions
follow: how large the target excess of a score with mixture excess
at most \(t\) can be, which Theorem~\ref{thm:affine_frontier}
answers for every fixed pair of distributions and every positive
target weight, and how
small \(t\) can be for the weighted SVM fitted on \(n\) target and
\(m\) source observations, which Theorem~\ref{thm:weighted_crerm}
bounds.  Their combination is the target guarantee
\(\mathcal E_T(\bar f_D)\le\Psi_\alpha(\mathcal E_\alpha(\bar f_D))
\le\Psi_\alpha(t_{\alpha,n,m})\).  Figure~\ref{fig:conversion_exact}
draws this combination: the sampling bound fixes the abscissa, the
laws and the weight fix the curve, and the guarantee is the height
of the curve at that abscissa.  The source weight enters both: it
changes the classification problem being learned, through how
often the mixture margin \(m_\alpha=s_\alpha/p_\alpha\) is small
and through the approximation and complexity of the kernel class,
and it changes the relation between mixture errors and target
errors, through the distribution of \(Q_\alpha=s_\alpha/s_T\)
under the target-margin measure.  Adding source observations at a
fixed weight moves the abscissa along a fixed curve; changing the
weight moves the abscissa and changes the curve.

\subsection{Exact Conversion}
\label{subsec:risk_comparison}

Suppose a clipped predictor has mixture excess at most \(t\).  How
large can its target excess still be?  For a fixed pair of laws
and a fixed weight, the answer is the following conversion
modulus.  Fix \(P_T\), \(P_S\), and \(\alpha\in(0,1]\), and write
\(\mathcal E_D(f):=\mathcal R_{L,P_D}(f)-\mathcal R^*_{L,P_D}\),
\(D\in\{T,\alpha\}\), for the hinge excess risk of any measurable
\(f\) under the target and under the mixture; for a clipped score
\(g\) it has the pointwise form \eqref{eq:hinge_pointwise_excess}.  For \(t\ge0\), define
\begin{equation}
  \Psi_\alpha(t)
  :=
  \sup\bigl\{\mathcal E_T(g):
  \mathcal E_\alpha(g)\le t,\
  g:\mathcal X\to[-1,1]\ \text{measurable}\bigr\}.
  \label{eq:modulus_def}
\end{equation}
The \emph{conversion modulus} \(\Psi_\alpha\) is a population
quantity: the laws and the weight are fixed, and the supremum runs
over clipped predictors.  A linear comparison is often easier to use once a
statistical bound on mixture excess is available.  For a chosen
slope \(c\), we therefore also ask for the smallest additive
term that makes the comparison valid for every clipped predictor:
for each slope \(c\ge0\), define
\begin{equation}
  \mathcal B_\alpha(c)
  :=
  \sup_{g:\mathcal X\to[-1,1]}
  \bigl\{\mathcal E_T(g)-c\,\mathcal E_\alpha(g)\bigr\},
  \label{eq:frontier_def}
\end{equation}
so that, by construction,
\(\mathcal E_T(g)\le c\,\mathcal E_\alpha(g)+\mathcal B_\alpha(c)\)
holds uniformly over clipped \(g\): the slope \(c\) fixes how much
the mixture excess is amplified, and \(\mathcal B_\alpha(c)\) is the
smallest uniform additive allowance that slope requires.  We call
\(c\mapsto\mathcal B_\alpha(c)\) the \emph{affine frontier}.  The
theorem below computes both objects exactly.

The worst predictor places its mixture excess where that excess
buys the most target excess.  At a tie or a strict flip, it can
take the target-opposed sign at no mixture cost, so the budget
\(t=0\) already allows the target excess \(H_{\rm flip}(u)\) of
\eqref{eq:flip_profile}.  On the stable set
\(\mathcal D_\alpha^+=\mathcal D_T\cap\{Q_\alpha>0\}\), one unit of
mixture excess buys \(A_\alpha=Q_\alpha^{-1}\) units of target
excess \eqref{eq:pointwise_stable_conversion}, and each location
accepts a limited amount: moving a clipped score across its full
range at a stable point creates at most \(2p_\alpha|m_\alpha|\)
units, so the budget available at each location is the measure
\(2\nu_\alpha\) of \eqref{eq:margin_measure} restricted to the
stable set.  Its total mass is the vote capacity
\begin{equation}
  2\,\nu_\alpha(\mathcal D_\alpha^+)
  =2\int_{\mathcal D_\alpha^+}p_\alpha|m_\alpha|\,d\mu .
  \label{eq:surviving_vote_measure}
\end{equation}
The worst predictor spends its budget at the largest multiplier
until that capacity is exhausted, then at the next largest, so
each further unit buys less target excess; once every location
carries its full target cost \(2p_T|m_T|\), the total is \(2J_T\).
The ordering is the decreasing rearrangement
\(A_\alpha^{\downarrow}(s)
:=\inf\{a\ge0:2\nu_\alpha(\mathcal D_\alpha^+\cap\{A_\alpha>a\})\le s\}\),
\(0<s<2\nu_\alpha(\mathcal D_\alpha^+)\), of the multiplier under \(2\nu_\alpha|_{\mathcal D_\alpha^+}\).  Its
essential supremum is
\(M_\alpha:=\operatorname*{ess\,sup}_{\nu_\alpha|_{\mathcal D_\alpha^+}}A_\alpha\),
with \(M_\alpha:=0\) when \(\nu_\alpha(\mathcal D_\alpha^+)=0\).

\begin{theorem}[Exact conversion]
\label{thm:affine_frontier}
\label{rem:vote_profile}
\label{prop:pareto_conversion}
\label{thm:vote_conversion}
Fix \(\alpha\in(0,1]\) and write \(u=1-\alpha\).
\begin{itemize}
\item[(i)] \emph{(Exact modulus.)}  For every \(t\ge0\),
\begin{equation}
  \Psi_\alpha(t)
  =
  H_{\rm flip}(u)
  +\int_0^{\,t\wedge2\nu_\alpha(\mathcal D_\alpha^+)}A_\alpha^{\downarrow}(s)\,ds,
  \label{eq:conversion_modulus}
\end{equation}
and the supremum in \eqref{eq:modulus_def} is attained by a
clipped score at every \(t\).
\item[(ii)] \emph{(Affine frontier and duality.)}  For
every \(c\ge0\),
\begin{align}
  \mathcal B_\alpha(c)
  &=
  \int\bigl(
  p_T|m_T|-c\,p_\alpha|m_\alpha|
  +\bigl|p_Tm_T-c\,p_\alpha m_\alpha\bigr|
  \bigr)\,d\mu
  \notag\\
  &=
  H_{\rm flip}(u)
  +2\int_{\mathcal D_\alpha^+}(A_\alpha-c)_+\,d\nu_\alpha
  =
  2\int_{\mathcal D_T}
       \bigl[1-c\,(Q_\alpha)_+\bigr]_+\,d\nu_T,
  \label{eq:frontier_closed}
\end{align}
where \((z)_+:=\max\{z,0\}\), and
\begin{equation}
  \Psi_\alpha(t)
  =
  \inf_{c\ge0}\bigl\{c\,t+\mathcal B_\alpha(c)\bigr\},
  \qquad t\ge0.
  \label{eq:conversion_duality}
\end{equation}
\end{itemize}
\end{theorem}

Figure~\ref{fig:conversion_exact}(a) draws
\eqref{eq:conversion_modulus}.  The curve starts at
\begin{equation}
  \Psi_\alpha(0)=H_{\rm flip}(u)
  =2\,\nu_T\{Q_\alpha\le0\},
  \label{eq:profile_representation}
\end{equation}
the cost of the tied and flipped locations; its initial slope is
the largest multiplier \(M_\alpha\); it bends as the budget passes
to locations with smaller multipliers; and it is flat at
\(2\nu_T(\mathcal X)=2J_T\) once
\(t\ge2\nu_\alpha(\mathcal D_\alpha^+)\).  Attainment at every
budget makes it the strongest uniform target guarantee available
from mixture excess alone.  The straight lines in the panel are
the bounds \(ct+\mathcal B_\alpha(c)\) at the same weight: the
slope \(c\) varies at a fixed weight, and by
\eqref{eq:conversion_duality} the curve is their lower envelope.
A different weight changes the mixture vote \(s_\alpha\) at every
location, and with it the intercept and every slope of the curve
(Example~\ref{ex:three_point} in the Supplement computes one
curve).

The three expressions in \eqref{eq:frontier_closed} identify what
determines this comparison: the first uses no assumption on the
domains, the second allocates the mixture budget on the stable
set under the mixture-margin measure, and the third reads the
same quantity under the target-margin measure, charging flipped
locations in full and stable locations according to
\(Q_\alpha\).  The unnormalized distribution of
\((Q_\alpha)_+\) under \(\nu_T\) therefore determines the entire
frontier and, by duality, the conversion curve.  Different spatial
arrangements with the same weighted profile have the same
conversion, although their approximation and learning costs can
differ.

The full proof, hard-classifier attainment,
Example~\ref{ex:three_point}, and the assumption-free comparison of
Remark~\ref{rem:finite_conversion} are in Supplement
Section~\ref{pf:frontier}.

\begin{remark}[Affine bounds and sharpness]
\label{rem:affine_sharpness}
A worst-case multiplier treats every error as though it occurred
at the most expensive location.  This multiplier can be infinite
near cancellation even when the total target cost of those
locations is small.  Every affine intercept is at least
\(H_{\rm flip}(u)\), and \(\inf_{c\ge0}\mathcal B_\alpha(c)=H_{\rm flip}(u)\)
(Figure~\ref{fig:conversion_exact}(a)).  When \(M_\alpha<\infty\),
every measurable \(g:\mathcal X\to[-1,1]\) satisfies
\begin{equation}
  \mathcal E_T(g)
  \le
  M_\alpha\,\mathcal E_\alpha(g)
  +H_{\rm flip}(u),
  \label{eq:global_vote_conversion}
\end{equation}
and \(M_\alpha\) is the smallest slope attaining the minimal
intercept.  If \(M_\alpha=\infty\), that is, if
\(\nu_\alpha(\mathcal D_\alpha^+\cap\{A_\alpha>K\})>0\) for every finite \(K\), no
finite slope attains that intercept, and a finite-slope comparison
needs a larger intercept, set by how much target cost the
near-canceled locations carry: Proposition~\ref{prop:radial_frontier}
computes it under boundary tilt, and
Corollary~\ref{cor:regime_flip}(iii) under a small-vote mass
condition.
\end{remark}

\begin{figH}[\figposconv]
\centering
\begin{tikzpicture}[font=\footnotesize,>=stealth,line cap=round,line join=round]
\definecolor{ink}{HTML}{111111}
\definecolor{bblue}{HTML}{1D4E89}
\definecolor{bred}{HTML}{D62828}
\definecolor{byellow}{HTML}{F4C20D}
\begin{scope}[x=3.85cm,y=2.8cm]
  \def\Hf{0.25}  
  \def\tb{0.35}  
  \fill[pattern=north east lines,pattern color=ink!45] (0,0) rectangle (1.25,\Hf);
  \draw[ink,dotted,thick] (0,\Hf) -- (1.25,\Hf);
  \node[ink,anchor=south west] at (0.17,\Hf) {tie-or-flip contribution};
  \draw[->,ink] (0,0) -- (1.62,0) node[anchor=north east,yshift=-10pt] {mixture excess $t$};
  \draw[->,ink] (0,0) -- (0,1.42) node[anchor=south west,xshift=-5pt] {target excess bound};
  \draw[ink,dotted,thick] (0,1) -- (1.25,1);
  \draw[ink,dashed] (1,0) -- (1,1);
  \draw[byellow,very thick] (0,\Hf) -- (0.633,1.2);
  \node[byellow!80!black,anchor=south] at (0.60,1.21) {$M_\alpha t+H_{\rm flip}(u)$};
  \draw[bred,very thick] (0,0.3419) -- (0.88,1.2);
  \node[bred,anchor=south west] at (0.91,1.17) {$ct+\mathcal B_\alpha(c)$};
  \draw[bblue,very thick] plot[domain=0:1,samples=80] (\x,{\Hf+(1-\Hf)*(1-(1-\x)^2)}) -- (1.25,1);
  \node[bblue,anchor=north] at (0.84,0.90) {$\Psi_\alpha(t)$};
  \fill[ink] (\tb,{\Hf+(1-\Hf)*(1-(1-\tb)^2)}) circle (1.4pt);
  \draw[ink] (-0.02,\Hf) -- (0,\Hf) node[anchor=east,xshift=-2pt] {$H_{\rm flip}(u)$};
  \draw[ink] (-0.02,0.3419) -- (0,0.3419) node[anchor=east,xshift=-2pt,yshift=4pt] {$\mathcal B_\alpha(c)$};
  \draw[ink] (-0.02,1) -- (0,1) node[anchor=east,xshift=-2pt] {$2J_T$};
  \node[anchor=east,xshift=-2pt] at (0,0) {0};
  \draw[ink] (1,-0.025) -- (1,0) node[anchor=north east,yshift=-2pt,xshift=2pt] {$2\nu_\alpha(\mathcal D_\alpha^+)$};
  \node[anchor=north,yshift=-2pt] at (0,0) {0};
  \node[anchor=north west] at (0.02,1.40) {\textbf{(a)}};
\end{scope}
\begin{scope}[xshift=7.05cm,x=4.15cm,y=2.8cm]
  \def\Hf{0.25}
  \def\ta{0.27}
  \def\tb{0.8}
  \fill[pattern=north east lines,pattern color=ink!45] (0,0) rectangle (1.1,\Hf);
  \draw[ink,dotted,thick] (0,\Hf) -- (1.1,\Hf);
  \draw[->,ink] (0,0) -- (0,1.42) node[anchor=south west,xshift=-5pt] {target excess bound};
  \draw[->,ink] (0,0) -- (1.2,0) node[anchor=north east,yshift=-10pt] {mixture excess $t$};
  \draw[ink,dashed,thin] (\ta,0) -- (\ta,1.28);
  \draw[ink,dashed,thin] (\tb,0) -- (\tb,1.28);
  \draw[byellow,very thick] (0,0) -- (0.6,1.2);
  \node[byellow!80!black,anchor=east] at (0.52,1.12) {$M_\alpha t$};
  \draw[bblue,very thick] (0,0) -- (1.05,0.525);
  \node[bblue,anchor=north west] at (0.90,0.42) {$t/q_0$};
  \draw[bred,very thick] plot[domain=0:1.05,samples=60] (\x,{\Hf+0.8*sqrt(\x)});
  \node[bred,anchor=north west] at (0.83,0.92) {$H_{\rm flip}(u)+Ct^{\chi}$};
  \fill[ink] (\ta,0.6657) circle (1.4pt);
  \fill[ink] (\ta,0.54) circle (1.4pt);
  \fill[ink] (\ta,0.135) circle (1.4pt);
  \fill[ink] (\tb,0.9655) circle (1.4pt);
  \draw[->,ink,thick] (\ta,0.632) -- (\ta,0.17);
  \node[ink,anchor=west,inner sep=1pt] at (0.30,0.46) {better conversion};
  \draw[->,ink,thick] (\tb,1.28) -- (0.30,1.28);
  \node[ink,anchor=south] at (0.55,1.30) {smaller mixture bound};
  \draw[ink] (-0.015,\Hf) -- (0,\Hf) node[anchor=east,xshift=-2pt] {$H_{\rm flip}(u)$};
  \node[anchor=east,xshift=-2pt] at (0,0) {0};
  \draw[ink] (\ta,-0.025) -- (\ta,0) node[anchor=north,yshift=-2pt] {$t_a$};
  \draw[ink] (\tb,-0.025) -- (\tb,0) node[anchor=north,yshift=-2pt] {$t_b$};
  \node[anchor=north,yshift=-2pt] at (0,0) {0};
  \node[anchor=north west] at (0.02,1.40) {\textbf{(b)}};
\end{scope}
\end{tikzpicture}

\caption{A mixture bound \(t\) is worth exactly \(\Psi_\alpha(t)\)
on the target, and the shape of \(\Psi_\alpha\) sets the rate
that survives the conversion.  (a) One law pair and one weight,
with \(M_\alpha<\infty\).  The curve \(\Psi_\alpha\) (blue,
Theorem~\ref{thm:affine_frontier}) starts at the flip cost
\(H_{\rm flip}(u)\), rises concavely, and saturates at \(2J_T\)
once the budget reaches the vote capacity
\(2\nu_\alpha(\mathcal D_\alpha^+)\).  Each sharp affine bound
\(ct+\mathcal B_\alpha(c)\) touches the curve at one budget (red,
tangent at the dot); the line of slope \(M_\alpha\) through
\(H_{\rm flip}(u)\) (yellow) has the smallest intercept any
affine bound can have; by \eqref{eq:conversion_duality} the curve
is the lower envelope of all these lines.  (b) The three
conversion shapes of Corollary~\ref{cor:regime_flip}.  Without
opposition, or with opposition held below the flip threshold, the
conversion is linear through the origin with slope \(M_\alpha\)
(yellow).  For the constant vote \(Q_\alpha\equiv q_0\) of
Example~\ref{ex:constant_Q} the exact curve is the line \(t/q_0\) up to its saturation level \(2J_T\)
(blue).  Near cancellation the bound is H\"older,
\(H_{\rm flip}(u)+Ct^{\chi}\) (red), and it starts at the flip
cost.  At one budget \(t_a\) the three shapes give three target
guarantees (dots); a smaller mixture bound (\(t_b\to t_a\)) moves
along one curve, a better source--target relation moves down to
another curve, and a different weight changes the curve.  In both panels the hatched band
below the level \(H_{\rm flip}(u)=\Psi_\alpha(0)\) is the largest
cost the tied and flipped locations can still carry at zero
mixture excess; the two linear relations have no such cost.  A
mixture bound of order \(N_{\rm eff}^{-\omega}\), with
\(N_{\rm eff}\) the effective sample size of
\eqref{eq:two_block_scales}, converts to \(O(N_{\rm eff}^{-\omega})\)
under the linear shapes and to
\(H_{\rm flip}(u)+O(N_{\rm eff}^{-\omega\chi})\) under the
H\"older shape.  Colors keep one meaning in both panels: blue
is the exact curve, yellow the slope-\(M_\alpha\) line, red a
bound that pays a flip intercept.  Schematic population
quantities.}
\label{fig:conversion_exact}
\label{fig:regime_bounds}
\end{figH}

\subsection{Target Excess for the Weighted SVM}
\label{sec:general_oracle}
\label{subsec:baseline_oracle}

The conversion says what mixture accuracy is worth; this
subsection bounds how much accuracy the samples provide, combining the sampling contributions of the two domains into an effective sample size.  Index the
combined sample by \(Z_1,\ldots,Z_N\), \(N=n+m\), assign
weights \(w_i=\alpha/n\) to target observations and
\(w_i=(1-\alpha)/m\) to source observations, and let
\(D_w=\sum_iw_i\delta_{Z_i}\), so that
\(\mathbb E\,\mathbb E_{D_w}g=\mathbb E_{P_\alpha}g\) and
\(\mathcal R_{L,D_w}\) is exactly the empirical risk in
\eqref{eq:transfer_svm_empirical}.  Second moments square the
observation weights, so the two blocks contribute at the scales
\(\alpha^2/n\) and \(u^2/m\), obtained by summing the squared
observation weights within each block: for every measurable \(h\),
\begin{equation}
  \sum_{i=1}^Nw_i^2\,\mathbb E\,h^2(Z_i)
  =
  \frac{\alpha^2}{n}\,\mathbb E_{P_T}h^2
  +\frac{u^2}{m}\,\mathbb E_{P_S}h^2.
  \label{eq:variance_two_block_exact}
\end{equation}
The largest observation weight gives a single scale for
controlling concentration; its reciprocal is the effective sample
size \(N_{\rm eff}\), and
the scalar reduction
\begin{equation}
  \frac{\alpha^2}{n}\,\mathbb E_{P_T}h^2
  +\frac{u^2}{m}\,\mathbb E_{P_S}h^2
  \;\le\;
  N_{\rm eff}^{-1}\,\mathbb E_{P_\alpha}h^2,
  \qquad
  N_{\rm eff}
  :=\Bigl(\max_iw_i\Bigr)^{-1}
  =\min\Bigl\{\frac{n}{\alpha},\frac{m}{u}\Bigr\},
  \label{eq:two_block_scales}
\end{equation}
uses the convention \(m/0:=+\infty\) and holds because \(\alpha^2/n\le\alpha/N_{\rm eff}\) and
\(u^2/m\le u/N_{\rm eff}\); target-only training (\(\alpha=1\))
gives \(N_{\rm eff}=n\), and proportional weighting
(\(\alpha=n/(n+m)\)) gives \(N_{\rm eff}=n+m\).
The concentration step applies the Klein--Rio inequality
\citep{klein2005concentration} to the independent weighted
observations from the two laws. Its scalar form uses
\(N_{\rm eff}\); the more detailed bound retains the two domainwise contributions
(Supplement, Section~\ref{pf:two_block_hinge_variance}).

Let \(H\) be a separable RKHS with measurable kernel bounded by
one, write \(\bar f\) for the clipping of \(f\) to \([-1,1]\), and
set \(h_f:=L\circ f-L\circ f^*_\alpha\), the excess hinge loss of
\(f\) relative to the mixture Bayes representative \(f^*_\alpha\) of
Section~\ref{sec:problem_setup}, which minimizes the hinge risk
(the Supplement writes \(f^*_{L,P_\alpha}\) in its general-loss
statements), so that
\(\mathbb E_{P_\alpha}h_{\bar f}=\mathcal E_\alpha(\bar f)\).
Functions with small excess risk may also have smaller loss
fluctuations; the following condition quantifies this relation and
allows concentration to improve as the candidate approaches the
mixture optimum: there are \(\vartheta_\alpha\in[0,1]\) and
\(V_\alpha\ge2^{2-\vartheta_\alpha}\) with
\begin{equation}
  \mathbb E_{P_\alpha}h_{\bar f}^2
  \le
  V_\alpha
  \bigl(\mathbb E_{P_\alpha}h_{\bar f}\bigr)^{\vartheta_\alpha},
  \qquad f\in H,
  \label{eq:variance_bound_alpha}
\end{equation}
always admissible with \(\vartheta_\alpha=0\), \(V_\alpha=4\), and
improved under source dominance in
Section~\ref{subsec:source_dominance}.

Sampling fluctuations are controlled only among functions whose
regularized mixture excess is small.  Let
\[
  r_0:=\inf_{f\in H}\bigl\{\lambda\|f\|_H^2+\mathcal E_\alpha(\bar f)\bigr\}
\]
be the attainable regularized excess, and let
\(\mathcal H_r:=\{h_{\bar f}:f\in H,\ \lambda\|f\|_H^2+\mathcal E_\alpha(\bar f)\le r\}\)
be the clipped excess-loss class at regularized mixture excess
\(r\).

\begin{itemize}
\item[{\rm(W)}] \emph{Localized weighted complexity:} there is
\(\varphi:[0,\infty)\to[0,\infty)\) with \(\varphi(4r)\le2\varphi(r)\)
such that, for every \(r>r_0\),
\[
  \mathbb E_D\,\mathbb E_\sigma\sup_{h\in\mathcal H_r}
  \Bigl|\sum_{i=1}^Nw_i\sigma_i\,h(Z_i)\Bigr|\le\varphi(r),
\]
where \(\sigma_1,\dots,\sigma_N\) are independent Rademacher signs
independent of \(D\).
\end{itemize}
Condition {\rm(W)} is the weighted form of the localized complexity bound used
for unweighted SVMs by \citet{steinwart2008support}; the Supplement
(Section~\ref{cmp:oracle-inputs}, equation
\eqref{eq:weighted_rad_assumption}) states it for general losses.
Fix a comparator
\(f_0\in H\) and \(B_0\ge2\) with \(\|L\circ f_0\|_\infty\le B_0\),
and a confidence parameter \(\varrho>0\).  The estimation radius
solves a localization fixed point: for a universal constant
\(C_0=216\), \(r_{\rm est}\) is any radius satisfying
\begin{equation}
  r_{\rm est}>C_0\max\Bigl\{\varphi(r_{\rm est}),\;
    \Bigl(\frac{V_\alpha\varrho}{N_{\rm eff}}\Bigr)^{
      \frac{1}{2-\vartheta_\alpha}},\;
    \frac{B_0\varrho}{N_{\rm eff}},\,r_0\Bigr\}.
  \label{eq:r_condition}
\end{equation}
The radius must dominate the local complexity, the variance
contribution, and the largest single-observation contribution,
while remaining above the attainable approximation level.  Given
an optimization slack \(\varepsilon_{\rm opt}\ge0\) and such a
radius, define the oracle bound on mixture excess,
\begin{equation}
  t_{\alpha,n,m}
  =
  t_{\alpha,n,m}(f_0,\lambda,\varrho,\varepsilon_{\rm opt};
  r_{\rm est})
  :=
  6\bigl(\lambda\|f_0\|_{H}^2
    +\mathcal E_\alpha(f_0)\bigr)
  +3r_{\rm est}
  +3\varepsilon_{\rm opt},
  \label{eq:oracle_budget}
\end{equation}
abbreviated \(t_{\alpha,n,m}\) hereafter.

\begin{theorem}[Target excess bound for the weighted SVM]
\label{thm:weighted_crerm}\label{thm:target_excess}
Let \(L\) be the hinge loss, \(\alpha\in(0,1]\), and fix
\(f_0\in H\) and \(B_0\ge2\) with
\(\|L\circ f_0\|_\infty\le B_0\).  Assume {\rm(W)} and the
variance condition \eqref{eq:variance_bound_alpha}, and let the
estimation radius \(r_{\rm est}\) satisfy the fixed-point condition
\eqref{eq:r_condition}.  Then every measurable
\(\varepsilon_{\rm opt}\)-approximate minimizer \(f_D\) of the
regularized empirical risk
\(f\mapsto\lambda\|f\|_{H}^2+\mathcal R_{L,D_w}(f)\) satisfies, with
probability at least \(1-3e^{-\varrho}\), the following two
bounds simultaneously:
\begin{equation}
  \lambda\|f_D\|_{H}^2+\mathcal E_\alpha(\bar f_D)
  \;\le\;
  t_{\alpha,n,m},
  \label{eq:oracle_mixture_bound}
\end{equation}
and
\begin{equation}
  \mathcal E_T(\bar f_D)
  \;\le\;
  \Psi_\alpha\bigl(\mathcal E_\alpha(\bar f_D)\bigr)
  \;\le\;
  \Psi_\alpha\bigl(t_{\alpha,n,m}\bigr).
  \label{eq:target_excess_master}
\end{equation}
The same bounds hold for the target classification excess of
\(\operatorname{sign}\bar f_D\), by
\eqref{eq:hinge_classification_calibration}.
\end{theorem}

The proof in the Supplement (Section~\ref{pf:weighted_crerm})
specializes Theorem~\ref{thm:weighted_crerm_general} to the hinge
loss and applies Theorem~\ref{thm:affine_frontier} on the same
event.  In Figure~\ref{fig:conversion_exact}(a), \(t_{\alpha,n,m}\)
is the abscissa and \(\Psi_\alpha\) the curve.  Two weights are
compared through \(\Psi_\alpha(t_{\alpha,n,m})\) at each, and the
weight with the smaller mixture budget can have the larger target
bound; at the same confidence level, borrowing is assessed by
comparing \(\Psi_\alpha(t_{\alpha,n,m})\) with the target-only
guarantee \(t_{1,n,m}\), the budget \eqref{eq:oracle_budget} at
\(\alpha=1\).

\begin{example}[The target value of the same mixture improvement]
\label{ex:constant_Q}
Let both designs be uniform on \([-1,1]\), take \(m_T(x)=x/4\) and
\(\alpha=1/2\), and compare a reinforcing source
\(m_S^{\rm rein}=3m_T\) with an opposing, sign-preserving source
\(m_S^{\rm opp}=-m_T/2\); by \eqref{eq:margin_identity} their
normalized mixture votes are \(Q_\alpha\equiv2\) and
\(Q_\alpha\equiv1/4\), respectively.  For a constant positive vote
\(Q_\alpha\equiv q_0\) on the target-active set, \(H_{\rm flip}(u)=0\),
\(A_\alpha\equiv1/q_0\), and \(2\nu_\alpha(\mathcal D_\alpha^+)=2q_0J_T\), so
Theorem~\ref{thm:affine_frontier}(i) gives
\begin{equation}
  \Psi_\alpha(t)=\min\{t/q_0,\;2J_T\},
  \qquad J_T=\nu_T(\mathcal X)=1/8.
  \label{eq:constant_Q_modulus}
\end{equation}
Take illustrative budgets \(t_\alpha=0.05\) and \(t_1=0.10\) for the
mixture and the target-only guarantee; then
\[
  \Psi_\alpha^{\rm rein}(0.05)=0.025<0.10,
  \qquad
  \Psi_\alpha^{\rm opp}(0.05)=0.20>0.10.
\]
Both mixtures have the target Bayes rule and zero tie-or-flip
cost, and the same mixture-excess guarantee produces target
guarantees that differ by a factor of eight.
\end{example}

\subsection{When Source Labels Make the Mixture Easier to Learn}
\label{subsec:source_dominance}
Source labels change how hard the mixture is to learn, and the
change begins at single locations: with a common design and
\(\alpha=1/2\), target and source probabilities \(0.55\) and \(0.80\) for
label \(+1\) give a mixture probability of \(0.675\), so the
mixture margin \(0.35\) exceeds the target margin \(0.10\).  A
learning guarantee depends on how often small mixture margins
occur under the mixture design, through the variance pair
\eqref{eq:variance_bound_alpha} of Theorem~\ref{thm:weighted_crerm}.
The following calculation turns the pointwise gain into a gain in
that frequency.  With a common design and
\(m_S=\operatorname{sign}(m_T)\,|m_T|^{\zeta}\) for some
\(0<\zeta<1\), raising \(|m_T|\le1\) to the power \(\zeta\) enlarges
every weak margin, and
\(|m_\alpha|=\alpha|m_T|+u|m_T|^{\zeta}\ge u|m_T|^{\zeta}\).  A
mixture margin below \(z\) therefore requires a target margin
below \((z/u)^{1/\zeta}\): if
\(P_{X,T}(|m_T|\le z)\le c_{q,T}z^{q_T}\) for all \(z>0\), then
\[
  P_{X,T}(|m_\alpha|\le z)
  \le c_{q,T}\,u^{-q_T/\zeta}\,z^{q_T/\zeta},
  \qquad z>0,
\]
and the mixture satisfies the target's margin condition with the
larger exponent \(q_\alpha=q_T/\zeta\).  This exponent enters
Theorem~\ref{thm:weighted_crerm} through the variance term of the
localization radius \eqref{eq:r_condition}, whose sample-size
exponent is \((q_\alpha+1)/(q_\alpha+2)\): with \(q_T=1\) and
\(\zeta=1/2\), \(q_\alpha=2\), and that term improves from order
\(N_{\rm eff}^{-2/3}\) to \(N_{\rm eff}^{-3/4}\).  The proposition
below states this gain for general designs and controls the small
mixture margins that opposition adds, under three conditions
\citep{mammen1999smooth,tsybakov2004optimal}, which describe the
frequency of ambiguous target labels, the strength of the source
signal on the agreement region, and its representation under the
mixture design; they enter only Proposition~\ref{lem:one_acceleration}
and the results built on it.

\begin{itemize}
\item[{\rm(T)}] \emph{Target margin condition:} there exist
\(q_T\ge0\) and \(c_{q,T}>0\) such that
\(P_{X,T}(|m_T|\le z)\le c_{q,T}\,z^{q_T}\) for all \(z>0\).
\item[{\rm(A1)}] \emph{Margin dominance on the agreement region:}
there exist \(c_\zeta>0\) and \(\zeta\in(0,1]\) such that
\(|m_S|\ge c_\zeta\,|m_T|^{\zeta}\) \(P_{X,T}\)-almost everywhere
on \(\Omega_+\).
\item[{\rm(A2)}] \emph{Bounded density ratio:} the marginals
\(P_{X,S}\) and \(P_{X,T}\) are equivalent, with
\(\phi^{-1}\le dP_{X,S}/dP_{X,T}\le\phi\)
\(P_{X,T}\)-almost everywhere, for some \(\phi\ge1\).
\end{itemize}
On the agreement region, {\rm(A1)}--{\rm(A2)} give
\(|m_\alpha|\gtrsim u|m_T|^\zeta\), and the calculation above
applies with the exponent \(q_\alpha=q_T/\zeta\).  Opposition can
create additional small mixture margins, and the proposition
controls that contribution in three ways: opposition is absent; a
positive target-margin buffer and a sufficiently small source
weight keep the mixture away from ties; or the
small-mixture-margin mass on the disagreement region is controlled
directly under the mixture design,
\begin{equation}
  P_{X,\alpha}\bigl(\Omega_-\cap\{|m_\alpha|<z\}\bigr)
  \le c_2z^{q_\alpha},
  \qquad z>0,
  \label{eq:tie_mass}
\end{equation}
where \(z\) is the margin threshold and \(c_2<\infty\) bounds the
mass of the low-signal part of \(\Omega_-\).

\begin{proposition}[Mixture noise exponent under source dominance]
\label{lem:one_acceleration}
Assume {\rm(T)}, {\rm(A1)}, {\rm(A2)}, fix
\(u=1-\alpha\in(0,1)\), and set \(q_\alpha:=q_T/\zeta\).
In case {\rm(ii)}, assume in addition that
\(|m_T|\ge c_T\) \(P_{X,T}\)-almost everywhere on \(\Omega_-\)
for some \(c_T>0\) {\rm(A3)}, and
\(u\le(1-\delta)\underline u_c\) for some \(\delta\in(0,1)\),
where \(\rho_-\) is the essential supremum of \(|m_S|/|m_T|\) on
\(\Omega_-\cap\mathcal D_T\) and
\(\underline u_c=(1+\phi\rho_-)^{-1}\); by {\rm(A2)},
\(r(x)\le\phi\rho_-\), so \(u_c(x)\ge\underline u_c\) on
\(\Omega_-\cap\mathcal D_T\) and the condition keeps every
location strictly below its threshold
\eqref{eq:pointwise_flip_threshold} (Supplement,
Lemma~\ref{lem:flip_extinction}).
In each of the three cases below,
\[
  P_{X,\alpha}(|m_\alpha|\le z)\le c_{q,\alpha}z^{q_\alpha},
  \qquad z>0,
\]
with
\begin{equation}
  c_{q,\alpha}\le
  \begin{cases}
    c_1u^{-q_\alpha},
    &\text{(i) no opposition, } P_{X,\alpha}(\Omega_-)=0;\\[2pt]
    c_1u^{-q_\alpha}
    +\bigl(\phi/(\delta c_T)\bigr)^{q_\alpha},
    &\text{(ii) buffered opposition;}\\[2pt]
    c_1u^{-q_\alpha}+c_2,
    &\text{(iii) controlled small mixture margins, \eqref{eq:tie_mass};}
  \end{cases}
  \label{eq:mixture_noise_constants}
\end{equation}
where \(c_1\) depends only on \(c_\zeta,\zeta,q_T,c_{q,T}\), and
\(\phi\).  In every case the pair
\begin{equation}
  \vartheta_\alpha=\frac{q_\alpha}{q_\alpha+1},
  \qquad
  V_\alpha=6\,c_{q,\alpha}^{1/(q_\alpha+1)}
  \label{eq:acceleration_variance_pair}
\end{equation}
is admissible in \eqref{eq:variance_bound_alpha}.
\end{proposition}

The three cases share the exponent \(q_\alpha\) and differ in how
opposition enters the constant: a smaller source weight or a
thinner buffer enlarges it, and it blows up as \(u\downarrow0\).
With the pair
\eqref{eq:acceleration_variance_pair}, the variance entry of the
localization radius \eqref{eq:r_condition} becomes
\begin{equation}
  \Bigl(\frac{V_\alpha\varrho}{N_{\rm eff}}\Bigr)^{
    \frac{q_\alpha+1}{q_\alpha+2}}
  =
  \Bigl(\frac{V_\alpha\varrho}{N_{\rm eff}}\Bigr)^{
    \frac{q_T+\zeta}{q_T+2\zeta}},
  \label{eq:noise_exponent_main}
\end{equation}
whose exponent exceeds the target-only exponent \((q_T+1)/(q_T+2)\)
exactly when \(q_T>0\) and \(0<\zeta<1\) (the difference is
computed at \eqref{eq:noise_exponent_gain} in the Supplement).

For \(q_\alpha>0\), condition \eqref{eq:tie_mass} excludes mixture
ties of positive mixture-design mass while permitting strict
flips: past its threshold a location's opposing mixture margin
grows again, so the mixture can be easy to learn where it is wrong
for the target.  What that accuracy is worth is the subject of the
next subsection.  The proof of
Proposition~\ref{lem:one_acceleration} and the without-dominance
calculation are in Supplement
Section~\ref{supp:mixture_exponents}.

\subsection{How Agreement and Cancellation Affect Target Guarantees}
\label{subsec:bernstein_transfer_oracle}
What a given mixture budget is worth on the target depends on the
local votes through \(Q_\alpha\).  Where the mixture retains a fixed fraction of the
target vote, one unit of mixture excess is worth a bounded number
of units of target excess, and the conversion is linear.  Where
the retained vote can be arbitrarily small, some target errors are
nearly free under the mixture, and their total target cost becomes
decisive.  We measure that cost with the target-margin measure
\(\nu_T\) of \eqref{eq:margin_measure}, under which the flip
profile is the mass \(H_{\rm flip}(u)=2\,\nu_T\{Q_\alpha\le0\}\) of
the tie-or-flip set \eqref{eq:profile_representation}.  The
conditions below control small positive \(Q_\alpha\) under
\(\nu_T\); the tie-mass condition \eqref{eq:tie_mass} of the
preceding subsection controls small \(|m_\alpha|\) under the
mixture design.  They are two distinct noise conditions and the
two sides of the target guarantee.  A source can leave every Bayes
label unchanged and still move both.  With equal designs and
\(m_S=\operatorname{sign}(m_T)\,|m_T|^{\zeta}\), the margin-power
path of Section~\ref{sec:numerical}, \(|m_\alpha|=\alpha|m_T|+u|m_T|^{\zeta}\)
and \(Q_\alpha=\alpha+u\,|m_T|^{\zeta-1}\) wherever \(m_T\ne0\), so
a smaller \(\zeta\) enlarges the mixture margin and raises
\(Q_\alpha\) at every location with \(0<|m_T|<1\) while
\(H_{\rm flip}(u)=0\) throughout.

\begin{corollary}[Target bounds in the three vote regimes]
\label{cor:regime_reinforcement}
\label{rem:no_flip_conversion}
\label{cor:regime_flip}
\label{lem:destab_exponent}
Fix \(u=1-\alpha\in(0,1)\) and consider three conditions:
(i) \emph{no opposition}, \(\nu_T(\Omega_-)=0\);
(ii) \emph{opposition below the flip threshold}, \(Q_\alpha\ge\delta\)
\(\nu_T\)-almost everywhere on \(\mathcal D_T\) for some
\(\delta\in(0,1]\); (iii) \emph{near cancellation, ties, and
flips}, under the small-vote mass condition
\begin{equation}
  \nu_T\{0<Q_\alpha<\varepsilon\}
  \le c_0\,\varepsilon^{\chi_0},
  \qquad \varepsilon>0,
  \label{eq:mixture_margin_profile}
\end{equation}
for some \(\chi_0>0\) and \(c_0<\infty\).  Under (i) and (ii),
\(H_{\rm flip}(u)=0\).  Each case assumes only its stated
condition, and for all \(t\ge0\),
\begin{equation}
  \Psi_\alpha(t)\le
  \begin{cases}
    M_\alpha\,t\le\alpha^{-1}t,
    &\text{(i) no opposition;}\\[2pt]
    M_\alpha\,t\le\delta^{-1}t,
    &\text{(ii) opposition below the flip threshold;}\\[2pt]
    H_{\rm flip}(u)+C\,t^{\chi},
    &\text{(iii) near cancellation, with }\chi:=\chi_0/(\chi_0+1),
  \end{cases}
  \label{eq:regime_bounds}
\end{equation}
where \(C\) depends only on \((c_0,\chi_0)\).
Every bound may in addition be capped at the saturation value
\(2J_T=2\nu_T(\mathcal X)\).  On the event of
Theorem~\ref{thm:weighted_crerm}, each applicable bound above
becomes a fitted-SVM target guarantee by taking
\(t=t_{\alpha,n,m}\) in \eqref{eq:target_excess_master}.
\end{corollary}

\begin{remark}[Where the H\"older power comes from]
\label{rem:holder_mechanism}
The power in case {\rm(iii)} comes from a tradeoff near
cancellation.  Fix a cutoff \(\varepsilon>0\) for the retained
vote.  Away from the locations with \(Q_\alpha<\varepsilon\), each
unit of mixture excess is amplified by at most \(1/\varepsilon\), so
these locations contribute at most \(t/\varepsilon\) to the target
excess.  The remaining stable locations contribute at most
\(2c_0\varepsilon^{\chi_0}\), by the target-margin mass condition
\eqref{eq:mixture_margin_profile}.  Hence
\(\Psi_\alpha(t)-H_{\rm flip}(u)\le t/\varepsilon+2c_0\varepsilon^{\chi_0}\),
and balancing the two costs at \(\varepsilon\asymp t^{1/(\chi_0+1)}\)
gives the power \(\chi_0/(\chi_0+1)\); ties and flips contribute
separately through \(H_{\rm flip}(u)\). For every \(\chi_0>0\),
the Supplement gives a model attaining this exponent and the
explicit constant from the sharper affine-frontier bound, up to
its saturation threshold.
\end{remark}

The three cases carry a mixture rate to the target in three ways.
If Theorem~\ref{thm:weighted_crerm} gives a mixture budget of
order \(N_{\rm eff}^{-\omega}\) at a fixed weight, for some
\(\omega>0\), then under (i) and (ii) the target bound has the
same order, multiplied by \(1/\alpha\) or \(1/\delta\), and under
(iii) it is \(H_{\rm flip}(u)+O(N_{\rm eff}^{-\omega\chi})\): the
exponent is
multiplied by \(\chi<1\), and the flip cost stays at every sample
size.  Case (iii) includes near cancellation without flips,
\(H_{\rm flip}(u)=0\) with the exponent still reduced, the global
form of the pre-flip weakening of Section~\ref{subsec:local_cost}.
Figure~\ref{fig:regime_bounds}(b) draws the three shapes at one
weight, with the budget entering through the abscissa.  The
proof, the explicit constant \(C\), the sufficient condition
\(u\le(1-\delta)\underline u_c\) for case {\rm(ii)}, and a model
in which the exponent of case {\rm(iii)} cannot be improved are in
the Supplement (Section~\ref{supp:mixing}).

Changing the weight moves the mixture margin, the effective sample
size, and the conversion at once.  For the opposing source of
Example~\ref{ex:constant_Q}, \(m_S=-m_T/2\) with a common uniform
design and \(m_T(x)=x/4\), \eqref{eq:margin_identity} gives
\(m_\alpha=(1-3u/2)\,m_T\) and \(Q_\alpha=1-3u/2\) at every
target-active location, so every location has the threshold
\(u_c=2/3\).  For \(u<2/3\), the mixture satisfies the target's
margin condition {\rm(T)} with the same exponent \(q_\alpha=1\)
and the constant \(4/|1-3u/2|\), which grows as \(u\) approaches
\(2/3\); the effective sample size is
\(N_{\rm eff}=\min\{n/(1-u),m/u\}\); and
\eqref{eq:constant_Q_modulus} gives
\(\Psi_\alpha(t)=\min\{t/(1-3u/2),\,1/4\}\).  At \(u=1/3\) the
multiplier is \(2\); at \(u=1/2\) it is \(4\); for \(u\ge2/3\)
every target-active location is tied or flipped and
\(\Psi_\alpha\equiv2J_T=1/4\): proportional pooling,
\(u=m/(n+m)\), maximizes \(N_{\rm eff}\) and, when \(m\ge2n\),
converts every mixture budget to the trivial target bound.
Between \(u=1/3\) and \(u=1/2\) no label changes, the multiplier
doubles, and with \(n=50\) and \(m=1000\) the effective sample
size rises from \(75\) to \(100\) while the margin constant rises
from \(8\) to \(16\).  The budget \(t_{\alpha,n,m}\) depends on
both, and on the approximation of the mixture decision by the
kernel class.  While both target bounds are below the saturation
value \(1/4\), the weight \(u=1/3\) gives the smaller one exactly
when \(2\,t_{2/3,n,m}<4\,t_{1/2,n,m}\), the budgets being those of
Theorem~\ref{thm:weighted_crerm} at \(\alpha=2/3\) and
\(\alpha=1/2\); the two budgets decide this comparison once they
are computed.

\begin{remark}[Reference cases]
\label{rem:benchmark_cases}
At the target-only endpoint, all source observation weights are zero.
Theorems~\ref{thm:affine_frontier} and~\ref{thm:weighted_crerm}
therefore include \(\alpha=1\). With \(J_T=\nu_T(\mathcal X)\), we have
\[
 Q_1\equiv1,\qquad N_{\rm eff}=n,\qquad
 \Psi_1(t)=\min\{t,2J_T\},\qquad
 \mathcal B_1(c)=2J_T(1-c)_+.
\]
These identities follow directly from the definitions; no limit of
the source-dominance constants is taken. If \(P_S=P_T\) then
\(Q_\alpha\equiv1\) for every \(\alpha\) and source data change
only the sample size: proportional pooling is then a pure
increase of \(N_{\rm eff}\), while for the opposing source above
the same rule flips every target-active location once \(m>2n\).  Under pure covariate shift \(m_S=m_T\)
the sign never flips but the multiplier still depends on the
density shift.  Even the fixed-slope frontier
\(\alpha\mapsto\mathcal B_\alpha(c)\), with \(c>0\) fixed, need
not be convex (Supplement, Remark~\ref{rem:alpha_nonconvex}, Section~\ref{supp:frontier}).
\end{remark}

\subsection{Explicit Target Guarantees with Gaussian Kernels}
\label{sec:gaussian}
\label{subsec:gaussian_est}
\label{subsec:target_excess}

Theorem~\ref{thm:weighted_crerm} leaves the approximation and
complexity of the RKHS to be specified; for Gaussian kernels both
become explicit, and the conversion then determines their target
value.  For the Gaussian kernel
\(k_\gamma(x,x'):=\exp(-\|x-x'\|_2^2/\gamma^2)\), \(\gamma>0\), with
RKHS \(H_\gamma\), both terms are explicit through established
bounds \citep{steinwart2007fast,steinwart2008support}; on the
radial model of Section~\ref{sec:adaptive}, Gaussian averaging also
preserves the direction of a linear signal, and the same
kernel family is used in the experiments.  Let
\(f_{D,\lambda,\gamma}\) be the weighted SVM
\eqref{eq:transfer_svm_empirical} over \(H_\gamma\) at a fixed
weight \(\alpha\in(0,1)\), and \(\bar f_{D,\lambda,\gamma}\) its
clipping.  Three properties of \(P_\alpha\) enter: the mixture
Tsybakov condition of Proposition~\ref{lem:one_acceleration}, with
parameters \((q_\alpha,c_{q,\alpha})\) and \(q_\alpha=0\) always
admissible; a margin-noise condition
\(\int_{\{\Delta_\alpha\le z\}}|m_\alpha|\,dP_{X,\alpha}
\le c_{\beta,\alpha}z^{\beta_\alpha}\), \(\beta_\alpha>0\), where
\(\Delta_\alpha(x)\) is the distance from \(x\) to the opposite
mixture decision class; and a marginal tail condition
\(P_{X,\alpha}\{\|X\|_2>R\}\le c_{\tau,\alpha}R^{-\tau_\alpha}\)
for \(R>0\) and \(0<\tau_\alpha<\infty\). The convention
\(\tau_\alpha=\infty\) means
\(P_{X,\alpha}\{\|X\|_2>R_\alpha\}=0\) for a fixed
\(R_\alpha<\infty\) (\(\tau_\alpha\) is unrelated to the tilt
\(\tau\) of Section~\ref{sec:adaptive}).  Without loss of generality, take \(c_{\beta,\alpha}\ge1\) and
\(c_{\tau,\alpha}\ge1\), by enlarging these upper-bound constants.
The Gaussian entropy
numbers enter through two free parameters \(\varepsilon>0\) and
\(d/(d+\tau_\alpha)<p<1\).  With \(N_{\rm eff}=\min\{n/\alpha,m/u\}\)
and the variance pair \((\vartheta_\alpha,V_\alpha)\) of
\eqref{eq:acceleration_variance_pair}, which the mixture Tsybakov
condition alone provides (Supplement,
Lemma~\ref{lem:variance-bound-hinge}, Section~\ref{supp:faster}), define the approximation and
estimation blocks
\begin{align}
  \mathcal A_\alpha(\lambda,\gamma)
  &:=K\Bigl[\lambda^{\frac{\tau_\alpha}{\tau_\alpha+d}}
  \gamma^{-\frac{d\tau_\alpha}{\tau_\alpha+d}}
  +c_{\beta,\alpha}\gamma^{\beta_\alpha}\Bigr],
  \label{eq:gaussian_approximation_block}\\
  \mathcal S_\alpha(\lambda,\gamma)
  &:=K\Bigl[\varrho\bigl(N_{\rm eff}
  \lambda^p\gamma^{(1-p)(1+\varepsilon)d}\bigr)^{
    -\frac{q_\alpha+1}{q_\alpha+2-p}}
  +\Bigl(\frac{\varrho V_\alpha}{N_{\rm eff}}\Bigr)^{
    \frac{q_\alpha+1}{q_\alpha+2}}
  +\frac{\varrho}{N_{\rm eff}}\Bigr],
  \label{eq:gaussian_estimation_block}
\end{align}
where \(K\) collects fixed constants and, for
\(\tau_\alpha=\infty\), the first term of \(\mathcal A_\alpha\) is
read as \(\lambda\gamma^{-d}\).  The approximation block bounds how
closely a regularized Gaussian-RKHS function can follow the
mixture decision (Supplement, Lemma~\ref{lem:gaussian_mixture_approximation},
Section~\ref{supp:gauss_detail}),
the estimation block bounds the additional error from fitting with
finite weighted samples, and a narrower kernel reduces the boundary-smoothing term while
increasing the regularization and estimation terms at fixed
\(\lambda\), so the bound makes the tuning tradeoff explicit.

\begin{corollary}[Gaussian target-risk bound]
\label{thm:main_target_bound}
Let \(\varepsilon>0\), \(d/(d+\tau_\alpha)<p<1\),
\(\lambda,\gamma\in(0,1]\), and \(\varrho\ge1\).  Under the three
mixture conditions above there is a constant \(K\), depending only
on the fixed weight $\alpha$, the displayed model and entropy
parameters, and a support radius $R_\alpha$ when
$\tau_\alpha=\infty$, such that, with probability at least
\(1-3e^{-\varrho}\),
\begin{equation}
  \mathcal E_T(\bar f_{D,\lambda,\gamma})
  \;\le\;
  \Psi_\alpha\bigl(
  \mathcal A_\alpha(\lambda,\gamma)
  +\mathcal S_\alpha(\lambda,\gamma)\bigr)
  =
  \inf_{c\ge0}
  \bigl\{c\,[\mathcal A_\alpha+\mathcal S_\alpha]
  +\mathcal B_\alpha(c)\bigr\},
  \label{eq:main_target_bound}
\end{equation}
and the same bound holds for the target classification excess of
\(\operatorname{sign}\bar f_{D,\lambda,\gamma}\).
\end{corollary}

The proof (Supplement, Section~\ref{pf:main_target_bound}) applies the
Gaussian mixture oracle inequality,
Theorem~\ref{thm:gaussian_transfer_estimation} (Section~\ref{supp:gauss}), and then
Theorem~\ref{thm:affine_frontier} pathwise as in
\eqref{eq:target_excess_master}; the scalar and assumption-free
evaluations of the modulus are those of Remark~\ref{rem:affine_sharpness}
and Remark~\ref{rem:finite_conversion}.  Because \(\Psi_\alpha\) is
nondecreasing, \((\lambda,\gamma)\) may be balanced inside its
argument.  On bounded support the balanced mixture bound has order
\(N_{\rm eff}^{-(r_G-\xi)}\) at fixed confidence, where
\[
  r_G=\frac{\beta_\alpha(q_\alpha+1)}
  {\beta_\alpha(q_\alpha+2)+d(q_\alpha+1)}
\]
and \(\xi>0\) is an arbitrarily small slack. The target-only
calculation at \(\alpha=1\), with \(N_{\rm eff}=n\), is given in
Supplement Remark~\ref{rem:gaussian_target_only} (Section~\ref{supp:gauss}).
The full tail-dependent
rates and the balancing choices of \((\lambda,\gamma)\) are in the
Supplement (Corollary~\ref{cor:transfer_rates}, Section~\ref{pf:main_target_bound}).
Their target consequences follow from Corollary~\ref{cor:regime_flip}.
At a fixed pair of laws and weight, kernel tuning changes the
error budget supplied to \(\Psi_\alpha\) and leaves the conversion
curve itself unchanged.  Two weights are compared through
\(\Psi_\alpha(\mathcal A_\alpha+\mathcal S_\alpha)\) at each: the
weight enters the budget through \(N_{\rm eff}\), the mixture
exponents and constants, and the approximation of the mixture
decision, and it enters the curve through \(Q_\alpha\).

\begin{remark}[A nonparametric minimax benchmark under reinforcement]
\label{rem:reinforcement_benchmark}
On the fixed nonparametric reinforcement class
\(\mathcal P_\Theta(\nu,q_T,\zeta)\) of Supplement
Section~\ref{supp:reinforcement}, under the sufficient lower-bound
conditions of Theorem~\ref{thm:gaussian_reinforcement_minimax}
(Section~\ref{supp:reinforcement_statement}), Gaussian SVMs with deterministic tuning for known \((\nu,\zeta)\)
attain the all-procedure minimax classification rate
\[
  \bigl[n+m^{(2\nu+d)/(2\zeta\nu+d)}\bigr]^{-\nu(q_T+1)/(2\nu+d)},
  \qquad \nu q_T\le1 .
\]
Here the scalar \(\nu\) is the H\"older smoothness of the target
margin (unrelated to the measures \(\nu_D\)),
\(q_T\) its Tsybakov exponent, and \(\zeta\in(0,1)\) the
reinforcement power of a source with \(m_Tm_S\ge0\) and
\(|m_S|\asymp|m_T|^\zeta\); the stronger source signal allows a
finer spatial resolution, which produces the effective source
contribution \(m^{(2\nu+d)/(2\zeta\nu+d)}\).  In the notation of
\citet{cai2021transfer}, \((\nu,q_T,\zeta,n,m)\) correspond to
\((\beta,\alpha,\gamma,n_Q,n_P)\); likewise, the exponent
\(\beta_\alpha\) above is the margin-noise exponent \(\beta\) of
\citet{steinwart2008support}, which corresponds to \(d\) times the geometric noise
exponent \(\alpha\) of \citet{steinwart2007fast}, and the
confidence parameter \(\varrho\) of
Section~\ref{sec:general_oracle} that of the parameter \(\tau\) of
\citet{steinwart2008support}; these are renamed here because
\(\alpha\) is the target weight and \(\tau\) the tilt of
Section~\ref{sec:adaptive}.
\citet[Theorems~1--2]{cai2021transfer} established this rate form
for weighted nearest-neighbor classification under posterior drift.
Theorem~\ref{thm:gaussian_reinforcement_minimax} proves Gaussian-SVM
attainment and the matching lower bound on the stated fixed class.
This benchmark uses a fixed regularization \(\lambda_0\ge1/2\), a different class
and regime from the adaptive analysis of Section~\ref{sec:adaptive}.
\end{remark}

We can now assess a proposed weight through the accuracy
attainable on its mixture and the target meaning of that accuracy.
Choosing the weight from data asks a new question: how much of the
source--target relation must be learned to make this assessment
useful?

\section{Data-Dependent Transfer: Adaptation and Limits}
\label{sec:adaptive}

Additional source observations do not resolve every uncertainty
relevant to borrowing.  We use a tilted-boundary class, in which
the source boundary is tilted away from the target's while keeping
some strength along it, to separate
source sampling uncertainty from uncertainty about the
source--target relation, identify the resulting minimax risk, and
then construct a selector that approaches it, while the output
remains one fitted Gaussian SVM.  The selector proposes a
source-assisted fit and checks its target decisions.  The analysis
controls classification error directly through the finite-slope
frontier of Theorem~\ref{thm:affine_frontier}.

\subsection{Target-Safe Minimax Adaptation}
\label{subsec:adaptive_minimax}

We use a common radial design to isolate two features of source
quality: signal strength along the target direction and
disagreement caused by a different direction.  The target decision
boundary is a hyperplane with unknown normal, and multiplying the
linear signal by a positive radial profile changes its strength
across the domain while preserving that boundary.  Fix a dimension
\(d\ge2\).  Let
\(\mathcal X=B^d=\{x\in\mathbb R^d:\|x\|_2\le1\}\), let both
covariate designs equal \(\mu_0:=\operatorname{Unif}(B^d)\), and
let the class \(\mathcal P_1\) consist of the ordered law
pairs whose posterior margins are
\begin{equation}
  m_D(x)=a\,\varphi(x)\,w_D^\top x,
  \qquad D\in\{T,S\},
  \qquad w_T\in\mathbb S^{d-1},
  \label{eq:sgn_class}
\end{equation}
where \(w_T\) is an unknown unit target normal, \(w_S\in\mathbb R^d\)
is the source normal,
\(\varphi(x)=\bar\varphi(\|x\|_2)\) is a positive radial profile
whose normalization \((d+2)\int\varphi(x)x_1^2\,d\mu_0(x)=1\)
fixes the scale of the amplitude \(a\) (this radial profile
\(\varphi\) is unrelated to the complexity function \(\varphi\) of
Section~\ref{sec:general_oracle}).  The source normal
decomposes along the target normal: its \emph{aligned strength}
\(s:=\langle w_S,w_T\rangle\) scales the target direction, and its
residual \(w_S-sw_T\), orthogonal to \(w_T\), tilts it; the
length \(\tau:=\|w_S-sw_T\|_2\) of the residual is the
\emph{tilt} (the scalar \(s\) is unrelated to the local votes
\(s_T,s_S\) of Section~\ref{sec:problem_setup}).\footnote{The
restriction \(d\ge2\) is needed for the lower bound of
Theorem~\ref{thm:target_safe_adaptive}: for \(d=1\) the
normal \(w_T\) takes two values, the residual is zero,
and the sign of \(n^{-1}\sum_iY_iX_i\) recovers \(w_T\) with
exponentially small error probability, so no \(n^{-1}\)-order
lower bound holds.}  The class parameters are fixed before the
sample sizes:
\begin{equation}
\begin{gathered}
  a\in[a_-,a_+],
  \qquad
  s\in[-s_H,s_U],
  \qquad
  \tau\le\tau_U,
  \qquad
  0<\varphi_-\le\varphi\le\varphi_+,
  \quad \varphi_-\le1\le\varphi_+,\\
  a_+\varphi_+\max\{1,s_{\max}+\tau_U\}
  \le\min\bigl\{1-\varepsilon_Y,\;2\lambda_0(1-\varepsilon_H)\bigr\},
\end{gathered}
  \label{eq:sgn_envelope}
\end{equation}
with \(0<a_-\le a_+<\infty\), \(s_H\ge0\), \(0<s_U<\infty\),
\(\tau_U\ge0\), \(s_{\max}=\max\{s_H,s_U\}\), positive
\(\varepsilon_Y,\varepsilon_H\), and one fixed Gaussian bandwidth
and regularization \((\gamma_0,\lambda_0)\), \(0<\lambda_0<1/2\),
valid throughout the class.  The lower amplitude bound \(a_->0\)
is a uniform signal floor used throughout the risk analysis; the
first envelope condition keeps the label probabilities away from
\(0\) and \(1\); the second places every population soft-margin
minimizer in the all-active regime; and the explicit constants are
in the Supplement (Section~\ref{supp:adaptive}).  These class
constants are fixed and available to the procedure below (its
truncation level and grid use \(s_U\), and its thresholds use
\(a_-\), \(\varphi_-\), and \(\lambda_0\)); what remains unknown,
and what the selector adapts to, are the two normals \(w_T\) and
\(w_S\), that is, the target direction, the aligned strength \(s\),
and the tilt, together with the amplitude \(a\).

Borrowing changes both the strength and the direction of the
training signal.  On this class
the vote identity \eqref{eq:margin_identity} reads on the normals:
with \(u=1-\alpha\),
\begin{equation}
\begin{gathered}
  w_\alpha:=\alpha w_T+uw_S=\kappa_\alpha w_T+u(w_S-sw_T),
  \qquad
  \kappa_\alpha:=\langle w_\alpha,w_T\rangle=\alpha+us,\\
  b_\alpha:=\|w_\alpha\|_2
  =\sqrt{\kappa_\alpha^2+u^2\tau^2},
  \qquad
  m_\alpha(x)=a\,\varphi(x)\,w_\alpha^\top x.
\end{gathered}
  \label{eq:mixture_vote_vector}
\end{equation}
The component \(\kappa_\alpha w_T\) of the mixture normal along the
target normal retains the target direction, with \(\kappa_\alpha\)
its retained strength, whereas \(u(w_S-sw_T)\) tilts the mixture
boundary, and \(b_\alpha\) is the total strength of the mixture
signal (the model has no intercept).  Positive aligned strength in
parameter space does not imply pointwise agreement in covariate
space: for \(0<\alpha<1\),
a single law with \(s>0\) and \(\tau>0\) carries positive design
measure of all three vote regimes of \eqref{eq:vote_cases}---agreement,
opposition without a crossing, and strict mixture flips.  We next
quantify the target cost of the tilt, including the case where a
single worst-case conversion multiplier is infinite.  Write
\(J_T=\nu_T(\mathcal X)=\mathbb E_{\mu_0}|m_T(X)|\) as in
\eqref{eq:margin_measure} and, when \(b_\alpha>0\),
\(\rho_\alpha:=\kappa_\alpha/b_\alpha\in[-1,1]\), the cosine of
the angle between the target normal \(w_T\) and the mixture normal
\(w_\alpha\).

\begin{proposition}[Radial affine frontier]
\label{prop:radial_frontier}
Let \(P\in\mathcal P_1\) and \(\alpha\in(0,1]\) with
\(b_\alpha>0\).  Then, for every \(c\ge0\),
\begin{equation}
  \mathcal B_\alpha(c)
  =J_T\Bigl[\,1-c\,b_\alpha
  +\sqrt{1-2c\,\kappa_\alpha+c^2b_\alpha^2}\;\Bigr],
  \qquad
  H_{\rm flip}(u)=J_T(1-\rho_\alpha),
  \label{eq:radial_frontier}
\end{equation}
and the target classification excess of the mixture Bayes rule
\(h^*_\alpha=\operatorname{sign}(w_\alpha^\top x)\) is
\(\mathcal E^{01}_T(h^*_\alpha)=\tfrac12J_T(1-\rho_\alpha)\).
On the target-oriented branch \(\kappa_\alpha>0\),
\begin{equation}
  \mathcal B_\alpha(2/b_\alpha)
  =J_T\bigl\{\sqrt{5-4\rho_\alpha}-1\bigr\}
  \le2J_T(1-\rho_\alpha^2)
  =2J_T\,\frac{u^2\tau^2}{b_\alpha^2}.
  \label{eq:radial_frontier_slope}
\end{equation}
\end{proposition}

The tilted-boundary class makes the distinction in
Remark~\ref{rem:affine_sharpness} concrete.  Every nonzero tilt
gives \(M_\alpha=\infty\): when \(u\tau>0\), \(|\rho_\alpha|<1\)
and \eqref{eq:radial_frontier} shows that \(\mathcal B_\alpha(c)\)
exceeds \(H_{\rm flip}(u)\) at every finite \(c\), so no finite
slope attains the smallest intercept and no bound of the form
\eqref{eq:global_vote_conversion} exists.  Yet on the
target-oriented branch a finite-slope comparison controls the
integrated target cost.  For a hard classifier \(h\), hinge excess
is twice classification excess, and hence
\(\mathcal E^{01}_T(h)\le c\,\mathcal E^{01}_\alpha(h)
+\tfrac12\mathcal B_\alpha(c)\).  Choosing \(c=2/b_\alpha\) bounds the additive term by
\(J_T\sin^2\theta_\alpha\), where \(\theta_\alpha\) is the angle
between the target and mixture normals.
Corollary~\ref{cor:fixed_weight_signed} uses exactly this term.
The frontier lets us account for boundary tilt, which a
bounded-ratio assumption would exclude.  The proof in the
Supplement (Section~\ref{supp:radial_frontier}) uses the radial
symmetry of the design and does not involve the kernel.

Data-driven weighting needs source quality expressed through
quantities determined by the two distributions.  The signed kernel
mean \(\mu_D^{\gamma_0}:=\mathbb E_D[Y\,k_{\gamma_0}(X,\cdot)]\in H_{\gamma_0}\),
\(D\in\{T,S\}\), averages the labeled feature representation; its
source projection onto the target mean measures aligned strength,
and the normalized residual measures tilt.  Define
\begin{equation}
  s(P)
  =\frac{\langle\mu_S^{\gamma_0},\mu_T^{\gamma_0}\rangle_{H_{\gamma_0}}}
        {\|\mu_T^{\gamma_0}\|_{H_{\gamma_0}}^2},
  \qquad
  \tau(P)
  =\frac{\|\mu_S^{\gamma_0}-s(P)\,\mu_T^{\gamma_0}\|_{H_{\gamma_0}}}
        {\|\mu_T^{\gamma_0}\|_{H_{\gamma_0}}},
  \qquad
  \load(P):=[s(P)]_+.
  \label{eq:quality_functionals}
\end{equation}
Radial Gaussian averaging maps each normal to its kernel mean,
so the projection and the residual computed in feature space agree
with those of \(w_S\) on \(w_T\): \(s(P)=s\) and \(\tau(P)=\tau\)
(Supplement, Section~\ref{supp:adaptive_calibration}).  We call \(\load\) the
\emph{transferable strength}; the subscript records that only the
positive part of the aligned strength is transferable.  Neither
\(\load\) nor \(\tau\) is a margin, noise, or smoothness exponent.
We compare procedures over pairs with the same transferable
strength \(\load\) and tilt at most \(\tau\), which keeps the useful
source strength fixed while allowing uncertainty in the
source-specific component:
\begin{equation}
  \mathcal P_1(\load,\tau)
  =\bigl\{P\in\mathcal P_1:
  \load(P)=\load,\ \tau(P)\le\tau\bigr\}.
  \label{eq:quality_ball}
\end{equation}
For \(\load>0\) the class \(\mathcal P_1(\load,\tau)\), a fiber of
\(\load(\cdot)\), fixes \(s=\load\).  The class with \(\load=0\)
includes all nonpositive aligned strengths, so its minimax risk is taken
over that combined class (the exact-shell refinement is in the
Supplement).

The classification problem does not require the soft score to
reach the hinge optimum.  Let \(f^P_\alpha\) minimize the
regularized population mixture hinge risk at the fixed
\((\gamma_0,\lambda_0)\).  Radial Gaussian averaging preserves the
boundary normal (the signed kernel mean of
\(a\varphi(x)\,w^\top x\) is \(a\,\psi_\varphi(\|x\|_2)\,w^\top x\)
with \(\psi_\varphi>0\), Lemma~\ref{lem:supp_radial} in
Section~\ref{supp:adaptive_calibration} of the Supplement), so
\(\operatorname{sign}f^P_\alpha\) is the mixture Bayes rule.  At the
same time, the all-active envelope gives
\(\|f^P_\alpha\|_\infty\le1-\varepsilon_H\), and hence, by
\eqref{eq:hinge_pointwise_excess},
\(\mathcal E_\alpha(f^P_\alpha)\ge\varepsilon_H\,
\mathbb E_{P_{X,\alpha}}|m_\alpha(X)|\), even though its mixture classification excess is zero.  We
therefore analyze how empirical fitting perturbs the sign and use
the frontier to evaluate the target cost of the resulting
decisions.  The next bound separates this sampling error from the
cost of the tilted population boundary, relative to the
target-oriented signal retained by the mixture.  Write
\(\mathcal E^{01}_T(\,\cdot\,)\) for the target classification
excess risk.

\begin{corollary}[Fixed-weight risk on the tilted-boundary class]
\label{cor:fixed_weight_signed}
Let \(P\in\mathcal P_1\) have amplitude \(a\), aligned strength
\(s\), and tilt \(\tau\).  Fix a weight
\(\alpha\in(0,1]\) measurable with respect to a sigma-field
\(\mathcal G_{\rm pre}\) independent of the fit blocks
\((T_I,S_I)\) of sizes \((n_f,m_f)\), let
\(\widehat f_\alpha\) be the \(\alpha\)-weighted no-intercept
Gaussian hinge-SVM fitted on them, and write \(u=1-\alpha\) for
the source weight and \(q=u/\alpha\) for the source-to-target odds.  There is a constant \(C\), depending only on the
fixed class constants, such that the following hold.
\begin{itemize}
\item[(i)] On the target-oriented branch
\(\kappa_\alpha=\alpha+us>0\),
\begin{equation}
  \mathbb E\bigl[\mathcal E^{01}_T
  \{\operatorname{sign}\widehat f_\alpha\}
  \mid\mathcal G_{\rm pre}\bigr]
  \;\le\;
  \frac Ca\,
  \frac{\alpha^2/n_f+u^2/m_f+a^2u^2\tau^2}{\kappa_\alpha^2}.
  \label{eq:fixed_weight_target}
\end{equation}
\item[(ii)] \emph{(Class-informed weight.)}  For \(\load>0\) and
\(\tau\ge0\), the deterministic odds
\(q^\star(\load,\tau):=m_f\load/\{n_f(1+m_f\tau^2)\}\), which use neither
\(a\) nor the sample, satisfy
\[
  \sup_{P\in\mathcal P_1(\load,\tau)}
  \mathbb E_P\,\mathcal E^{01}_T
  \{\operatorname{sign}\widehat f_{\alpha^\star}\}
  \le C\Bigl[n_f+\frac{m_f\load^2}{1+m_f\tau^2}\Bigr]^{-1},
  \qquad \alpha^\star=(1+q^\star)^{-1},
\]
with no logarithmic factor.
\end{itemize}
\end{corollary}

The same joint assessment now appears directly in the
fitted-classifier bound: the retained target-oriented signal
normalizes both sampling error and the target cost of boundary
tilt.  The terms \(\alpha^2/n_f\) and \(u^2/m_f\) in
\eqref{eq:fixed_weight_target} are the sampling contributions of
the two fitted blocks; after normalization by \(a\kappa_\alpha^2\),
the tilt term \(a^2u^2\tau^2\) bounds the finite-slope frontier
contribution in \eqref{eq:radial_frontier_slope}, using
\(J_T\asymp a\) and \(b_\alpha^2\ge\kappa_\alpha^2\); a small
retained target component \(\kappa_\alpha\) makes both more costly
(proof in the Supplement, Section~\ref{supp:adaptive}).  Part (ii)
shows that a weight informed by the class parameters \((\load,\tau)\)
attains, with fixed split fractions and \(a\in[a_-,a_+]\), the
risk scale \([n+I_S]^{-1}\) (the minimization of
\eqref{eq:fixed_weight_target} over the odds \(q\) is in the
Supplement proof of Corollary~\ref{cor:fixed_weight_signed}), where
\begin{equation}
  I_S:=\frac{m\load^2}{1+m\tau^2}=\frac{\load^2}{m^{-1}+\tau^2}
  \label{eq:effective_information}
\end{equation}
is the effective source information.

\begin{remark}[Why source information saturates]
\label{rem:saturation}
Knowing the source law exactly need not determine the target
decision.  When tilt is allowed within
\(\mathcal P_1(\load,\tau)\), distinct target directions can remain
compatible with the same source law; Example~\ref{ex:supp_witness}
in Section~\ref{supp:adaptive_positive} of the Supplement gives an
explicit pair.  Source observations alone cannot distinguish them,
even with an arbitrarily large source sample.  The expression
\(I_S=\load^2/(m^{-1}+\tau^2)\) in \eqref{eq:effective_information}
separates the two uncertainties: \(m^{-1}\) decreases as the
source sample grows, whereas \(\tau^2\) records the ambiguity
allowed by the source--target relation.  For fixed positive
tilt, the source contribution approaches \(\load^2/\tau^2\).  The
all-procedure lower bound in
Theorem~\ref{thm:target_safe_adaptive}(ii) establishes this limit
for every learning procedure, the weighted SVM included.
\end{remark}

Adaptation aims to attain the risk scale allowed by these
uncertainties, up to the factor in
Theorem~\ref{thm:target_safe_adaptive}.  The class-informed weight
\(q^\star\) uses both strength and tilt; our selector does not
reproduce it.  An estimated strength proposes one source-assisted
fit, and independent target data check the candidate's decisions.

Calibration blocks \(T_C,S_C\) estimate the aligned strength of
the source through the empirical kernel means, producing the
positivity certificate, the event \(\Omega_{\rm cert}\) that
the estimated strength clears its confidence margin, and a clipped,
discretized estimate \(\widetilde s>0\), and propose the odds
\begin{equation}
  q_I=\frac{m_I\widetilde s}{n_I},
  \qquad
  \widehat\alpha_I=(1+q_I)^{-1}\in(0,1).
  \label{eq:qI_def}
\end{equation}
Independent fit blocks \(T_I,S_I\) produce the
\(\widehat\alpha_I\)-weighted no-intercept Gaussian hinge-SVM
\(\widehat f_I\), so that Corollary~\ref{cor:fixed_weight_signed}
applies with \((n_f,m_f)=(n_I,m_I)\); a further target block \(T_0\) produces the
target-only fallback \(\widehat f_0\) and takes no part in the
acceptance decision; and a pilot target block \(T_P\) of size
\(n_P\) estimates the target normal by the moment
\(\widehat w_T\propto n_P^{-1}\sum_{i\in T_P}Y_iX_i\) and checks
whether to retain the source-assisted fit.  The blocks are formed by
a fixed deterministic rule with fixed positive fractions, so the
proposed weight is fixed relative to the fit blocks in
Corollary~\ref{cor:fixed_weight_signed}.

The target check protects decisions that target data can already
resolve.  Outside a band of width \(b_N\asymp\sqrt{\ell_N/n_P}\)
around the pilot boundary, it requires the candidate to agree in
sign with \(\widehat w_T^\top x\) on a finite net, with a positive
score buffer; the buffer together with the SVM's Lipschitz bound
extends the check to nearby points.  Inside the band, the check
does not force the candidate to copy the pilot; source information
can refine those unresolved decisions.

On the pilot confidence event, accepted errors are confined to the
band.  Its probability mass and its target margin are both
proportional to its width, so its total excess-risk cost is
quadratic in the width, of order \(b_N^2\asymp\ell_N/n\); this is
the cost of the errors the check permits.  Write
\(\Omega_{\rm gate}\) for the event that the target check accepts
and \(A:=\Omega_{\rm cert}\cap\Omega_{\rm gate}\) for acceptance.
The selector is
\begin{equation}
  \widehat h
  =
  \begin{cases}
    \operatorname{sign}(\widehat f_I),
    &\text{on \(A\)},\\[1mm]
    \operatorname{sign}(\widehat f_0),
    &\text{otherwise}.
  \end{cases}
  \label{eq:selector_def}
\end{equation}
The output is one of the two stored fitted SVMs, with no refit,
label reversal, or negative weight; the complete procedure, with
meshes, radii, and the exact acceptance rule, is in the Supplement
(Section~\ref{supp:adaptive_positive}).  Let
\[
  \mathcal T
  :=\bigl\{P_T:\ m_T(x)=a\,\varphi(x)\,w_T^\top x
  \text{ under the fixed target envelope}\bigr\}
\]
be the target marginal class.

\begin{theorem}[Target-safe adaptive minimaxity]
\label{thm:target_safe_adaptive}
Let \(N=n+m\) and \(\ell_N=\log(eN)\), and let \(\widehat h\)
be the split-sample selector \eqref{eq:selector_def}.  Assume
\(n,m\ge n_{\min}\) and \(\ell_N\le c_{\rm imb}\,n\), where
\(n_{\min}\) and \(c_{\rm imb}\) depend only on the fixed dimension
\(d\), the fixed class and check constants. There exist constants
\(C,c>0\), depending only on these fixed constants, the fixed Gaussian
kernel and regularization, and the fixed split proportions, such that the
following hold.
\begin{itemize}
\item[(i)] \emph{(Arbitrary-source target safety.)}  For every
\(P_T\in\mathcal T\) and every independent source
law \(P_S\) on \(B^d\times\{-1,1\}\), with arbitrary source
design and regression,
\[
  \mathbb E\,\mathcal E^{01}_T(\widehat h)
  \le C\,\frac{\ell_N}{n}:
\]
the selector retains the target-only rate up to the factor
\(\ell_N\).
\item[(ii)] \emph{(Minimax rate and adaptive attainment.)}  For
every transferable strength \(0\le\load\le s_U\) and tilt \(0\le\tau\le\tau_U\),
\[
  c\,\Bigl[n+\tfrac{m\load^2}{1+m\tau^2}\Bigr]^{-1}
  \le
  \inf_{\widetilde h}
  \sup_{P\in\mathcal P_1(\load,\tau)}
  \mathbb E_P\,\mathcal E^{01}_T(\widetilde h)
  \le
  C\,\Bigl[n+\tfrac{m\load^2}{1+m\tau^2}\Bigr]^{-1},
\]
the infimum running over all randomized measurable classifiers
using the full domain-labeled sample, and the selector, which
uses neither the strength nor the tilt, satisfies
\[
  \sup_{P\in\mathcal P_1(\load,\tau)}
  \mathbb E_P\,\mathcal E^{01}_T(\widehat h)
  \le
  C\ell_N
  \Bigl[n+\tfrac{m\load^2}{1+m\tau^2}\Bigr]^{-1}.
\]
\end{itemize}
Moreover, the logarithmic factor drops on the clean fibers
\(\tau=0\), and the selector attains the target-only order
\(n^{-1}\) on the coarse fiber \(\load=0\):
\[
  \sup_{P\in\mathcal P_1(\load,0)}
  \mathbb E_P\,\mathcal E^{01}_T(\widehat h)
  \le\frac{C}{n+m\load^2},
  \qquad
  \sup_{P\in\mathcal P_1(0,\tau)}
  \mathbb E_P\,\mathcal E^{01}_T(\widehat h)
  \le\frac{C}{n}.
\]
\end{theorem}

Theorem~\ref{thm:target_safe_adaptive} gives both the limit and a
way to approach it.  A safety guarantee is not an error floor.  The target check bounds
what an accepted fit may get wrong; the informative-source
analysis gives a sharper risk bound for the same fit.  We use that
sharper bound and do not add the safety allowance to it.  Recall the acceptance event
\(A=\Omega_{\rm cert}\cap\Omega_{\rm gate}\), and write
\(\Delta_I=\mathcal E^{01}_T(\operatorname{sign}\widehat f_I)\)
and \(\Delta_0=\mathcal E^{01}_T(\operatorname{sign}\widehat f_0)\)
for the excess risks of the candidate and of the fallback.  Since
the fallback is fitted on a block independent of the blocks that
determine \(A\), for every fixed law pair
\[
  \mathbb E\,\mathcal E^{01}_T(\widehat h)
  =\mathbb E[\mathbf 1_A\Delta_I]+P(A^c)\,\mathbb E\Delta_0 ,
\]
where \(\mathbf 1_A\) is the indicator of acceptance.  The
fallback is paid only when used: its target-only risk of order
\(n^{-1}\) is multiplied by the rejection probability, and on
informative classes the proof keeps the candidate's sharper risk
on \(A\) and controls this product at the source-assisted scale
of part (ii).  Thus target data can protect the fit without
dictating its final accuracy.  Proofs, the exact construction, and
the exact-shell refinement of the lower bound are in the
Supplement (Section~\ref{supp:adaptive}).

\begin{remark}[Adaptation to source quality]
\label{rem:adapt_scope}
\label{rem:no_weight_recovery}
The selector adapts to the source's aligned strength and tilt, the two
quantities that determine the quality class, at fixed class-valid
bandwidth and regularization; simultaneous data-driven selection
of \((\alpha,\gamma,\lambda)\) and adaptation to unknown smoothness
and margin parameters lie outside the theorem.  The guarantee
concerns the risk of the selected fit.  The target check is a finite-net
certificate over a covering of \(B^d\) at resolution \(N^{-3}\),
\(N=n+m\); it provides an explicit construction establishing
statistical attainability, and the present analysis does not
establish a computationally efficient implementation of this
certificate.  The numerical studies of Section~\ref{sec:numerical}
use a practical target-validation rule.
\end{remark}

\subsection{A Minimax Limit of Nonnegative Reweighting}
\label{subsec:representation_limit}

The preceding construction shows how target evidence can guide a
useful choice.  Is better selection the whole story?  So far, reweighting has been
the given way of using the source.  We now ask whether choosing
its weights sufficiently well can always realize the benefit
available from an informative source.

To separate a limitation of selection from a limitation of
fitting, we allow the weights and tuning parameters to be
arbitrary functions of the full sample and consider a class with a
known opposition relation, on which the source strongly favors the
opposite of the target label and reversing the source labels gives
an accurate target rule.  Here source labels become useful through
reversal, so the distinction between reweighting the observed
labels and correcting their relation to the target becomes
essential.

Let both covariate designs equal
\(\mu_1:=\operatorname{Unif}([-1,1]^d)\).  For
\(\sigma\in\{-1,+1\}\) define the posterior margins
\(m_{T,\sigma}(x)=\sigma(\theta_0+\theta_1x_1)\) and
\(m_{S,\sigma}(x)=-\sigma\theta_2\), and set
\begin{equation}
  \mathcal P_2
  =\Bigl\{
  P_{\sigma;\theta_0,\theta_1,\theta_2}:
  \sigma\in\{-1,+1\},\
  \theta_0\in[\tfrac12,\tfrac58],\
  \theta_1\in[\tfrac18,\tfrac14],\
  \theta_2\in[\tfrac9{10},\tfrac{19}{20}]
  \Bigr\}.
  \label{eq:opp_class}
\end{equation}
Here \(\theta_0\) and \(\theta_1\) are the target margin offset
and slope and \(\theta_2\) is the constant source margin level;
since \(\theta_0>\theta_1\), the target Bayes sign is \(\sigma\)
everywhere and the source Bayes sign is \(-\sigma\) everywhere.
There is no boundary to estimate inside the support; the unknown
is the target sign \(\sigma\).  This removes boundary-estimation
difficulty and isolates how the labels may be used.  The class
still has a full-dimensional common design, smooth margins,
bounded strong density and a strong target signal.  The limitation
below comes from the final fitting operation, on an event on which
every admissible final fit receives only the target-opposed label
(Remark~\ref{rem:opp_firewall} in the Supplement records why
opposition alone is not such an obstruction).

\begin{defn}[Label-preserving nonnegative-weight
Gaussian-SVM procedures]
\label{def:glob_class}
Let \(\mathcal A\) consist of all measurable randomized wrappers
whose returned classifier is \(\operatorname{sign}(\widehat f+\widehat b)\),
where \((\widehat f,\widehat b)\) is an exact minimizer of
\[
  \widehat\lambda\,\|f\|_{H_{\widehat\gamma}}^2
  +\sum_{i=1}^{n}\widehat w_i^T
  \bigl(1-Y_i^T[f(X_i^T)+b]\bigr)_+
  +\sum_{j=1}^{m}\widehat w_j^S
  \bigl(1-Y_j^S[f(X_j^S)+b]\bigr)_+
\]
over \(f\in H_{\widehat\gamma}\) and \(b\in\mathcal I\), with
\(0<\widehat\gamma<\infty\), \(0<\widehat\lambda<\infty\), the
original observed labels, and nonnegative final weights
\(\widehat w_i^T,\widehat w_j^S\ge0\) of strictly positive total
mass.  The tuning parameters \(\widehat\gamma,\widehat\lambda\) and all
final observation weights may be arbitrary measurable functions of
the full domain-labeled sample and an auxiliary random seed; the
weights may be zero or label-dependent, which permits
data-dependent deletion, validation, tuning, and refitting.  The
intercept set \(\mathcal I\) is either \(\{0\}\) (no intercept; the
minimizer is then unique) or \(\mathbb R\) (an unpenalized
intercept; any measurable selection of an exact minimizer is
allowed).
\end{defn}

The restriction concerns what the final optimization can return.
On an all-common-label sample, every admissible no-intercept fit
is a positive combination of the Gaussian kernel with that label's
sign, and with an unpenalized intercept fitting the common label
has zero objective value and forces the same constant prediction
(Lemma~\ref{lem:sign_lock} in the Supplement,
Section~\ref{supp:adaptive_negative}, pathwise for data-dependent
subsets, weights, bandwidths, and regularization).  Thus even a
procedure that recognizes the source orientation cannot return the
opposite decision while retaining the prescribed final fit.  Exact
minimization is essential when the total observation weight is
unrestricted; the tolerance dependence is explained in
Remark~\ref{cmp:exact-fitting-scope} of the Supplement (Section~\ref{supp:adaptive_negative}).  Write \(\mathcal E^{01}_{T}\) for
the target classification excess under the sampled law, define
\[
  \mathfrak R_{\mathcal A}(n,m)
  :=\inf_{\widehat h\in\mathcal A}\,
  \sup_{P\in\mathcal P_2}
  \mathbb E_P\,\mathcal E^{01}_T(\widehat h),
\]
and let \(\mathfrak R(n,m)\) be the same
minimax risk with the infimum over all measurable randomized
domain-aware classifiers.

\begin{theorem}[Exponential minimax separation under opposition]
\label{thm:opposition_separation}
For all \(n,m\ge1\),
\begin{equation}
  \mathfrak R_{\mathcal A}(n,m)
  \ge\frac12\Bigl(\frac14\Bigr)^{n}
  \Bigl(\frac{19}{20}\Bigr)^{m},
  \qquad
  \mathfrak R(n,m)
  \le\frac58\Bigl(\frac{19}{100}\Bigr)^{m/2}.
  \label{eq:opp_bounds}
\end{equation}
Consequently,
\begin{equation}
  \begin{aligned}
  \frac{\mathfrak R_{\mathcal A}(n,m)}
       {\mathfrak R(n,m)}
  &\ge\frac45\,4^{-n}\Bigl(\frac{19}{4}\Bigr)^{m/2}
  &&\text{for all } n,m\ge1,\\
  \frac{\mathfrak R_{\mathcal A}(n,m)}
       {\mathfrak R(n,m)}
  &\ge\frac45\Bigl(\frac{19}{16}\Bigr)^{m/2}
  &&\text{whenever } n\le m/2.
  \end{aligned}
  \label{eq:opp_ratio}
\end{equation}
\end{theorem}

On the all-label event, whose probability under the member
\(P_{\sigma;\theta_0,\theta_1,\theta_2}\) is
\((\frac{1-\theta_0}{2})^n(\frac{1+\theta_2}{2})^m\), every
admissible final fit predicts the target-opposed label by
Lemma~\ref{lem:sign_lock} and pays the full target margin.  The
rule \(-\operatorname{sign}(\sum_jY_j^S)\) can instead use the
known opposition relation by reversing the source labels.  Both
minimax risks tend to zero, but the error forced by the all-label
event decays more slowly than the bound for the reversing rule
whenever the target
sample is at most half the source sample, producing the exponential
ratio in \eqref{eq:opp_bounds}--\eqref{eq:opp_ratio}.

The distinction is between a decision the model can represent and
one the prescribed fit can return.  The gap remains after
optimizing over data-dependent weights and tuning, and it persists
when the class is known, since the minimax comparison already
optimizes over the whole output family.  On this class, attaining
the unrestricted minimax scale requires changing the permitted
operation; choosing its weights more carefully does not suffice.

\begin{remark}[A limit to simultaneous minimax adaptation]
\label{rem:simultaneous_adaptation}
Theorem~\ref{thm:opposition_separation} rules out simultaneous
minimax adaptation, even up to constant or logarithmic risk
factors, over any collection of transfer classes containing
\(\mathcal P_2\) by a procedure whose output remains in
\(\mathcal A\): for such a procedure, the worst ratio of its risk on
a class to that class's own all-procedure minimax risk is at least
\(\mathfrak R_{\mathcal A}(n,m)/\mathfrak R(n,m)\), which grows
exponentially in \(m\) whenever
\(n\le m/2\) by
\eqref{eq:opp_ratio}.  The conclusion concerns classwise adaptation within the output
family of Definition~\ref{def:glob_class}; it says nothing about
the minimax risk on the union of the classes, which may be
governed by a harder class.
\end{remark}

\section{Numerical Studies}
\label{sec:numerical}

The experiments vary three features of the source separately in one
nonlinear classification model: source sample size
(Section~\ref{sec:general_oracle}), source-label strength with the
Bayes decision fixed (Section~\ref{subsec:source_dominance}), and
source-boundary orientation
(Sections~\ref{subsec:local_global_dictionary}
and~\ref{subsec:bernstein_transfer_oracle}).  Comparing mixture and target excess risks examines
the joint assessment developed above: improved learning of the
training problem and improved target prediction need not coincide.
Spatial maps further distinguish predictions that vary across
training samples from predictions that repeatedly take the wrong
target sign.  We compare fixed and data-selected weights along
these paths while keeping the covariate design unchanged; a
supplementary experiment examines mixture-to-target conversion at
a fixed weight.  The spiral
design below provides a nonlinear boundary whose orientation and
label strength can be varied separately while both domains keep
the same covariate distribution.

The two-dimensional covariate distribution \(P_X\) is generated
around two centrally symmetric spiral arms.  Let
\[
  \mathbf c(t)=\frac{t}{\pi}(\cos t,\sin t)^\top,
  \qquad
  X=Z\,\mathbf c(U)+0.1\,\xi,
\]
where \(U\sim\operatorname{Unif}[0,3\pi]\), \(Z\) takes the
values \(-1\) and \(+1\) with equal probability, the
coordinates of \(\xi\) are independent Student-\(t\) variables
with two degrees of freedom, and all these variables are
independent; \(0.1\) is the noise scale, and the arm indicator
\(Z\) generates feature locations only.  Let \(d_\theta(x)\) be
the signed distance contrast to the two spiral arms after rotation
by \(\theta\), positive when \(x\) is closer to the positive arm,
so that its zero set is the decision boundary; the discretized
reference arms and the exact distance definition are given in the
Supplement (Section~\ref{cmp:numerical-dgp}).  The target and
source posterior margins are
\begin{equation}
  m_T(x)=\tanh\Bigl\{\frac{\beta_T\,d_0(x)}{2}\Bigr\},
  \qquad
  m_S(x)
  =\operatorname{sign}\{d_\theta(x)\}\,
  \Bigl|\tanh\Bigl\{\frac{\beta_S\,d_\theta(x)}{2}\Bigr\}
  \Bigr|^{\zeta},
  \label{eq:sim_dgp}
\end{equation}
with labels drawn from \(P_D(Y=1\mid X=x)=\{1+m_D(x)\}/2\),
\(D\in\{T,S\}\), so that \(m_\alpha=\alpha m_T+(1-\alpha)m_S\) is
the plain average of \eqref{eq:margin_identity}.  We fix
\(n=50\) and \(\beta_T=\beta_S=5\).  The source-size path varies
\(m\) with \(\theta=0\) and \(\zeta=1\), so that \(P_S=P_T\); the
margin and rotation paths use \(m=1{,}000\), varying \(\zeta\) at
\(\theta=0\) and varying \(\theta\) at \(\zeta=1\); the path grids
are listed in the Supplement (Section~\ref{cmp:numerical-dgp}).
Each path changes one feature of the source, and the local vote
\eqref{eq:margin_identity} says what the change does to the
mixture.  On the source-size path \(P_S=P_T\), so
\(P_\alpha=P_T\) and \(Q_\alpha\equiv1\) at every weight,
including the pooling weight \(u=m/(n+m)\) that grows with \(m\);
more source observations change only how precisely this fixed
mixture is estimated.  At zero rotation, changing \(\zeta\)
preserves which label is more likely at every location: values
below one make that preferred label more likely, while values
above one move its probability toward one half (at \(\theta=0\),
\(|m_S|=|m_T|^{\zeta}\), so \(\zeta<1\) is the dominance regime of
{\rm(A1)}); this moves the mixture margin and \(Q_\alpha\) while
\(H_{\rm flip}\) stays zero
(Section~\ref{subsec:bernstein_transfer_oracle}).  Rotation
changes the sign of \(m_S\) over regions and so moves locations
through opposition without a crossing, near cancellation, ties,
and strict flips.

Three fixed rules serve as references: target-only training
\((\alpha=1)\); proportional pooling \((\alpha=n/(n+m))\), which
gives every observation equal weight; and equal-domain weighting
\((\alpha=1/2)\), which gives the two samples equal total weight.
We compare them with a CV-adaptive rule that selects the mixture
weight by target clipped-hinge validation with a target-only
fallback (a practical rule, distinct from the split-sample selector
of Theorem~\ref{thm:target_safe_adaptive}).  A replicate-wise
oracle selects, in each replication, the lowest-risk candidate
among those already fitted.  All fits use
the intercept-augmented version of
\eqref{eq:transfer_svm_empirical} in the standard RBF-SVC
implementation, with kernel \(k(x,x')=\exp\{-2\|x-x'\|_2^2\}\),
regularization \(\lambda=3\times10^{-4}\), an unpenalized
intercept, and observation weights \(\alpha/n\) and
\((1-\alpha)/m\).  The kernel and regularization are fixed across
methods and settings so that the comparisons isolate source
weighting at a common tuning choice, selected in an independent
pilot with matched source and target laws (Supplement,
Section~\ref{supp:numerical}).  Within each replication, methods
share the same underlying samples.  The
design lies outside the opposition class of
Theorem~\ref{thm:opposition_separation}; the study examines the
conversion and borrowing mechanisms.

Figure~\ref{fig:sim_adaptive} and the supplementary fixed-weight
experiment summarize \(1{,}000\) independent Monte Carlo
replications.  In
each replication, an independent reference sample of \(20{,}000\)
covariates from \(P_X\), shared across methods and settings, is
used to evaluate the excess risks and the conversion envelope with
the known margins; this removes the variability of generated test
labels and puts all methods on the same feature locations.
The fixed-weight experiment reports clipped hinge excess to match
the conversion theory, whereas Figures~\ref{fig:sim_adaptive}
and~\ref{fig:sim_votemap} report target classification excess.
Risk curves report Monte Carlo means with error bars of one
standard error.  Remaining implementation details are in the
Supplement (Section~\ref{supp:numerical}).

\begin{figH}[\figposrisk]
\centering
\includegraphics[width=\textwidth]{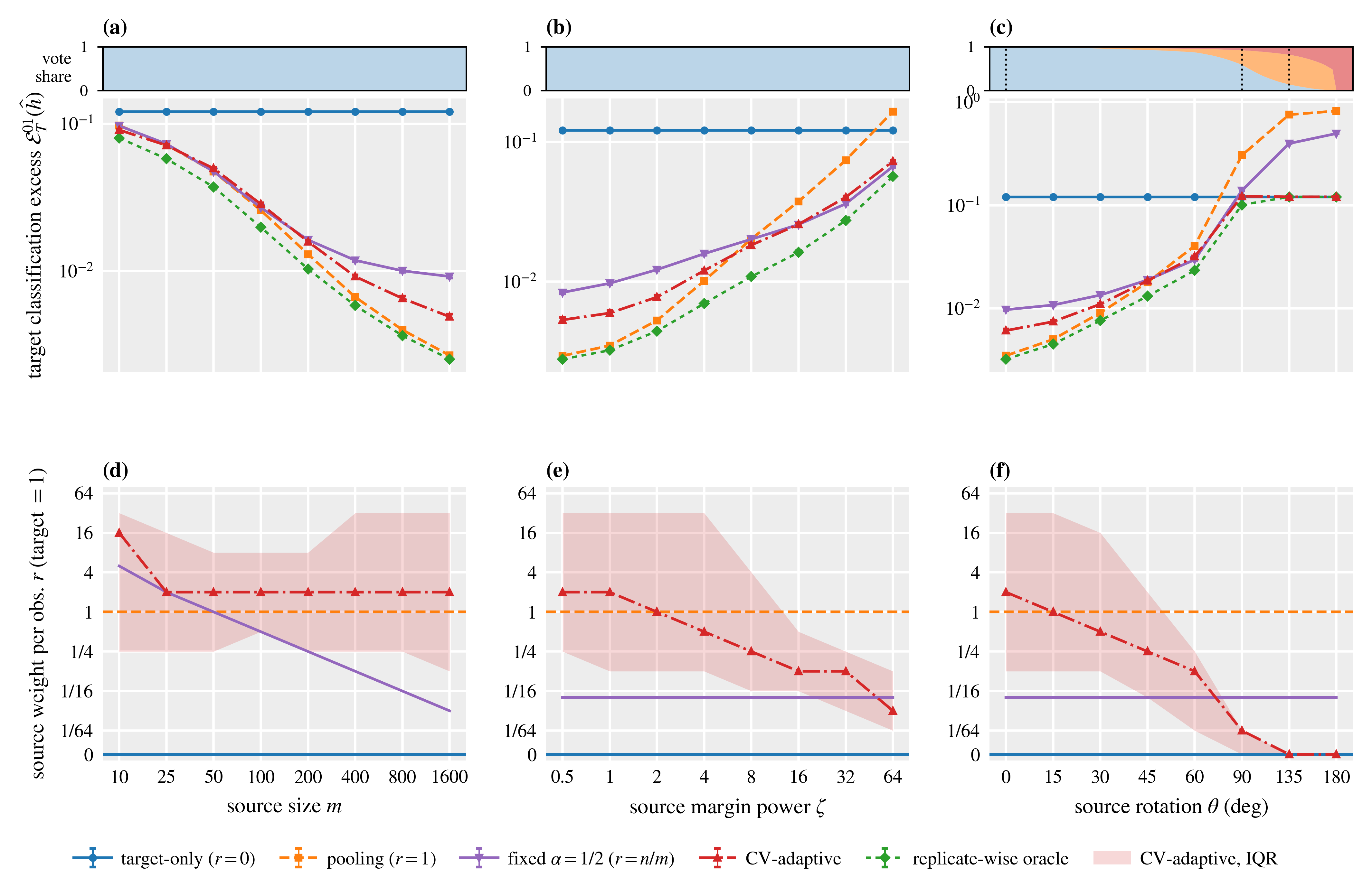}
\caption{Target risk and selected borrowing along the three
perturbation paths.  The top row reports mean target
classification excess \(\mathcal E^{01}_T(\widehat h)\) (log scale)
with pointwise Monte Carlo standard errors for target-only
(\(r=0\)), proportional pooling (\(r=1\)), fixed \(\alpha=1/2\)
(\(r=n/m\)), the CV-adaptive rule, and the replicate-wise oracle.
The bottom row reports the median selected source-to-target
per-observation weight ratio \(r\) with its interquartile range.
Ribbons show population vote regions at the reference weight
\(\alpha=1/2\): agreement (blue), surviving opposition (orange),
and thresholded tie-or-flip (red).  Dotted vertical lines in panel
(c) mark the angles examined in Figure~\ref{fig:sim_votemap}.  The
path grids are specified in the Supplement
(Section~\ref{cmp:numerical-dgp}).}
\label{fig:sim_adaptive}
\end{figH}

\label{subsec:sim_adaptive}\label{subsec:sim_fixed}%
The supplementary fixed-weight experiment
(Section~\ref{cmp:fixed-weight-experiment}, Figure~\ref{fig:sim_bridge})
reports the mixture excess, target excess, and conversion
envelope, the value of \(\Psi_\alpha\) at the fitted mixture
excess, along all three perturbation paths at the equal-domain
weight.
There, source rotation separates mixture accuracy from target
performance: the fitted SVM can remain accurate relative to the
mixture optimum while making increasingly costly target errors, and
the target-excess curve and the conversion envelope both cross the
target-only reference between \(60^\circ\) and \(90^\circ\).  At
\(180^\circ\) the two margins cancel exactly, so every clipped score
has zero mixture excess even though each domain retains label
information.  Here we turn to selecting the source contribution
from target data.

We let target validation select the source weight and compare
the resulting risk with the fixed rules
(Figure~\ref{fig:sim_adaptive}).  To make selected
weights comparable across source sample sizes, we report the
weight of one source observation relative to one target
observation, the scalar ratio \(r=n(1-\alpha)/\{m\alpha\}\)
(distinct from the local vote ratio \(r(x)\) of
Section~\ref{sec:problem_setup} and the residual vector \(r\) of
Supplement Section~\ref{supp:adaptive}):
\(r=0\) is target-only training, \(r=1\) is proportional
pooling, and equal-domain weighting \(\alpha=1/2\) corresponds to
\(r=n/m\), which is \(0.05\) when \(m=1{,}000\).  The candidate
weights form a geometric grid in
\(r\) together with \(r=0\) and \(\alpha=1/2\) (Supplement,
Section~\ref{supp:numerical}), compared by repeated target
cross-validation (five folds, five repetitions, each fold
retaining the full source sample).  A
source-using candidate is retained only if its mean held-out
clipped-hinge loss improves on target-only training by more than
\(0.01\); otherwise the rule selects the target-only candidate.
The reported classifier is refitted on the full training data at
the selected weight.

\begin{figH}[\figposmap]
\centering
\includegraphics[width=\votemapwidth\textwidth]{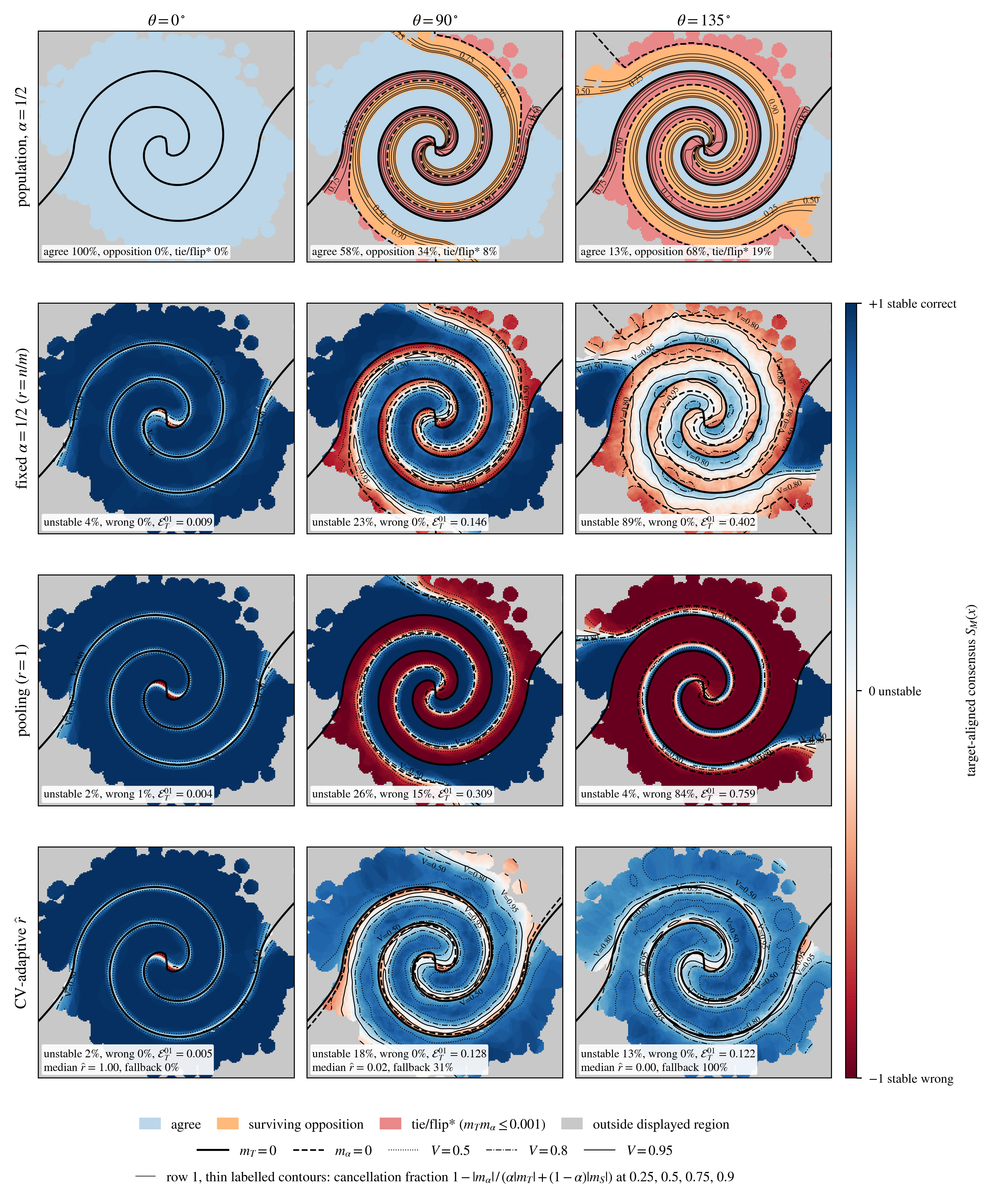}
\caption{Prediction maps behind Figure~\ref{fig:sim_adaptive}(c),
at \(0^\circ\), \(90^\circ\), and \(135^\circ\), from \(100\)
additional training replications per angle (\(n=50\),
\(m=1{,}000\)) on a shared reference sample of \(20{,}000\)
covariates.  The first row shows population vote regions at
\(\alpha=1/2\); tie/flip\(^*\) uses \(m_Tm_\alpha\le10^{-3}\) on
opposing points, including small positive products, and therefore
differs from the exact theoretical set.  The thin labeled contours
in this row show the cancellation fraction
\(1-|m_\alpha|/(\alpha|m_T|+u|m_S|)\), which equals
\(1-|s_\alpha|/(\alpha|s_T|+u|s_S|)\) under the common design, at
\(0.25\), \(0.50\), \(0.75\), and \(0.90\); it is \(0\) where the
domains agree and \(1\) at a tie.  The remaining rows show
equal-domain, pooling, and CV-adaptive fits, colored by the
target-aligned consensus \(S_M\) of \eqref{eq:sim_consensus}: blue
indicates repeated target-sign agreement, red persistent
disagreement, and near-white unstable predictions.  Contours show
\(V_M=1-S_M^2\) at \(0.5\), \(0.8\), and \(0.95\); reported
unstable and stable-wrong shares use \(V_M>0.8\) and \(S_M<-0.8\).
Solid and dashed curves mark the target and population mixture
boundaries; for the CV-adaptive row the mixture boundary is drawn
at the median selected \(\alpha\) over the replications, and the
dashed curve is omitted when this median equals \(1\), as at
\(135^\circ\), where every replication selects target-only
training and the population mixture is the target law.  Gray masks
delimit the displayed region.
Definitions and the shared-reference protocol are in the
Supplement (Section~\ref{cmp:spatial-experiment}).}
\label{fig:sim_votemap}
\end{figH}

When the source and target laws coincide (panels (a) and (d)),
proportional pooling closely tracks the replicate-wise oracle.
The CV-adaptive rule also improves substantially over target-only
training, although a gap from pooling remains at the largest
source sizes.  The value of selection becomes clearer as source
quality deteriorates: the CV rule reduces borrowing as the source
labels weaken (panels (b) and (e)) or disagree more strongly with
the target (panels (c) and (f)), limiting the deterioration seen
under continued pooling.  On the rotation path, where the ribbon
records the expansion of disagreement and flip regions at the
reference weight \(\alpha=1/2\), the CV rule retains substantial
source weight at small rotations, sharply reduces \(r\) as
disagreement grows, and under severe rotation returns to
target-only training, so its target risk recovers target-only
performance; pooling remains fixed at \(r=1\) and incurs large
negative transfer.

Figure~\ref{fig:sim_votemap} examines the predictions behind
the rotation curves in Figure~\ref{fig:sim_adaptive}(c), at the
three marked angles \(0^\circ\), \(90^\circ\), and \(135^\circ\),
using \(100\) additional training replications at each angle.
The first row maps the population vote regions summarized by the
ribbon above Figure~\ref{fig:sim_adaptive}(c), at the reference
weight \(\alpha=1/2\).  The remaining rows correspond to the
equal-domain, pooling, and CV-adaptive curves, respectively; the
replication statistics behind the colors are defined in the
Supplement (Section~\ref{cmp:spatial-experiment}).

At \(0^\circ\), source and target agree, and all three methods
produce predominantly target-aligned predictions, consistent with
the low source-assisted risks in Figure~\ref{fig:sim_adaptive}(c).
At \(90^\circ\), disagreement between the domains expands.
Pooling produces broad red regions of persistent target-sign
error, accompanying its rise above the target-only risk.  The CV
rule reduces borrowing, as shown in
Figure~\ref{fig:sim_adaptive}(f), and substantially reduces these
regions.

At \(135^\circ\), the two fixed borrowing rules fail differently,
and the difference follows their weights.  Equal-domain weighting
(\(u=1/2\)) gives opposing votes of nearly equal magnitude a
mixture vote near zero, so where the rotated source margin is close to
\(-m_T\) the mixture has little preference and the fitted sign
varies across training samples: the near-white regions.
Proportional pooling (\(u\approx0.952\)) gives the source twenty
times the target's total weight, so the same locations lie far past
their thresholds (at a location with \(m_S=-m_T\),
\(Q_\alpha=1-2u\approx-0.9\)); the mixture favors the
target-opposed label with a large margin, and a fit that learns
this mixture well is repeatedly wrong for the target: the red
regions.  High target error thus has two prediction patterns in
this experiment, unstable decisions under equal-domain weighting
and persistent target-sign errors under pooling, and these are
the prediction patterns behind the elevated risks of both rules
in Figure~\ref{fig:sim_adaptive}(c).  The population vote regions
in the first row are drawn at \(\alpha=1/2\) and describe the
equal-domain row; at its own weight, pooling has a larger
tie-or-flip set.  In all \(100\) spatial-study replications at
this angle, the CV rule selects target-only training, which
removes the source's influence on the population decision while
leaving the target's own classification difficulty in place, and
its predominantly blue prediction map accompanies the return to
target-only performance along the rotation path.  The maps show
finite-sample patterns consistent with the population regimes at
each rule's weight.

\FloatBarrier
\section{Discussion}\label{sec:discussion}

The target value of borrowing depends jointly on the learning
problem created by the source and on what the attainable accuracy
on that problem means for the target.  The conversion and learning
analyses separate these effects so that they can be evaluated
together.  More observations improve precision at a fixed mixture;
changing the weight changes the mixture as well as its sampling
bound.  The full frontier also shows why an infinite worst-case
multiplier need not end the analysis: the integrated target cost can
still admit a useful finite-slope comparison.  In the experiments,
rotating the source boundary separates mixture accuracy from target
performance: the fitted SVM
stays accurate for the mixture while its target error grows, and
target validation withdraws borrowing where the fixed rules
continue it.  Neither an easier training problem nor agreement of
the optimal decisions alone determines the target guarantee.

More precise source knowledge does not resolve every uncertainty
relevant to transfer.  On the tilted-boundary class, the
same source law can remain compatible with different target
directions, producing an information limit that applies to every
learner.  A proposed fit can be checked against target evidence; when the
source is informative the check rarely rejects it, so the fit keeps
its faster risk, and the independent fallback is paid only on
rejection.  The procedure learns enough to support the choice
without first recovering the entire relation.

The opposition result adds a further requirement.  Even a known,
informative relation may call for a correction that the prescribed
fit cannot implement.  The gap remains after optimizing over
weights and tuning, so improving the selection within that family
is not enough.  The value of source data depends on the training
problem they create, what the achievable accuracy on that problem
means for the target, and whether the permitted way of using the
source can implement the needed correction.

These distinctions suggest concrete next problems.  A useful
conversion bound might be estimated through a finite-slope,
target-cost-weighted functional, which avoids uniform recovery of
the local multiplier near cancellation.  The target check of
Section~\ref{sec:adaptive} might be made computationally
practical while preserving the candidate's faster risk, and the
necessity of its logarithmic adaptation factor in
Theorem~\ref{thm:target_safe_adaptive} remains open.  When
weighting cannot implement the needed correction, the next object
to choose is the correction itself, as illustrated by the
decision-set and link-scale approaches of \citet{wang2025decision}
and \citet{maity2024linear}.

\section*{Supplementary Material}
\noindent The Supplementary Material follows the references: proofs of Theorems~\ref{thm:affine_frontier}--\ref{thm:opposition_separation} and of all supporting lemmas, the closed-form radial frontier, the nonparametric reinforcement benchmark, the numerical protocol, and auxiliary concentration, entropy-number, and approximation results.

\bibliographystyle{plainnat}
\bibliography{main}

\clearpage
\setcounter{section}{0}
\renewcommand{\thesection}{S.\arabic{section}}
\setcounter{figure}{0}\renewcommand{\thefigure}{S\arabic{figure}}
\setcounter{table}{0}\renewcommand{\thetable}{S\arabic{table}}
\setcounter{theorem}{0}\setcounter{proposition}{0}\setcounter{lemma}{0}
\setcounter{corollary}{0}\setcounter{defn}{0}\setcounter{example}{0}\setcounter{remark}{0}
\renewcommand{\theHsection}{S\arabic{section}}\renewcommand{\theHfigure}{S\arabic{figure}}\renewcommand{\theHtable}{S\arabic{table}}
\renewcommand{\theHtheorem}{S\arabic{theorem}}\renewcommand{\theHproposition}{S\arabic{proposition}}\renewcommand{\theHlemma}{S\arabic{lemma}}
\renewcommand{\theHcorollary}{S\arabic{corollary}}\renewcommand{\theHdefn}{S\arabic{defn}}\renewcommand{\theHexample}{S\arabic{example}}\renewcommand{\theHremark}{S\arabic{remark}}
\renewcommand{\figposrisk}{H}
\renewcommand{\figposconv}{H}
\begin{center}
{\LARGE\bfseries Supplementary Material}\\[10pt]
{\large Supplement to ``Reweighting for Transfer Classification: Exact Risk Conversion and Minimax Adaptation''}
\end{center}
\medskip
\noindent The first four sections give the proofs of the main-text results, organized by exact conversion, the weighted statistical oracle, the Gaussian corollary, and adaptive source use together with the limit of the exact-fitting output family. Kernel-specific inputs appear next to their applications. The numerical protocol and the independent reinforcement benchmark follow these core proofs.  Theorems, propositions, lemmas, corollaries, examples and remarks are numbered consecutively, separately in the main text and in this supplement: statements of the main text carry no prefix, and statements of this supplement carry the prefix ``S''.  Sections, equations and figures of this supplement also carry the prefix ``S''.

\renewcommand{\thetheorem}{S\arabic{theorem}}
\renewcommand{\theproposition}{S\arabic{proposition}}
\renewcommand{\thelemma}{S\arabic{lemma}}
\renewcommand{\thecorollary}{S\arabic{corollary}}
\renewcommand{\thedefn}{S\arabic{defn}}
\renewcommand{\theexample}{S\arabic{example}}
\renewcommand{\theremark}{S\arabic{remark}}

\section*{Contents}\mbox{}\par\medskip
\begingroup
\setlength{\parindent}{0pt}
\setlength{\parskip}{3pt}
\newcommand{\tocsec}[2]{\par\noindent\textbf{\ref{#1}\enspace \hyperref[#1]{#2}}\nobreak\ \dotfill\nobreak\ \pageref{#1}\par}
\newcommand{\tocsub}[2]{\par\noindent\hangindent 3.2em\hspace*{1.2em}\ref{#1}\enspace \hyperref[#1]{#2}\nobreak\ \dotfill\nobreak\ \pageref{#1}\par}
\newcommand{\tocsubsub}[2]{\par\noindent\hangindent 5.4em\hspace*{3.0em}\ref{#1}\enspace \hyperref[#1]{#2}\nobreak\ \dotfill\nobreak\ \pageref{#1}\par}
\tocsec{sec:proofs}{Exact conversion: proofs and threshold tools}
\tocsub{supp:frontier}{The affine frontier and conversion modulus (proof of Theorem~\ref*{thm:affine_frontier}; Lemma~\ref*{lem:flip_extinction}, Example~\ref*{ex:three_point})}
\tocsub{supp:mixing}{Regime bounds and their sharpness (proof of Corollary~\ref*{cor:regime_flip})}
\tocsec{supp:oracle}{The weighted statistical oracle}
\tocsub{supp:crerm}{Weighted-sample oracle inequalities}
\tocsubsub{cmp:oracle-inputs}{Localized inputs for the main-text oracle (Lemmas~\ref*{lem:weighted_bernstein_scalar}--\ref*{lem:sym_peel})}
\tocsubsub{cmp:general-oracle-sec}{The general oracle inequality (Theorem~\ref*{thm:weighted_crerm_general})}
\tocsubsub{pf:weighted_crerm}{Proof of Theorem~\ref*{thm:weighted_crerm}}
\tocsubsub{pf:benign}{Weighted chaining and the entropy oracle (Lemma~\ref*{lem:weighted_chaining}, Theorem~\ref*{thm:transfer_bernstein_svm})}
\tocsub{supp:faster}{Hinge variance (Lemma~\ref*{lem:variance-bound-hinge})}
\tocsub{supp:mixture_exponents}{Proof of Proposition~\ref*{lem:one_acceleration}}
\tocsub{pf:two_block_hinge_variance}{A domainwise variance profile under reinforcement}
\tocsec{supp:gaussian}{Gaussian inputs and the proof of Corollary~\ref*{thm:main_target_bound}}
\tocsub{supp:gauss_detail}{Mixture conditions and a Gaussian comparator (Lemma~\ref*{lem:gaussian_mixture_approximation})}
\tocsub{supp:apalpha}{Proof of Lemma~\ref*{lem:gaussian_mixture_approximation}}
\tocsub{sec:appB}{Weighted empirical Gaussian entropy (Lemma~\ref*{lem:weighted_empirical_entropy})}
\tocsub{supp:gauss}{The Gaussian mixture oracle (Theorem~\ref*{thm:gaussian_transfer_estimation})}
\tocsub{pf:main_target_bound}{Proof of Corollary~\ref*{thm:main_target_bound} and balanced tuning (Corollary~\ref*{cor:transfer_rates})}
\tocsec{supp:adaptive}{Adaptation and the exact-fitting limit}
\tocsub{supp:adaptive_calibration}{Radial geometry and source-quality estimation (Lemmas~\ref*{lem:supp_radial}--\ref*{lem:signed_projection_confidence})}
\tocsub{supp:radial_frontier}{The radial affine frontier (proof of Proposition~\ref*{prop:radial_frontier})}
\tocsub{supp:adaptive_positive}{The target-checked selector (Lemmas~\ref*{lem:supp_bridge}--\ref*{lem:supp_lower}, Example~\ref*{ex:supp_witness}; proofs of Corollary~\ref*{cor:fixed_weight_signed} and Theorem~\ref*{thm:target_safe_adaptive})}
\tocsub{supp:adaptive_negative}{The exact-fitting limit (Lemma~\ref*{lem:sign_lock}; proof of Theorem~\ref*{thm:opposition_separation})}
\tocsec{supp:numerical}{Numerical protocol}
\tocsub{cmp:fixed-weight-experiment}{Fixed-weight experiment: mixture learning and target performance}
\tocsub{cmp:spatial-experiment}{Spatial diagnostics for the prediction maps}
\tocsub{supp:numerical_protocol_details}{Cross-validation rules and implementation details}
\tocsec{supp:reinforcement}{A nonparametric reinforcement benchmark}
\tocsub{supp:reinforcement_statement}{The reinforcement class and the benchmark theorem (Lemma~\ref*{lem:exact_gaussian_vote}, Theorem~\ref*{thm:gaussian_reinforcement_minimax})}
\tocsub{supp:gaussian_minimax}{Proof of Theorem~\ref*{thm:gaussian_reinforcement_minimax} (Lemmas~\ref*{lem:gaussian_regular_moments}--\ref*{lem:gmm_hard_membership})}
\tocsubsub{pf:gaussian_reinforcement_upper}{Direct classification upper bound}
\tocsubsub{pf:gmm_fixed_class_lower}{Matching lower bound with fixed class constants}
\tocsubsub{pf:unrestricted_source}{No acceleration without a relation}
\endgroup
\bigskip

\section{Exact Conversion: Proofs and Threshold Tools}\label{sec:proofs}

\subsection{The affine frontier and conversion modulus (Theorem~\ref{thm:affine_frontier})}
\label{supp:frontier}
\phantomsection\label{pf:frontier}%

The proof separates target excess that can occur at zero mixture
cost from excess created on the stable region.  On the stable
region, the extremal score allocates its mixture-error budget to
locations with the largest target-to-mixture ratio.  Throughout,
write \(a=p_T|m_T|\), \(b=p_\alpha|m_\alpha|\), recall
\(s_T=p_Tm_T\), \(s_\alpha=p_\alpha m_\alpha\), abbreviate the
stable set by
\(S_\alpha:=\mathcal D_\alpha^+=\mathcal D_T\cap\{Q_\alpha>0\}\),
and write \(\varpi_\alpha:=2\nu_\alpha|_{S_\alpha}\) for the
restricted capacity measure of \eqref{eq:surviving_vote_measure},
so that \(d\varpi_\alpha=2b\,\mathbf 1_{S_\alpha}d\mu\), and
\(L_\alpha:=\varpi_\alpha(\mathcal X)=2\nu_\alpha(\mathcal D_\alpha^+)\)
for the vote capacity; the tie-or-flip
target-active set is
\(\mathcal C(u)\) of \eqref{eq:charged_as_threshold_sublevel}, and
target ties \(\{p_T|m_T|=0\}\) carry no target excess.

\begin{proof}[Proof of Theorem~\ref{thm:affine_frontier}]
We also prove the sharpness claims of Remark~\ref{rem:affine_sharpness}
and the attainment claims for randomized and deterministic rules.
At \(\alpha=1\), \(\mathcal E_1(g)=\mathcal E_T(g)\) for every
clipped score. Recall \(J_T=\nu_T(\mathcal X)=\mathbb E_{P_{X,T}}|m_T|\). The scores
\(g=(1-r)f_T^*\), \(0\le r\le2\), attain all target excesses
\(rJ_T\) between zero and \(2J_T\). Hence
\(\Psi_1(t)=\min\{t,2J_T\}\) and
\(\mathcal B_1(c)=2J_T(1-c)_+\), including when \(J_T=0\).
Here \(Q_1=1\) on \(\mathcal D_T\), \(H_{\rm flip}(0)=0\), and
\(L_1=2J_T\), so these are exactly the stated formulas and their
duality. The argument below applies to every \(\alpha\in(0,1]\)
and also establishes the hard-classifier attainment claims.

\emph{Step 1: the frontier formula, the first display of part
(ii).}
For clipped scores both excess risks are linear:
by \eqref{eq:hinge_pointwise_excess},
\(\mathcal E_T(g)=\int(a-s_Tg)\,d\mu\) and
\(\mathcal E_\alpha(g)=\int(b-s_\alpha g)\,d\mu\)
after absorbing the Bayes terms, so
\[
  \mathcal E_T(g)-c\,\mathcal E_\alpha(g)
  =\int\bigl\{(a-cb)-(s_T-c\,s_\alpha)g\bigr\}\,d\mu .
\]
The integrand is bounded by \(2a+2cb\), and the supremum over
\(|g|\le1\) decouples pointwise: taking
\(g_c=-\operatorname{sign}(s_T-c\,s_\alpha)
=\operatorname{sign}(c\,s_\alpha-s_T)\) turns the last term into
\(+|s_T-c\,s_\alpha|\), with equality for a given \(g\) iff
\((c\,s_\alpha-s_T)g=|c\,s_\alpha-s_T|\) a.e.  This proves the
first equality in \eqref{eq:frontier_closed}.  For the regional
form, evaluate the integrand \(a-cb+|s_T-cs_\alpha|\) on the four
regions: on \(S_\alpha\) the two votes share a sign, so
\(|s_T-cs_\alpha|=|a-cb|\) and the integrand is
\(a-cb+|a-cb|=2(a-cb)_+\); on the flip part the signs are opposite
and the integrand is \(a-cb+(a+cb)=2a\); on a mixture tie \(b=0\)
and the integrand is \(2a\); on target ties \(a=0\) and it
vanishes.  Summing, and recalling
\(H_{\rm flip}(u)=2\int_{\mathcal C(u)}a\,d\mu\), gives
\[
  \mathcal B_\alpha(c)
  =H_{\rm flip}(u)+2\int_{S_\alpha}(a-cb)_+\,d\mu
  =H_{\rm flip}(u)+\int(A_\alpha-c)_+\,d\varpi_\alpha,
\]
the last step factoring \(2(a-cb)_+=2b(a/b-c)_+\) and using
\(A_\alpha=a/b\) on \(S_\alpha\) with
\(d\varpi_\alpha=2b\,\mathbf 1_{S_\alpha}d\mu\).  The layer-cake
identity \((x-c)_+=\int_c^\infty\mathbf 1\{x>t\}\,dt\) and Fubini
give
\(\mathcal B_\alpha(c)-H_{\rm flip}(u)
=\int_c^\infty\varpi_\alpha\{A_\alpha>t\}\,dt\),
a survival-function integral of a nonnegative decreasing function
of \(c\); convexity, monotonicity, continuity, and the
one-sided derivatives (by dominated convergence on difference
quotients) follow, as do the endpoints:
\(\mathcal B_\alpha(0)
=H_{\rm flip}(u)+\int A_\alpha\,d\varpi_\alpha
=2\int_{S_\alpha\cup\mathcal C(u)}a\,d\mu
=2\,\mathbb E_{P_{X,T}}|m_T|\)
(target ties contribute zero), and
\(\mathcal B_\alpha(c)\downarrow H_{\rm flip}(u)\) by monotone
convergence.  A kink at \(c\) means the two one-sided derivatives
differ, i.e.\ \(\varpi_\alpha\{A_\alpha=c\}>0\).

\emph{Step 2: the modulus \eqref{eq:conversion_modulus} of part
(i) and the duality \eqref{eq:conversion_duality} of part (ii).}
Fix the baseline score
\[
  g_0(x):=
  \begin{cases}
    \operatorname{sign}s_\alpha(x), & s_\alpha(x)\ne0,\\
    -\operatorname{sign}s_T(x), & s_\alpha(x)=0,\ s_T(x)\ne0,\\
    0, & s_\alpha(x)=s_T(x)=0:
  \end{cases}
\]
the mixture Bayes sign wherever the mixture vote is nonzero
(this is the common Bayes sign on \(S_\alpha\), the target-opposed
sign on flips, and the mixture sign on target-inactive points with
\(b>0\), including points with \(p_T=0\)), the target-opposed
sign on mixture ties with a nonzero target vote, and zero where
both votes vanish.  Its mixture excess is zero, since
\(b-s_\alpha g_0=0\) pointwise, and its target excess is exactly
\(H_{\rm flip}(u)\): the integrand \(a-s_Tg_0\) equals \(0\) on
\(S_\alpha\), \(2a\) on flips and on mixture ties with \(a>0\),
and \(0\) on target ties, whatever the value of \(b\) there.  (A
score chosen arbitrarily on target ties need not have zero mixture
excess: at a point with \(m_T=0\), \(m_S=1\), and common design,
\(m_\alpha=1-\alpha\), and the choice \(g=-1\) costs mixture
excess \(2(1-\alpha)>0\) while costing no target excess.)  Any
deviation from the baseline on the flip part spends mixture budget
and strictly lowers target excess; on mixture ties deviations are
free and the baseline already collects the full charge; on target
ties deviations spend mixture budget and buy no target excess, so
the extremal allocation never uses them.  On
\(S_\alpha\), parametrize \(g\) by
\(y=(1-\sigma g)/2\in[0,1]\), \(\sigma\) the common Bayes sign; the
pointwise mixture and target costs of \(y\) are \(2b\,y\) and
\(2a\,y\).  Hence
\begin{equation}
  \Psi_\alpha(t)
  =H_{\rm flip}(u)+
  \sup\Bigl\{\int A_\alpha\,y\,d\varpi_\alpha:\
  0\le y\le1,\ \int y\,d\varpi_\alpha\le t\Bigr\}.
  \label{eq:knapsack_form}
\end{equation}
For the upper bound, layer cake in the value variable gives
\[
\begin{aligned}
\int A_\alpha y\,d\varpi_\alpha
&=\int_0^\infty\int_{\{A_\alpha>z\}}y\,d\varpi_\alpha\,dz\\
&\le\int_0^\infty\min\{t,\varpi_\alpha(A_\alpha>z)\}\,dz\\
&=\int_0^{t\wedge L_\alpha}A_\alpha^{\downarrow}(s)\,ds.
\end{aligned}
\]
For attainment at \(0<t<L_\alpha\), let
\(c^*=A_\alpha^{\downarrow}(t)\), so
\(\varpi_\alpha\{A_\alpha>c^*\}\le t
\le\varpi_\alpha\{A_\alpha\ge c^*\}\) by right continuity of the
tail, and take the threshold rule
\(y^*=\mathbf 1\{A_\alpha>c^*\}+\lambda\mathbf 1\{A_\alpha=c^*\}\)
with the fractional filling
\(\lambda=(t-\varpi_\alpha\{A_\alpha>c^*\})/
\varpi_\alpha\{A_\alpha=c^*\}\) when the level set has positive
mass, \(\lambda=0\) otherwise; it spends exactly \(t\) and attains
the layer-cake value.  At \(t=0\) take \(y=0\), and at
\(t\ge L_\alpha\) take \(y=1\), which yields the saturation value
\(H_{\rm flip}(u)+\int A_\alpha d\varpi_\alpha
=2\,\mathbb E_{P_{X,T}}|m_T|\).  For the dual form, the pointwise
inequality \(A_\alpha y\le c\,y+(A_\alpha-c)_+\) for
\(0\le y\le1\), integrated, gives
\(\Psi_\alpha(t)\le c\,t+\mathcal B_\alpha(c)\) for every
\(c\ge0\); at \(c=c^*\) the threshold rule turns both inequalities
into equalities, closing the duality gap, and the endpoint cases
\(t=0\) (infimum approached as \(c\to\infty\), attained iff
\(M_\alpha<\infty\)) and \(t\ge L_\alpha\) (optimal \(c=0\)) are
checked directly.  Concavity and monotonicity are clear from the
quantile form.  Geometrically, the attainable pairs
\((\mathcal E_\alpha(g),\mathcal E_T(g))\) form a convex set, since
both coordinates are affine in the clipped score; \(\Psi_\alpha(t)\)
is the largest target coordinate among pairs whose mixture
coordinate is at most \(t\), and \(\mathcal B_\alpha(c)\) is the
smallest intercept of a line of slope \(c\) dominating all pairs,
which is the picture behind Remark~\ref{rem:affine_sharpness}.

\emph{Step 3: attainment, the last assertion of part (i), and its
scope.}  Every score
\(g\in[-1,1]\) is the conditional
mean of a randomized hard rule with
\(P(H=1\mid X=x)=(1+g(x))/2\), at identical linear risks; hence
exactness for scores transfers to externally randomized rules.  For
deterministic rules, \(y\) in \eqref{eq:knapsack_form} is
restricted to \(\{0,1\}\): an indivisible knapsack.  If
\(\varpi_\alpha\) is nonatomic, a measurable-splitting argument
(for a finite nonatomic measure \(\lambda\), every measurable
\(E\) and every \(s\in[0,\lambda(E)]\) admit a measurable
\(E_s\subseteq E\) with \(\lambda(E_s)=s\); recursively bisect
positive-mass pieces and pass to the limit) selects, inside the
level set \(\{A_\alpha=c^*\}\), a deterministic set of exactly the
required mass, so the deterministic value equals
\eqref{eq:conversion_modulus}.  A one-point example shows failure
on an atom:
there \(\mu(\{x_0\})=2\), \(s_T=s_\alpha=1/2\), the two
deterministic rules give the risk pairs \((0,0)\) and \((2,2)\),
so the deterministic modulus vanishes on \(0<t<2\), while the
score \(g=1-t\) attains \((t,t)\).

\emph{Step 4: sharpness and the flip profile.}
The minimality claims of Remark~\ref{rem:affine_sharpness} follow
from the frontier formula: an affine intercept \(\beta\) is valid
exactly when \(\beta\ge\mathcal B_\alpha(c)\), by
\eqref{eq:frontier_def};
\(\int(A_\alpha-c)_+\,d\varpi_\alpha=0\) exactly when
\(c\ge M_\alpha\) by definition of the essential supremum, and
\(\mathcal B_\alpha(c)\downarrow H_{\rm flip}(u)\) was shown by
monotone convergence; primal attainment was established in Steps
2--3.  Finally, \(\mathcal D_T\cap\{Q_\alpha\le0\}=\mathcal C(u)\),
so
\[
 H_{\rm flip}(u)=2\int_{\{Q_\alpha\le0\}}|m_T|\,dP_{X,T},
\]
which gives \eqref{eq:profile_representation} together with
\(\Psi_\alpha(0)=H_{\rm flip}(u)\).
\end{proof}

We next record the uniform no-flip range used in
Proposition~\ref{lem:one_acceleration}, and then a
three-point example that exhibits two different slopes and a
zero-cost contribution in Theorem~\ref{thm:affine_frontier}.

Take essential suprema and infima with respect to \(P_{X,T}\),
with the conventions \(\operatorname*{ess\,sup}:=0\) and
\(\operatorname*{ess\,inf}:=1\) over a \(P_{X,T}\)-null index set,
and let
\(\bar r_-:=\operatorname*{ess\,sup}_{\Omega_-\cap\mathcal D_T}r\).

\begin{lemma}[A uniform no-flip range]
\label{lem:flip_extinction}
The essential onset of the tie-or-flip region is
\begin{equation}
  u_c^{\rm ess}
  :=
  \operatorname*{ess\,inf}_{
    \Omega_-\cap\mathcal D_T}u_c(x)
  =
  \frac{1}{1+\bar r_-},
  \label{eq:global_flip_threshold}
\end{equation}
and for \(u<u_c^{\rm ess}\), no positive-target-margin point is tied or
flipped \(P_{X,T}\)-almost surely: \(H_{\rm flip}(u)=0\).  If,
moreover, \(P_{X,S}\ll P_{X,T}\) with
\(dP_{X,S}/dP_{X,T}\le\phi\) for some \(1\le\phi<\infty\), and
\(\rho_-:=\operatorname*{ess\,sup}_{\Omega_-\cap\mathcal
D_T\cap\{p_S>0\}}|m_S|/|m_T|<\infty\), then
\[
  u_c^{\rm ess}
  \ge
  \underline u_c:=\frac{1}{1+\phi\rho_-},
\]
and \(u<\underline u_c\), equivalently
\((u/(1-u))\,\phi\rho_-<1\), is a sufficient uniform no-flip
condition.
\end{lemma}

\begin{proof}
A point of \(\Omega_-\cap\mathcal D_T\) ties or flips at
source weight \(u\) exactly when \(u\ge u_c(x)\) by
\eqref{eq:charged_as_threshold_sublevel},
\(u_c(x)=(1+r(x))^{-1}\) is decreasing in \(r(x)\), and the
domination assumptions give \(r\le\phi\rho_-\) on the
disagreement region.  Crossing \(\underline u_c\) does \emph{not} imply
that a flip occurs: the two worst cases \(\phi\) and \(\rho_-\)
need not be attained at the same point, whereas the exact onset
depends on the pointwise product \(r(x)\), and a point contributes
to \(H_{\rm flip}\) only after its own threshold is crossed.
\end{proof}

\begin{remark}[The source-only endpoint]
\label{cmp:source-only-endpoint}
The flip profile \eqref{eq:flip_profile} and the modulus of
Theorem~\ref{thm:affine_frontier} are stated for \(0\le u<1\).
When \(u=1\) the mixture vote is the source vote alone, so every
target-active location with a zero source vote (\(s_S=0\),
\(s_T\ne0\)) is a mixture tie and can carry the full target cost
\(2|m_T|\) at zero mixture excess, although it lies outside
\(\Omega_-\): with \(m_T=1\) and \(m_S=0\), the integral in
\eqref{eq:flip_profile} is zero, while the score \(g=-1\) has zero
source hinge excess and target hinge excess \(2\).  The zero-budget
target cost at \(u=1\) is therefore
\[
  2\int_{\mathcal D_T\cap\{s_Ts_S\le0\}}|m_T|\,dP_{X,T}
  =\lim_{u\uparrow1}H_{\rm flip}(u)
  +2\int_{\mathcal D_T\cap\{s_S=0\}}|m_T|\,dP_{X,T},
\]
which exceeds the limit of \eqref{eq:flip_profile} by the
target-margin mass of \(\{s_S=0\}\).  Consequently the interior
identity \(\Psi_\alpha(0)=H_{\rm flip}(u)\) cannot be extended to
\(\alpha=0\) without this correction; a source-only fit, where it
appears, is treated separately, and the target-only endpoint
\(\alpha=1\) (\(Q_1\equiv1\), \(\Psi_1(t)=t\) before saturation)
is unaffected.
\end{remark}

\begin{example}[A three-point conversion]
\label{ex:three_point}
Let \(\alpha=\tfrac12\), let \(\mu\) put mass \(\tfrac13\) on
each of three points \(x_0,x_1,x_2\), let the target vote be
\(s_T\equiv1\) there (\(p_T\equiv1\), \(m_T\equiv+1\)), and let
the source votes be \(s_S=(-\tfrac32,-\tfrac12,1)\).  Then
\(Q_\alpha=(-\tfrac14,\tfrac14,1)\): the source flips the mixture
at \(x_0\), destabilizes it at \(x_1\), and reinforces it at
\(x_2\).  Hence \(H_{\rm flip}(u)=\tfrac23\), the multipliers are
\(A_\alpha(x_1)=4\) and \(A_\alpha(x_2)=1\) with
\(\varpi_\alpha\)-masses \(\tfrac16\) and \(\tfrac23\), and
\[
  \Psi_\alpha(t)=
  \begin{cases}
    \tfrac23+4t, & 0\le t\le\tfrac16,\\[2pt]
    \tfrac43+(t-\tfrac16), & \tfrac16\le t\le\tfrac56,\\[2pt]
    2, & t\ge\tfrac56,
  \end{cases}
  \qquad
  \mathcal B_\alpha(c)=\tfrac23+\tfrac16(4-c)_+
  +\tfrac23(1-c)_+ .
\]
The extremal score errs at \(x_0\) at no mixture cost, then buys
target excess at rate \(4\) until the atom \(x_1\) is exhausted,
then at rate \(1\).  The single worst-case slope
\(M_\alpha=4\) with intercept \(\tfrac23\) is one supporting
line; the envelope also records the second slope \(1\), the kink
at \(t=\tfrac16\), and the saturation at \(t=\tfrac56\).
\end{example}

The next statement records the assumption-free affine comparison
at the slope \(1/\alpha\), referred to in the main text.

\begin{remark}[Cancellation, and a finite conversion without any
stability assumption]
\label{rem:finite_conversion}
Substituting \(c=1/\alpha\) into \eqref{eq:frontier_closed} gives
the cancellation identity
\begin{equation}
  \alpha\,\mathcal B_\alpha(1/\alpha)
  =
  2\int
  \min\bigl\{\alpha\,p_T|m_T|,\;u\,p_S|m_S|\bigr\}
  \mathbf 1\{s_Ts_S<0\}\,d\mu :
  \label{eq:cancellation_identity}
\end{equation}
at the assumption-free slope, the intercept is exactly the
\(L_1\)-defect of the triangle inequality for the two weighted
votes---the total overlap of conflicting votes, not the full
opposing source mass.  Consequently, for every
\(\alpha\in(0,1]\), every measurable
\(g:\mathcal X\to[-1,1]\), and
\(\Gamma_{-}:=2\int_{\Omega_-}|m_S|\,dP_{X,S}\),
\begin{equation}
  \mathcal E_T(g)
  \le
  \frac1\alpha\mathcal E_\alpha(g)
  +\frac2\alpha\int_{\Omega_-}
  \min\bigl\{\alpha\,p_T|m_T|,\;u\,p_S|m_S|\bigr\}\,d\mu
  \le
  \frac1\alpha\mathcal E_\alpha(g)
  +\frac{1-\alpha}{\alpha}\Gamma_{-},
  \label{eq:finite_vote_conversion}
\end{equation}
where the second inequality bounds the minimum by its second argument, which gives a coarser bound.  This comparison is always
finite---\(M_\alpha\) may be infinite and no tail condition is
needed---but it does not replace the exact conversion.
\end{remark}

\begin{proof}[Proof of \eqref{eq:cancellation_identity}]
We also prove the first inequality in Remark~\ref{rem:finite_conversion}.
Substitute \(c=1/\alpha\) into \eqref{eq:frontier_closed}:
\(s_T-\alpha^{-1}s_\alpha=-(u/\alpha)s_S\), so
\[
\begin{aligned}
  \alpha\,\mathcal B_\alpha(1/\alpha)
  &=\int\bigl(\alpha|s_T|-|s_\alpha|+u|s_S|\bigr)\,d\mu\\
  &=\alpha\|s_T\|_1+u\|s_S\|_1-\|\alpha s_T+u s_S\|_1 ,
\end{aligned}
\]
and the elementary identity
\(|x|+|y|-|x+y|=2\min\{|x|,|y|\}\mathbf 1\{xy<0\}\) with
\(x=\alpha s_T\), \(y=us_S\) gives
\eqref{eq:cancellation_identity}.  The first inequality of
\eqref{eq:finite_vote_conversion} is
\(\mathcal E_T\le c\,\mathcal E_\alpha+\mathcal B_\alpha(c)\) at
\(c=1/\alpha\); the second drops the minimum against \(u|s_S|\) and
recognizes \(\Gamma_-\).
\end{proof}

\begin{remark}[Weight dependence of the affine frontier]
\label{rem:alpha_nonconvex}
For fixed finite \(c>0\), the map
\(\alpha\mapsto\mathcal B_\alpha(c)\) is Lipschitz but in general
not convex: on a one-point feature space with \(Y=+1\) almost
surely under the target and \(Y=-1\) almost surely under the
source,
\[
  \mathcal B_\alpha(1)=
  \begin{cases}
    2, & \alpha\le1/2,\\
    4(1-\alpha), & \alpha\ge1/2,
  \end{cases}
\]
a downward kink at \(\alpha=1/2\)
(proved below).  \end{remark}

\begin{proof}[Proof of the nonconvexity example in
Remark~\ref{rem:alpha_nonconvex}]
On the one-point space with \(Y\equiv+1\) under \(P_T\) and
\(Y\equiv-1\) under \(P_S\): \(\mu(\{x_0\})=2\), \(p_T=p_S=1/2\),
\(s_T=1/2\), \(s_S=-1/2\), and \(s_\alpha=(2\alpha-1)/2\).  For
\(\alpha\le1/2\) the pair is tied or flipped, so
\(\mathcal B_\alpha(1)=2\int|s_T|d\mu=2\); for \(\alpha\ge1/2\) the
point is stable with \(a=1/2\), \(b=(2\alpha-1)/2\), so
\(\mathcal B_\alpha(1)=2\cdot2\cdot(a-b)_+=4(1-\alpha)\).  The map
is Lipschitz in \(\alpha\) for each fixed finite \(c\) (the
integrand of \eqref{eq:frontier_closed} is \(2c\)-Lipschitz in
\(s_\alpha\), which is affine in \(\alpha\)), and the displayed
downward kink at \(\alpha=1/2\) rules out convexity.
\end{proof}

\subsection{Regime bounds and their sharpness (Corollary~\ref{cor:regime_flip})}
\label{supp:mixing}

This subsection proves the three regime bounds and their sufficient no-flip condition.
As in Section~\ref{supp:frontier}, \(\varpi_\alpha:=2\nu_\alpha|_{\mathcal D_\alpha^+}\)
denotes the capacity measure restricted to the stable set, so that
\(d\varpi_\alpha=2\,p_\alpha|m_\alpha|\mathbf 1_{\mathcal D_\alpha^+}d\mu\),
and \(L_\alpha:=\varpi_\alpha(\mathcal X)\) is the vote capacity.

\phantomsection\label{pf:regime_bounds}%
\begin{proof}[Proof of Corollary~\ref{cor:regime_flip}]
On \(\mathcal D_T\), \(p_\alpha|m_\alpha|=Q_\alpha\,p_T|m_T|\), so
a \(\nu_T\)-null subset of \(\mathcal D_T\) is
\(\varpi_\alpha\)-null.

\emph{Case (i).} \(\{Q_\alpha\le0\}\subset
\Omega_-\cap\mathcal D_T\) is \(\nu_T\)-null, so
\(H_{\rm flip}(u)=2\,\nu_T\{Q_\alpha\le0\}=0\); off the
\(\nu_T\)-null set \(\Omega_-\), the stable multiplier is
\(A_\alpha=(\alpha+ur)^{-1}\le\alpha^{-1}\), so
\(M_\alpha\le\alpha^{-1}\).

\emph{Case (ii).} \(Q_\alpha\ge\delta\)
\(\nu_T\)-a.e.\ makes \(\{Q_\alpha\le0\}\) \(\nu_T\)-null and
gives \(A_\alpha=Q_\alpha^{-1}\le\delta^{-1}\)
\(\varpi_\alpha\)-a.e., so \(H_{\rm flip}(u)=0\) and
\(M_\alpha\le\delta^{-1}\).  In both cases the duality
\eqref{eq:conversion_duality} at \(c=M_\alpha\) gives
\(\Psi_\alpha(t)\le M_\alpha t+\mathcal B_\alpha(M_\alpha)
=M_\alpha t\), since
\(\mathcal B_\alpha(M_\alpha)=H_{\rm flip}(u)=0\) by
\eqref{eq:frontier_closed}.  The closed-form sufficient condition
for the buffer of case (ii) is the next proof.

\emph{Case (iii).} The condition
\eqref{eq:mixture_margin_profile} bounds the affine frontier at every finite slope, and the
duality \eqref{eq:conversion_duality} gives the H\"older bound,
with the explicit constant
\(C=(2c_0)^{1/(\chi_0+1)}\{(\chi_0+1)/\chi_0\}^{\chi_0/(\chi_0+1)}\)
(the detailed calculation follows).

\phantomsection\label{pf:regime_case_iii}%
On the stable target-active set,
\(A_\alpha=Q_\alpha^{-1}\) and
\(d\varpi_\alpha
=2\,Q_\alpha\,|m_T|\,dP_{X,T}\),
since \(p_\alpha|m_\alpha|=Q_\alpha\,p_T|m_T|\) there.  Write
\(k:=\chi_0>0\) and
\(F(z):=\nu_T\{0<Q_\alpha<z\}\), \(z>0\), the target-margin mass
at which the retained vote is positive but below \(z\); condition
\eqref{eq:mixture_margin_profile} reads \(F(z)\le c_0z^{k}\).  We
integrate the affine frontier directly.  For every slope \(c>0\),
by \eqref{eq:frontier_closed},
\begin{equation}
\begin{aligned}
  \mathcal B_\alpha(c)-H_{\rm flip}(u)
  &=\int(A_\alpha-c)_+\,d\varpi_\alpha\\
  &=2\int_{\{Q_\alpha>0\}}(1-cQ_\alpha)_+\,d\nu_T\\
  &=2c\int_0^{1/c}F(z)\,dz
  \le\frac{2c_0}{k+1}\,c^{-k},
\end{aligned}
  \label{eq:destab_tail}
\end{equation}
where the third equality uses
\((1-cq)_+=c\int_0^{1/c}\mathbf 1\{q<z\}\,dz\) for \(q>0\) and
Tonelli's theorem (an atom of \(\nu_T\circ Q_\alpha^{-1}\) at the
cutoff \(1/c\) does not affect the integral).  By the duality
\eqref{eq:conversion_duality}, for every \(t>0\),
\[
  \Psi_\alpha(t)-H_{\rm flip}(u)
  \le\inf_{c>0}\Bigl\{ct+\frac{2c_0}{k+1}\,c^{-k}\Bigr\}
  =(2c_0)^{\frac1{k+1}}
  \Bigl(\frac{k+1}{k}\Bigr)^{\frac{k}{k+1}}
  t^{\frac{k}{k+1}},
\]
the infimum being attained at
\(c_*=\{2kc_0/((k+1)t)\}^{1/(k+1)}\); the first term is the
charge of the slope \(c\) on the mixture budget and the second is
the additive allowance that slope requires, and balancing them
gives the H\"older power.  This is case (iii) of
Corollary~\ref{cor:regime_flip} with the explicit constant
\[
  C
  =(2c_0)^{\frac{1}{\chi_0+1}}
  \Bigl(\frac{\chi_0+1}{\chi_0}\Bigr)^{\frac{\chi_0}{\chi_0+1}},
  \qquad
  \chi=\frac{\chi_0}{\chi_0+1}.
\]
At \(t=0\) the bound is \(\Psi_\alpha(0)=H_{\rm flip}(u)\) by
\eqref{eq:profile_representation}; if \(c_0=0\) the stable
target-margin mass vanishes and \(\Psi_\alpha\equiv H_{\rm flip}(u)\).
When \(\chi_0=1\) the two displays reduce to
\(\mathcal B_\alpha(c)-H_{\rm flip}(u)\le c_0/c\) and
\(\Psi_\alpha(t)-H_{\rm flip}(u)\le2\sqrt{c_0t}\).  (Bounding the
multiplier tail instead, \(\varpi_\alpha\{A_\alpha>1/\varepsilon\}
\le2\varepsilon F(\varepsilon)\le2c_0\varepsilon^{k+1}\), and
integrating the decreasing rearrangement, gives the same exponent
with the larger constant \(((k+1)/k)(2c_0)^{1/(k+1)}\); the
frontier integral avoids the loss incurred by replacing
\(Q_\alpha\) by the cutoff on the whole near-cancellation set.)

\emph{Saturation and fitted-SVM consequences.}
The saturation cap is Theorem~\ref{thm:affine_frontier}(i).  The
fitted-SVM target bounds stated after
Corollary~\ref{cor:regime_flip} follow by evaluating
\eqref{eq:target_excess_master} with the displayed bound on
\(\Psi_\alpha\).
\end{proof}

\phantomsection\label{pf:safe_conversion}%
\begin{proof}[Proof of the case-(ii) buffer condition]
We show that \(u\le(1-\delta)\underline u_c\) implies the buffer
condition in Corollary~\ref{rem:no_flip_conversion}.
Write \(u=1-\alpha\).  Lemma~\ref{lem:flip_extinction}, under its hypotheses \(dP_{X,S}/dP_{X,T}\le\phi\) and \(\rho_-<\infty\), gives
\(u_c(x)\ge \underline u_c\) for
\(P_{X,T}\)-almost every \(x\in\Omega_-\cap\mathcal D_T\).
Since
\[
 u\le(1-\delta)\underline u_c<\underline u_c,
\]
the tie-or-flip set is \(P_{X,T}\)-null.  It remains to bound the
normalized vote from below.  On \(\Omega_+\cap\mathcal D_T\),
\[
 Q_\alpha=\alpha+u r(x)\ge\alpha=1-u
 \ge 1-\frac{u}{\underline u_c},
\]
where the last inequality uses \(\underline u_c\le1\).  On the
no-flip part of \(\Omega_-\cap\mathcal D_T\), the bounds
\(p_S/p_T\le\phi\) and
\(|m_S|/|m_T|\le\rho_-\) yield
\[
 Q_\alpha=\alpha-u r(x)
 \ge 1-u-u\phi\rho_-
 =1-\frac{u}{\underline u_c}.
\]
Therefore
\(Q_\alpha\ge1-u/\underline u_c\ge\delta\)
\(\nu_T\)-almost everywhere on \(\mathcal D_T\), which is the
buffer condition of case (ii).  At \(\alpha=1\), both
distributions in the comparison are \(P_T\), so the conclusion is
immediate.
\end{proof}

\paragraph*{Sharpness for every near-cancellation exponent}
Fix $k=\chi_0>0$ and $u=1-\alpha\in(0,1)$. Give both domains
the uniform design on $[-1,1]\times[0,1]$, and set
\[
 m_T(x_1,z)=a_0x_1,\qquad
 m_S(x_1,z)=\frac{a_0x_1}{u}\{uz^{1/k}-\alpha\},
 \qquad 0<a_0\le\min\{1,u/\alpha\}.
\]
These are valid margins: $|m_T|\le a_0\le1$, and
$|m_S|\le a_0\max\{1,\alpha/u\}\le1$.
The mixture margin and normalized vote are
$m_\alpha=a_0x_1uz^{1/k}$ and $Q_\alpha=uz^{1/k}$ on
$\mathcal D_T$. The source opposes the target on the positive-mass
set $\{x_1\ne0,\ z<\min\{1,(\alpha/u)^k\}\}$, but the mixture
preserves the target sign almost everywhere. Thus
$H_{\rm flip}(u)=0$ and $M_\alpha=\infty$.

Put $J_T:=\mathbb E_{P_{X,T}}|m_T|=a_0/2$ and
$c_0:=J_T/u^k$. Independence of the two uniform coordinates gives,
for every $\varepsilon>0$,
\[
 \nu_T\{0<Q_\alpha<\varepsilon\}
 =J_T\min\{(\varepsilon/u)^k,1\}
 \le c_0\varepsilon^k,
\]
with equality when $0<\varepsilon\le u$. Hence the condition
\eqref{eq:mixture_margin_profile} holds over its full stated range.
For $0\le r\le u$, let $g_r$ equal
$-\operatorname{sign}m_T$ on $\{0<Q_\alpha<r\}$ and
$\operatorname{sign}m_T$ elsewhere. Since the target and mixture Bayes signs agree
almost everywhere, \eqref{eq:hinge_pointwise_excess} yields
\[
 \mathcal E_T(g_r)=2c_0r^k,\qquad
 \mathcal E_\alpha(g_r)
 =2\int_{\{0<Q_\alpha<r\}}Q_\alpha\,d\nu_T
 =\frac{2kc_0}{k+1}r^{k+1}.
\]
In particular, the total mixture capacity is
$L_\alpha=2kJ_Tu/(k+1)$. For $0\le t\le L_\alpha$, choose
$r=\{(k+1)t/(2kc_0)\}^{1/(k+1)}$. This score attains the
case-(iii) upper bound proved above. For $t\ge L_\alpha$, $g_u$
attains the saturation cap $2J_T$. Therefore, for every $t\ge0$,
\begin{equation}
 \Psi_\alpha(t)
 =\min\left\{
 (2c_0)^{\frac1{k+1}}
 \left(\frac{k+1}{k}\right)^{\frac{k}{k+1}}
 t^{\frac{k}{k+1}},\;2J_T\right\}.
 \label{eq:sqrt_modulus_example}
\end{equation}
Thus, for each $k>0$, both the exponent and the displayed
coefficient are attained for the corresponding $c_0=J_T/u^k$;
no larger exponent with a finite constant can hold as $t\downarrow0$.
At $k=1$ the formula reduces to $2\sqrt{c_0t}\wedge2J_T$.
The target boundary remains the hyperplane $\{x_1=0\}$: the
nonlinear conversion comes from the source--target relation.

\section{The Weighted Statistical Oracle}\label{supp:oracle}
The general oracle first isolates the concentration and localization
arguments used in Theorem~\ref{thm:weighted_crerm}. Its entropy
specialization is the interface used by the Gaussian application.
The hinge-variance lemma and the proof of
Proposition~\ref{lem:one_acceleration} then verify the mixture-noise
inputs. A domainwise variance profile gives a further input to the
general oracle.

\subsection{Weighted-sample oracle inequalities}
\label{supp:crerm}
\label{subsec:proof_bernstein_transfer}

\emph{Citation convention for this subsection.} The argument adapts
the single-sample localization proof of
\citet{steinwart2008support} to the weighted two-sample setting;
plain references of the form Theorem~7.x, Proposition~7.x,
Lemma~7.x, Corollary~7.x, or (7.x) are to that monograph.

We adapt the localized empirical-process argument to two
independent samples with different laws and weights.  The main
change is to retain the separate target and source fluctuation
terms throughout concentration and localization.  The mixture
bound of Theorem~\ref{thm:weighted_crerm} was stated for the hinge
loss in Section~\ref{sec:general_oracle}; this subsection states
and proves it in full generality---for an arbitrary clippable
loss, a general regularizer, and a domainwise variance profile
(Theorem~\ref{thm:weighted_crerm_general} below)---and also states
and proves Theorem~\ref{thm:transfer_bernstein_svm}, its
instantiation for kernel classes, which is the input to the
Gaussian guarantees of Section~\ref{sec:gaussian}.  The auxiliary concentration results and the general oracle
come first; weighted chaining then yields the entropy specialization.
A short reading note follows the proofs.

For the proofs, write \(N=n+m\) and index the combined independent
sample as \(Z_i=(X_i^T,Y_i^T)\sim P_T\), \(i\le n\), and
\(Z_i=(X_{i-n}^S,Y_{i-n}^S)\sim P_S\), \(i>n\).  Let \(P_i\) denote
the law of \(Z_i\), and set
\[
  w_i=\frac{\alpha}{n}\quad(i\le n),
  \qquad
  w_i=\frac{1-\alpha}{m}\quad(i>n).
\]
Then \(\sum_iw_i=1\),
\(D_w=\sum_iw_i\delta_{Z_i}\),
\(\mathbb E_{D_w}g=\sum_iw_ig(Z_i)\), and
\(\mathbb E\mathbb E_{D_w}g=\mathbb E_{P_\alpha}g\).  Moreover,
\[
  w_{\max}:=\max_iw_i=N_{\rm eff}^{-1}=:b_{\alpha,n,m}.
\]
The exact weighted second moment is
\begin{equation}
  \sum_{i=1}^N w_i^2\,\mathbb{E}_{P_i}\,g^2
  =\frac{\alpha^2}{n}\mathbb E_{P_T}g^2
   +\frac{(1-\alpha)^2}{m}\mathbb E_{P_S}g^2.
  \label{eq:collapse_proof}
\end{equation}
Only when the exact block moments are unavailable do we use the
fallback upper bound
\[
  \sum_iw_i^2\mathbb E_{P_i}g^2
  \le b_{\alpha,n,m}\mathbb E_{P_\alpha}g^2.
\]

\subsubsection{Localized inputs for the main-text oracle}
\label{cmp:oracle-inputs}

For a loss \(L\) that can be clipped at \(M>0\), write
\(\bar f:=\max\{-M,\min\{f,M\}\}\). Fix a measurable Bayes
decision \(f^*_{L,P_\alpha}:\mathcal X\to[-M,M]\), and set
\[
 h_f:=L\circ f-L\circ f^*_{L,P_\alpha},\qquad
 \mathcal E_\alpha(f):=\mathbb E_{P_\alpha}h_f.
\]
Let \(\mathcal F\subset\mathcal L_0(\mathcal X)\) be equipped
with a complete, separable metric dominating pointwise convergence,
and let \(\Upsilon:\mathcal F\to[0,\infty)\) be a continuous
regularizer. Define
\[
\begin{aligned}
 r_0&:=\inf_{f\in\mathcal F}
   \{\Upsilon(f)+\mathcal E_\alpha(\bar f)\},\\
 \mathcal F_r&:=\{f\in\mathcal F:
   \Upsilon(f)+\mathcal E_\alpha(\bar f)\le r\},\qquad
 \mathcal H_r:=\{h_{\bar f}:f\in\mathcal F_r\}.
\end{aligned}
\]
The expected weighted Rademacher complexity of these localized
classes is bounded by a function
\(\varphi:[0,\infty)\to[0,\infty)\) with
\(\varphi(4r)\le2\varphi(r)\):
\begin{equation}
  \mathbb{E}_{D}\,\mathrm{Rad}_w(\mathcal{H}_r)
  := \mathbb{E}_{D}\,\mathbb{E}_{\varepsilon}
     \sup_{h \in \mathcal{H}_r}
     \Bigl|\sum_{i=1}^N w_i\,\varepsilon_i\,h(Z_i)\Bigr|
  \;\leq\; \varphi(r),
  \qquad r > r_0,
  \label{eq:weighted_rad_assumption}
\end{equation}
where the \(\varepsilon_i\) are independent Rademacher signs,
independent of \(D\). For Theorem~\ref{thm:weighted_crerm}, take
\(L\) to be the hinge loss, \(M=1\), \(\mathcal F=H\), and
\(\Upsilon(f)=\lambda\|f\|_H^2\), where \(H\) is the separable
RKHS of a measurable kernel bounded by one. These definitions then
agree with those in Section~\ref{sec:general_oracle}. The choice
\(C_0=216\) in \eqref{eq:r_condition} is verified in the proof of
Theorem~\ref{thm:weighted_crerm} below.

Write \(u=1-\alpha\) and
\(b_{\alpha,n,m}:=N_{\rm eff}^{-1}=\max_iw_i\) for the
largest-increment scale of the weighted process.  An admissible
two-block variance profile is a nondecreasing
\(\mathcal V_{\alpha,n,m}:[0,\infty)\to[0,\infty)\) with
\(r\mapsto\mathcal V_{\alpha,n,m}(r)/r^2\) nonincreasing on
\((0,\infty)\) and
\begin{equation}
  \frac{\alpha^2}{n}\mathbb E_{P_T}h_{\bar f}^2
  +\frac{u^2}{m}\mathbb E_{P_S}h_{\bar f}^2
  \le
  \mathcal V_{\alpha,n,m}
  \bigl(\mathbb E_{P_\alpha}h_{\bar f}\bigr),
  \qquad f\in\mathcal F,
  \label{eq:variance_bound_two_block}
\end{equation}
the threshold \eqref{eq:r_condition} being replaced by
\begin{equation}
  r>\max\bigl\{30\varphi(r),\,
    216B_0\varrho b_{\alpha,n,m},\,r_0\bigr\},
  \qquad
  216\varrho\,\mathcal V_{\alpha,n,m}(r)<r^2.
  \label{eq:r_condition_general}
\end{equation}
The shape conditions hold for every finite sum
\(\sum_jA_jr^{\vartheta_j}\) with \(A_j\ge0\) and
\(\vartheta_j\in[0,1]\).  The scalar bound of the main text is the
choice \(\mathcal V_{\alpha,n,m}(r)=V_\alpha
b_{\alpha,n,m}r^{\vartheta_\alpha}\), admissible by
\eqref{eq:two_block_scales} and \eqref{eq:variance_bound_alpha}.

We first record the three auxiliary results that replace their
i.i.d.\ counterparts in \citet{steinwart2008support}.

\begin{lemma}[Weighted scalar Bernstein inequality]
\label{lem:weighted_bernstein_scalar}
Let $g$ be measurable with $\|g\|_\infty \leq B'$ and
$\bar\sigma^2 := \sum_i w_i^2 \mathbb{E}_{P_i}(g - \mathbb{E}_{P_i}g)^2$.
Then, for all $\varrho > 0$, with probability not less than $1 - e^{-\varrho}$,
\[
  \mathbb{E}_{D_w} g - \mathbb{E}_{P_\alpha} g
  \;\leq\; \sqrt{2\varrho\bar\sigma^2} + \frac{2 w_{\max} B' \varrho}{3}.
\]
\end{lemma}

\begin{proof}[Proof of Lemma~\ref{lem:weighted_bernstein_scalar}]
Put \(\xi_i=w_i\{g(Z_i)-\mathbb E_{P_i}g\}\),
\(S=\sum_i\xi_i\), \(V=\bar\sigma^2\), and
\(M=2w_{\max}B'\). The variables \(\xi_i\) are independent and
centered, with \(|\xi_i|\le M\). If \(V=0\), then \(S=0\) almost
surely, so the non-strict conclusion is immediate. Otherwise put
\(c=M/3\).  We apply the moment form of Bernstein's inequality,
\citet[Theorem~2.10]{boucheron2013concentration}: if
\(\sum_i\mathbb E\xi_i^2\le v\) and
\(\sum_i\mathbb E(\xi_i)_+^q\le\tfrac{q!}{2}\,v\,c^{q-2}\) for all
integers \(q\ge3\), then
\(\mathbb P\{S\ge\sqrt{2v\varrho}+c\varrho\}\le e^{-\varrho}\) for
every \(\varrho>0\).  Both conditions hold with \(v=V\) and
\(c=M/3\): the first is the definition of \(V\), and for
\(q\ge3\),
\[
 \sum_i\mathbb E(\xi_i)_+^q
 \le M^{q-2}\sum_i\mathbb E\xi_i^2
 \le M^{q-2}V
 \le\frac{q!}{2}\Bigl(\frac M3\Bigr)^{q-2}V,
\]
because \(q!\ge2\cdot3^{q-2}\).  Hence
\[
 \mathbb P\{S\ge\sqrt{2V\varrho}+c\varrho\}
 \le e^{-\varrho},
\]
and \(c=2w_{\max}B'/3\) gives the assertion.  (For the reader's
convenience, the mechanism is the bound
\(\log\mathbb E e^{tS}\le Vt^2/\{2(1-ct)\}\) for
\(0<t<1/c\), obtained from the exponential series and
\(k!\ge2\cdot3^{k-2}\), followed by the exponential Markov
inequality at \(t=\sqrt{2\varrho}/\{\sqrt V+c\sqrt{2\varrho}\}\).
The more familiar tail bound
\(\mathbb P\{S\ge x\}\le\exp[-x^2/\{2(v+cx)\}]\),
\citet[Corollary~2.11]{boucheron2013concentration}, is a relaxation
of the same moment-generating-function bound; inverting it gives the
weaker deviation \(\sqrt{2v\varrho}+2c\varrho\), that is,
\(4w_{\max}B'\varrho/3\) in place of \(2w_{\max}B'\varrho/3\).)
The bound is one-sided; an absolute-deviation bound has failure
probability at most \(2e^{-\varrho}\).
\end{proof}

\begin{lemma}[High-probability form of the Klein--Rio inequality]
\label{lem:kr_highprob}
Let $\mathcal{G}$ be a separable Carath\'eodory set of measurable
functions with $\|g\|_\infty \leq b_{\mathcal{G}}$ and
$\mathbb{E}_{P_\alpha} g = 0$ for all $g \in \mathcal{G}$, and set
\[
  Z_w := \sup_{g\in\mathcal{G}}\ \mathbb{E}_{D_w} g,
  \qquad
  W := \sup_{g\in\mathcal{G}} \sum_{i=1}^N w_i^2
       \operatorname{Var}\bigl(g(Z_i)\bigr).
\]
Then, for all $\varrho > 0$, with probability not less than $1-e^{-\varrho}$,
\[
  Z_w \;\leq\; \frac{5}{4}\,\mathbb{E} Z_w + \sqrt{2\varrho W}
        + 14\, w_{\max}\, b_{\mathcal{G}}\,\varrho .
\]
\end{lemma}

\begin{proof}
If \(b_{\mathcal G}=0\), every \(g\in\mathcal G\) vanishes, so
\(Z_w=W=0\) and the assertion is immediate. Assume henceforth
\(b_{\mathcal G}>0\).
Since $\sum_i w_i \mathbb{E}_{P_i} g = \mathbb{E}_{P_\alpha}g = 0$, we
may write $\mathbb{E}_{D_w}g = \sum_i w_i\bigl(g(Z_i) -
\mathbb{E}_{P_i}g\bigr)$, a sum of independent, \emph{individually
centered} terms. Set $b := 2 w_{\max} b_{\mathcal{G}}$ and $s^i_g :=
b^{-1} w_i (g - \mathbb{E}_{P_i}g)$, so that $s^i_g(Z_i) \in [-1,1]$
and $\mathbb{E}\,s^i_g(Z_i) = 0$. Applying
\citet[Theorem~1.1(c)]{klein2005concentration} to $Z' := Z_w/b =
\sup_g \sum_i s^i_g(Z_i)$ gives, for all $x > 0$,
\[
  \mathbb{P}\bigl(Z' \geq \mathbb{E}Z' + x\bigr)
  \leq \exp\Bigl(-\frac{x^2}{2v + 3x}\Bigr),
  \qquad
  v := 2\,\mathbb{E}Z' + V'_N,
  \quad
  V'_N := \sup_g \sum_i \operatorname{Var}\bigl(s^i_g(Z_i)\bigr).
\]
Solving $x^2 = 2v\varrho + 3x\varrho$ and using $\sqrt{x+y} \leq \sqrt{x} + \sqrt{y}$ shows that with probability not less than $1 - e^{-\varrho}$,
\[
  Z' \leq \mathbb{E}Z' + \sqrt{2v\varrho} + 3\varrho
     \leq \mathbb{E}Z' + \sqrt{4\varrho\,\mathbb{E}Z'} + \sqrt{2\varrho V'_N}
       + 3\varrho
     \leq \frac{5}{4}\,\mathbb{E}Z' + \sqrt{2\varrho V'_N} + 7\varrho,
\]
where the last step uses $2\sqrt{\varrho\,\mathbb{E}Z'} \leq
\frac{1}{4}\mathbb{E}Z' + 4\varrho$. Multiplying by $b$ and noting
$b^2 V'_N = W$ and $7b\varrho = 14 w_{\max} b_{\mathcal{G}} \varrho$ yields
the assertion.
\end{proof}

\begin{lemma}[Symmetrization and peeling for the weighted process]
\label{lem:sym_peel}
Let $\varphi$ satisfy \eqref{eq:weighted_rad_assumption} with
$\varphi(4r) \leq 2\varphi(r)$, and for $f \in \mathcal{F}$ and
$r > r_0$ define
\[
  g_{f,r} := \frac{\mathbb{E}_{P_\alpha} h_{\bar f} - h_{\bar f}}
                  {\Upsilon(f) + \mathbb{E}_{P_\alpha} h_{\bar f} + r}.
\]
Then
$\displaystyle \mathbb{E}_{D}\,\sup_{f\in\mathcal{F}}\
\mathbb{E}_{D_w}\,g_{f,r} \;\leq\; \frac{8\varphi(r)}{r}$.
\end{lemma}

\begin{proof}
Put
\(s(f):=\Upsilon(f)+\mathbb E_{P_\alpha}h_{\bar f}\), and form the
\(4\)-adic shells
\[
\mathcal F^{(0)}=\{f:s(f)\le r\},\qquad
\mathcal F^{(k)}=\{f:4^{k-1}r<s(f)\le4^kr\},\quad k\ge1.
\]
For each index, take an independent ghost observation with law
\(P_i\). The symmetrization argument of Proposition~7.10 gives,
for every \(R>r_0\),
\[
  \mathbb{E}_{D}\sup_{f\in\mathcal{F}_R}
  \Bigl|\mathbb E_{P_\alpha}h_{\bar f}
        -\mathbb E_{D_w}h_{\bar f}\Bigr|
  \;\leq\; 2\,\mathbb{E}_D\,\mathrm{Rad}_w(\mathcal{H}_R)
  \;\leq\; 2\varphi(R).
\]
Each shell supremum's positive part is bounded by the corresponding
absolute supremum. Applying the display with \(R=4^kr\), and
iterating \(\varphi(4r)\le2\varphi(r)\), bounds the contribution of
\(\mathcal F^{(0)}\) by \(2\varphi(r)/r\), and that of
\(\mathcal F^{(k)}\), \(k\ge1\), by
\[
\frac{2\varphi(4^kr)}{(4^{k-1}+1)r}
\le
\frac{2^{k+1}}{4^{k-1}+1}\frac{\varphi(r)}r.
\]
The numerical series is bounded by
\[
2+\frac4{2}+\frac8{5}
  +\sum_{k\ge3}2^{3-k}
=7.6<8.
\]
Summing the shell bounds proves the assertion.
\end{proof}

\subsubsection{The general oracle inequality (Theorem~\ref{thm:weighted_crerm_general})}
\label{cmp:general-oracle-sec}\label{pf:crerm}

We now prove the mixture oracle using the localized inputs and
admissible variance profiles defined above.  The main-text theorem
is its hinge-loss specialization with the scalar profile
\(\mathcal V_{\alpha,n,m}(r)=V_\alpha b_{\alpha,n,m}r^{\vartheta_\alpha}\).

\begin{theorem}[General-loss mixture oracle inequality]
\label{thm:weighted_crerm_general}
Let \(L:\mathcal Y\times\mathbb R\to[0,\infty)\) be a convex,
locally Lipschitz loss, clippable at \(M>0\) and bounded by \(B\)
on \([-M,M]\). Assume that the Bayes risk under \(P_\alpha\) is
attained by a measurable decision with values in \([-M,M]\),
and use the localization quantities defined above. Let
\(\mathcal F\) and \(\Upsilon\) be as above, and fix
\(f_0\in\mathcal F\) and \(B_0\ge B\) with
\(\|L\circ f_0\|_\infty\le B_0\).  Assume the capacity bound
\eqref{eq:weighted_rad_assumption} for \(\mathcal H_r\), and let
\(\mathcal V_{\alpha,n,m}\) be an admissible variance profile
satisfying \eqref{eq:variance_bound_two_block}.  Fix
\(\epsilon\ge0\), \(\varrho>0\), and any \(r\) satisfying the
threshold \eqref{eq:r_condition_general}.  Then every measurable
\(\epsilon\)-approximate minimizer \(f_D\) of
\(f\mapsto\Upsilon(f)+\mathcal R_{L,D_w}(f)\) satisfies, with
probability not less than \(1-3e^{-\varrho}\),
\[
  \Upsilon(f_D)+\mathcal E_\alpha(\bar f_D)
  \;\le\;
  6\bigl(\Upsilon(f_0)+\mathcal E_\alpha(f_0)\bigr)+3r+3\epsilon.
\]
\end{theorem}

The proof adapts the classical localization argument of
\citet[Theorem~7.20]{steinwart2008support} to the weighted
two-sample setting. In the classical proof, the concentration step
applies Bousquet's inequality to the empirical process of localized
excess-loss functions of an i.i.d.\ sample; here the combined sample
is independent but not identically distributed, and the
corresponding step is supplied by the concentration inequality of
\citet{klein2005concentration}, which has the same shape with the
variance parameter $\sum_i w_i^2\sigma_i^2$ in place of $\sigma^2/n$
and the scale $\max_i w_i$ in place of $1/n$.  We retain these two
quantities separately: the profile in
\eqref{eq:variance_bound_two_block} controls the first, while
\(b_{\alpha,n,m}\) controls the second.

\begin{proof}[Proof of Theorem~\ref{thm:weighted_crerm_general}]
Write \(P:=P_\alpha\), \(b:=b_{\alpha,n,m}\),
\(\mathcal V:=\mathcal V_{\alpha,n,m}\), and
\(h_f:=L\circ f-L\circ f^*_{L,P_\alpha}\). It suffices that \(f_D\)
satisfies
\begin{equation}
  \Upsilon(f_D) + \mathbb{E}_{D_w} h_{\bar f_D}
  \;\leq\; \Upsilon(f_0) + \mathbb{E}_{D_w} h_{f_0} + \epsilon,
  \label{eq:crerm_property}
\end{equation}
which holds for every $\epsilon$-approximate clipped regularized empirical risk minimizer, and in particular for the
transfer SVM with $\epsilon = 0$ since clipping does not increase the
empirical risk. Exactly as in
(7.39), \eqref{eq:crerm_property} gives
\begin{equation}
  \Upsilon(f_D) + \mathbb{E}_{P} h_{\bar f_D}
  \leq \bigl(\Upsilon(f_0) + \mathbb{E}_{P} h_{f_0}\bigr)
  + \bigl(\mathbb{E}_{D_w} h_{f_0} - \mathbb{E}_{P} h_{f_0}\bigr)
  + \bigl(\mathbb{E}_{P} h_{\bar f_D} - \mathbb{E}_{D_w} h_{\bar f_D}\bigr)
  + \epsilon.
  \label{eq:basic_decomp}
\end{equation}

\emph{Step 1: the $f_0$-terms.} Split $\mathbb{E}_{D_w}h_{f_0} -
\mathbb{E}_P h_{f_0} = [\mathbb{E}_{D_w}(h_{f_0}-h_{\bar f_0}) -
\mathbb{E}_P(h_{f_0}-h_{\bar f_0})] + [\mathbb{E}_{D_w}h_{\bar f_0} -
\mathbb{E}_P h_{\bar f_0}]$. Since $0 \leq h_{f_0} - h_{\bar f_0} =
L\circ f_0 - L\circ \bar f_0 \leq B_0$, we have, for every $i$,
$\mathbb{E}_{P_i}(h_{f_0}-h_{\bar f_0})^2 \leq B_0\,
\mathbb{E}_{P_i}(h_{f_0}-h_{\bar f_0})$, so the scalar fallback after
\eqref{eq:collapse_proof} bounds the weighted variance by
$b B_0\, \mathbb{E}_{P}(h_{f_0}-h_{\bar f_0})$.
Lemma~\ref{lem:weighted_bernstein_scalar} and the elementary estimate
$\sqrt{2\varrho w_{\max} B_0\, x} \leq x + \varrho w_{\max} B_0/2$ then give,
with probability not less than $1-e^{-\varrho}$,
\begin{equation}
  \mathbb{E}_{D_w}(h_{f_0}-h_{\bar f_0})
  - \mathbb{E}_P(h_{f_0}-h_{\bar f_0})
  \;\leq\; \mathbb{E}_P(h_{f_0}-h_{\bar f_0}) + \frac{7 B_0\varrho b}{6}.
  \label{eq:f0_clip_term}
\end{equation}
For the second part, put
\(s_0:=\Upsilon(f_0)+\mathbb E_Ph_{\bar f_0}\).
The exact variance bound \eqref{eq:variance_bound_two_block} gives
weighted variance at most
\(\mathcal V(\mathbb E_Ph_{\bar f_0})\).  Monotonicity of
\(\mathcal V\), monotonicity of
\(\mathcal V(t)/t^2\), and the profile fixed point in
\eqref{eq:r_condition_general} imply
\begin{equation}
 \sqrt{2\varrho\,
   \mathcal V(\mathbb E_Ph_{\bar f_0})}
 \le \frac{s_0+r}{\sqrt{108}}
 <\frac{s_0+r}{10}.
 \label{eq:profile_comparator_root}
\end{equation}
Indeed, if \(s_0\le r\), use
\(\mathcal V(s_0)\le\mathcal V(r)\); if \(s_0>r\), use
\(\mathcal V(s_0)/s_0^2\le\mathcal V(r)/r^2\).
Lemma~\ref{lem:weighted_bernstein_scalar} therefore yields, with
probability not less than $1-e^{-\varrho}$,
\begin{equation}
  \mathbb{E}_{D_w}h_{\bar f_0} - \mathbb{E}_P h_{\bar f_0}
  \;\leq\; \frac{s_0+r}{10}
  + \frac{2B\varrho b}{3}
  \label{eq:f0_var_term}
\end{equation}
Combining, using \(h_{\bar f_0}\le h_{f_0}\),
\(s_0\le\Upsilon(f_0)+\mathbb E_Ph_{f_0}\), and
\(\frac76+\frac23=\frac{11}{6}\): with probability not less than
$1 - 2e^{-\varrho}$,
\begin{equation}
  \mathbb{E}_{D_w} h_{f_0} - \mathbb{E}_{P} h_{f_0}
  \;\leq\; \frac{11}{10}
  \bigl\{\Upsilon(f_0)+\mathbb E_Ph_{f_0}\bigr\}
  +\frac r{10}
  +\frac{11B_0\varrho b}{6}.
  \label{eq:f0_combined}
\end{equation}

\emph{Step 2: the normalized process.} Define \(g_{f,r}\) as in
Lemma~\ref{lem:sym_peel}. As in
(7.43), $\|g_{f,r}\|_\infty \leq 2B/r$
and \(\mathbb{E}_{P} g_{f,r}=0\) by construction.  Since centering
does not change the within-block variance,
\eqref{eq:variance_bound_two_block} and the two shape properties of
\(\mathcal V\) give
\[
  W := \sup_{f} \sum_i w_i^2 \operatorname{Var}\bigl(g_{f,r}(Z_i)\bigr)
  \;\leq\; \frac{\mathcal V(r)}{r^2}.
\]
Indeed, write
\(s(f)=\Upsilon(f)+\mathbb E_Ph_{\bar f}\).  The contribution of one
\(f\) is at most
\(\mathcal V(\mathbb E_Ph_{\bar f})/\{s(f)+r\}^2\).
If \(s(f)\le r\), use monotonicity of \(\mathcal V\); if \(s(f)>r\),
use monotonicity of \(\mathcal V(t)/t^2\).
Lemmas~\ref{lem:kr_highprob} and~\ref{lem:sym_peel} (with
$b_{\mathcal{G}} = 2B/r$) therefore give, with probability not less
than $1 - e^{-\varrho}$,
\begin{equation}
  \sup_{f\in\mathcal{F}} \mathbb{E}_{D_w} g_{f,r}
  \;\leq\; \frac{10\varphi(r)}{r}
  + \sqrt{\frac{2\varrho\mathcal V(r)}{r^2}}
  + \frac{28 B\varrho b}{r}
  \;=:\; \widehat S.
  \label{eq:sup_bound}
\end{equation}
We estimate $\widehat S$ using \eqref{eq:r_condition_general}.
First, $r>30\varphi(r)$ gives
$10\varphi(r)/r \leq 1/3$. Second, the profile fixed point gives
\(\sqrt{2\varrho\mathcal V(r)/r^2}<1/10\).
Third, using the bounded-increment threshold directly,
\[
  \frac{28B\varrho b}{r}
  \leq \frac{28}{216}.
\]
Hence $\widehat S \leq \frac13 + \frac1{10} + \frac{7}{54} \leq 0.564$.
Unraveling the definition of $g_{f_D,r}$, \eqref{eq:sup_bound} implies
\[
  \mathbb{E}_{P} h_{\bar f_D} - \mathbb{E}_{D_w} h_{\bar f_D}
  \;\leq\; \widehat S\,\bigl(\Upsilon(f_D)
        + \mathbb{E}_{P}h_{\bar f_D}\bigr) + \widehat S\, r .
\]
Substituting this and \eqref{eq:f0_combined} into
\eqref{eq:basic_decomp} and rearranging,
\[
  (1-\widehat S)\bigl(\Upsilon(f_D)+\mathbb{E}_{P}h_{\bar f_D}\bigr)
  \;\leq\; \frac{21}{10}
  \bigl\{\Upsilon(f_0)+\mathbb E_Ph_{f_0}\bigr\}
  +\frac r{10}
  +\frac{11B_0\varrho b}{6}+\widehat S\,r+\epsilon
\]
with probability not less than $1-3e^{-\varrho}$. Since $(1-\widehat
S)^{-1} \leq 1/0.436 \leq 2.3$, and using
\(11B_0\varrho b/6<11r/1296\) from
\eqref{eq:r_condition_general},
\begin{align*}
  &\Upsilon(f_D)+\mathbb{E}_{P}h_{\bar f_D}\\
  &\quad\leq 2.3\cdot\frac{21}{10}
    \bigl\{\Upsilon(f_0)+\mathbb E_Ph_{f_0}\bigr\}\\
  &\qquad+2.3\Bigl(\frac r{10}+\frac{11r}{1296}
    +0.564\,r\Bigr)+2.3\,\epsilon\\
  &\quad\leq 6\bigl(\Upsilon(f_0)+\mathbb{E}_{P}h_{f_0}\bigr)
    +3r+3\epsilon,
\end{align*}
because \(2.3(0.1+11/1296+0.564)<1.55<3\). This is the assertion.
\end{proof}

\subsubsection{Proof of Theorem~\ref{thm:weighted_crerm}}
\label{pf:weighted_crerm}

\begin{proof}[Proof of Theorem~\ref{thm:weighted_crerm}]
Apply Theorem~\ref{thm:weighted_crerm_general} with hinge loss,
\(M=1\), \(B=2\), \(B_0\ge2\), \(\mathcal F=H\),
\(\Upsilon(f)=\lambda\|f\|_H^2\), and
\(\epsilon=\varepsilon_{\rm opt}\). A clipped Bayes decision is
\(f^*_{L,P_\alpha}=\operatorname{sign}(m_\alpha)\).
By \eqref{eq:two_block_scales} and \eqref{eq:variance_bound_alpha}, the variance profile is
\(\mathcal V(r)=V_\alpha b r^{\vartheta_\alpha}\), where
\(b=N_{\rm eff}^{-1}\). Put \(r=r_{\rm est}\). With
\(C_0=216\),
\eqref{eq:r_condition} implies
\[
\begin{aligned}
 r&>\max\{30\varphi(r),216B_0\varrho b,r_0\},\\
 r^{2-\vartheta_\alpha}
 &>216^{2-\vartheta_\alpha}V_\alpha\varrho b
 \ge216V_\alpha\varrho b.
\end{aligned}
\]
Thus \eqref{eq:r_condition_general} holds, and the general oracle
gives \eqref{eq:oracle_mixture_bound}. On the same event,
\(\mathcal E_\alpha(\bar f_D)\le t_{\alpha,n,m}\).
Theorem~\ref{thm:affine_frontier} and monotonicity of \(\Psi_\alpha\)
give \eqref{eq:target_excess_master};
\eqref{eq:hinge_classification_calibration} gives the
classification bounds.
\end{proof}

\subsubsection{Weighted chaining and the entropy oracle (Lemma~\ref{lem:weighted_chaining}, Theorem~\ref{thm:transfer_bernstein_svm})}
\label{pf:benign}

We use a weighted version of the chaining bound in
\citet[Theorem~7.16]{steinwart2008support}.

\begin{lemma}[Anchored weighted chaining bound]
\label{lem:weighted_chaining}
Let $\mathcal{G}$ be a nonempty separable Carath\'eodory set, and
let \(B',\sigma\ge0\) satisfy $\|g\|_\infty \leq B'$ and
$\mathbb{E}_{P_\alpha} g^2 \leq \sigma^2$ for all
$g \in \mathcal{G}$. Suppose that for constants
$p \in (0,1)$ and $a'\geq0$,
\[
  \mathbb{E}_{D}\, e_i\bigl(\mathcal{G}, L_2(D_w)\bigr)
  \leq a'\, i^{-\frac{1}{2p}}, \qquad i \geq 1 .
\]
Then there are constants $C_1(p), C_2(p)$ depending only on \(p\)
such that
\[
  \mathbb{E}_{D}\,\mathrm{Rad}_w(\mathcal{G})
  \;\leq\; \sigma b_{\alpha,n,m}^{1/2}
  +\max\Bigl\{
    C_1(p)\, (a')^{p}\, \sigma^{1-p}\, b_{\alpha,n,m}^{\frac12},\;
    C_2(p)\, (a')^{\frac{2p}{1+p}}\, (B')^{\frac{1-p}{1+p}}\,
      b_{\alpha,n,m}^{\frac{1}{1+p}}
  \Bigr\}.
\]
\end{lemma}

\begin{proof}
Write \(b=b_{\alpha,n,m}\), fix \(g_0\in\mathcal G\), and set
\(\mathcal A=\mathcal G-g_0\). By the triangle inequality,
\[
 \operatorname{Rad}_w(\mathcal G)
 \le\mathbb E_\varepsilon\left|\sum_iw_i\varepsilon_i g_0(Z_i)\right|
   +\operatorname{Rad}_w(\mathcal A).
\]
The expectation of the first term is at most
\[
 \left\{\sum_iw_i^2\mathbb E_{P_i}g_0^2\right\}^{1/2}
 \le\sigma\sqrt b.
\]
Translation preserves the empirical entropy numbers, and
\(0\in\mathcal A\), \(\|h\|_\infty\le2B'\), and
\(\mathbb E_{P_\alpha}h^2\le4\sigma^2\) for \(h\in\mathcal A\).
Put
\[
 \delta_D:=\sup_{h\in\mathcal A}\|h\|_{L_2(D_w)},\qquad
 R:=\mathbb E_D\operatorname{Rad}_w(\mathcal A).
\]
Any single ball covering \(\mathcal A\) also contains zero, so
\(\delta_D\le2e_1(\mathcal A,L_2(D_w))\). If \(a'=0\), then
\(\delta_D=0\) almost surely and \(R=0\). Assume \(a'>0\), and set
\[
 s_1:=\mathbb E_D\delta_D,\qquad
 s_j:=\mathbb E_D e_j(\mathcal A,L_2(D_w)),\quad j\ge2.
\]
This is a nonincreasing sequence, since zero can be used as a
covering center. It satisfies \(s_1\le2a'\) and
\(s_j\le a'j^{-1/(2p)}\) for \(j\ge2\).

Conditionally on \(D\), the Rademacher increments have the bound
\[
\begin{aligned}
 &\mathbb E_\varepsilon\exp\!\left(
 t\sum_iw_i\varepsilon_i\{h-h'\}(Z_i)\right)\\
 &\qquad\le
 \exp\!\left(\frac{t^2}{2}\sum_iw_i^2\{h-h'\}^2(Z_i)\right)
 \le\exp\!\left(\frac{t^2b}{2}
       \|h-h'\|_{L_2(D_w)}^2\right).
\end{aligned}
\]
Dudley's bound \citep[Theorem~7.12]{steinwart2008support},
anchored at zero, therefore gives
\[
 R\le\sqrt{b\log16}\sum_{k=0}^{\infty}2^{k/2}s_{2^k}.
\]
Apply \citet[Lemma~7.15]{steinwart2008support} with entropy constant
\(A_p:=2^{1+1/(2p)}a'\): its initial condition follows from
\(s_1\le2a'=A_p2^{-1/(2p)}\). Consequently,
\begin{equation}
 R\le C_p\sqrt b\,(a')^p(\mathbb E_D\delta_D)^{1-p},
 \label{eq:weighted_chaining_radius}
\end{equation}
where \(C_p\) depends only on \(p\).

Let \(\mathcal A^2:=\{h^2:h\in\mathcal A\}\).
Weighted symmetrization, using independent ghost observations
with laws \(P_i\), and the contraction principle give
\begin{align*}
 \mathbb E_D\delta_D^2
 &\le4\sigma^2+
   \mathbb E_D\sup_{h\in\mathcal A}
      |(D_w-P_\alpha)h^2|\\
 &\le4\sigma^2+2\mathbb E_D\operatorname{Rad}_w(\mathcal A^2)\\
 &\le4\sigma^2+16B'R.
\end{align*}
Indeed, squaring is \(4B'\)-Lipschitz on \([-2B',2B']\);
for each positive weight, the map
\(t\mapsto w_i(t/w_i)^2\) is \(4B'\)-Lipschitz on
\([-2w_iB',2w_iB']\) and vanishes at zero. Applying the
coordinatewise contraction principle to \(w_i h(Z_i)\) gives
\(\operatorname{Rad}_w(\mathcal A^2)
\le8B'\operatorname{Rad}_w(\mathcal A)\).
Zero-weight terms vanish. This is the weighted form of the
argument in \citet[Theorem~A.8.4 and Corollary~A.8.5]{steinwart2008support}.
By \eqref{eq:weighted_chaining_radius} and Cauchy--Schwarz,
\[
 R\le C_p\sqrt b\,(a')^p
       (4\sigma^2+16B'R)^{(1-p)/2}.
\]
If \(4\sigma^2\ge16B'R\), then
\[
 R\le C_p8^{(1-p)/2}(a')^p\sigma^{1-p}\sqrt b.
\]
Otherwise \(R>0\), and division by \(R^{(1-p)/2}\) yields
\[
 R\le
 \bigl(C_p32^{(1-p)/2}\bigr)^{2/(1+p)}
 (a')^{2p/(1+p)}(B')^{(1-p)/(1+p)}b^{1/(1+p)}.
\]
These two bounds, together with the singleton term
\(\sigma\sqrt b\), prove the assertion.
\end{proof}

\begin{theorem}[Transfer oracle inequality for SVMs under
entropy-number decay]
\label{thm:transfer_bernstein_svm}
Let $L$ be a convex, locally Lipschitz continuous loss that can be clipped at
$M > 0$ and satisfies the supremum bound \(L(y,t)\le B\) for all \(y\) and \(t\in[-M,M]\) with $B > 0$. Let
$H$ be a separable RKHS of a measurable kernel over $\mathcal{X}$, and
assume the surrogate Bayes risk under \(P_\alpha\) is attained by a
measurable decision with values in \([-M,M]\).  Further,
  assume $P_\alpha$ satisfies the mixture variance bound
\eqref{eq:variance_bound_alpha} with $\vartheta \in [0,1]$ and $V \geq
B^{2-\vartheta}$. Assume further that there exist $p \in (0,1)$ and
$a \geq B$ such that the weighted empirical marginal
$D_{w,X} := \sum_i w_i \delta_{x_i}$ satisfies
\begin{equation}
  \mathbb{E}_{D}\, e_i\bigl(\mathrm{id} : H \to
  L_2(D_{w,X})\bigr) \;\leq\; a\, i^{-\frac{1}{2p}},
  \qquad i \geq 1.
  \label{eq:weighted_entropy}
\end{equation}
Fix $f_0 \in H$ and $B_0 \geq B$ with
$\|L \circ f_0\|_\infty \leq B_0$. Then there exists a constant
$K \geq 1$ depending only on $p$, $M$, $B$, \(|L|_{M,1}\) (the Lipschitz constant of \(L(y,\cdot)\) on \([-M,M]\)),
$\vartheta$, and $V$ such
that for all $\varrho \geq 1$ and $\lambda > 0$, the transfer SVM solution
$f_{D,\lambda} := \arg\min_{f \in H}
\lambda\|f\|_{H}^2 + \mathcal{R}_{L,D_w}(f)$ satisfies
\begin{align*}
  &\lambda\|f_{D,\lambda}\|_{H}^2
  + \mathcal{R}_{L,P_\alpha}(\bar f_{D,\lambda})
  - \mathcal{R}^*_{L,P_\alpha}\\
  &\quad\leq 12\bigl(\lambda\|f_0\|_{H}^2
    + \mathcal{R}_{L,P_\alpha}(f_0) - \mathcal{R}^*_{L,P_\alpha}\bigr)\\
  &\qquad
  + K\Bigl(\frac{a^{2p}b_{\alpha,n,m}}{\lambda^p}\Bigr)^{\!
      \frac{1}{2-p-\vartheta+\vartheta p}}\\
  &\qquad
  + 6\bigl(3600\,V\varrho b_{\alpha,n,m}\bigr)^{\!
      \frac{1}{2-\vartheta}}
  + 1296 B_0 \varrho b_{\alpha,n,m}
\end{align*}
with probability not less than $1 - 3e^{-\varrho}$.
\end{theorem}

\begin{proof}[Proof of Theorem~\ref{thm:transfer_bernstein_svm}]
Write \(b:=b_{\alpha,n,m}\). We first dispose of the case
\(a^{2p}b>\lambda^p\). The supremum bound and minimality of
\(f_{D,\lambda}\) give
\begin{align*}
  &\lambda\|f_{D,\lambda}\|_{H}^2
  +\mathcal{R}_{L,P_\alpha}(\bar f_{D,\lambda})
  -\mathcal{R}^*_{L,P_\alpha}\\
  &\quad\leq\lambda\|f_{D,\lambda}\|_{H}^2
  +\mathcal{R}_{L,D_w}(f_{D,\lambda})+B\\
  &\quad\leq\mathcal{R}_{L,D_w}(0)+B\leq2B,
\end{align*}
and since $2 - p - \vartheta + \vartheta p \geq 1$ implies
$(a^{2p}b/\lambda^p)^{1/(2-p-\vartheta+\vartheta p)} \geq 1$, the
assertion holds whenever $K \geq 2B$.

Assume now \(a^{2p}b\leq\lambda^p\). We apply
Theorem~\ref{thm:weighted_crerm_general} with $\Upsilon(f) :=
\lambda\|f\|^2_{H}$, $\mathcal{F} := H$, and
$\epsilon = 0$, using the scalar variance bound
\eqref{eq:variance_bound_alpha} (equivalently the profile
\(\mathcal V_{\alpha,n,m}(r)=Vbr^\vartheta\), admissible by
\eqref{eq:two_block_scales});
\eqref{eq:crerm_property} holds since clipping does not
increase the empirical risk. It remains to produce an admissible
$\varphi$. For $r > r_0$ we have $\mathcal{F}_r \subseteq
(r/\lambda)^{1/2} B_{H}$, so by the local Lipschitz
continuity of $L$, the contraction property of clipping, and the
scaling of entropy numbers,
\[
  \mathbb{E}_{D}\, e_i\bigl(\mathcal{H}_r, L_2(D_w)\bigr)
  \;\leq\; 2\,|L|_{M,1}\Bigl(\frac{r}{\lambda}\Bigr)^{1/2}
  a\, i^{-\frac{1}{2p}},
\]
exactly as in the proof of
Theorem~7.23, with $L_2(D)$ replaced by
$L_2(D_w)$ and \eqref{eq:weighted_entropy} in place of (7.48) there.
Moreover, $\|h_{\bar f}\|_\infty \leq B$ and, by the variance bound,
$\mathbb{E}_{P_\alpha} h_{\bar f}^2 \leq V
(\mathbb{E}_{P_\alpha}h_{\bar f})^{\vartheta} \leq V r^{\vartheta}$ for
$f \in \mathcal{F}_r$. Lemma~\ref{lem:weighted_chaining} applied with
$a' := 2|L|_{M,1}(r/\lambda)^{1/2}a$, $\sigma^2 := Vr^{\vartheta}$,
$B' := B$ therefore shows that
\eqref{eq:weighted_rad_assumption} holds with
\begin{align*}
  \varphi(r)&:=V^{1/2}r^{\vartheta/2}b^{1/2}\\
  &\quad+\max\left\{\begin{aligned}
    &C_1(p)\,2^p|L|^p_{M,1}\,a^p
      \Bigl(\frac{r}{\lambda}\Bigr)^{\frac p2}
      \bigl(Vr^{\vartheta}\bigr)^{\frac{1-p}{2}}
      b^{\frac12},\\
    &C_2(p)\,\bigl(2^p|L|^p_{M,1}a^p\bigr)^{\frac{2}{1+p}}
      \Bigl(\frac{r}{\lambda}\Bigr)^{\frac{p}{1+p}}
      B^{\frac{1-p}{1+p}}\,b^{\frac{1}{1+p}}
  \end{aligned}\right\},
\end{align*}
which satisfies $\varphi(4r) \leq 2\varphi(r)$ since the exponents of
$r$ are $\frac{p+\vartheta(1-p)}{2} \leq \frac12$ and
$\frac{p}{1+p} \leq \frac12$.  We first make the maximum in
\(\varphi\) at most \(r/60\).  For the first branch, rearranging shows
this holds whenever
\[
  r \;\geq\;
  \Bigl(60\,C_1(p)\,2^p|L|^p_{M,1}V^{\frac{1-p}{2}}\Bigr)^{
    \frac{2}{2-p-\vartheta+\vartheta p}}
  \Bigl(\frac{a^{2p}b}{\lambda^p}\Bigr)^{
    \frac{1}{2-p-\vartheta+\vartheta p}} ,
\]
using $2 - p - \vartheta(1-p) = 2-p-\vartheta+\vartheta p$. For the
second branch it holds whenever
\[
  r\geq60^{1+p}C_2(p)^{1+p}2^{2p}|L|^{2p}_{M,1}
  B^{1-p}\frac{a^{2p}b}{\lambda^p}
\]
\[
  \left(\leq60^{1+p}C_2(p)^{1+p}2^{2p}|L|^{2p}_{M,1}B^{1-p}
  \Bigl(\frac{a^{2p}b}{\lambda^p}\Bigr)^{
    \frac{1}{2-p-\vartheta+\vartheta p}}\right),
\]
where the parenthetical inequality uses \(a^{2p}b\leq\lambda^p\) and
$2-p-\vartheta+\vartheta p \geq 1$.  Hence the maximum is at most
\(r/60\) whenever
\[
  r\geq\widetilde K
  \Bigl(\frac{a^{2p}b}{\lambda^p}\Bigr)^{
    \frac{1}{2-p-\vartheta+\vartheta p}},
\]
\[
  \widetilde K:=\max\left\{\begin{aligned}
  &\bigl(60C_1(p)2^p|L|^p_{M,1}V^{\frac{1-p}{2}}\bigr)^{
    \frac{2}{2-p-\vartheta+\vartheta p}},\\
  &60^{1+p}C_2(p)^{1+p}2^{2p}|L|^{2p}_{M,1}B^{1-p}
  \end{aligned}\right\}.
\]
The anchored term is at most \(r/60\) whenever
\[
  r\geq\left(3600Vb\right)^{1/(2-\vartheta)}.
\]
Under the last two displays, \(30\varphi(r)\le r\).
Let
\[
  r_1 := 2\max\Bigl\{
    \widetilde K\Bigl(\frac{a^{2p}b}{\lambda^p}\Bigr)^{
      \frac{1}{2-p-\vartheta+\vartheta p}},\;
    \bigl(3600V\varrho b\bigr)^{\frac{1}{2-\vartheta}},\;
    216B_0\varrho b,\; r_0\Bigr\}
\]
so \(r_1\) strictly exceeds every threshold in
\eqref{eq:r_condition_general}.  Using
\[
3r_1\leq
 6\widetilde K(a^{2p}b/\lambda^p)^{1/(2-p-\vartheta+\vartheta p)}
 +6(3600V\varrho b)^{1/(2-\vartheta)}
 +1296B_0\varrho b+6r_0
\]
together with
\[
r_0 \leq \lambda\|f_0\|^2_{H}
+\mathcal{R}_{L,P_\alpha}(f_0)-\mathcal{R}^*_{L,P_\alpha},
\]
yields the assertion with
\(K := \max\{1,\,6\widetilde K,\,2B\}\).
\end{proof}

\paragraph*{Reading the oracle bounds}
In the RKHS setting, the comparator contribution in
Theorem~\ref{thm:weighted_crerm_general} is
\(6\{\lambda\|f_0\|_H^2+\mathcal E_\alpha(f_0)\}\), and
clipping gives
\(r_0\le\lambda\|f_0\|_H^2+\mathcal E_\alpha(f_0)\).
For Gaussian kernels, Lemma~\ref{lem:gaussian_mixture_approximation}
supplies a bounded comparator, while
Lemma~\ref{lem:weighted_empirical_entropy} verifies
\eqref{eq:weighted_entropy} under the stated marginal tail
conditions. The resulting rates are recorded in
Corollary~\ref{thm:main_target_bound} of the main text and
Corollary~\ref{cor:transfer_rates} of this supplement. The
domainwise variance profile in
Section~\ref{pf:two_block_hinge_variance} provides another input
to Theorem~\ref{thm:weighted_crerm_general}.

\subsection{Hinge variance (Lemma~\ref{lem:variance-bound-hinge})}\label{supp:faster}

The following distribution-level bound converts a small-margin
condition into a variance bound for the clipped hinge excess.

\begin{lemma}[Variance bound for the hinge loss]
\label{lem:variance-bound-hinge}
Let $P$ be a distribution on $\mathcal{X} \times \{-1, +1\}$ with
finite noise exponent $q \in [0, \infty)$, i.e.\ there exists a
  constant $c_q > 0$ such that for all $t > 0$,
\[
P_X\bigl(\{x \in \mathcal{X} : |2\eta(x) - 1| < t\}\bigr) \leq c_q\, t^q,
\]
where $\eta(x) := P(Y = 1 \mid X = x)$. Let $L$ denote the hinge
loss, $L(y, t) = \max\{0, 1 - yt\}$, let
$f^*_{L,P} : \mathcal{X} \to [-1, 1]$ be a fixed Bayes decision
function for $L$, and write
$\bar f:=\max\{-1,\min\{f,1\}\}$. Then for every measurable
$f : \mathcal{X} \to \mathbb{R}$,
\begin{equation}
\mathbb{E}_P\bigl(L \circ \bar f - L \circ f^*_{L,P}\bigr)^2
\;\leq\; 6\, c_q^{1/(q+1)}
\bigl(\mathbb{E}_P(L \circ \bar f
- L \circ f^*_{L,P})\bigr)^{q/(q+1)}.
\label{eq:variance-bound}
\end{equation}
If instead $P_X(|2\eta-1|<t_0)=0$ for some $t_0>0$, then
\[
  \mathbb E_P\bigl(L\circ\bar f-L\circ f^*_{L,P}\bigr)^2
  \leq \frac{2}{t_0}\,
  \mathbb E_P\bigl(L\circ\bar f-L\circ f^*_{L,P}\bigr).
\]
\end{lemma}

\begin{proof}[Proof of Lemma~\ref{lem:variance-bound-hinge}]
Write \(m=2\eta-1\), \(d=|\bar f-f^*_{L,P}|\le2\), and
\(E=\mathbb E_P(L\circ\bar f-L\circ f^*_{L,P})\).
Linearity of the hinge loss on \([-1,1]\) gives
\(L\circ\bar f-L\circ f^*_{L,P}=Y(f^*_{L,P}-\bar f)\) and
\(E=\int |m|d\,dP_X\); hence
\begin{equation}
 \mathbb E_P(L\circ\bar f-L\circ f^*_{L,P})^2
 =\int d^2\,dP_X.
 \label{eq:variance-as-L2}
\end{equation}
For every \(s>0\), the region \(\{|m|\ge s\}\) contributes at
most \(2E/s\), since \(d^2\le2d\le2|m|d/s\), whereas
\begin{equation}
 \int_{\{|m|<s\}}d^2\,dP_X
 \le4P_X(|m|<s)\le4c_qs^q.
 \label{eq:noisy-bound}
\end{equation}
Thus
\begin{equation}
 \mathbb E_P(L\circ\bar f-L\circ f^*_{L,P})^2
 \le\frac{2E}{s}+4c_qs^q.
 \label{eq:combined}
\end{equation}
For \(q>0\) and \(E>0\), minimizing at
\(s^*=(E/(2qc_q))^{1/(q+1)}\) gives the factor
\[
 2(q+1)2^{1/(q+1)}q^{-q/(q+1)}
 c_q^{1/(q+1)}E^{q/(q+1)}.
\]
The numerical coefficient is at most \(6\), with equality at
\(q=1/2\); if \(E=0\), let \(s\downarrow0\) in
\eqref{eq:combined}. For \(q=0\), the noise condition forces
\(c_q\ge1\), and the left-hand side is at most \(4\le6c_q\),
using the convention \(E^0=1\), including at \(E=0\).
Finally, under the hard-gap assumption, choose \(s=t_0\) in
the same split: the small-margin region has zero mass, leaving
only \(2E/t_0\).
\end{proof}

\subsection{Proof of Proposition~\ref{lem:one_acceleration}}
\label{supp:mixture_exponents}\label{pf:prop31}

\begin{proof}[Proof of Proposition~\ref{lem:one_acceleration}]
Put \(u=1-\alpha>0\).  We first establish the agreement-region
estimate shared by all three cases.

A direct calculation from the definition of $P_\alpha$ as a mixture yields
\[
m_\alpha(x) = \frac{\alpha p_T(x) m_T(x) + (1 - \alpha) p_S(x) m_S(x)}{\alpha p_T(x) + (1 - \alpha) p_S(x)},
\]
so $m_\alpha$ is the density-weighted average of $m_T$ and $m_S$.

We begin on $\Omega_+$, where the two signed votes do not oppose one
another.  Consequently absolute values pass through their sum
\(\mu\)-almost everywhere, including when one vote is zero:
\[
|m_\alpha(x)| = \frac{\alpha p_T(x) |m_T(x)| + (1 - \alpha) p_S(x) |m_S(x)|}{\alpha p_T(x) + (1 - \alpha) p_S(x)}.
\]
Discarding the non-negative term $\alpha p_T(x) |m_T(x)|$ and applying (A1) to the remaining one,
\[
|m_\alpha(x)| \geq \frac{(1 - \alpha) p_S(x) \cdot c_\zeta \, |m_T(x)|^\zeta}{\alpha p_T(x) + (1 - \alpha) p_S(x)} = c_\zeta \cdot \frac{(1 - \alpha) (p_S/p_T)(x)}{\alpha + (1 - \alpha) (p_S/p_T)(x)} \cdot |m_T(x)|^\zeta.
\]
The function $r \mapsto (1 - \alpha) r / (\alpha + (1 - \alpha) r)$ is increasing in $r > 0$, so the lower bound (A2) on $p_S / p_T$ gives
\[
\frac{(1 - \alpha) (p_S/p_T)(x)}{\alpha + (1 - \alpha) (p_S/p_T)(x)} \geq \frac{(1 - \alpha) \phi^{-1}}{\alpha + (1 - \alpha) \phi^{-1}} = \frac{1 - \alpha}{\alpha\phi + (1 - \alpha)}.
\]
Setting
\(
\tilde{c}_\alpha := \frac{(1 - \alpha)\, c_\zeta}{\alpha\phi + (1 - \alpha)},
\)
for all $x \in \Omega_+$ we obtain the pointwise bound
$|m_\alpha(x)| \geq \tilde{c}_\alpha \, |m_T(x)|^\zeta$,
which translates into the level-set inclusion, for every $t > 0$,
\begin{equation}
\{x \in \Omega_+ : |m_\alpha(x)| < t\} \subseteq \{x \in \Omega_+ : |m_T(x)| < (t/\tilde{c}_\alpha)^{1/\zeta}\}. \label{eq:level-set-inclusion}
\end{equation}

The agreement-region contribution now follows. By the upper half of
(A2), $p_S \leq \phi p_T$ on the common support, so $P_{X,\alpha}(B)
= \alpha P_{X,T}(B) + (1 - \alpha) P_{X,S}(B) \leq \bigl(\alpha + (1
- \alpha) \phi\bigr) P_{X,T}(B)$ for every measurable $B$. Applying
this, \eqref{eq:level-set-inclusion}, and the noise exponent of
$P_T$,
\[
P_{X,\alpha}\bigl(\Omega_+ \cap \{|m_\alpha| < t\}\bigr)
\leq \bigl(\alpha + (1 - \alpha) \phi\bigr)\, c_{q,T} \,
(t/\tilde{c}_\alpha)^{q_T / \zeta}
\leq c_1\, (1-\alpha)^{-q_T/\zeta}\, t^{q_T/\zeta},
\]
with $c_1 := (1+\phi)\, c_{q,T}\, (\phi/c_\zeta)^{q_T/\zeta}$, using
$\tilde{c}_\alpha \geq (1-\alpha)\, c_\zeta/\phi$.

\emph{Case (i): pure reinforcement.}
When \(P_{X,\alpha}(\Omega_-)=0\), the agreement estimate is
the whole mixture small-margin bound, with constant
\(c_1u^{-q_T/\zeta}\).

\emph{Case (ii): buffered opposition.}
Let
\(u=1-\alpha\leq(1-\delta)\underline u_c\).  On \(\Omega_-\) the
two votes are opposed, and by (A2)
and the definition of $\rho_-$ (which (A3) bounds by $1/c_T$), $(1 -
\alpha)\, p_S |m_S| \leq u\,\phi\rho_-\, p_T |m_T|$, so the target
wins with margin:
\[
p_\alpha |m_\alpha|
= \bigl|\alpha p_T m_T + (1-\alpha) p_S m_S\bigr|
\;\geq\; \alpha\, p_T|m_T| - u\,\phi\rho_-\, p_T|m_T|
= \bigl(1 - u/\underline u_c\bigr)\, p_T |m_T| .
\]
Together with $p_T/p_\alpha \geq \phi^{-1}$ and (A3), $|m_\alpha|
\geq \phi^{-1}(1-u/\underline u_c)\,c_T
\geq \phi^{-1}\delta\, c_T =: \kappa_T$
on $\Omega_-$. Hence $\Omega_- \cap \{|m_\alpha| < t\}$ is empty for
$t < \kappa_T$, while for $t \geq \kappa_T$ its $P_{X,\alpha}$-mass is at most
$1 \leq (t/\kappa_T)^{q_T/\zeta}$; either way the contribution is at most
$(\phi/(\delta c_T))^{q_T/\zeta}\, t^{q_T/\zeta}$.
Adding this to the agreement contribution gives the second
constant in \eqref{eq:mixture_noise_constants}.

\emph{Case (iii): tie-mass control.}
The disagreement-region bound is precisely the hypothesis
\eqref{eq:tie_mass}, contributing \(c_2t^{q_T/\zeta}\).
Adding the agreement contribution gives the third constant in
\eqref{eq:mixture_noise_constants}. In all three cases, letting
\(t\downarrow z\) gives the stated bound with \(\lvert m_\alpha\rvert\le z\).

\phantomsection\label{cmp:mixture-noise-constants}%
\emph{Variance pair and exponent.}
In each of the three cases, apply
Lemma~\ref{lem:variance-bound-hinge} under \(P_\alpha\), with
\(q_\alpha=q_T/\zeta\) and the constant just proved.  This gives
\[
  \vartheta_\alpha=\frac{q_\alpha}{q_\alpha+1},
  \qquad
  V_\alpha=6c_{q,\alpha}^{1/(q_\alpha+1)},
\]
which is \eqref{eq:acceleration_variance_pair}.  Substitution into
\eqref{eq:r_condition} gives \eqref{eq:noise_exponent_main}, and
\begin{equation}
  \frac{q_\alpha+1}{q_\alpha+2}
  -\frac{q_T+1}{q_T+2}
  =\frac{q_T(1-\zeta)}{(q_T+2\zeta)(q_T+2)}.
  \label{eq:noise_exponent_gain}
\end{equation}
Thus the noise-exponent improvement is strict exactly when
\(q_T>0\) and \(0<\zeta<1\).  For \(q_\alpha>0\), letting
\(t\downarrow0\) in the mixture noise bound excludes exact ties
of positive mixture-design mass, while case (iii) permits strict
flips.  The constants retain their stated dependence on \(u\)
and the buffer \(\delta\).
\end{proof}

\paragraph*{Without dominance}
\phantomsection\label{pf:noise_floor}%
Suppose \(P_T\) satisfies {\rm(T)} with exponent \(q_T\) and constant
\(c_{q,T}\), and \(dP_{X,S}/dP_{X,T}\le\phi\) holds
\(P_{X,T}\)-a.e.\ (the upper bound in {\rm(A2)}).  Fix
\(u=1-\alpha\in(0,1)\).  If \(P_{X,\alpha}(\Omega_-)=0\), then for
all \(t>0\),
\[
  P_{X,\alpha}(|m_\alpha|<t)
  \le
  c_{q,T}\,
  \frac{(\alpha+u\phi)^{q_T+1}}{\alpha^{q_T}}\,
  t^{q_T},
\]
so \(q_T\) remains an admissible mixture Tsybakov exponent; on
\(\alpha\in[1/2,1]\) the constant is at most
\(2^{q_T}\phi^{q_T+1}c_{q,T}\).  Under {\rm(A3)} and
\(0<u\le(1-\delta)\underline u_c\), the disagreement region adds the
constant \((\phi/(\delta c_T))^{q_T}\), as in
case (ii) of the proof above.

On \(\Omega_+\) the votes add:
\(p_\alpha|m_\alpha|=\alpha p_T|m_T|+u\,p_S|m_S|\ge\alpha p_T|m_T|\),
while \(p_\alpha\le(\alpha+u\phi)p_T\); hence
\(\Omega_+\cap\{|m_\alpha|<t\}\subseteq
\{|m_T|<\tfrac{\alpha+u\phi}{\alpha}\,t\}\)
up to a \(P_{X,\alpha}\)-null set.  As measures,
\(P_{X,\alpha}\le(\alpha+u\phi)P_{X,T}\), so {\rm(T)} gives the
display.  The buffered disagreement case reuses the margin-gap calculation in
case (ii) above with \(q_T\) in
place of \(q_T/\zeta\).

\paragraph*{Interpretation and posterior drift}
For \(q_T>0\), decreasing \(\zeta\) to zero gives
\(q_\alpha\to\infty\) and \(\vartheta_\alpha\to1\), whereas
\(\zeta\to1\) removes the exponent gain. The limit
\(|m_T|^\zeta\to1\) is pointwise away from target ties and uniform
on each \(\Omega_+\cap\{|m_T|\ge\tau\}\), \(\tau>0\).
Under \(m=2\eta-1\), (A1) on \(\Omega_+\) is precisely
\(|\eta_S-1/2|\ge2^{\zeta-1}c_\zeta|\eta_T-1/2|^\zeta\),
the relative-signal condition of
\citet[Definition~1]{cai2021transfer} with \(P=S\), \(Q=T\),
\(\gamma=\zeta\), and \(C_\gamma=2^{\zeta-1}c_\zeta\).
Writing their smoothness exponent as \(\nu\), their effective
source contribution \(m^{(2\nu+d)/(2\zeta\nu+d)}\) increases as
\(\zeta\) decreases, in the same qualitative direction, under
different estimator and model assumptions. Here (A1) controls only
relative magnitudes on the agreement region: it is incomparable
with Bayes-identical transfer at sublinear exponents and weaker
than uniform informativeness \(|m_S|\ge c\). Its role is to bound
the mixture's small-margin mass; the classifier benchmark in
Section~\ref{supp:reinforcement} is analyzed separately.

\subsection{A domainwise variance profile under reinforcement}
\label{pf:two_block_hinge_variance}

Let \(L\) be the hinge loss clipped at \(M=1\).  Suppose
\(P_{X,\alpha}(\Omega_-)=0\), \(P_T\) satisfies {\rm(T)} with
exponent \(q_T\) and constant \(c_{q,T}\), and {\rm(A1)} and
{\rm(A2)} hold.

Fix a common
\(\{-1,1\}\)-valued hinge Bayes representative \(f^*\) for
\(P_T\), \(P_S\), and \(P_\alpha\), write
\[
  \mathcal E^*_D(f):=\mathbb E_{P_D}
  (L\circ\bar f-L\circ f^*),
  \qquad D\in\{T,S\}
\]
(the star records the common representative), and set
\[
  q_S:=\frac{q_T}{\zeta},
  \quad
  c_{q,S}:=\phi c_{q,T}c_\zeta^{-q_T/\zeta},
  \quad
  \vartheta_T:=\frac{q_T}{q_T+1},
  \quad
  \vartheta_S:=\frac{q_T}{q_T+\zeta},
\]
\[
  V_T:=6c_{q,T}^{1/(q_T+1)},
  \qquad
  V_S:=6c_{q,S}^{1/(q_S+1)}.
\]
The source satisfies {\rm(T)} with exponent \(q_S\) and constant
\(c_{q,S}\).  For every measurable \(f\), with
\(h_{\bar f}:=L\circ\bar f-L\circ f^*\) and
\(r_f:=\alpha \mathcal E^*_T(f)+u\mathcal E^*_S(f)\),
\begin{align}
 &\frac{\alpha^2}{n}\mathbb E_{P_T}h_{\bar f}^2+
  \frac{u^2}{m}\mathbb E_{P_S}h_{\bar f}^2
  \notag\\
 &\quad\le
  \frac{\alpha^2V_T}{n}\mathcal E^*_T(f)^{\vartheta_T}
  +\frac{u^2V_S}{m}\mathcal E^*_S(f)^{\vartheta_S}
  \notag\\
 &\quad\le
  \underbrace{\frac{V_T\alpha^{2-\vartheta_T}}{n}
    r_f^{\vartheta_T}
  +\frac{V_Su^{2-\vartheta_S}}{m}
    r_f^{\vartheta_S}}_{
    \mathcal V^{+}_{\alpha,n,m}(r_f)}.
  \label{eq:two_block_hinge_variance}
\end{align}
Thus \(\mathcal V^{+}_{\alpha,n,m}\) is an admissible profile for
Theorem~\ref{thm:weighted_crerm_general}.

\begin{proof}[Verification of the profile]
With no opposition, choose the common Bayes sign preferred by
any nonzero margin, and either sign when both vanish. By
{\rm(A1)}, {\rm(A2)}, and {\rm(T)}, for every \(t>0\),
\[
 P_{X,S}(|m_S|<t)
 \le\phi P_{X,T}\bigl(|m_T|<(t/c_\zeta)^{1/\zeta}\bigr)
 \le c_{q,S}t^{q_S}.
\]
Lemma~\ref{lem:variance-bound-hinge}, applied separately under
\(P_T\) and \(P_S\), now gives the first inequality in
\eqref{eq:two_block_hinge_variance}. For the second, both domain
excesses are nonnegative, so their mixture decomposition yields
\(\mathcal E_T^*(f)\le r_f/\alpha\) and
\(\mathcal E_S^*(f)\le r_f/u\) for each present block; an absent
block contributes zero. Each exponent lies in \([0,1]\), so the
profile is nondecreasing and its ratio to \(r^2\) is
nonincreasing on \((0,\infty)\). At exponent zero, use
\(r^0=1\), including at \(r=0\).
\end{proof}

For \(q_T>0\), the two margin conditions include the zero-margin
sets and imply \(P_{X,T}(m_T=0)=P_{X,S}(m_S=0)=0\).

\section{Gaussian Inputs and the Proof of Corollary~\ref{thm:main_target_bound}}
\label{supp:gaussian}
A bounded comparator and weighted empirical entropy specialize
Theorem~\ref{thm:transfer_bernstein_svm} to the Gaussian mixture bound in
Theorem~\ref{thm:gaussian_transfer_estimation}. Combining that bound with
Theorem~\ref{thm:affine_frontier} proves
Corollary~\ref{thm:main_target_bound}.

\subsection{Mixture conditions and a Gaussian comparator (Lemma~\ref{lem:gaussian_mixture_approximation})}
\label{supp:gauss_detail}

Throughout this subsection,
$k_\gamma(x,x') := \exp\bigl(-\|x-x'\|_2^2/\gamma^2\bigr)$ denotes
the Gaussian RBF kernel with width parameter $\gamma > 0$, and
$H_\gamma$ its RKHS on $\mathcal X\subset\mathbb R^d$.
We use the full-space Gaussian RKHS for the analytic bounds.
By the kernel-restriction property, every function of
$H_\gamma(\mathcal X)$ has a full-space extension of the same norm,
and restricting a full-space function cannot increase its norm
\citep[Section~4.2]{steinwart2008support}. Thus the unit balls
have identical images on samples in $\mathcal X$.
This section supplies the
kernel-specific inputs behind the Gaussian guarantees of
Section~\ref{sec:gaussian}: the mixture conditions in detail, the
entropy and approximation estimates, and the balancing exponents
of Corollary~\ref{cor:transfer_rates}.

For fixed \(\alpha\), let
\[
  f_{D,\lambda,\gamma}
  \in\arg\min_{f\in H_\gamma}
  \left\{
  \lambda\|f\|_{H_\gamma}^2
  +
  \alpha\mathcal R_{L,D_T^n}(f)
  +(1-\alpha)\mathcal R_{L,D_S^m}(f)
  \right\},
\]
and write \(\bar f_{D,\lambda,\gamma}\) for its clipping to
\([-1,1]\).  We suppress \(\alpha\) from the subscript only to keep
the notation readable.

Corollaries~\ref{thm:main_target_bound} and~\ref{cor:transfer_rates}
use three properties of the training mixture
\(P_\alpha\).  First, \(P_\alpha\) satisfies the mixture analog of
condition {\rm(T)} of Section~\ref{subsec:source_dominance},
with exponent \(q_\alpha\in[0,\infty)\) and constant
\(c_{q,\alpha}\).
Second, under the convention \(\operatorname{sign}(0)=1\), define
\[
  \mathcal X_{\alpha,1}:=\{m_\alpha\ge0\},
  \qquad
  \mathcal X_{\alpha,-1}:=\{m_\alpha<0\},
\]
and let \(\Delta_\alpha(x)\) be the distance from \(x\) to the
opposite decision class, with distance to the empty set defined as
\(\infty\).  The margin-noise condition is
\[
  \int_{\{\Delta_\alpha\le t\}}|m_\alpha|\,dP_{X,\alpha}
  \le c_{\beta,\alpha}t^{\beta_\alpha},
  \qquad t>0,
\]
for some \(\beta_\alpha>0\).  Third, the marginal tail condition is
\[
  P_{X,\alpha}\{\|X\|_2>R\}
  \le c_{\tau,\alpha}R^{-\tau_\alpha},
  \qquad R>0,
\]
for some \(0<\tau_\alpha<\infty\). The convention
\(\tau_\alpha=\infty\) means
\(P_{X,\alpha}\{\|X\|_2>R_\alpha\}=0\) for a fixed
\(R_\alpha<\infty\). Expressions containing \(\tau_\alpha\) are
then read by continuity, for example
\(\lambda^{\tau_\alpha/(\tau_\alpha+d)}
\gamma^{-d\tau_\alpha/(\tau_\alpha+d)}=\lambda\gamma^{-d}\).  The
tail condition is inherited by the mixture with
\(\tau_\alpha=\tau_T\wedge\tau_S\), after enlarging the constant;
the margin-noise condition is a genuine condition on the mixture
decision boundary.

Without loss of generality, take \(c_{\beta,\alpha}\ge1\) and
\(c_{\tau,\alpha}\ge1\), by enlarging these upper-bound constants.

For the capacity input, let $0<\gamma\le1$. For every \(\varepsilon>0\) and
\(d/(d+\tau_\alpha)<p<1\), the weighted empirical entropy condition
\eqref{eq:weighted_entropy} holds for \(H_\gamma\) with
\begin{equation}
  a_\gamma
  =
  c_{\varepsilon,p}\,
  \gamma^{-\frac{(1-p)(1+\varepsilon)d}{2p}},
  \qquad c_{\varepsilon,p}\ge2.
  \label{eq:gaussian_entropy_constant}
\end{equation}
The entropy constant can always be enlarged to ensure
$a_\gamma\ge B=2$, as required by the hinge-loss oracle.
For finite \(\tau_\alpha\) and a fixed \(\alpha\in(0,1)\), the mixture
tail condition implies
\[
 P_{X,T}(\|X\|_2>R)\le c_{\tau,\alpha}\alpha^{-1}R^{-\tau_\alpha},
 \qquad
 P_{X,S}(\|X\|_2>R)\le c_{\tau,\alpha}(1-\alpha)^{-1}R^{-\tau_\alpha}.
\]
\begin{sloppypar}
For bounded mixture support, both marginals have the same support
bound. Thus Lemma~\ref{lem:weighted_empirical_entropy} applies to both blocks;
its constants may depend on this fixed weight. At \(\alpha=1\),
only the target block is used.
This is the only place where the Gaussian width enters the capacity
calculation; the weighted two-block entropy proof is given in
Section~\ref{sec:appB}.
\end{sloppypar}
For approximation, the oracle requires a bounded comparator in
$H_\gamma$, supplied directly by the following lemma.

\begin{lemma}[A bounded Gaussian comparator for the training mixture]
\label{lem:gaussian_mixture_approximation}
Suppose $P_\alpha$ satisfies the margin-noise and tail conditions
above. For every $\lambda,\gamma>0$, there exists $g^*\in H_\gamma$
with $\|g^*\|_\infty\le1$ such that
\begin{equation}
\begin{aligned}
 &\lambda\|g^*\|_{H_\gamma}^2+\mathcal E_\alpha(g^*)\\
 &\quad\le C_{d,\tau_\alpha}
 \lambda^{\frac{\tau_\alpha}{\tau_\alpha+d}}
 \gamma^{-\frac{d\tau_\alpha}{\tau_\alpha+d}}
 +C_{d,\beta_\alpha}c_{\beta,\alpha}\gamma^{\beta_\alpha}.
\end{aligned}
\label{eq:apalpha_bound}
\end{equation}
The first constant may depend on $c_{\tau,\alpha}$, while the
second depends only on $d,\beta_\alpha$. For $\tau_\alpha=\infty$,
the first term is $C_{d,R_\alpha}\lambda\gamma^{-d}$, where
$R_\alpha$ is a support radius. The proof in
Section~\ref{supp:apalpha} gives explicit constants and covers
one-class distributions as well.
\end{lemma}

Theorem~\ref{thm:gaussian_transfer_estimation} combines this
comparator and \eqref{eq:gaussian_entropy_constant} with the
entropy-based oracle of Theorem~\ref{thm:transfer_bernstein_svm}.
Its right side is the sum of the approximation and estimation
blocks \eqref{eq:gaussian_approximation_block}--\eqref{eq:gaussian_estimation_block}.
Theorem~\ref{thm:affine_frontier} then gives the target bound;
a separate proof of Corollary~\ref{thm:main_target_bound} appears
in Section~\ref{pf:main_target_bound}.

\subsection{Proof of Lemma~\ref{lem:gaussian_mixture_approximation}}
\label{supp:apalpha}
\phantomsection\label{pf:apalpha_bound}

\begin{proof}[Proof of Lemma~\ref{lem:gaussian_mixture_approximation}]
Write $P=P_\alpha$, $m=m_\alpha$, $s=\operatorname{sign}(m)$,
$\beta=\beta_\alpha$, $\tau=\tau_\alpha$, and use the corresponding
margin and tail constants without the subscript $\alpha$.
Write $\Delta=\Delta_\alpha$, let $B$ be the closed unit Euclidean
ball, and write $R=R_\alpha$ in the bounded-support case.
We construct a function in the full-space Gaussian RKHS
$H_\gamma(\mathbb R^d)$ and then restrict it to $\mathcal X$;
kernel restriction cannot increase its RKHS norm.
Let $\mathcal E_P$ denote hinge excess under $P$ and put
\[
 C_d=\frac{2^d}{\Gamma(d/2+1)},\qquad
 M_{d,\beta}=2^{-\beta/2}
       \frac{\Gamma((d+\beta)/2)}{\Gamma(d/2)}.
\]

The posterior sign $s$ is initially defined only on $\mathcal X$.
Write $A_+=\{x\in\mathcal X:m(x)\ge0\}$ and
$A_-=\{x\in\mathcal X:m(x)<0\}$. If $A_-$ is empty, set
$\widetilde s=1$ on $\mathbb R^d$; if $A_+$ is empty, set
$\widetilde s=-1$. Otherwise define the nearest-class continuation
\[
 \widetilde s(z)=
 \begin{cases}
  1,&\operatorname{dist}(z,A_+)\le\operatorname{dist}(z,A_-),\\
  -1,&\operatorname{dist}(z,A_+)>\operatorname{dist}(z,A_-).
 \end{cases}
\]
Distances to nonempty subsets of $\mathbb R^d$ are finite
$1$-Lipschitz functions, so $\widetilde s$ is Borel measurable.
For $x\in\mathcal X$ and $\|z-x\|_2<\Delta(x)/2$, the distance of
$z$ to the class containing $x$ is at most $\|z-x\|_2$, whereas its
distance to the opposite class is at least
$\Delta(x)-\|z-x\|_2$; hence $\widetilde s(z)=s(x)$.
The same assertion holds for the constant continuation when one
class is empty and $\Delta(x)=\infty$. At points with
$\Delta(x)=0$ it is vacuous. In particular, the tie convention can
change the sign on $\mathcal X$ only in this latter set, whose
$|m|\,dP_X$ measure is zero by the margin-noise condition and
$\beta>0$.

Define the probability density
\[
 Q_\gamma(z):=
 \left(\frac{2}{\pi\gamma^2}\right)^{d/2}
 e^{-2\|z\|_2^2/\gamma^2},
\]
and let $Z_\gamma$ have this density. Polar integration gives
\[
 \mathbb E\|Z_\gamma\|_2^\beta=M_{d,\beta}\gamma^\beta.
\]
First fix $\rho\ge1$ and put
\[
 u_\rho=\widetilde s\,\mathbf1_{2\rho B},\qquad
 g_\rho=Q_\gamma*u_\rho.
\]
Since $|u_\rho|\le1$ and $Q_\gamma$ integrates to one,
$\|g_\rho\|_\infty\le1$.
To verify RKHS membership and control its norm, use the feature map
\[
 \Phi_\gamma(x)(z)=
 \frac{2^{d/2}}{\pi^{d/4}\gamma^{d/2}}
 e^{-2\|x-z\|_2^2/\gamma^2}.
\]
Completing the square gives
$\langle\Phi_\gamma(x),\Phi_\gamma(x')\rangle_{L_2}=k_\gamma(x,x')$.
With $w_\rho=(\pi\gamma^2)^{-d/4}u_\rho$, we have
$g_\rho(x)=\langle w_\rho,\Phi_\gamma(x)\rangle_{L_2}$.
The feature-space representation of an RKHS
\citep[Section~4.2]{steinwart2008support} therefore yields
\[
 \|g_\rho\|_{H_\gamma}^2
 \le\|w_\rho\|_{L_2}^2
 \le(\pi\gamma^2)^{-d/2}\operatorname{vol}_d(2\rho B)
 =C_d\rho^d\gamma^{-d}.
\]

For the error bound, introduce the untruncated convolution
$h_\gamma=Q_\gamma*\widetilde s$ only as an auxiliary bounded function;
no RKHS membership is claimed for $h_\gamma$.
The nearest-class construction gives, for every $x\in\mathcal X$,
\[
 |h_\gamma(x)-s(x)|
 \le2\mathbb P\{\|Z_\gamma\|_2\ge\Delta(x)/2\}.
\]
When $\Delta(x)=\infty$ the right side is zero; when
$\Delta(x)=0$ the bound is the trivial inequality $|h_\gamma-s|\le2$.
For $x\in\mathcal X\cap\rho B$, points omitted by the truncation to $2\rho B$
are farther than $\rho$ from $x$, so
\[
 |g_\rho(x)-h_\gamma(x)|
 \le\mathbb P\{\|Z_\gamma\|_2>\rho\}.
\]
For $|f|\le1$, the hinge loss is $1-yf$; hence
\eqref{eq:hinge_pointwise_excess} gives
$\mathcal E_P(f)=\int |f-s|\,|m|\,dP_X$.
Combining the last two displays and using Tonelli and the
margin-noise assumption gives
\begin{align*}
 \int_{\rho B}|g_\rho-s|\,|m|\,dP_X
 &\le2\mathbb E\left[
   \int_{\{\Delta\le2\|Z_\gamma\|_2\}}|m|\,dP_X\right]
       +\mathbb P\{\|Z_\gamma\|_2>\rho\}\\
 &\le(2^{\beta+1}c_\beta+\rho^{-\beta})
          \mathbb E\|Z_\gamma\|_2^\beta\\
 &\le(2^{\beta+1}+1)M_{d,\beta}c_\beta\gamma^\beta.
\end{align*}
Here $\int_{\rho B}|m|\,dP_X\le1$, Markov's inequality,
$\rho\ge1$, and $c_\beta\ge1$ were used.
No condition on how many decision classes meet $\rho B$ was needed.
In particular, the argument covers a globally constant Bayes sign.

Suppose first that $\tau<\infty$.
Outside $\rho B$, the integrand is bounded by two, and hence
\[
\begin{aligned}
 &\lambda\|g_\rho\|_{H_\gamma}^2+\mathcal E_P(g_\rho)\\
 &\quad\le C_d\lambda\gamma^{-d}\rho^d
      +2c_\tau\rho^{-\tau}
      +(2^{\beta+1}+1)M_{d,\beta}c_\beta\gamma^\beta,
 \qquad \rho\ge1.
\end{aligned}
\]
For $0<\rho<1$, take $g_\rho=0$ instead. Its norm is zero and
$\mathcal E_P(0)=\int|m|\,dP_X\le1$, whereas
$2c_\tau\rho^{-\tau}\ge2$. Thus the same total bound is valid
for every $\rho>0$, with the indicated choice of comparator.
The first two terms are minimized at
\[
 \rho_*=
 \left(\frac{2\tau c_\tau}
 {dC_d\lambda\gamma^{-d}}\right)^{1/(d+\tau)}.
\]
Use the convolution comparator if $\rho_*\ge1$, and the zero
comparator if $\rho_*<1$. Substitution gives the first term in
\eqref{eq:apalpha_bound}, with coefficient
\[
\begin{aligned}
 C_{d,\tau}
 &=C_d^{\tau/(d+\tau)}(2c_\tau)^{d/(d+\tau)}\\
 &\quad\times\left[
  \left(\frac{\tau}{d}\right)^{d/(d+\tau)}
 +\left(\frac{d}{\tau}\right)^{\tau/(d+\tau)}
 \right].
\end{aligned}
\]
Taking $C_{d,\beta}=(2^{\beta+1}+1)M_{d,\beta}$ proves the finite-tail assertion
for all $\lambda,\gamma>0$.

If $\tau=\infty$, set $\rho=\max\{1,R\}$ and use the convolution
comparator. The exterior integral is zero and the two previously
proved norm and interior bounds give the asserted bound with $C_{d,R}=C_d\max\{1,R\}^d$.

\end{proof}

\subsection{Weighted empirical Gaussian entropy (Lemma~\ref{lem:weighted_empirical_entropy})}\label{sec:appB}
Within this subsection the index $a\in\{T,S\}$ runs over the two
domains, $\mu_a$ and $\nu_a$ are local names for covariate laws, and
$\tau_\wedge:=\tau_T\wedge\tau_S$; these symbols are unrelated to the
amplitude $a$, the design measure $\mu$, and the margin measures
$\nu_D$ of the main text.

\begin{lemma}[Weighted Gaussian empirical entropy]
\label{lem:weighted_empirical_entropy}
Let $\mu_T=P_{X,T}$ and $\mu_S=P_{X,S}$ be distributions on
$\mathbb R^d$. For $a\in\{T,S\}$, assume either that
$\mu_a(\|X\|>r)\le c_a r^{-\tau_a}$ for all $r>0$, with
$0<\tau_a<\infty$ and $c_a>0$, or that $\tau_a=\infty$ and
$\mu_a$ is supported on a ball of finite radius $R_a$ centered
at the origin. Let $D_{T,X}$ and $D_{S,X}$ be the empirical
marginals of independent i.i.d.\ blocks of sizes $n,m\ge1$, and set
\[
 D_{w,X}=\alpha D_{T,X}+(1-\alpha)D_{S,X},
 \qquad 0<\alpha<1,\qquad \tau_\wedge=\tau_T\wedge\tau_S.
\]
For every $\varepsilon>0$ and $d/(d+\tau_\wedge)<p<1$, there is a
constant $C\ge1$, depending only on $\varepsilon,p,d$ and the
stated tail or support data, such that
\[
 \mathbb E_D\,
 e_i\bigl(\mathrm{id}:H_\gamma(\mathbb R^d)\to L_2(D_{w,X})\bigr)
 \le C\gamma^{-\kappa}i^{-1/(2p)},
 \qquad
 \kappa:=\frac{(1-p)(1+\varepsilon)d}{2p},
\]
for all $0<\gamma\le1$ and $i\ge1$. The constant is independent
of $n,m,\gamma,i,\alpha$ when the marginal tail data are fixed.
The same bound holds for $H_\gamma(\mathcal X)$ when the
marginals are supported on $\mathcal X\subset\mathbb R^d$.
For $\alpha=1$, the same conclusion holds using only the target
block and $\tau_\wedge=\tau_T$.
\end{lemma}

\begin{proof}
The RKHS of a restricted kernel consists of restrictions of
full-space functions with the norm of the extension of smallest
norm \citep[Section~4.2]{steinwart2008support}. Consequently, the
unit balls of $H_\gamma(\mathcal X)$ and
$H_\gamma(\mathbb R^d)$ have identical images on every sample in
$\mathcal X$. It suffices to prove the full-space bound.

\emph{Standard inputs and tail normalization.}
We use two precise results of \citet{steinwart2008support}.
First, Theorem~7.34 (p.~278, proof pp.~278--279) states that if
$\nu(\|X\|>r)\le r^{-\tau_0}$ for every $r>0$, then
\[
 e_i\bigl(\mathrm{id}:H_\delta(\mathbb R^d)\to L_2(\nu)\bigr)
 \le C_{\varepsilon,p,d,\tau_0}\,
       \delta^{-\kappa}i^{-1/(2p)},
 \qquad 0<\delta\le1,
\]
provided $\varepsilon>0$ and $d/(d+\tau_0)<p<1$.
Here $\tau_0=\infty$ means support in the closed unit ball,
as in Definition~7.32 (p.~277), and the constant can be taken
at least one. Second, Corollary~7.31 (p.~276) says that for a
measurable kernel with separable RKHS and square-integrable
diagonal envelope, a population bound
$e_i(H\to L_2(\nu))\le A i^{-1/(2p)}$, $A\ge1$, implies
\[
 \mathbb E_{D\sim\nu^k}e_i(H\to L_2(D))
 \le c_p A\min\{i,k\}^{1/(2p)}i^{-1/p}
 \le c_p A i^{-1/(2p)},\qquad i,k\ge1,
\]
where $c_p$ depends only on $p$.

For each marginal with finite $\tau_a$, put
$L_a=\max\{1,c_a^{1/\tau_a}\}$; in the bounded-support case put
$L_a=\max\{1,R_a\}$. The law $\nu_a$ of $X/L_a$ satisfies
exactly the normalized tail or unit-support condition above.
The identity
\[
 k_\gamma(L_ax,L_ay)=k_{\gamma/L_a}(x,y)
\]
shows that $f\mapsto f(L_a\,\cdot)$ is an isometry from
$H_\gamma(\mathbb R^d)$ onto $H_{\gamma/L_a}(\mathbb R^d)$:
it preserves inner products on finite kernel spans, and then
extends by completion. The same map is an isometry of the
corresponding $L_2(\mu_a)$ and $L_2(\nu_a)$ spaces.
Since $\gamma/L_a\le1$ and
$p>d/(d+\tau_\wedge)\ge d/(d+\tau_a)$, the first input gives
\[
 e_i\bigl(\mathrm{id}:H_\gamma\to L_2(\mu_a)\bigr)
 \le A_a\gamma^{-\kappa}i^{-1/(2p)},
 \qquad
 A_a:=C_{\varepsilon,p,d,\tau_a}L_a^\kappa\ge1.
\]
The Gaussian kernel is measurable, its RKHS is separable, and
$k_\gamma(x,x)=1$, so the integrability assumptions of the
second input hold. Applying it to each block, of size
$k_T=n$ or $k_S=m$, yields
\[
 \mathbb E_D e_i^{(a)}
 \le c_p A_a\gamma^{-\kappa}i^{-1/(2p)},
 \qquad
 e_i^{(a)}:=e_i\bigl(\mathrm{id}:H_\gamma\to L_2(D_{a,X})\bigr).
\]

\emph{Two-block product cover.}
For every function $f$ on the combined sample,
\begin{equation}
 \|f\|_{L_2(D_{w,X})}^2
 =\alpha\|f\|_{L_2(D_{T,X})}^2
 +(1-\alpha)\|f\|_{L_2(D_{S,X})}^2
 \le\max_{a\in\{T,S\}}\|f\|_{L_2(D_{a,X})}^2.
 \label{eq:dw_interpolation}
\end{equation}
Fix the sample and $i\ge1$. For each block, cover the image of
the unit ball of $H_\gamma$ by $2^{i-1}$ balls with radius
$\eta_a>e_i^{(a)}$. Intersect the inverse images of these covers.
There are at most $2^{2i-2}$ nonempty cells, each with diameter
at most $2\eta_a$ in its block norm, and hence at most
$2\max\{\eta_T,\eta_S\}$ in the weighted norm.
Selecting one function in each cell as its center and letting
$\eta_a\downarrow e_i^{(a)}$ proves the deterministic bound
\[
 e_{2i-1}\bigl(\mathrm{id}:H_\gamma\to L_2(D_{w,X})\bigr)
 \le2\max\{e_i^{(T)},e_i^{(S)}\}
 \le2(e_i^{(T)}+e_i^{(S)}).
\]
For arbitrary $j\ge1$, take $i=\lceil j/2\rceil$.
Monotonicity gives $e_j\le e_{2i-1}$, while
$i^{-1/(2p)}\le2^{1/(2p)}j^{-1/(2p)}$.
Taking expectations proves the assertion with
$C=\max\{1,2^{1+1/(2p)}c_p(A_T+A_S)\}$.
For $\alpha=1$, the per-target-block bound alone proves the
stated endpoint, without any condition on the source marginal.

\emph{Application to a fixed training mixture.}
If the available assumption is instead
$P_{X,\alpha}(\|X\|>r)\le c_{\tau,\alpha}r^{-\tau_\alpha}$,
then positivity of the mixture weights gives
\begin{align*}
 P_{X,T}(\|X\|>r)&\le(c_{\tau,\alpha}/\alpha)r^{-\tau_\alpha},\\
 P_{X,S}(\|X\|>r)&\le(c_{\tau,\alpha}/(1-\alpha))r^{-\tau_\alpha}.
\end{align*}
Use these two marginal constants and
$\tau_T=\tau_S=\tau_\alpha$ above. Bounded mixture support likewise
implies the same support bound for both marginals.
Thus the entropy constant used in the Gaussian oracle may depend
on the fixed interior weight $\alpha$. The combined sample need
not be i.i.d.\ from the mixture: only the two individual blocks
are used in the empirical-entropy input.
\end{proof}

\subsection{The Gaussian mixture oracle (Theorem~\ref{thm:gaussian_transfer_estimation})}\label{supp:gauss}

\begin{theorem}[Gaussian oracle inequality on the training mixture]
\label{thm:gaussian_transfer_estimation}
Fix \(\alpha\in(0,1)\).  Suppose that \(P_\alpha\) satisfies the
Tsybakov, margin-noise, and tail conditions stated in
Section~\ref{subsec:gaussian_est}, with parameters
\((q_\alpha,c_{q,\alpha})\), \((\beta_\alpha,c_{\beta,\alpha})\),
and \((\tau_\alpha,c_{\tau,\alpha})\), respectively, and use the
variance pair \((\vartheta_\alpha,V_\alpha)\) in
\eqref{eq:acceleration_variance_pair}. For every
\(\varepsilon>0\) and
\[
  \frac{d}{d+\tau_\alpha}<p<1,
\]
there is a constant \(K<\infty\), depending only on the displayed
distributional parameters, the fixed weight \(\alpha\), \(d,p,\varepsilon\), and the fixed
constants in the Gaussian entropy bound (including a support
radius $R_\alpha$ when $\tau_\alpha=\infty$), such that, for all
\(\lambda,\gamma\in(0,1]\) and \(\varrho\ge1\), with probability at
least \(1-3e^{-\varrho}\),
\begin{align}
  &\lambda\|f_{D,\lambda,\gamma}\|_{H_\gamma}^2
  +\mathcal R_{L,P_\alpha}(\bar f_{D,\lambda,\gamma})
  -\mathcal R^*_{L,P_\alpha}
  \notag\\
  &\quad\le K\Bigg[
  \lambda^{\frac{\tau_\alpha}{\tau_\alpha+d}}
  \gamma^{-\frac{d\tau_\alpha}{\tau_\alpha+d}}
  +c_{\beta,\alpha}\gamma^{\beta_\alpha}
  \notag\\
  &\qquad\qquad
  +\varrho\Bigl(
    \frac{\lambda^p\gamma^{(1-p)(1+\varepsilon)d}}
         {b_{\alpha,n,m}}
  \Bigr)^{-\frac{q_\alpha+1}{q_\alpha+2-p}}
  +\left(\varrho V_\alpha b_{\alpha,n,m}\right)^{
    \frac{q_\alpha+1}{q_\alpha+2}}
  +\varrho b_{\alpha,n,m}
  \Bigg].
  \label{eq:gaussian_estimation_bound}
\end{align}
For \(\tau_\alpha=\infty\), the first product in the brackets is
read as \(\lambda\gamma^{-d}\).  This theorem controls the training
mixture only; target-risk conversion is performed separately in
Corollary~\ref{thm:main_target_bound}.
\end{theorem}

\phantomsection\label{pf:gauss_est}%
\begin{proof}[Proof of Theorem~\ref{thm:gaussian_transfer_estimation}]
We apply Theorem~\ref{thm:transfer_bernstein_svm} with $M = 1$, $B =
2$, the variance pair \eqref{eq:acceleration_variance_pair},
justified by Lemma~\ref{lem:variance-bound-hinge}, the entropy constant
\eqref{eq:gaussian_entropy_constant}, and the comparator $f_0 := g^*$,
where $g^*\in H_\gamma$ is the bounded comparator supplied by
Lemma~\ref{lem:gaussian_mixture_approximation}. Since
$\|g^*\|_\infty\le1$, we may take
$B_0 := 2$, and the lemma gives
\[
  \lambda\|g^*\|^2_{H_\gamma}
  + \mathcal{R}_{L,P_\alpha}(g^*) - \mathcal{R}^*_{L,P_\alpha}
  \;\leq\;
  C_{d,\tau_\alpha}\,\lambda^{\frac{\tau_\alpha}{\tau_\alpha+d}}
    \gamma^{-\frac{d\tau_\alpha}{\tau_\alpha+d}}
  + C_{d,\beta_\alpha}\,c_{\beta,\alpha}\,\gamma^{\beta_\alpha}.
\]
With $\vartheta_\alpha = q_\alpha/(q_\alpha+1)$ we have
\[
\begin{aligned}
\frac{1}{2-p-\vartheta_\alpha+\vartheta_\alpha p}
 &= \frac{q_\alpha+1}{q_\alpha+2-p}, &
\frac{1}{2-\vartheta_\alpha} &= \frac{q_\alpha+1}{q_\alpha+2},\\
a^{2p} &= c_{\varepsilon,p}^{2p}\gamma^{-(1-p)(1+\varepsilon)d}.
\end{aligned}
\]
The capacity term of Theorem~\ref{thm:transfer_bernstein_svm} becomes
\[
 K'\varrho\,(\lambda^p
 \gamma^{(1-p)(1+\varepsilon)d}/b_{\alpha,n,m})^{
 -(q_\alpha+1)/(q_\alpha+2-p)}
\]
(the confidence factor $\varrho \geq 1$ is kept explicitly,
enlarging the constant), the variance term becomes
\[
 6(3600V_\alpha\varrho b_{\alpha,n,m})^{
 (q_\alpha+1)/(q_\alpha+2)},
\]
and the remainder is
\[
 1296B_0\varrho b_{\alpha,n,m}=2592\varrho b_{\alpha,n,m}.
\]
  Collecting constants into $K$ yields
\eqref{eq:gaussian_estimation_bound}.
\end{proof}

\begin{remark}[Target-only Gaussian rate]
\label{rem:gaussian_target_only}
For \(\alpha=1\), suppose \(P_T\) satisfies the Tsybakov,
margin-noise, and tail conditions with target parameters.
Apply Theorem~\ref{thm:transfer_bernstein_svm} directly with
\(b_{1,n,m}=1/n\),
\(\vartheta_T=q_T/(q_T+1)\), and
\(V_T=6c_{q,T}^{1/(q_T+1)}\). The target-block case of
Lemma~\ref{lem:weighted_empirical_entropy} and the comparator of
Lemma~\ref{lem:gaussian_mixture_approximation} for \(P_T\) give
\eqref{eq:gaussian_estimation_bound} with target parameters and
\(N_{\rm eff}=n\). For every \(\xi>0\), choose the tuning exponents
in \eqref{eq:gaussian_rate_exponent} at the target parameters, with
\(p\) and \(\varepsilon\) chosen as in the balancing proof below.
Then the target hinge excess is at most
\(C\varrho n^{-(r_G-\xi)}\) with probability at least
\(1-3e^{-\varrho}\), where \(r_G=\beta_T\sigma_\gamma\), using
the bounded-support limit when \(\tau_T=\infty\). Constants may
depend on the fixed target and entropy parameters and on \(\xi\).
\end{remark}

\subsection{Proof of Corollary~\ref{thm:main_target_bound} and balanced tuning (Corollary~\ref{cor:transfer_rates})}
\label{pf:main_target_bound}

\begin{proof}[Proof of Corollary~\ref{thm:main_target_bound}]
On the event supplied by
Theorem~\ref{thm:gaussian_transfer_estimation}, whose probability
is at least $1-3e^{-\varrho}$, use
$b_{\alpha,n,m}=N_{\rm eff}^{-1}$ and enlarge the fixed constant
in $\mathcal A_\alpha$ and $\mathcal S_\alpha$ if necessary to obtain
\[
  \mathcal E_\alpha(\bar f_{D,\lambda,\gamma})
  \le \lambda\|f_{D,\lambda,\gamma}\|_{H_\gamma}^2
       +\mathcal E_\alpha(\bar f_{D,\lambda,\gamma})
  \le \mathcal A_\alpha(\lambda,\gamma)
       +\mathcal S_\alpha(\lambda,\gamma).
\]
The deterministic conversion of
Theorem~\ref{thm:affine_frontier}, and monotonicity of
$\Psi_\alpha$, therefore give on the same event
\[
  \mathcal E_T(\bar f_{D,\lambda,\gamma})
  \le \Psi_\alpha\!\left(
       \mathcal E_\alpha(\bar f_{D,\lambda,\gamma})\right)
  \le \Psi_\alpha(\mathcal A_\alpha+\mathcal S_\alpha).
\]
The affine representation follows from
\eqref{eq:conversion_duality}. Finally,
\eqref{eq:hinge_classification_calibration} bounds the target
classification excess of $\operatorname{sign}\bar f_{D,\lambda,\gamma}$
by its clipped hinge excess. No additional probability event is
needed.
\end{proof}

\phantomsection\label{pf:gauss_rates}%
\paragraph*{Balanced tuning}
The balancing exponents of Corollary~\ref{cor:transfer_rates}
are
\begin{equation}
\begin{aligned}
  \sigma_\gamma&:=\frac{\tau_\alpha(d+\tau_\alpha)(q_\alpha+1)}{
  \beta_\alpha\tau_\alpha^2(q_\alpha+2)
  +d\{d\tau_\alpha+\beta_\alpha d+2\beta_\alpha\tau_\alpha
  +\tau_\alpha^2\}(q_\alpha+1)},\\
  \sigma_\lambda&:=
  \frac{\beta_\alpha(\tau_\alpha+d)+d\tau_\alpha}{\tau_\alpha}
  \sigma_\gamma,
\end{aligned}
  \label{eq:gaussian_rate_exponent}
\end{equation}
and \(r_G:=\beta_\alpha\sigma_\gamma\)
(for \(\tau_\alpha=\infty\):
\(\sigma_\gamma=\frac{q_\alpha+1}{\beta_\alpha(q_\alpha+2)
+d(q_\alpha+1)}\), \(\sigma_\lambda=(\beta_\alpha+d)\sigma_\gamma\),
and \(r_G=\frac{\beta_\alpha(q_\alpha+1)}{\beta_\alpha(q_\alpha+2)
+d(q_\alpha+1)}\)).
Corollary~\ref{cor:transfer_rates} chooses
\(\gamma=N_{\rm eff}^{-\sigma_\gamma}\) and
\(\lambda=N_{\rm eff}^{-\sigma_\lambda}\) and balances the
argument of \(\Psi_\alpha\) to \(C\varrho N_{\rm eff}^{-(r_G-\xi)}\);
the calculation is in the proof of Corollary~\ref{cor:transfer_rates} below.

The following corollary is the balanced-tuning form of
Corollary~\ref{thm:main_target_bound} referred to in the main text.

\begin{corollary}[Balanced Gaussian rate]
\label{cor:transfer_rates}\label{cmp:balanced-gaussian-rate}
Let \(\xi>0\) be arbitrary, and choose
\(\gamma_{\alpha,n,m}=N_{\rm eff}^{-\sigma_\gamma}\) and
\(\lambda_{\alpha,n,m}=N_{\rm eff}^{-\sigma_\lambda}\).
Under the mixture conditions of Corollary~\ref{thm:main_target_bound},
choose \(p>d/(d+\tau_\alpha)\) and \(\varepsilon>0\) sufficiently
close to \(d/(d+\tau_\alpha)\) and zero, respectively, as a function
of \(\xi\). Then the balanced oracle bound satisfies
\(\mathcal A_\alpha+\mathcal S_\alpha\le C\varrho\,N_{\rm eff}^{-(r_G-\xi)}\),
where \(C\) may depend on \(\xi\) and the fixed model parameters.
Consequently, with probability at least \(1-3e^{-\varrho}\),
\begin{equation}
  \mathcal E_T(\bar f_{D,\lambda,\gamma})
  \;\le\;
  \Psi_\alpha\bigl(C\varrho\,N_{\rm eff}^{-(r_G-\xi)}\bigr);
  \label{eq:transfer_rate_universal}
\end{equation}
for fixed \(\varrho\), the respective assumptions of
Corollary~\ref{cor:regime_flip} give upper bounds of orders
\(M_\alpha N_{\rm eff}^{-(r_G-\xi)}\),
\(\delta^{-1}N_{\rm eff}^{-(r_G-\xi)}\), and
\(H_{\rm flip}(u)+C'N_{\rm eff}^{-(r_G-\xi)\chi}\).
\end{corollary}

These are fixed-confidence upper-bound orders.  Linear conversion
preserves the order of the mixture bound, the H\"older regime
multiplies the vanishing exponent by \(\chi\), and the tie-or-flip
contribution remains additive. Remark~\ref{rem:gaussian_target_only}
gives the target-only rate at \(\alpha=1\). Kernel tuning shrinks only the
argument of \(\Psi_\alpha\); the modulus is fixed by the pair and
the weight.

\begin{proof}[Proof of Corollary~\ref{cor:transfer_rates}]
Write \(b:=b_{\alpha,n,m}\) and consider the four
dominant terms of the mixture-side sum
\(\mathcal A_\alpha+\mathcal S_\alpha\) at
the limiting parameters $\varepsilon \to 0$, $p \to d/(d+\tau_\alpha)$:
\[
  T_1 := \lambda^{\frac{\tau_\alpha}{\tau_\alpha+d}}
         \gamma^{-\frac{d\tau_\alpha}{\tau_\alpha+d}},
  \quad
  T_2 := \gamma^{\beta_\alpha},
  \quad
  T_3 := \bigl(\lambda^p\gamma^{(1-p)d}/b\bigr)^{
         -\frac{q_\alpha+1}{q_\alpha+2-p}},
  \quad
  T_4 := b^{\frac{q_\alpha+1}{q_\alpha+2}} .
\]
Balancing $T_1 = T_2$ gives $\lambda =
\gamma^{[\beta_\alpha(\tau_\alpha+d)+d\tau_\alpha]/\tau_\alpha}$,
i.e.\ $\sigma_\lambda = \sigma_\gamma\,
[\beta_\alpha(\tau_\alpha+d)+d\tau_\alpha]/\tau_\alpha$. Substituting
\(\gamma=b^{\sigma_\gamma}\) and balancing \(T_2=T_3\) gives
\[
  \beta_\alpha\sigma_\gamma\,(q_\alpha+2-p)
  = (q_\alpha+1)\bigl[1 - p\,\sigma_\lambda - (1-p)\,d\,\sigma_\gamma\bigr];
\]
at $p = d/(d+\tau_\alpha)$ one computes
$p\,\sigma_\lambda/\sigma_\gamma + (1-p)d =
d(\tau_\alpha+\beta_\alpha)/\tau_\alpha$, and solving for
$\sigma_\gamma$ and simplifying the denominator (using
$\beta_\alpha\tau_\alpha(q_\alpha+2)(d+\tau_\alpha)
- \beta_\alpha\tau_\alpha d
+ (q_\alpha+1)d(\tau_\alpha+\beta_\alpha)(d+\tau_\alpha)
= \beta_\alpha\tau_\alpha^2(q_\alpha+2)
+ d(d\tau_\alpha+\beta_\alpha d+2\beta_\alpha\tau_\alpha
+\tau_\alpha^2)(q_\alpha+1)$) yields the stated $\sigma_\gamma$ with
common value \(T_1=T_2=T_3=b^{\beta_\alpha\sigma_\gamma}\).
Keep these tuning exponents fixed and put
\(p_0:=d/(d+\tau_\alpha)\), with \(p_0=0\) for bounded support.
For admissible \(p>p_0\) and \(\varepsilon>0\), the exponent of
\(b\) in the capacity term is
\[
 s(p,\varepsilon):=
 \frac{(q_\alpha+1)
 [1-p\sigma_\lambda-(1-p)(1+\varepsilon)d\sigma_\gamma]}
 {q_\alpha+2-p}.
\]
Since \(s(p_0,0)=r_G\), continuity permits a choice of
\(p,\varepsilon\), depending on \(\xi\), such that
\(s(p,\varepsilon)\ge r_G-\xi\). The two approximation terms
retain exponent \(r_G\).

It remains to absorb \(T_4\) and the remainder \(\varrho b\). From the
balancing equation at \(p=p_0\),
\[
  \beta_\alpha\sigma_\gamma
  = \frac{\beta_\alpha(q_\alpha+1)}
  {\beta_\alpha(q_\alpha+2-p_0)
   + (q_\alpha+1)\bigl[p_0\,\sigma_\lambda/\sigma_\gamma + (1-p_0)d\bigr]}
  \leq \frac{q_\alpha+1}{q_\alpha+2},
\]
since $p_0\,\sigma_\lambda/\sigma_\gamma + (1-p_0)d \geq
p_0\beta_\alpha(\tau_\alpha+d)/\tau_\alpha + d - p_0d \geq p_0\beta_\alpha$,
whence the denominator is at least
$\beta_\alpha(q_\alpha+2-p_0) + (q_\alpha+1)p_0\beta_\alpha \geq
\beta_\alpha(q_\alpha+2)$. Hence \(T_4\le
b^{\beta_\alpha\sigma_\gamma}\). Since \(0<b\le1\) and
\(\varrho\ge1\), writing \(a=(q_\alpha+1)/(q_\alpha+2)\le1\)
gives \((\varrho b)^a\le\varrho b^a\) and
\(\varrho b\le\varrho b^a=\varrho T_4\). Thus
\(\mathcal A_\alpha+\mathcal S_\alpha
\leq C\varrho b^{\beta_\alpha\sigma_\gamma-\xi}\).
Since \(\Psi_\alpha\) is nondecreasing
(Theorem~\ref{thm:affine_frontier}), substituting this balanced
bound for the argument of the modulus in
\eqref{eq:main_target_bound} of Corollary~\ref{thm:main_target_bound}
proves \eqref{eq:transfer_rate_universal} on the same event, and
evaluating the three cases of Corollary~\ref{cor:regime_flip} at
the balanced bound gives the regime statement of
Corollary~\ref{cor:transfer_rates}; the scalar forms
\(M_\alpha t\) and \(\delta^{-1}t\) are evaluations of the modulus
under the respective case assumptions, and no finiteness of
\(M_\alpha\) is needed for \eqref{eq:transfer_rate_universal}
itself.
\end{proof}

\section{Adaptation and the Exact-Fitting Limit}\label{supp:adaptive}

\subsection{Radial geometry and source-quality estimation}
\label{supp:adaptive_calibration}

The radial kernel identities identify the source's aligned strength
and tilt from signed kernel means.  We then control their empirical
estimates and construct the positivity certificate used to
propose a source-using fit; the fixed class constants and their
order are collected in a separate paragraph below.  Throughout,
the class \(\mathcal P_1\), the envelope
\eqref{eq:sgn_envelope}, and the fixed pair
\((\gamma_0,\lambda_0)\) are as in
Section~\ref{subsec:adaptive_minimax}; \(d\ge2\), \(N=n+m\),
\(\ell_N=\log(eN)\), \(\mu_0=\operatorname{Unif}(B^d)\),
the kernel is \(k_{\gamma_0}\) with \(k_{\gamma_0}(x,z)=e^{-\|x-z\|_2^2/\gamma_0^2}\),
\(H=H_{\gamma_0}\), and \(\operatorname{sign}(0)=1\).  Throughout
this section we write \(v:=w_T\) for the unit target normal and
\(r:=w_S-sw_T\) for the residual of the source normal of
\eqref{eq:sgn_class}, so that \(w_S=sv+r\), \(r\perp v\),
\(\tau=\|r\|_2\), and \(w_\alpha=\kappa_\alpha v+ur\) as in
\eqref{eq:mixture_vote_vector}; the pilot estimate
\(\widehat w_T\) of the main text is written \(\widehat v\).  The profile
class \(\mathcal Q_{\rm rad}\) consists of the measurable radial
\(\varphi(x)=\bar\varphi(\|x\|_2)\) with
\(\varphi_-\le\varphi\le\varphi_+\) and
\((d+2)\int\varphi(x)x_1^2\,d\mu_0(x)=1\); since
\(\varphi_-\le1\le\varphi_+\) and
\(\mathbb E_{\mu_0}X_1^2=(d+2)^{-1}\), the constant profile
\(\varphi\equiv1\) belongs to it.  Of the two envelope conditions
in \eqref{eq:sgn_envelope}, the first keeps the label
probabilities in \([\varepsilon_Y/2,1-\varepsilon_Y/2]\), and the
second places every population soft-margin minimizer in the
all-active regime (Lemma~\ref{lem:supp_bridge} below).  All
constants \(c,C\) depend only on the fixed class constants,
\((d,\gamma_0,\lambda_0)\), the split fractions, and the gate
constants, and may change from line to line.  Two elementary
facts are used repeatedly.  Every empirical SVM of this section
has objective value at most one at \(f=0\), hence
\(\|\widehat f\|_H\le\lambda_0^{-1/2}\); and since
\(\|k_{\gamma_0}(x,\cdot)-k_{\gamma_0}(z,\cdot)\|_H^2=2\{1-e^{-\|x-z\|^2/\gamma_0^2}\}
\le2\|x-z\|^2/\gamma_0^2\), every such \(\widehat f\) satisfies
\begin{equation}
  |\widehat f(x)-\widehat f(z)|\le L_{\max}\|x-z\|_2,
  \qquad
  L_{\max}:=\frac{\sqrt2}{\gamma_0\sqrt{\lambda_0}}.
  \label{eq:supp_lipschitz}
\end{equation}

\paragraph*{Radial Gaussian geometry}
For a profile \(\varphi\) and a vector \(w\in\mathbb R^d\), put
\begin{equation}
\begin{gathered}
  G_{\varphi,w}
  :=\int_{B^d}\varphi(z)\,(w^\top z)\,k_{\gamma_0}(z,\cdot)\,d\mu_0(z)\in H,\\
  \psi_1(t)\,t:=G_{1,e_1}(te_1),
  \qquad
  \kappa_1:=\|G_{1,e_1}\|_H^2 .
\end{gathered}
  \label{eq:supp_G_def}
\end{equation}

\begin{lemma}[Radial Gaussian geometry]
\label{lem:supp_radial}
Uniformly over \(\varphi\in\mathcal Q_{\rm rad}\):
\begin{itemize}
\item[(a)] there are a continuous
\(\psi_\varphi:[0,1]\to(0,\infty)\) and a constant
\(\kappa_\varphi>0\) with
\(G_{\varphi,w}(x)=\psi_\varphi(\|x\|_2)\,w^\top x\) and
\(\langle G_{\varphi,w},G_{\varphi,w'}\rangle_H=\kappa_\varphi
w^\top w'\), and
\[
  \varphi_-\psi_1\le\psi_\varphi\le\varphi_+\psi_1,
  \qquad
  \varphi_-^2\kappa_1\le\kappa_\varphi\le\varphi_+^2\kappa_1,
  \qquad
  \psi_\varphi\ge\varphi_-\underline\psi_1,
\]
where
\[
\underline\psi_1:=4\gamma_0^{-2}e^{-1/\gamma_0^2}
\int_{\{z_1>0\}}z_1^2e^{-\|z\|^2/\gamma_0^2}\,d\mu_0(z)>0;
\]
moreover \(\|G_{\varphi,w}\|_\infty\le\varphi_+\|w\|_2\);
\item[(b)] for \(P\in\mathcal P_1\) with parameters
\((a,v,s,r,\varphi)\), the signed kernel means are
\(\mu_T^{\gamma_0}=aG_{\varphi,v}\) and \(\mu_S^{\gamma_0}=aG_{\varphi,sv+r}\), the
quality functionals \eqref{eq:quality_functionals} satisfy
\(s(P)=s\) and \(\tau(P)=\|r\|_2\), and
\(\|\mu_T^{\gamma_0}\|_H=a\sqrt{\kappa_\varphi}\ge
P_-:=a_-\varphi_-\sqrt{\kappa_1}>0\);
\item[(c)] \(\mathbb E_{P_T}[YX]=a\,v/(d+2)\).
\end{itemize}
\end{lemma}

\begin{proof}
(a) Radiality of \(\varphi\), of \(\mu_0\), and of the kernel
gives orthogonal equivariance,
\(G_{\varphi,Uw}(Ux)=G_{\varphi,w}(x)\) for every orthogonal
\(U\), so \(G_{\varphi,w}(x)\) is a function of \(\|x\|\) times
\(w^\top x\); call it \(\psi_\varphi(\|x\|)w^\top x\).  To bound
\(\psi_\varphi\), take \(w=e_1\), \(x=te_1\) with \(0<t\le1\), and
pair \(z=(z_1,z')\) with \((-z_1,z')\): the paired integrand is
\[
  \varphi(z)\,z_1\,e^{-(t^2+\|z\|^2)/\gamma_0^2}
  \bigl(e^{2tz_1/\gamma_0^2}-e^{-2tz_1/\gamma_0^2}\bigr)\ge0
  \qquad(z_1>0),
\]
so \(G_{\varphi,e_1}(te_1)\) is an integral of a nonnegative
function times \(\varphi\in[\varphi_-,\varphi_+]\), which gives
\(\varphi_-\psi_1(t)\le\psi_\varphi(t)\le\varphi_+\psi_1(t)\).
Dividing the paired integrand by \(t\) and applying dominated
convergence defines the continuous extension at zero:
\[
 \psi_\varphi(0)=\frac4{\gamma_0^2}
 \int_{\{z_1>0\}}\varphi(z)z_1^2
 e^{-\|z\|^2/\gamma_0^2}\,d\mu_0(z)>0.
\]
The same bounds hold at \(t=0\).  Since
\(2\sinh(2tz_1/\gamma_0^2)\ge4tz_1/\gamma_0^2\) and
\(e^{-t^2/\gamma_0^2}\ge e^{-1/\gamma_0^2}\), the same pairing
gives \(\psi_1(t)\ge\underline\psi_1\).  Next,
\[
\begin{aligned}
\langle G_{\varphi,w},G_{\varphi,w'}\rangle_H
&=\int\varphi(x)(w^\top x)G_{\varphi,w'}(x)\,d\mu_0(x)\\
&=\int\varphi(x)\psi_\varphi(\|x\|)(w^\top x)(w'^\top x)\,d\mu_0(x)\\
&=\kappa_\varphi w^\top w'
\end{aligned}
\]
by radial symmetry, with
\[
\kappa_\varphi=\int\varphi\psi_\varphi x_1^2\,d\mu_0\in
[\varphi_-^2\kappa_1,\varphi_+^2\kappa_1].
\]
Finally,
\[
\begin{aligned}
|G_{\varphi,w}(x)|
&\le\varphi_+\|w\|\int\|z\|k_{\gamma_0}(z,x)\,d\mu_0(z)\\
&\le\varphi_+\|w\|.
\end{aligned}
\]

(b) Since \(\mathbb E[Y\mid X=x]=m_D(x)\) and both designs equal
\(\mu_0\), \(\mu_T^{\gamma_0}=\int m_T(x)k_{\gamma_0}(x,\cdot)\,d\mu_0=aG_{\varphi,v}\)
and \(\mu_S^{\gamma_0}=aG_{\varphi,sv+r}\).  By (a) and \(r\perp v\),
\(\langle\mu_S^{\gamma_0},\mu_T^{\gamma_0}\rangle=a^2\kappa_\varphi s\),
\(\|\mu_T^{\gamma_0}\|^2=a^2\kappa_\varphi\), and
\(\|\mu_S^{\gamma_0}-s\mu_T^{\gamma_0}\|=a\|G_{\varphi,r}\|=a\sqrt{\kappa_\varphi}\|r\|\),
which is \eqref{eq:quality_functionals} with \(s(P)=s\),
\(\tau(P)=\|r\|\).

(c) \(\mathbb E[YX]=\mathbb E[m_T(X)X]
=a\,\mathbb E[\varphi(X)XX^\top]v=a\,v\,(d+2)^{-1}\) by radial
symmetry and the normalization of \(\varphi\).
\end{proof}

Part (b) is the radial-kernel calculation invoked after
\eqref{eq:quality_functionals}; part (a) also identifies the
population soft-margin centers used below, and part (c) is the
moment on which the target gate is built.

\paragraph*{The fixed constants and their order}
The constants are chosen in the following order so that each later
requirement depends only on quantities already fixed.
Fix a split constant \(\kappa_P\in(0,1/5]\), an integer
\(n_{\min}\), and deterministic index splits such that, whenever
\(n,m\ge n_{\min}\), every target block has size at least
\(\kappa_Pn\) and every source block at least \(\kappa_Pm\); the
concrete rule \((n_C,n_I,n_P,n_0)=(b_n,b_n,b_n,n-3b_n)\) with
\(b_n=\lfloor n/4\rfloor\) and \((m_C,m_I)=(b_m,m-b_m)\) with
\(b_m=\lfloor m/2\rfloor\) qualifies with \(\kappa_P=1/5\) once
\(n_{\min}\ge20\).  Fix \(A>3\).  The remaining constants are
fixed in the following noncircular order.  First the
concentration constants \(C_b\) (pilot) and \(C_A\)
(calibration), chosen so that the two concentration events below
have probability at least \(1-CN^{-A}\); then the perturbation
constant \(C_{\rm pert}\) of Lemma~\ref{lem:signed_projection_confidence}
and the derived quantities
\[
  \varepsilon_0(C_{\rm pert}):=\min\{1/4,(16C_{\rm pert})^{-1}\},
  \qquad
  c_{\rm gate}:=\frac{\varphi_-\underline\psi_1}{4\lambda_0};
\]
then \(C_{\rm cal}\ge2C_{\rm pert}\), and the absorption
constant \(C_{\rm abs}\) of Lemma~\ref{lem:signed_projection_confidence},
which depends only on \(C_{\rm pert}\) and \(\varepsilon_0(C_{\rm pert})\);
then
\[
  0<\varepsilon_*\le\min\Bigl\{\varepsilon_0(C_{\rm pert}),\;
  \frac1{2C_{\rm abs}C_{\rm cal}}\Bigr\};
\]
and finally the imbalance constant
\begin{equation}
  0<c_{\rm imb}\le\kappa_P\min\Bigl\{
  \Bigl[\frac{a_-}{33(d+2)C_b}\Bigr]^2,\;
  \Bigl(\frac{\varepsilon_*P_-}{2C_A}\Bigr)^2\Bigr\}.
  \label{eq:supp_cimb}
\end{equation}
The sample range is \(n,m\ge n_{\min}\) and
\(\ell_N\le c_{\rm imb}n\); no polynomial relation between \(n\)
and \(m\) is imposed.  Where the argument needs \(N\) large in
a way that depends only on the fixed constants (for instance
\(N^{-3}\le1/4\), or the mesh \(\mathcal M_N\) below nonempty),
\(n_{\min}\) is enlarged accordingly.

\paragraph*{Calibration and the positivity certificate}
On the calibration blocks put
\[
  \widehat p=\widehat\mu_T^{\gamma_0}=\frac1{n_C}\sum_{T_C}Y_ik_{\gamma_0}(X_i,\cdot),
  \quad
  \widehat q=\widehat\mu_S^{\gamma_0}=\frac1{m_C}\sum_{S_C}Y_jk_{\gamma_0}(X_j,\cdot),
  \quad
  \widehat D=\|\widehat p\|_H,
\]
\[
  e_T=\min\{2,C_A\sqrt{\ell_N/n_C}\},
  \qquad
  e_S=\min\{2,C_A\sqrt{\ell_N/m_C}\},
\]
and define the calibration event
\[
  \Omega_C=
  \bigl\{\|\widehat p-\mu_T^{\gamma_0}\|_H\le e_T,\
  \|\widehat q-\mu_S^{\gamma_0}\|_H\le e_S\bigr\}.
\]
The required dimension-free concentration follows directly from
bounded differences.
For independent Hilbert-space variables \(Z_i\) with \(\|Z_i\|\le1\),
write \(F=\|k^{-1}\sum_i(Z_i-\mathbb EZ_i)\|\).
Independence gives \(\mathbb EF\le k^{-1/2}\); changing one
observation changes \(F\) by at most \(2/k\). The
bounded-differences inequality yields
\begin{equation}
 \mathbb P\{F>k^{-1/2}+\sqrt{2t/k}\}\le e^{-t},\qquad t>0.
 \label{eq:supp_hilbert_mean}
\end{equation}
Apply this to \(Z_i=Y_ik_{\gamma_0}(X_i,\cdot)\) and choose
\(t=A\log(eN)\). A sufficiently large \(C_A\) gives
\(P(\Omega_C^c)\le CN^{-A}\); the cap at two makes the
corresponding event deterministic. The same dimension-free
argument applies to the pilot vectors \(Y_iX_i\).
Throughout, \(d\) remains fixed and the geometric constants may
depend on \(d\); this is not a dimension-uniform risk assertion.
On \(\Omega_C\), \eqref{eq:supp_cimb} gives
\(e_T\le\varepsilon_*P_-/2\), hence \(\widehat D\ge P_-/2\) and
\(e_T/\widehat D\le\varepsilon_*\).  When \(\widehat D>0\) and
\(e_T/\widehat D\le\varepsilon_*\), set
\[
  \widehat s_0
  =\frac{\langle\widehat q,\widehat p\rangle}{\widehat D^2},
  \qquad
  \widehat r_0=\widehat q-\widehat s_0\widehat p,
  \qquad
  \widehat\tau_0=\frac{\|\widehat r_0\|}{\widehat D},
  \qquad
  \widehat\Delta
  =C_{\rm cal}\Bigl\{\frac{e_S}{\widehat D}
  +\bigl(|\widehat s_0|+\widehat\tau_0\bigr)\frac{e_T}{\widehat D}
  +N^{-2}\Bigr\},
\]
and define the positivity certificate as the event
\(\Omega_{\rm cert}:=\{\widehat D>0,\ e_T/\widehat D\le\varepsilon_*,\
\widehat s_0>4\widehat\Delta\}\).  On \(\Omega_{\rm cert}\), the estimate
\(\widetilde s\) of Section~\ref{subsec:adaptive_minimax} is
\(\widehat s_0\) clipped above at \(2s_U\) and rounded to the
nearest point of the mesh
\(\mathcal M_N=\{kN^{-2}:1\le k\le\lfloor2s_UN^2\rfloor\}\)
(the smaller point at ties); on \(\Omega_{\rm cert}^c\),
\(\widetilde s=N^{-2}\).  In either case \(\widetilde s\ge N^{-2}\).  On
\(\Omega_{\rm cert}\), \(\widehat s_0>4\widehat\Delta\ge4C_{\rm cal}N^{-2}\),
so rounding moves \(\widehat s_0\) by at most
\(N^{-2}/2\le\widehat s_0/(8C_{\rm cal})\), and on
\(\Omega_{\rm cert}\cap\Omega_C\) the clip is inactive because
\(\widehat s_0<\tfrac43s\le\tfrac43s_U\) by
Lemma~\ref{lem:signed_projection_confidence}; hence
\(\widetilde s/\widehat s_0\) lies in a fixed interval around one.

\begin{lemma}[Signed-projection confidence]
\label{lem:signed_projection_confidence}
On \(\Omega_C\), \(|\widehat s_0-s|\le\widehat\Delta\).  On
\(\Omega_{\rm cert}\cap\Omega_C\), \(s>0\) and
\(\tfrac34\widehat s_0<s<\tfrac54\widehat s_0\), so that
\(c_-\le\widetilde s/s\le c_+\) for fixed
\(0<c_-<c_+<\infty\).  If \(s>0\), then on
\(\Omega_C\cap\Omega_{\rm cert}^c\)
\(s^2\le C\ell_N(m_C^{-1}+\tau^2n_C^{-1})+CN^{-4}\), and hence
\(ms^2/(1+m\tau^2)\le C\ell_N\); if in addition \(\tau=0\), then
\(ms^2\le C\ell_N\).  If \(s\le0\), then
\(\Omega_C\subset\Omega_{\rm cert}^c\).
\end{lemma}

\begin{proof}
By Lemma~\ref{lem:supp_radial}(b), write \(p=\mu_T^{\gamma_0}\) and
\(q=sp+r_H\), where \(\langle r_H,p\rangle=0\) and
\(\|r_H\|=\tau\|p\|\). Let
\(\widehat p=p+\delta_T\), \(\widehat q=q+\delta_S\), so that
\(\|\delta_T\|\le e_T\) and \(\|\delta_S\|\le e_S\) on
\(\Omega_C\); also
\(\widehat q=\widehat s_0\widehat p+\widehat r_0\) with
\(\widehat r_0\perp\widehat p\).  Comparing the two
decompositions and pairing with \(p\) gives the exact identity
\[
  (s-\widehat s_0)\|p\|^2
  =\widehat s_0\langle\delta_T,p\rangle
  -\langle\widehat r_0,\delta_T\rangle
  -\langle\delta_S,p\rangle .
\]
Since \(e_T/\widehat D\le\varepsilon_*\le1/4\), the norms
\(\widehat D\) and \(\|p\|\) are uniformly comparable, whence
\(|\widehat s_0-s|\le C_{\rm pert}\{e_S/\widehat D
+(|\widehat s_0|+\widehat\tau_0)e_T/\widehat D\}\); the choice
\(C_{\rm cal}\ge2C_{\rm pert}\) gives the first claim, and the
definition of \(\Omega_{\rm cert}\) then brackets \(s\) between
\(\widehat s_0\mp\widehat\Delta\), which proves the second; the
factor-two comparison of \(\widetilde s\) and \(\widehat s_0\)
gives the bracket for \(\widetilde s/s\).  For
the residual, projection onto \(\widehat p^\perp\) kills the
\(s\widehat p\) component, so
\(\widehat r_0=\Pi_{\widehat p^\perp}(r_H+\delta_S-s\delta_T)\),
and \(\|(\Pi_{\widehat p^\perp}-I)r_H\|
=|\langle r_H,\delta_T\rangle|/\widehat D\le\|r_H\|e_T/\widehat D\);
with \(|\widehat D-\|p\||\le e_T\) this yields
\(|\widehat\tau_0-\tau|\le C_{\rm pert}
\{e_S/\widehat D+(|s|+\tau)e_T/\widehat D\}\), and, after absorption
for \(e_T/\widehat D\le\varepsilon_0(C_{\rm pert})\),
\(|\widehat s_0|+\widehat\tau_0\le C\{s+\tau+e_S/\widehat D\}\)
whenever \(s>0\).  On \(\Omega_C\cap\Omega_{\rm cert}^c\), the
denominator conditions cannot be the reason, and
\(s\le|\widehat s_0-s|+\max(\widehat s_0,0)\le5\widehat\Delta\),
so \(s\le C_{\rm abs}C_{\rm cal}\{e_S/\widehat D+\tau e_T/\widehat D
+se_T/\widehat D+N^{-2}\}\), where \(C_{\rm abs}\) collects the
fixed factors of the two perturbation bounds and of the absorption
step, and depends only on \(C_{\rm pert}\) and
\(\varepsilon_0(C_{\rm pert})\).  The coefficient of \(s\) on the
right is \(C_{\rm abs}C_{\rm cal}e_T/\widehat D
\le C_{\rm abs}C_{\rm cal}\varepsilon_*\le1/2\) by the last
constraint on \(\varepsilon_*\), so it can be moved to the left, and
squaring with
\(\widehat D\asymp1\), \(e_S^2\le C\ell_N/m_C\), and
\(e_T^2\le C\ell_N/n_C\) proves the third claim.  Block sizes are
fixed fractions of \(n,m\), so multiplying by \(m/(1+m\tau^2)\)
bounds the leading terms by \(C\ell_N\), and \(mN^{-4}\le N^{-3}\);
at \(\tau=0\) the same multiplication by \(m\) gives
\(ms^2\le C\ell_N\).  Finally, if \(s\le0\) then
\(\widehat s_0\le s+\widehat\Delta\le\widehat\Delta<4\widehat\Delta\).
\end{proof}

\subsection{The radial affine frontier (Proposition~\ref{prop:radial_frontier})}
\label{supp:radial_frontier}

We evaluate the affine frontier explicitly on the
radial class, including the two parallel endpoints and the
nonparallel case.  Throughout this subsection \(P\in\mathcal P_1\) has parameters
\((a,v,s,r,\varphi)\), \(\alpha\in(0,1]\), \(u=1-\alpha\),
\(w_\alpha=\kappa_\alpha v+ur\) with \(b_\alpha=\|w_\alpha\|>0\),
\(\rho_\alpha=\kappa_\alpha/b_\alpha\in[-1,1]\),
\(\bar w_\alpha=w_\alpha/b_\alpha\), and
\(J_T=\mathbb E_{\mu_0}|m_T(X)|=a\,\mathbb E_{\mu_0}\{\varphi(X)|v^\top X|\}\).
On the common design, the frontier may be computed in the
\(\mu_0\)-normalization,
\(\mathcal B_\alpha(c)=\int\{|m_T|-c|m_\alpha|+|m_T-cm_\alpha|\}\,d\mu_0\).
The one geometric input is the radial identity
\begin{equation}
  \mathbb E_{\mu_0}\{\varphi(X)\,|z^\top X|\}
  =\|z\|_2\,\mathbb E_{\mu_0}\{\varphi(X)\,|v^\top X|\}
  =\|z\|_2\,J_T/a,
  \qquad z\in\mathbb R^d,
  \label{eq:supp_radial_identity}
\end{equation}
which holds because \(\varphi(x)\,d\mu_0(x)\) is invariant under
orthogonal maps.

\begin{proof}[Proof of Proposition~\ref{prop:radial_frontier}]
Since \(m_T=a\varphi\,v^\top x\) and \(m_\alpha=a\varphi\,w_\alpha^\top x\)
with \(a\varphi>0\),
\[
\begin{aligned}
  \mathcal B_\alpha(c)
  &=a\int\varphi(x)\bigl\{|v^\top x|-c|w_\alpha^\top x|
  +|(v-cw_\alpha)^\top x|\bigr\}\,d\mu_0(x)\\
  &=J_T\bigl[1-c\,b_\alpha+\|v-cw_\alpha\|_2\bigr]
\end{aligned}
\]
by \eqref{eq:supp_radial_identity} applied to \(z=v\),
\(z=w_\alpha\), and \(z=v-cw_\alpha\), and
\(\|v-cw_\alpha\|^2=1-2c\,v^\top w_\alpha+c^2b_\alpha^2
=1-2c\kappa_\alpha+c^2b_\alpha^2\) because \(v^\top w_\alpha=\kappa_\alpha\)
(\(r\perp v\), \(\|v\|=1\)).  This is the first display of
\eqref{eq:radial_frontier}.  By the proof of
Theorem~\ref{thm:affine_frontier}, \(\mathcal B_\alpha(c)\downarrow
H_{\rm flip}(u)\) as \(c\to\infty\); expanding
\(\sqrt{c^2b_\alpha^2-2c\kappa_\alpha+1}
=cb_\alpha-\kappa_\alpha/b_\alpha+O(c^{-1})\) gives
\(H_{\rm flip}(u)=J_T(1-\kappa_\alpha/b_\alpha)=J_T(1-\rho_\alpha)\).
For the mixture Bayes rule \(h^*_\alpha=\operatorname{sign}(w_\alpha^\top x)\),
the excess-risk identity gives
\[
\begin{aligned}
\mathcal E^{01}_T(h^*_\alpha)
&=\tfrac12\int\{|m_T|-m_Th^*_\alpha\}\,d\mu_0\\
&=\tfrac12\{J_T-a\,\mathbb E[\varphi(X)(v^\top X)
\operatorname{sign}(\bar w_\alpha^\top X)]\}.
\end{aligned}
\]
Write \(v=\rho_\alpha\bar w_\alpha+\sqrt{1-\rho_\alpha^2}\,w_\perp\)
with \(w_\perp\) a unit vector orthogonal to \(\bar w_\alpha\) (absent
if \(|\rho_\alpha|=1\)).  Then
\[
\begin{aligned}
\mathbb E[\varphi(X)(\bar w_\alpha^\top X)\operatorname{sign}(\bar w_\alpha^\top X)]
&=\mathbb E[\varphi|\bar w_\alpha^\top X|]\\
&=J_T/a
\end{aligned}
\]
by \eqref{eq:supp_radial_identity}, while
\[
\mathbb E[\varphi(X)(w_\perp^\top X)\operatorname{sign}(\bar w_\alpha^\top X)]=0,
\]
because the reflection \(x\mapsto x-2(w_\perp^\top x)w_\perp\)
preserves \(\varphi\,d\mu_0\) and \(\bar w_\alpha^\top x\) and reverses
\(w_\perp^\top x\).  Hence
\(\mathcal E^{01}_T(h^*_\alpha)=\tfrac12J_T(1-\rho_\alpha)\), which is
also \(\tfrac12H_{\rm flip}(u)\), as it must be: the mixture Bayes
rule pays exactly the target-margin mass of the flip set.  Finally,
at \(c=2/b_\alpha\), \(cb_\alpha=2\) and \(c\kappa_\alpha=2\rho_\alpha\),
so \(\mathcal B_\alpha(2/b_\alpha)=J_T[\sqrt{5-4\rho_\alpha}-1]\), and
for \(\rho_\alpha\in[0,1]\),
\[
\begin{aligned}
\sqrt{5-4\rho_\alpha}-1
&=4(1-\rho_\alpha)/(\sqrt{5-4\rho_\alpha}+1)\\
&\le2(1-\rho_\alpha)\le2(1-\rho_\alpha^2),
\end{aligned}
\]
with
\(1-\rho_\alpha^2=u^2\|r\|^2/b_\alpha^2\) by
\eqref{eq:mixture_vote_vector}.
\end{proof}

The closed form separates three geometric cases.  When
\(|\rho_\alpha|<1\) (nonparallel normals),
\(\mathcal B_\alpha(c)-H_{\rm flip}(u)
=J_T[\rho_\alpha-cb_\alpha+\sqrt{(cb_\alpha-\rho_\alpha)^2+1-\rho_\alpha^2}]>0\)
for every finite \(c\), so the minimal intercept is not attained
and \(M_\alpha=\infty\).  When \(\rho_\alpha=1\) (parallel normals,
\(u\|r\|=0\) and \(\kappa_\alpha>0\)), \(b_\alpha=\kappa_\alpha\)
and the square root equals \(|1-c\kappa_\alpha|\), so
\(\mathcal B_\alpha(c)=2J_T(1-c\kappa_\alpha)_+\), which vanishes
from \(c=1/\kappa_\alpha\) on: \(M_\alpha=1/\kappa_\alpha\) and
\(H_{\rm flip}(u)=0\).  When \(\rho_\alpha=-1\) (antiparallel
normals, \(u\|r\|=0\) and \(\kappa_\alpha<0\)), the square root
equals \(1+cb_\alpha\), so \(\mathcal B_\alpha(c)\equiv2J_T
=H_{\rm flip}(u)\) for every \(c\ge0\): the minimal intercept is
attained already at \(c=0\), and the modulus is saturated from
\(t=0\) on.  The closed form also gives the constant in
Lemma~\ref{lem:supp_bridge}(d): \(\mathcal B_\alpha(2/b_\alpha)
\le2J_T\sin^2\theta_\alpha\le2a\varphi_+\sin^2\theta_\alpha\), since
\(J_T\le a\varphi_+\mathbb E_{\mu_0}|X_1|\le a\varphi_+\).

\subsection{The target-checked selector (Corollary~\ref{cor:fixed_weight_signed}, Theorem~\ref{thm:target_safe_adaptive})}
\label{supp:adaptive_positive}

We first bound the risk of a fitted candidate at a weight chosen
independently of its training sample.  We then show that the
target check protects against an arbitrary source and retains the
useful candidate with sufficient probability under the structured
model.  The lower bound and the treatment of exceptional events
complete the minimax statement.

\paragraph*{The target gate}
On the pilot block compute
\(\widehat z=n_P^{-1}\sum_{i\in T_P}Y_iX_i\) and
\(b_N=C_b\sqrt{\ell_N/n_P}\).  Reject if
\(\|\widehat z\|\le4b_N\); otherwise set
\(\widehat v=\widehat z/\|\widehat z\|\) and
\(\widehat u=8b_N/\|\widehat z\|\), and reject if
\(\widehat u>1/4\).  With a lexicographically fixed
\(N^{-3}\)-net \(\mathcal N_N\) of \(B^d\), test the set
\(\mathcal Z_N=\{z\in\mathcal N_N:
|\widehat v^\top z|\ge3\widehat u\}\) (reject if empty), and
accept---this is the event \(\Omega_{\rm gate}\)---exactly when every \(z\in\mathcal Z_N\) satisfies
\(\operatorname{sign}(\widehat v^\top z)\,\widehat f_I(z)
\ge2L_{\max}N^{-3}\).  All maps are measurable, minimizers are
unique by strong convexity, ties and empty sets have deterministic
rules, and there is no refit after selection.  Write
\(\mathcal F_{\rm cal}=\sigma(T_C,S_C)\), \(\mathcal F_P=\sigma(T_P)\),
and \(\mathcal F=\mathcal F_{\rm cal}\vee\mathcal F_P\): the
weight \(\widehat\alpha_I\) is \(\mathcal F_{\rm cal}\)-measurable,
the fit blocks \((T_I,S_I)\) are independent of \(\mathcal F\), the
gate decision is \(\sigma(T_C,S_C,T_P,T_I,S_I)\)-measurable, and
the fallback block \(T_0\) is independent of that whole
sigma-field.  The gate compares the fitted candidate with an
independent estimate of the target orientation outside a slab of
half-width \(3\widehat u\asymp b_N\).

The proof of Theorem~\ref{thm:target_safe_adaptive} rests on
four lemmas.  Throughout, \(\mathcal G_{\rm pre}\) is a sigma-field independent
of the fit blocks \((T_I,S_I)\) of sizes \((n_f,m_f)\), \(\alpha\)
is \(\mathcal G_{\rm pre}\)-measurable, \(u=1-\alpha\),
\(\kappa_\alpha=\alpha+us\), \(w_\alpha=\kappa_\alpha v+ur\),
\(b_\alpha=\|w_\alpha\|=\sqrt{\kappa_\alpha^2+u^2\tau^2}\) as in
\eqref{eq:mixture_vote_vector}, and
\(\sigma_{\alpha,f}^2:=\alpha^2/n_f+u^2/m_f\) is the two-block
fluctuation scale.  Since both designs equal \(\mu_0\), the hinge
excess risks of Section~\ref{sec:problem_setup} read
\(\mathcal E_D(g)=\int|m_D|\,|g-\operatorname{sign}m_D|\,d\mu_0
=\int\{|m_D|-m_Dg\}\,d\mu_0\) for \(D\in\{T,\alpha\}\), with
\(m_\alpha=\alpha m_T+um_S\); the dominating measure
\(\mu=2\mu_0\) of Section~\ref{subsec:mixture_geometry} carries
\(p_T=p_S=1/2\), and the factors cancel in every affine integral,
so \(\mathcal B_\alpha\) of \eqref{eq:frontier_def} may be
computed in the \(\mu_0\)-normalization.  For a hard classifier
\(h\), \(\mathcal E_D(h)=2\mathcal E^{01}_D(h)\), so
Theorem~\ref{thm:affine_frontier}(ii) gives, for every \(c\ge0\),
\begin{equation}
  \mathcal E^{01}_T(h)
  \le c\,\mathcal E^{01}_\alpha(h)+\tfrac12\mathcal B_\alpha(c).
  \label{eq:supp_hard_conversion}
\end{equation}

\begin{lemma}[Fixed-weight SVM risk bound]
\label{lem:supp_bridge}
Let \(P\in\mathcal P_1\) have amplitude \(a\), aligned strength
\(s\), and tilt \(\tau\), and let \(\widehat f_\alpha\) be the
\(\alpha\)-weighted no-intercept Gaussian hinge-SVM on the fit
blocks.  Then:
\begin{itemize}
\item[(a)] \emph{(Population center.)}  The population objective
\(f\mapsto\lambda_0\|f\|_H^2+\mathbb E_{P_\alpha}(1-Yf(X))_+\) has
the unique minimizer
\(f^P_\alpha=\frac{a}{2\lambda_0}G_{\varphi,w_\alpha}\), with
\(\|f^P_\alpha\|_\infty\le1-\varepsilon_H\) and
\(\operatorname{sign}f^P_\alpha=\operatorname{sign}(w_\alpha^\top x)
=\operatorname{sign}m_\alpha\) off the mixture hyperplane.
\item[(b)] \emph{(Two-block mean square.)}
\(\mathbb E[\|\widehat f_\alpha-f^P_\alpha\|_H^2\mid\mathcal G_{\rm pre}]
\le\sigma_{\alpha,f}^2/(4\lambda_0^2)\).
\item[(c)] \emph{(Mixture risk.)}  If \(b_\alpha>0\),
\begin{equation}
  \mathbb E\bigl[\mathcal E^{01}_\alpha
  \{\operatorname{sign}\widehat f_\alpha\}\mid
  \mathcal G_{\rm pre}\bigr]
  \le\frac{C}{a\,b_\alpha}\,\sigma_{\alpha,f}^2 .
  \label{eq:supp_bridge_mixture}
\end{equation}
\item[(d)] \emph{(Frontier at the geometry-adapted slope.)}  If
\(\kappa_\alpha>0\), with \(\sin\theta_\alpha=u\tau/b_\alpha\)
the tilt angle between the mixture and target hyperplanes,
\begin{equation}
  \mathcal B_\alpha(2/b_\alpha)
  \le C(d,\varphi_+)\,a\sin^2\theta_\alpha .
  \label{eq:supp_bridge}
\end{equation}
When \(u\tau>0\) the target and mixture normals are nonparallel,
the strict-flip wedge has positive design measure, both tie
hyperplanes are \(\mu_0\)-null, and \(M_\alpha=\infty\), while
\(\mathcal B_\alpha(c)\le2a\varphi_+\) for every finite \(c\).
\end{itemize}
\end{lemma}

\begin{proof}
\emph{(a)}  On the class, \(m_\alpha=a\varphi\,w_\alpha^\top x\)
by \eqref{eq:mixture_vote_vector}, and the mixed signed kernel
mean is \(\mathbb E_{P_\alpha}[Yk_{\gamma_0}(X,\cdot)]=aG_{\varphi,w_\alpha}\).
By Lemma~\ref{lem:supp_radial}(a),
\[
\begin{aligned}
\|aG_{\varphi,w_\alpha}\|_\infty
&\le a\varphi_+\|w_\alpha\|\\
&\le a_+\varphi_+\max\{1,s_{\max}+\tau_U\},\\
\|w_\alpha\|&\le\alpha+u|s|+u\tau
\le\max\{1,s_{\max}+\tau_U\};
\end{aligned}
\]
hence the second envelope condition in \eqref{eq:sgn_envelope}
gives \(\|f^P_\alpha\|_\infty\le1-\varepsilon_H<1\) for
\(f^P_\alpha:=\frac a{2\lambda_0}G_{\varphi,w_\alpha}\).  Every
population hinge is then in its linear region on a neighborhood
of \(f^P_\alpha\), the Fr\'echet derivative of the population
objective is \(2\lambda_0f-aG_{\varphi,w_\alpha}\), which vanishes
at \(f^P_\alpha\), and strict convexity makes the minimizer
unique.  Its sign is that of \(\psi_\varphi(\|x\|)w_\alpha^\top x\),
with \(\psi_\varphi>0\).

\emph{(b)}  First freeze \(\alpha\).  Define the empirical signed
kernel mean
\(\widehat Z_\alpha=\frac\alpha{n_f}\sum_{T_I}Y_ik_{\gamma_0}(X_i,\cdot)
+\frac u{m_f}\sum_{S_I}Y_jk_{\gamma_0}(X_j,\cdot)\), so that
\(Z_\alpha:=\mathbb E\widehat Z_\alpha=2\lambda_0f^P_\alpha\).  Since
every population hinge is active at \(f^P_\alpha\),
\(-\widehat Z_\alpha\) is a subgradient of the empirical hinge
part at \(f^P_\alpha\), while Fermat's rule at the empirical
minimizer supplies a subgradient \(-\widetilde Z_\alpha\) of the
empirical hinge part with
\(2\lambda_0\widehat f_\alpha=\widetilde Z_\alpha\).
Monotonicity of the subdifferential of the convex empirical hinge
part gives
\(\langle\widetilde Z_\alpha-\widehat Z_\alpha,
\widehat f_\alpha-f^P_\alpha\rangle\le0\), that is,
\(2\lambda_0\|\widehat f_\alpha-f^P_\alpha\|_H^2
\le\langle\widehat Z_\alpha-Z_\alpha,
\widehat f_\alpha-f^P_\alpha\rangle_H\), whence
\(\|\widehat f_\alpha-f^P_\alpha\|_H
\le\|\widehat Z_\alpha-Z_\alpha\|_H/(2\lambda_0)\).  The two
blocks are independent, the summands are centered after
subtracting their means, cross terms vanish, and
\(\|Yk_{\gamma_0}(X,\cdot)\|_H^2=k_{\gamma_0}(X,X)=1\), so
\(\mathbb E\|\widehat Z_\alpha-Z_\alpha\|_H^2
\le\alpha^2/n_f+u^2/m_f=\sigma_{\alpha,f}^2\).  For a
\(\mathcal G_{\rm pre}\)-measurable \(\alpha\), note that for a
fixed fit sample the map \(\beta\mapsto f_\beta\) from weights to
empirical minimizers is Lipschitz: adding the two strong-convexity
optimality inequalities for \(f_\beta,f_{\beta'}\) and using that
each empirical hinge risk is one-Lipschitz on \(H\) gives
\(\|f_\beta-f_{\beta'}\|_H\le|\beta-\beta'|/\lambda_0\), and the
same holds for the population minimizers.  Approximate \(\alpha\)
by the dyadic \(\alpha_k=2^{-k}\lfloor2^k\alpha\rfloor\); conditional
on \(\mathcal G_{\rm pre}\) the fit blocks retain their product
law, so the fixed-weight bound applies to each \(\alpha_k\), and
the Lipschitz continuity, conditional Fatou, and continuity of
\(\sigma_{\alpha,f}^2\) in \(\alpha\) give the claim.

\emph{(c)}  Put \(\bar w_\alpha=w_\alpha/b_\alpha\),
\(T=\bar w_\alpha^\top X\) for \(X\sim\mu_0\), and
\(D=\|f-f^P_\alpha\|_H\) for a fixed \(f\in H\).  Since
\(|f(x)-f^P_\alpha(x)|\le D\) (the kernel is bounded by one), a
sign disagreement between \(f\) and \(f^P_\alpha\) at \(x\)
forces \(|f^P_\alpha(x)|\le D\), and by
Lemma~\ref{lem:supp_radial}(a),
\(|f^P_\alpha(x)|\ge\frac{a\varphi_-\underline\psi_1b_\alpha}
{2\lambda_0}|\bar w_\alpha^\top x|\).  Hence
\(\operatorname{sign}f\ne\operatorname{sign}m_\alpha\) only on the
slab \(\{|T|\le q\}\), \(q:=2\lambda_0D/(a\varphi_-\underline\psi_1
b_\alpha)\).  The projection \(T\) has the density
\(p_d(t)=\frac{\operatorname{Vol}(B^{d-1})}{\operatorname{Vol}(B^d)}
(1-t^2)^{(d-1)/2}\) on \([-1,1]\), which is bounded, so
\(\mathbb E[|T|\mathbf 1\{|T|\le q\}]\le C_d^\circ q^2\) with
\(C_d^\circ=\max\{1,\operatorname{Vol}(B^{d-1})/\operatorname{Vol}(B^d)\}\).
Since \(|m_\alpha|\le a\varphi_+b_\alpha|T|\), the
excess-risk identity gives
\[
  \mathcal E^{01}_\alpha(\operatorname{sign}f)
  \le a\varphi_+b_\alpha\,\mathbb E[|T|\mathbf 1\{|T|\le q\}]
  \le\frac{4C_d^\circ\varphi_+\lambda_0^2}
  {a\varphi_-^2\underline\psi_1^2\,b_\alpha}\,
  \|f-f^P_\alpha\|_H^2 .
\]
Apply this with \(f=\widehat f_\alpha\), condition on
\(\mathcal G_{\rm pre}\), and use (b); the factors
\(4\lambda_0^2\) cancel.

\emph{(d)}  If \(u\tau>0\), then \(r\ne0\) and
\(w_\alpha=\kappa_\alpha v+ur\) is nonparallel to \(v\), so
\(b_\alpha>0\) and \(|\rho_\alpha|<1\).  Since \(a\varphi>0\),
the strict-flip set is the open wedge where
\((v^\top x)(w_\alpha^\top x)<0\), which has positive
\(\mu_0\)-measure in \(B^d\).  Both tie sets are proper
hyperplanes and hence \(\mu_0\)-null.  The nonparallel case of
Proposition~\ref{prop:radial_frontier} gives
\(\mathcal B_\alpha(c)>H_{\rm flip}(u)\) for every finite
\(c\ge0\); the minimal intercept is therefore unattainable at a
finite slope, so \(M_\alpha=\infty\).  In contrast,
\eqref{eq:frontier_closed} gives
\(0\le\mathcal B_\alpha(c)\le2\int|m_T|\,d\mu_0\le2a\varphi_+\)
for every finite \(c\ge0\).

\emph{The angular estimate.}  On the target-oriented branch,
\(\kappa_\alpha>0\) implies \(b_\alpha>0\) and
\(\rho_\alpha=\kappa_\alpha/b_\alpha\in(0,1]\).
Proposition~\ref{prop:radial_frontier}, specifically
\eqref{eq:radial_frontier_slope}, gives
\[
  \mathcal B_\alpha(2/b_\alpha)
  \le2J_T(1-\rho_\alpha^2)
  =2J_T\sin^2\theta_\alpha
  \le2a\varphi_+\sin^2\theta_\alpha,
\]
where \(J_T=\mathbb E_{\mu_0}|m_T(X)|\le a\varphi_+\).
This proves \eqref{eq:supp_bridge}.  (At \(\alpha=1\), or when \(\tau=0\) and
\(\kappa_\alpha>0\), the two normals are parallel,
\(M_\alpha=1/\kappa_\alpha<\infty\), and
\(\mathcal B_\alpha(2/\kappa_\alpha)=0\); when \(\kappa_\alpha<0\)
the collinear case already carries a full-orientation flip with
\(\sin\theta_\alpha=0\), which is why the estimate is confined to
the target-oriented branch.)
\end{proof}

Write \(q=u/\alpha\) for the borrowing odds.
\par
The fixed-weight bound \eqref{eq:fixed_weight_target} reads
\((C/a)R_a(q;s,\tau)\), where
\begin{equation}
  R_a(q;s,\tau)
  =\frac{n_f^{-1}+q^2\bigl(m_f^{-1}+a^2\tau^2\bigr)}
        {(1+qs)^2},
  \label{eq:Ra_q}
\end{equation}
since \(\alpha=(1+q)^{-1}\), \(u=q/(1+q)\), and
\(\kappa_\alpha=\alpha(1+qs)\); this parameterization is used only
for the optimization below and in the proof of
Corollary~\ref{cor:fixed_weight_signed}(ii).

\begin{proof}[Proof of
Corollary~\ref{cor:fixed_weight_signed}]
Apply \eqref{eq:supp_hard_conversion} with \(c=2/b_\alpha\) to
\(h=\operatorname{sign}\widehat f_\alpha\), take conditional
expectations, and use
\eqref{eq:supp_bridge_mixture}--\eqref{eq:supp_bridge} together
with \(\sin^2\theta_\alpha=u^2\tau^2/b_\alpha^2\):
\[
  \mathbb E[\mathcal E^{01}_T\{\operatorname{sign}
  \widehat f_\alpha\}\mid\mathcal G_{\rm pre}]
  \le\frac2{b_\alpha}\cdot\frac{C\sigma_{\alpha,f}^2}{a\,b_\alpha}
  +\frac C2\,a\,\frac{u^2\tau^2}{b_\alpha^2}
  =\frac Ca\,\frac{\sigma_{\alpha,f}^2+a^2u^2\tau^2}{b_\alpha^2}
  \le\frac Ca\,\frac{\sigma_{\alpha,f}^2+a^2u^2\tau^2}
  {\kappa_\alpha^2},
\]
the last step using \(b_\alpha^2=\kappa_\alpha^2+u^2\tau^2\ge
\kappa_\alpha^2>0\) on the target-oriented branch.  With
\(\alpha=(1+q)^{-1}\), \(u=q/(1+q)\), and
\(\kappa_\alpha=\alpha(1+qs)\), the numerator equals
\(\alpha^2\{n_f^{-1}+q^2(m_f^{-1}+a^2\tau^2)\}\), the factor
\(\alpha^2\) cancels, and the condition \(\kappa_\alpha>0\) is
\(1+qs>0\); this is part (i) in the form
\((C/a)R_a(q;s,\tau)\) of \eqref{eq:Ra_q}.

\emph{Optimization of the borrowing odds.}
\phantomsection\label{lem:supp_optimizer}%
For \(R_a(q;s,\tau)\) as in \eqref{eq:Ra_q},
\[
  \inf_{q\ge0:\,1+qs>0} R_a(q;s,\tau)
  =
  \begin{cases}
    \Bigl[\,n_f+\dfrac{m_fs^2}{1+m_fa^2\tau^2}\Bigr]^{-1},
    & s>0,\\[3mm]
    n_f^{-1}, & s\le0,
  \end{cases}
\]
attained at \(q_*=m_fs/\{n_f(1+m_fa^2\tau^2)\}\) and at \(q=0\),
respectively.

Write \(R_a(q)=(A+Bq^2)/(1+qs)^2\) with \(A=n_f^{-1}\) and
\(B=m_f^{-1}+a^2\tau^2\), so that
\(R_a'(q)=2(Bq-sA)/(1+qs)^3\).  For \(s>0\) the derivative changes
sign once, at \(q_*=sA/B=m_fs/\{n_f(1+m_fa^2\tau^2)\}\), and the
identity \(R_a(q)-AB/(B+s^2A)=(Bq-sA)^2/\{(B+s^2A)(1+qs)^2\}\)
certifies global minimality with
\(R_a(q_*)=[n_f+m_fs^2/(1+m_fa^2\tau^2)]^{-1}\).  For \(s=0\),
\(R_a=A+Bq^2\); for \(s<0\), on \(0\le q<-1/s\),
\(Bq-sA>0\); in both cases \(q=0\) is the unique minimizer with
value \(A=n_f^{-1}\).

For part (ii), fix \(\load>0\), \(\tau\ge0\), and
\(P\in\mathcal P_1(\load,\tau)\), so that \(s=\load\) and
\(\tau_P:=\tau(P)\le\tau\).  Put
\(\widetilde B:=m_f^{-1}+\tau^2\) and
\(\widetilde R(q):=(A+\widetilde Bq^2)/(1+q\load)^2\).  Since
\(a\le a_+\) and \(\tau_P\le\tau\),
\(B=m_f^{-1}+a^2\tau_P^2\le\max\{1,a_+^2\}\widetilde B\), hence
\(R_a(q;\load,\tau_P)\le\max\{1,a_+^2\}\,\widetilde R(q)\) for every
\(q\ge0\).  The odds \(q^\star=\load A/\widetilde B
=m_f\load/\{n_f(1+m_f\tau^2)\}\) are deterministic, satisfy
\(\kappa_{\alpha^\star}=\alpha^\star(1+q^\star\load)>0\), and minimize
\(\widetilde R\) by the optimization above with \(B\)
replaced by \(\widetilde B\), so
\(\widetilde R(q^\star)=[n_f+m_f\load^2/(1+m_f\tau^2)]^{-1}\).  Part
(i) at \(q=q^\star\), with \(a\ge a_-\), gives the claim with
\(C\) depending on \(a_-\), \(a_+\), and the constant of part
(i).
\end{proof}

The implemented odds \(q_I=m_I\widetilde s/n_I\) and the
continuous optimizer \(q_*\) differ in general: with
\(m_I/n_I=s^{-2}\) and \(\tau^2=(m_Ia^2)^{-1}\), one has
\(q_I\approx1/s\) but \(q_*=1/(2s)\).  The next lemma shows that
the selector nevertheless attains the risk scale.

\begin{lemma}[Risk control for the target check]
\label{lem:supp_gate_audit}
Let \(P_T\in\mathcal T\) with parameters \((a,v,\varphi)\), let
\(A_0=a/(d+2)\), and define the pilot event
\(\Omega_P=\{\|\widehat z-A_0v\|_2\le b_N\}\).  Then
\(P(\Omega_P^c)\le2N^{-A}\) for \(C_b\) large enough, and on
\(\Omega_P\), under the sample range \eqref{eq:supp_cimb}:
\begin{itemize}
\item[(a)] \emph{(Pilot geometry.)}  \(b_N\le A_0/33\),
\(\|\widehat z\|\in[A_0-b_N,A_0+b_N]\),
\(\|\widehat v-v\|_2\le\widehat u/3\), neither
\(\|\widehat z\|\le4b_N\) nor \(\widehat u>1/4\) occurs,
\(\mathcal Z_N\ne\emptyset\), and
\begin{equation}
  c\,b_N/a\le\widehat u\le C\,b_N/a,
  \qquad
  c\,a^2\widehat u^2\le\ell_N/n\le C\,a^2\widehat u^2 .
  \label{eq:supp_pilot_scale}
\end{equation}
\item[(b)] \emph{(Pathwise safety of acceptance.)}  For an
arbitrary independent source law, if the gate accepts, then
\(\operatorname{sign}\widehat f_I(x)=\operatorname{sign}(v^\top x)\)
for every \(x\in B^d\) with \(|v^\top x|>4\widehat u\), so the accepted
candidate disagrees with the target Bayes sign only inside the
slab \(\{|v^\top x|\le4\widehat u\}\), and
\begin{equation}
  \mathcal E^{01}_T(\operatorname{sign}\widehat f_I)
  \le\int a\varphi(x)|v^\top x|\,
  \mathbf 1\{|v^\top x|\le4\widehat u\}\,d\mu_0(x)
  \le Ca\,\widehat u^2\le\frac{C\ell_N}{a\,n}.
  \label{eq:supp_slab_mass}
\end{equation}
\item[(c)] \emph{(Structured branch.)}  Let \(P\in\mathcal P_1\)
have aligned strength \(s\) and tilt \(\tau\), with \(s>0\) and \(\tau=\tau(P)\).
On \(\Omega_C\cap\Omega_P\cap\Omega_{\rm cert}\),
conditional on \(\mathcal F\), with
\(A_I=\widehat\alpha_I+(1-\widehat\alpha_I)s\),
\(V_I=\widehat\alpha_I^2/n_I+(1-\widehat\alpha_I)^2/m_I\), and
\(Q_I=(V_I+a^2(1-\widehat\alpha_I)^2\tau^2)/A_I^2\),
\begin{equation}
\begin{gathered}
  \mathbb E\bigl[\mathcal E^{01}_T(\operatorname{sign}
  \widehat f_I)\mid\mathcal F\bigr]
  \le\frac Ca\,Q_I,
  \qquad
  P(\text{gate rejects}\mid\mathcal F)
  \le\min\Bigl\{1,\frac{CQ_I}{a^2\widehat u^2}\Bigr\},\\
  \mathbb E\bigl[\mathcal E^{01}_T(\widehat h)\mid
  \mathcal F\bigr]
  \le\frac Ca\min\Bigl\{Q_I,\frac{\ell_N}{n}\Bigr\},\\
  \min\{n_I^{-1},Q_I\}
  \le C\Bigl[n_I+\frac{m_Is^2}{1+m_Ia^2\tau^2}\Bigr]^{-1},
  \qquad
  Q_I\big|_{\tau=0}\le\frac{C}{n_I+m_Is^2},
\end{gathered}
  \label{eq:supp_gate_audit}
\end{equation}
uniformly down to arbitrarily small positive \(s\).
\end{itemize}
\end{lemma}

The third inequality in \eqref{eq:supp_gate_audit} is the
statement that adaptation does not require estimating \(q_*\):
acceptance is safe because an accepted candidate can err only on
the slab, and rejection is harmless because the target-only
fallback already attains the displayed scale whenever the source-assisted
candidate does not.  The slab has design measure of order
\(\widehat u\asymp b_N\), but the excess risk charges the target
margin \(|m_T|\), which is itself of order \(b_N\) on the slab, so
the relevant mass is of order \(b_N^2\asymp\ell_N/n\).

\begin{proof}
\emph{Concentration.}  By Lemma~\ref{lem:supp_radial}(c),
\(\mathbb E(YX)=A_0v\), and \(\|YX\|_2\le1\), so the
dimension-free mean bound \eqref{eq:supp_hilbert_mean}, with
\(t=A\log(eN)\), gives
\(P(\Omega_P^c)\le2N^{-A}\) once \(C_b\) is large enough,
uniformly in the target law and irrespective of the source.

\emph{(a)}  Since \(n_P\ge\kappa_Pn\) and \(\ell_N\le c_{\rm imb}n\),
\eqref{eq:supp_cimb} gives
\(b_N\le C_b\sqrt{c_{\rm imb}/\kappa_P}\le a_-/(33(d+2))\le A_0/33\).
On \(\Omega_P\), \(\|\widehat z\|\in[A_0-b_N,A_0+b_N]\), so
\(\|\widehat z\|>4b_N\), and
\(\|\widehat v-v\|\le\|\widehat z-A_0v\|/\|\widehat z\|
+\bigl|\|\widehat z\|-A_0\bigr|/\|\widehat z\|
\le2b_N/(A_0-b_N)\le\frac{33}{16}\,\frac{b_N}{A_0}\);
also \(\widehat u=8b_N/\|\widehat z\|\in
[8b_N/(A_0+b_N),8b_N/(A_0-b_N)]\), which gives
\(\widehat u\le\frac{8\cdot33}{32}\,\frac{b_N}{A_0}\le1/4\) and
\(\widehat u/3\ge\frac{8\cdot33}{3\cdot34}\,\frac{b_N}{A_0}
\ge\frac{33}{16}\,\frac{b_N}{A_0}\ge\|\widehat v-v\|\).  Both displays in
\eqref{eq:supp_pilot_scale} follow from
\(\widehat u\asymp b_N/a\) and
\(b_N^2=C_b^2\ell_N/n_P\asymp\ell_N/n\).  A net point \(z_v\)
with \(\|z_v-\widehat v\|\le N^{-3}\le1/4\) satisfies
\(|\widehat v^\top z_v|\ge3/4\ge3\widehat u\), so
\(\mathcal Z_N\ne\emptyset\).

\emph{(b)}  Suppose the gate accepts on \(\Omega_P\), and let
\(x\in B^d\) with \(|v^\top x|>4\widehat u\).  Choose
\(z\in\mathcal N_N\) with \(\|z-x\|\le N^{-3}\).  Then
\(|\widehat v^\top z|\ge|v^\top x|-\|\widehat v-v\|-N^{-3}
\ge4\widehat u-\widehat u/3-N^{-3}\ge3\widehat u\), because
\(N^{-3}\le\widehat u/2\) on the sample range after enlarging
\(n_{\min}\) (\(\widehat u\ge c\sqrt{\ell_N/n}/a_+\ge cN^{-1/2}\)),
and the same estimate shows
\(\operatorname{sign}(\widehat v^\top z)=\operatorname{sign}(v^\top x)\).
Thus \(z\in\mathcal Z_N\), and acceptance gives
\(\operatorname{sign}(v^\top x)\widehat f_I(z)\ge2L_{\max}N^{-3}\),
while \eqref{eq:supp_lipschitz} gives
\(|\widehat f_I(x)-\widehat f_I(z)|\le L_{\max}N^{-3}\); hence
\(\operatorname{sign}\widehat f_I(x)=\operatorname{sign}(v^\top x)
=\operatorname{sign}m_T(x)\).  The target sign errors of the
accepted candidate therefore lie in the slab
\(\{|v^\top x|\le4\widehat u\}\), and the excess-risk identity
\eqref{eq:hinge_classification_calibration} charges each such
\(x\) with \(|m_T(x)|=a\varphi(x)|v^\top x|\le a\varphi_+|v^\top x|\).
The projection \(v^\top X\) has a bounded density under
\(\mu_0\), so the integral is at most \(Ca\widehat u^2\), and
\eqref{eq:supp_pilot_scale} converts this to \(C\ell_N/(an)\).
This uses from the source only independence from the target; its
design and regression are arbitrary.

\emph{(c) Gate margin and false rejection.}  On
\(E:=\Omega_C\cap\Omega_P\cap\Omega_{\rm cert}\),
Lemma~\ref{lem:signed_projection_confidence} gives \(s>0\) and
\(c_-\le\widetilde s/s\le c_+\), so
\(A_I=(1+q_Is)/(1+q_I)>0\), and \(Q_I\) is defined exactly on
the branch where it is used.  Let
\(g_I:=\frac{aA_I}{2\lambda_0}G_{\varphi,v}\).  By
Lemma~\ref{lem:supp_bridge}(a), the population center at the
weight \(\widehat\alpha_I\) is
\(f^P_I=\frac a{2\lambda_0}G_{\varphi,w_I}\) with
\(w_I=A_Iv+(1-\widehat\alpha_I)r\), so
\(f^P_I-g_I=\frac{a(1-\widehat\alpha_I)}{2\lambda_0}G_{\varphi,r}\),
whose squared norm is
\(a^2(1-\widehat\alpha_I)^2\kappa_\varphi\tau^2/(4\lambda_0^2)\)
by Lemma~\ref{lem:supp_radial}(a).  With
Lemma~\ref{lem:supp_bridge}(b) at \(\mathcal G_{\rm pre}=\mathcal F\)
(the fit blocks are independent of \(\mathcal F\)),
\begin{equation}
  \mathbb E[\|\widehat f_I-g_I\|_H^2\mid\mathcal F]
  \le C\{V_I+a^2(1-\widehat\alpha_I)^2\tau^2\}
  =CA_I^2Q_I .
  \label{eq:supp_gate_mse}
\end{equation}
For \(z\in\mathcal Z_N\), part (a) gives
\(\|\widehat v-v\|\le\widehat u/3\), hence
\(\operatorname{sign}(v^\top z)=\operatorname{sign}(\widehat v^\top z)\)
and \(|v^\top z|\ge(8/3)\widehat u\); the lower envelope
\(\psi_\varphi\ge\varphi_-\underline\psi_1\) then gives
\(\operatorname{sign}(\widehat v^\top z)\,g_I(z)
\ge\frac{aA_I}{2\lambda_0}\varphi_-\underline\psi_1\cdot\tfrac83\widehat u
\ge4c_{\rm gate}\,aA_I\widehat u\).  Since
\(\widetilde s\ge N^{-2}\) and \(s\ge cN^{-2}\) on \(E\), one has
\(A_I\ge1/2\) if \(q_I\le1\) and \(A_I\ge s/2\) if \(q_I>1\), so
\(A_I\ge cN^{-2}\); with \(a\ge a_-\) and
\(\widehat u\ge c\sqrt{\ell_N/n}/a\) this gives
\(c_{\rm gate}aA_I\widehat u\ge cN^{-2}\sqrt{\ell_N/n}\ge L_{\max}N^{-3}\)
after enlarging \(n_{\min}\).  Consequently, on the event
\(\{\|\widehat f_I-g_I\|_H\le c_{\rm gate}aA_I\widehat u\}\), every
tested node satisfies
\(\operatorname{sign}(\widehat v^\top z)\widehat f_I(z)
\ge4c_{\rm gate}aA_I\widehat u-c_{\rm gate}aA_I\widehat u
\ge2L_{\max}N^{-3}\), the tested set is nonempty, and the other
rejection modes are absent by (a): the gate accepts.  The
conditional Markov inequality applied to \eqref{eq:supp_gate_mse}
gives the second inequality of \eqref{eq:supp_gate_audit}; the
factor \(A_I^2\) in the mean square cancels against the squared
gate margin, so no lower bound on the unknown strength enters.

\emph{Conditional risk.}  The first inequality of
\eqref{eq:supp_gate_audit} is Corollary~\ref{cor:fixed_weight_signed}(i)
at the \(\mathcal F_{\rm cal}\)-measurable weight
\(\widehat\alpha_I\), since \(\kappa_{\widehat\alpha_I}=A_I>0\)
and the numerator of \eqref{eq:fixed_weight_target} is
\(V_I+a^2(1-\widehat\alpha_I)^2\tau^2\).  The fallback
\(\widehat f_0\) is independent of the gate decision and has
expected risk at most \(C/(an_0)\le C/(an)\) by the same
corollary at \(\alpha=1\).  Splitting on acceptance,
\[
  \mathbb E[\mathcal E^{01}_T(\widehat h)\mid\mathcal F]
  \le\frac Ca\Bigl[Q_I+\frac1n\min\Bigl\{1,
  \frac{CQ_I}{a^2\widehat u^2}\Bigr\}\Bigr]
  \le\frac Ca\,Q_I,
\]
because \(a^2\widehat u^2\asymp\ell_N/n\ge1/n\) by
\eqref{eq:supp_pilot_scale}.  Independently, part (b) on
acceptance and the fallback bound on rejection give
\(\mathbb E[\mathcal E^{01}_T(\widehat h)\mid\mathcal F]
\le C\ell_N/(an)\) on \(\Omega_P\); taking the better bound proves
the third inequality.

\emph{Binary comparison.}  Put \(c=\widetilde s/s\in[c_-,c_+]\),
\(x=m_Is^2/n_I\), and \(z=1+m_Ia^2\tau^2\).  Since
\(q_I=c\,m_Is/n_I\), cancellation of \(\widehat\alpha_I^2\) gives
the exact identities
\[
  Q_I=\frac1{n_I}\,\frac{1+c^2xz}{(1+cx)^2},
  \qquad
  \Bigl[n_I+\frac{m_Is^2}{1+m_Ia^2\tau^2}\Bigr]^{-1}
  =\frac1{n_I}\,\frac z{x+z}=:R^I_* .
\]
If \(x\le z\), then \(R^I_*\ge(2n_I)^{-1}\) and
\(\min\{n_I^{-1},Q_I\}\le n_I^{-1}\le2R^I_*\).  If \(x>z\ge1\),
then \(Q_I\le\frac1{n_I}\frac zx\bigl(1+\frac1{c^2xz}\bigr)
\le Cz/(n_Ix)\), while \(R^I_*\ge z/(2n_Ix)\).  This proves the
fourth inequality.  At \(\tau=0\), \(z=1\): if \(x\le1\) then
\(Q_I\le C/n_I\le2C/(n_I+m_Is^2)\), and if \(x>1\) then
\(Q_I\le C/(n_Ix)=C/(m_Is^2)\le2C/(n_I+m_Is^2)\), which is the
last inequality.
\end{proof}

\begin{example}[Two target directions with the same source law]
\label{ex:supp_witness}
Fix \(\load>0\), \(0<\theta<\pi/2\), and two coordinate unit vectors
\(e_1,e_2\), and put \(v_\pm=\cos\theta\,e_1\pm\sin\theta\,e_2\),
\(w_S=(\load/\cos\theta)e_1\), and \(r_\pm=w_S-\load v_\pm\).  Then
\(r_\pm\perp v_\pm\), \(\|r_\pm\|_2=\load\tan\theta\), and
\(\load v_++r_+=\load v_-+r_-=w_S\), so the two pairs with target directions
\(v_\pm\) and the common source margin \(a\varphi(x)w_S^\top x\) lie
in \(\mathcal P_1(\load,\tau)\) once \(\load\tan\theta\le\tau\).  They have
identical source laws and different target directions, so no
amount of source data separates them; this pair, with \(\theta=\delta\), is the one used in the
lower bound below whenever \(\load\tan\delta\le\tau\).  The same limit can be read from the minimax
scale \((n+I_S)^{-1}\) of
Theorem~\ref{thm:target_safe_adaptive}: relative to the
target-only scale \(n^{-1}\), the ratio is
\((1+I_S/n)^{-1}\), and
\(I_S^{-1}=(m\load^2)^{-1}+\tau^2/\load^2\) for \(\load>0\), so an
order-of-magnitude improvement over target-only training requires
\(I_S/n\to\infty\), that is, enough sufficiently strong source
observations together with a tilt that is small relative to
the target sample's precision; enlarging \(m\) alone cannot
compensate for a fixed relation uncertainty \(\tau>0\).
\end{example}

\begin{lemma}[Quality-class lower bound]
\label{lem:supp_lower}
Let \(d\ge2\).  For every \(0\le \load\le s_U\) and
\(0\le\tau\le\tau_U\),
\[
  \inf_{\widetilde h}\sup_{P\in\mathcal P_1(\load,\tau)}
  \mathbb E_P\,\mathcal E^{01}_T(\widetilde h)
  \ge c\Bigl[n+\frac{m\load^2}{1+m\tau^2}\Bigr]^{-1},
\]
the infimum running over all randomized measurable classifiers
using the full domain-labeled sample.  For \(d\ge3\) the same
bound holds on every exact shell
\(\{P:\load(P)=\load,\ \tau(P)=\tau\}\), and for \(d=2\) on the exact
shells with \(\tau=0\).
\end{lemma}

\begin{proof}
\emph{Two-point reduction.}  For two laws
\(P_\pm=(P_{T,\pm},P_{S,\pm})\) in the class, let
\(Q_\pm=P_{T,\pm}^{\otimes n}\otimes P_{S,\pm}^{\otimes m}\) be
the laws of the complete data and
\[
  \Gamma:=\int_{\{\operatorname{sign}m_{T,+}\ne
  \operatorname{sign}m_{T,-}\}}
  \min\{|m_{T,+}(x)|,|m_{T,-}(x)|\}\,d\mu_0(x).
\]
Couple \(Q_+\) and \(Q_-\) maximally and run a randomized
procedure \(\widetilde h\) with the same auxiliary seed under
both.  On the coincidence event, which has probability
\(1-\operatorname{TV}(Q_+,Q_-)\), the two outputs coincide, and
at every \(x\) in the disagreement set this common classifier
disagrees with at least one of the two Bayes signs, paying at
least \(\min\{|m_{T,+}(x)|,|m_{T,-}(x)|\}\) in the corresponding
experiment by \eqref{eq:hinge_classification_calibration}.
Integrating in \(x\) and then over the coupling,
\begin{equation}
  \max_{\sigma\in\{+,-\}}
  \mathbb E_{Q_\sigma}\mathcal E^{01}_{T,\sigma}(\widetilde h)
  \ge\frac\Gamma4\{1-\operatorname{TV}(Q_+,Q_-)\}.
  \label{eq:supp_lecam}
\end{equation}
For the divergence, if \(|u|,|w|\le1-\varepsilon_Y\) then
\(\operatorname{KL}(\operatorname{Ber}\tfrac{1+u}2\,\|\,
\operatorname{Ber}\tfrac{1+w}2)\le C_{\varepsilon_Y}(u-w)^2\)
by Taylor's theorem, since the second derivative in the first
probability is bounded on
\([\varepsilon_Y/2,1-\varepsilon_Y/2]\); the common design, the
chain rule for conditional divergences, and product additivity
turn this into
\(\operatorname{KL}(Q_+,Q_-)\le C\{n\|m_{T,+}-m_{T,-}\|_{L_2(\mu_0)}^2
+m\|m_{S,+}-m_{S,-}\|_{L_2(\mu_0)}^2\}\), and
\(\mathbb E_{\mu_0}(b^\top X)^2=\|b\|^2/(d+2)\).

\emph{The ball with \(\load>0\).}  Take the constant profile
\(\varphi\equiv1\) (admissible since \(\varphi_-\le1\le\varphi_+\)),
\(s=\load\), a fixed amplitude \(a_0\in[a_-,a_+]\), a small angle
\(\delta>0\), and
\begin{gather*}
  v_\pm=\cos\delta\,e_1\pm\sin\delta\,e_2,\qquad
  t_+=-\sin\delta\,e_1+\cos\delta\,e_2,\qquad
  t_-=\sin\delta\,e_1+\cos\delta\,e_2,\\
  \varsigma=\min\{\tau,\load\tan\delta\},\qquad
  r_+=-\varsigma t_+,\qquad r_-=\varsigma t_- .
\end{gather*}
Then \(r_\pm\perp v_\pm\) and \(\|r_\pm\|=\varsigma\le\tau\), so
both laws lie in \(\mathcal P_1(\load,\tau)\); this uses \(d\ge2\).
The target normals differ by \(2\sin\delta\,e_2\) and the source
normals \(\load v_\pm+r_\pm\) by \(2(\load\sin\delta-\varsigma\cos\delta)e_2\),
so
\[
  \operatorname{KL}(Q_+,Q_-)
  \le Ca_0^2\bigl\{n\sin^2\delta
  +m(\load\sin\delta-\varsigma\cos\delta)^2\bigr\}.
\]
Set \(D_{a_0}=n+m\load^2/(1+ma_0^2\tau^2)\) and
\(\tan^2\delta=c_0/(a_0^2D_{a_0})\) with a small fixed
\(c_0\in(0,1]\).  The target term is at most \(Cc_0\).  For the
source term, \(|\load\sin\delta-\varsigma\cos\delta|
=\cos\delta\,(\load\tan\delta-\varsigma)_+=\cos\delta\,(\load\tan\delta-\tau)_+\),
and with \(y=\sqrt m\,a_0\tau\),
\(\sqrt m\,a_0\load\tan\delta=\sqrt{c_0m\load^2/D_{a_0}}
\le\sqrt{c_0(1+y^2)}\le\sqrt{c_0}(1+y)\), so
\[
\begin{aligned}
  \sqrt m\,a_0(\load\tan\delta-\tau)_+
  &=\bigl[\sqrt m\,a_0\load\tan\delta-y\bigr]_+\\
  &\le\bigl[\sqrt{c_0(1+y^2)}-y\bigr]_+\\
  &\le\bigl[\sqrt{c_0}+(\sqrt{c_0}-1)y\bigr]_+
  \le\sqrt{c_0},
\end{aligned}
\]
the last step because \(c_0\le1\) and \(y\ge0\); the source term
is at most \(Cc_0\) as well.
Choose \(c_0\) so small that Pinsker's inequality gives
\(\operatorname{TV}(Q_+,Q_-)\le1/2\).  For the separation, on the
cylinder \(\{x_2\in[1/4,1/2],\ (x_3,\ldots,x_d)\in\text{a fixed
small box},\ |x_1|\le(\tan\delta)x_2\}\subset B^d\), the two target
margins \(m_{T,\pm}(x)=a_0(\cos\delta\,x_1\pm\sin\delta\,x_2)\)
have opposite signs and
\(\min\{|m_{T,+}|,|m_{T,-}|\}=a_0(\sin\delta\,x_2-\cos\delta|x_1|)\),
whose \(x_1\)-integral over \(|x_1|\le(\tan\delta)x_2\) equals
\(a_0x_2^2\sin^2\delta/\cos\delta\).  The remaining coordinates
contribute a fixed positive volume, so for small \(\delta\),
\(\Gamma\ge ca_0\delta^2\ge c/(a_0D_{a_0})\ge c/D_{a_0}\).
Combining with \eqref{eq:supp_lecam} and
\(D_{a_0}\asymp n+m\load^2/(1+m\tau^2)\) (fixed \(a_0\)) proves the
bound for \(\load>0\).

\emph{The coarse fiber \(\load=0\).}  Use \(s=0\), \(r=0\),
\(\varphi\equiv1\) under both alternatives, and rotate \(v_\pm\)
as above with \(\delta^2=c_0/(a_0^2n)\).  Both laws lie in
\(\mathcal P_1(0,\tau)\) for every \(\tau\ge0\), the two source
experiments coincide, so the source divergence vanishes, the
target term is at most \(Cc_0\), and the same wedge integral gives
\(\Gamma\ge c/n\).  Hence the minimax risk over
\(\mathcal P_1(0,\tau)\) is at least \(c/n=c[n+m\cdot0]^{-1}\);
that the fiber also contains negative aligned strengths only enlarges the
supremum.

\emph{Exact shells.}  For \(\load>0\) and \(d\ge3\), keep
\(v_\pm,t_\pm,\varsigma\) and put \(k=(\tau^2-\varsigma^2)^{1/2}\),
\(r_+=-\varsigma t_++ke_3\), \(r_-=\varsigma t_-+ke_3\).  Then
\(r_\pm\perp v_\pm\), \(\|r_\pm\|=\tau\), and the common
\(ke_3\) component cancels from the source difference, which is
again \(2(\load\sin\delta-\varsigma\cos\delta)e_2\); the divergence
and wedge calculations are unchanged on the shell
\(\{\load(P)=\load,\tau(P)=\tau\}\).  For \(\load=0\) and \(d\ge3\), use the
common residual \(r=\tau e_3\), orthogonal to both normals.  When
\(\tau=0\) no third direction is needed and \(d=2\) suffices.  For
\(d=2\) and \(\tau>0\) there is no third orthogonal direction in
which to place a common residual, and no exact-shell bound is
claimed.

\emph{Necessity of \(d\ge2\).}  For \(d=1\), \(v\in\{-1,+1\}\) and
\(r=0\), the target Bayes rule is \(\operatorname{sign}(vx)\), and
\(\widehat v=\operatorname{sign}(n^{-1}\sum_iY_iX_i)\) satisfies
\(v\,\mathbb E(YX)=a\,\mathbb E\{\varphi(X)X^2\}\ge a_-\varphi_-/3
>0\) with \(YX\in[-1,1]\), so Hoeffding's inequality makes
\(P(\widehat v\ne v)\) exponentially small in \(n\), and the
target-only risk decays faster than \(n^{-1}\): no
\(n^{-1}\)-order lower bound can hold.
\end{proof}

\begin{proof}[Proof of
Theorem~\ref{thm:target_safe_adaptive}]
\emph{Exceptional events.}
\phantomsection\label{lem:supp_exceptional}%
\(P(\Omega_C^c)\le CN^{-A}\) uniformly over \(\mathcal P_1\), and
\(P(\Omega_P^c)\le2N^{-A}\) uniformly over \(P_T\in\mathcal T\)
and arbitrary independent source laws, with \(A>3\).  Since
\(0\le\mathcal E^{01}_T\le\mathbb E_{\mu_0}|m_T|\le1\), the
contribution of \(\Omega_C^c\cup\Omega_P^c\) to
\(\mathbb E\,\mathcal E^{01}_T(\widehat h)\) is at most
\(CN^{-A}\le CN^{-3}\).  Under the fixed envelopes,
\(n+m\load^2/(1+m\tau^2)\le n+ms_U^2\le CN\), so
\(CN^{-A}\le C[n+m\load^2/(1+m\tau^2)]^{-1}\), and likewise
\(CN^{-A}\le C/(n+m\load^2)\) and \(CN^{-A}\le C/n\): the exceptional
events are dominated by every rate displayed in
Theorem~\ref{thm:target_safe_adaptive}, including the
logarithm-free rates on the clean fibers \(\tau=0\) and on the
coarse fiber \(\load=0\).

The probability bounds follow from Section~\ref{supp:adaptive_calibration}
and Lemma~\ref{lem:supp_gate_audit}. The remaining assertions follow
from the displayed comparisons, using \(A>3\), \(N\ge n\), and
\(\ell_N\ge1\).

\emph{Part (i).}  Fix \(P_T\in\mathcal T\) and an arbitrary
independent source law.  On \(\Omega_P\), if the gate accepts,
Lemma~\ref{lem:supp_gate_audit}(b) bounds the risk of the accepted
candidate by \(C\ell_N/(an)\); if it rejects, the selector returns
the independent target-only SVM \(\widehat f_0\), whose risk is at
most \(C/(an_0)\le C\ell_N/(an)\) by
Corollary~\ref{cor:fixed_weight_signed}(i) at \(\alpha=1\), where
only the target law enters (\(\kappa_\alpha=1\),
\(\theta_\alpha=0\), \(\sigma_{\alpha,f}^2=1/n_0\) at \(\alpha=1\)).
The event \(\Omega_P^c\) is absorbed
by the exceptional-event estimate above, and \(a\ge a_-\).  No
assumption on \(P_S\) has been used.

\emph{Part (ii), the minimax rate.}  The lower half is
Lemma~\ref{lem:supp_lower}.  For the upper half, apply
Corollary~\ref{cor:fixed_weight_signed}(ii) with the full
samples as fit blocks, \(n_f=n\), \(m_f=m\), and the
deterministic weight \(\alpha^\star(\load,\tau)\) when \(\load>0\), and
target-only training when \(\load=0\); both bounds are
\(C[n+m\load^2/(1+m\tau^2)]^{-1}\) with no logarithm.

\emph{Part (ii), the selector.}  Fix \(P\in\mathcal P_1(\load,\tau)\)
with \(\load>0\), so \(s=\load\) and \(\tau(P)\le\tau\).  On
\(E=\Omega_C\cap\Omega_P\cap\Omega_{\rm cert}\),
Lemma~\ref{lem:supp_gate_audit}(c) gives
\(\mathbb E[\mathcal E^{01}_T(\widehat h)\mid\mathcal F]
\le(C/a)\min\{Q_I,\ell_N/n\}\le(C/a)\,\ell_N\min\{n_I^{-1},Q_I\}
\le C\ell_N[n_I+m_Is^2/(1+m_Ia^2\tau(P)^2)]^{-1}\), the middle
step by considering \(Q_I\le n_I^{-1}\) and \(Q_I>n_I^{-1}\)
separately with fixed split fractions.  On
\(\Omega_C\cap\Omega_{\rm cert}^c\) the selector returns
\(\widehat f_0\), with risk \(C/n\), and
Lemma~\ref{lem:signed_projection_confidence} gives
\(ms^2/(1+m\tau(P)^2)\le C\ell_N\le Cn\), so
\(C/n\le C[n+m\load^2/(1+m\tau^2)]^{-1}\) already.  Block sizes are
fixed fractions of \(n\) and \(m\), \(a\in[a_-,a_+]\), and
\(\tau(P)\le\tau\), so
\(n_I+m_Is^2/(1+m_Ia^2\tau(P)^2)\ge c\,[n+m\load^2/(1+m\tau^2)]\).
The exceptional events are absorbed by
the exceptional-event estimate above.  For \(\load=0\), \(s\le0\), the
certificate is false on \(\Omega_C\) by
Lemma~\ref{lem:signed_projection_confidence}, the fallback gives
\(C/n\), and \(\Omega_C^c\) costs \(CN^{-A}\le C/n\).

\emph{The two refinements.}  The coarse-fiber bound \(C/n\) at
\(\load=0\) was just obtained.  On the clean fibers \(\tau=0\), return
to the pre-logarithmic inequality
\(\mathbb E[\mathcal E^{01}_T(\widehat h)\mid\mathcal F]\le(C/a)Q_I\)
on \(E\) and use \(Q_I|_{\tau=0}\le C/(n_I+m_Is^2)\) from
\eqref{eq:supp_gate_audit}; on \(\Omega_C\cap\Omega_{\rm cert}^c\),
Lemma~\ref{lem:signed_projection_confidence} gives
\(ms^2\le C\ell_N\le Cc_{\rm imb}n\), so \(n+ms^2\le Cn\) and the
fallback rate \(C/n\) is of the clean order; the exceptional
events are dominated by \(C/(n+m\load^2)\) by
the exceptional-event estimate above.  The exact-shell refinement of
the lower bound is the second statement of
Lemma~\ref{lem:supp_lower}.
\end{proof}

\subsection{The exact-fitting limit (Theorem~\ref{thm:opposition_separation})}
\label{supp:adaptive_negative}

The argument finds an event on which every admissible output
predicts the same target-opposed sign, and then compares the
restricted family with a rule that may reverse the source labels.
The sign lock stated informally before
Theorem~\ref{thm:opposition_separation} is the following lemma.

\begin{lemma}[Sign lock]
\label{lem:sign_lock}
Fix covariates \(x_1,\ldots,x_M\), a common label
\(y_0\in\{-1,+1\}\), finite weights \(w_j\ge0\) with
\(\sum_{j=1}^Mw_j>0\), and \(0<\gamma,\lambda<\infty\).
\begin{itemize}
\item[(a)] If \(\widehat f\) is the exact minimizer of
\(\lambda\|f\|_{H_\gamma}^2+\sum_{j=1}^Mw_j(1-y_0f(x_j))_+\)
over \(H_\gamma\), then \(\operatorname{sign}\widehat f(x)=y_0\)
for every \(x\).
\item[(b)] If \((\widehat f,\widehat b)\) is any exact minimizer of
\(\lambda\|f\|_{H_\gamma}^2+\sum_{j=1}^Mw_j(1-y_0[f(x_j)+b])_+\)
over \(H_\gamma\times\mathbb R\), then \(\widehat f=0\),
\(y_0\widehat b\ge1\), and
\(\operatorname{sign}(\widehat f+\widehat b)(x)=y_0\) for every
\(x\).
\end{itemize}
Both conclusions hold pathwise even when the subset, the
weights, the bandwidth, and the regularization are arbitrary
measurable functions of the full data and an auxiliary seed.  The
positive-total-mass condition is necessary, and (b) concerns exact
minimizers; approximate minimizers and the single-class behavior
of a particular software implementation are outside its scope.
\end{lemma}

\begin{proof}
(a) Hinge subgradient optimality gives
\(2\lambda\widehat f
=\sum_{j}w_ja_jy_0k_\gamma(x_j,\cdot)\) with
\(a_j\in[0,1]\).  At least one product \(w_ja_j\) is positive:
otherwise \(\widehat f=0\), and every observation with
\(w_j>0\) would have margin zero, forcing \(a_j=1\), a
contradiction.  Since \(k_\gamma(x_j,x)>0\) for all \(j,x\), the
sum has the sign \(y_0\) everywhere.
(b) The pair \((f,b)=(0,y_0)\) makes every hinge term
\((1-y_0\cdot y_0)_+=0\) and the penalty zero, so the minimum
value is zero.  Every term is nonnegative, so an exact minimizer
has \(\lambda\|\widehat f\|^2=0\), whence \(\widehat f=0\), and
\((1-y_0\widehat b)_+=0\) at every observation with \(w_j>0\),
whence \(y_0\widehat b\ge1\); then \(\widehat f+\widehat b\equiv\widehat b\)
has the sign \(y_0\).  The all-zero weight vector gives
\(\widehat f\equiv0\) in (a) and leaves \(\widehat b\)
unconstrained in (b), so the positive-total-mass condition cannot
be dropped.  Measurability of the subset, weights, bandwidth, and
regularization in the data is irrelevant to the pathwise
conclusion.
\end{proof}

\begin{remark}[Exact fitting and optimization tolerance]
\label{cmp:exact-fitting-scope}
Definition~\ref{def:glob_class} requires an exact minimizer, and
this is essential when the total observation weight is
unrestricted.  On an all-common-label sample with an unpenalized
intercept, the optimal objective value is zero, attained by
\(f=0\) and \(b=y_0\); the reversed constant output \(-y_0\) has
objective value \(2W\), where \(W=\sum_i\widehat w_i^T+\sum_j\widehat w_j^S\)
is the total weight, which the family allows to be arbitrarily
small.  A fixed positive optimization tolerance therefore need not
exclude the reversed output, and approximate minimizers are not
covered by Theorem~\ref{thm:opposition_separation}.  The theorem
concerns the training-and-output operation of exact
label-preserving nonnegative-weight fitting, not the expressive
power of the Gaussian RKHS.
\end{remark}

\begin{proof}[Proof of Theorem~\ref{thm:opposition_separation}]
\emph{Step 1: all-label lower bound.}
\phantomsection\label{prop:all_label_lower}%
For the opposition experiment \eqref{eq:opp_class} and every
\(\widehat h\in\mathcal A\),
\[
  \mathbb E_{P_{\sigma;\theta_0,\theta_1,\theta_2}}\,
  \mathcal E^{01}_T(\widehat h)
  \ge
  \theta_0\Bigl(\frac{1-\theta_0}{2}\Bigr)^{n}
  \Bigl(\frac{1+\theta_2}{2}\Bigr)^{m}.
\]

Consider the all-label event
\(\mathcal E_\sigma=\{Y_1^T=\cdots=Y_n^T=-\sigma\}\cap
\{Y_1^S=\cdots=Y_m^S=-\sigma\}\).  Since
\(\mathbb E_{\mu_1}X_1=0\),
\(P_\sigma(Y^T=-\sigma)
=\mathbb E_{\mu_1}\frac{1-(\theta_0+\theta_1X_1)}2
=\frac{1-\theta_0}2\), and the constant source regression gives
\(P_\sigma(Y^S=-\sigma)=\frac{1+\theta_2}2\); independence makes
\(P(\mathcal E_\sigma)\) the product of the two powers.  On
\(\mathcal E_\sigma\), every legal final fit of
\(\mathcal A\) uses observations that all carry the label
\(-\sigma\), with nonnegative weights of positive total mass, so
Lemma~\ref{lem:sign_lock}---part (a) on the no-intercept branch,
part (b) on the unpenalized-intercept branch---forces the returned
classifier to equal \(-\sigma\) everywhere, while the target
Bayes sign is \(\sigma\); the conditional target excess is
\(\int|m_{T,\sigma}|\,d\mu_1
=\frac12\int_{-1}^1(\theta_0+\theta_1z)\,dz=\theta_0\).  Multiplying by
\(P(\mathcal E_\sigma)\) proves the claim.

\emph{Step 2: label-reversal upper bound.}
\phantomsection\label{prop:reversal_upper}%
The domain-aware classifier
\[\widetilde h(x):=-\operatorname{sign}\bigl(\sum_{j=1}^mY_j^S\bigr)\]
satisfies
\[
  \mathbb E_{P_{\sigma;\theta_0,\theta_1,\theta_2}}\,
  \mathcal E^{01}_T(\widetilde h)
  \le\theta_0(1-\theta_2^2)^{m/2}
  \le\theta_0\,e^{-m\theta_2^2/2}.
\]

Let \(Z_j=-\sigma Y_j^S\); the \(Z_j\) are i.i.d.\ signs with
\(\mathbb EZ_j=\theta_2\), and, apart from the deterministic tie
convention, an orientation error is contained in
\(\{\sum_jZ_j\le0\}\).  Chernoff's bound at
\(\omega=\operatorname{arctanh}(\theta_2)\) gives
\(\mathbb Ee^{-\omega Z_1}
=\frac{1+\theta_2}2e^{-\omega}+\frac{1-\theta_2}2e^{\omega}
=\sqrt{1-\theta_2^2}\), so
\(P(\sum_jZ_j\le0)\le(1-\theta_2^2)^{m/2}\).  A wrong constant
orientation costs \(\theta_0\), and
\(\log(1-\theta_2^2)\le-\theta_2^2\) gives the second inequality.

\emph{Step 3: class bounds and separation.}
The all-label bound in Step 1 at the class point
\(\theta_0=\frac12\), \(\theta_2=\frac9{10}\), with any admissible
\(\theta_1\), gives, for every member of
\(\mathcal A\), risk at least
\(\frac12(\frac14)^n(\frac{19}{20})^m\), the first display.
The reversal bound in Step 2 at the worst class values
\(\theta_0=\frac58\), \(\theta_2=\frac9{10}\) gives
\(\mathfrak R
\le\frac58(1-\frac{81}{100})^{m/2}=\frac58(\frac{19}{100})^{m/2}\),
the second, using the first (Chernoff) inequality of Step 2
rather than its exponential relaxation.  Dividing,
\[
  \frac{\mathfrak R_{\mathcal A}}
       {\mathfrak R}
  \ge\frac{\frac12(\frac14)^{n}(\frac{19}{20})^m}
          {\frac58(\frac{19}{100})^{m/2}}
  =\frac45\,4^{-n}
  \Bigl(\frac{361}{400}\cdot\frac{100}{19}\Bigr)^{m/2}
  =\frac45\,4^{-n}\Bigl(\frac{19}{4}\Bigr)^{m/2},
\]
and \(4^{-n}\ge4^{-m/2}\) when \(n\le m/2\) gives the
second bound in \eqref{eq:opp_ratio}.  More generally the ratio
grows exponentially along \(n=\lfloor cm\rfloor\) for every
\(c<\log(19/4)/\log16\approx0.56\); for example, along
\(n=\lfloor m/10\rfloor\) it is at least
\(\frac45\exp\{(\frac12\log\frac{19}{4}-\frac1{10}\log4)m\}
\ge\frac45e^{0.64m}\).
\end{proof}

\begin{remark}[Scope of the impossibility conclusion]
\label{rem:opp_firewall}
None of opposition, label preservation, or nonnegative weights
is by itself a universal obstruction.  On a fixed-negative-strength
subclass of \(\mathcal P_1\), opposition laws still produce
target-aligned labels with positive probability per observation,
and a wrapper with observation-specific nonnegative weights may
retain only those observations; whether such filtering attains
the unrestricted order on that subclass is not proved in this
paper, and no rate is asserted for it here.  What the lower bound
of Theorem~\ref{thm:opposition_separation} uses is narrower: on
the all-label event, every retained observation carries the
target-opposed label, so this filtering channel is unavailable
there, whatever the wrapper does.  The obstruction is common-sign
locking of the exact final fit on the all-label event: kernel
positivity yields the lock without an intercept
(Lemma~\ref{lem:sign_lock}(a)), and the zero-objective argument
yields it with an unpenalized intercept
(Lemma~\ref{lem:sign_lock}(b)).  In particular,
Theorems~\ref{thm:target_safe_adaptive}
and~\ref{thm:opposition_separation} use two separate canonical
Euclidean subclasses, and neither is the negative half of the
other on a common class.
\end{remark}

\section{Numerical Protocol}\label{supp:numerical}

\subsection{Fixed-weight experiment: mixture learning and target performance}
\label{cmp:fixed-weight-experiment}

Holding the weight fixed at \(\alpha=1/2\) separates changes in
source characteristics from changes in the weighting rule.  For
each fitted score \(\widehat f\), Figure~\ref{fig:sim_bridge}
compares the mixture excess \(\widehat{\mathcal E}_\alpha(\widehat f)\),
the target excess \(\widehat{\mathcal E}_T(\widehat f)\), and the
conversion envelope
\(\Psi_\alpha\{\widehat{\mathcal E}_\alpha(\widehat f)\}\), the
largest target excess compatible with that same mixture excess.

\begin{figH}[t]
\centering
\includegraphics[width=\textwidth]{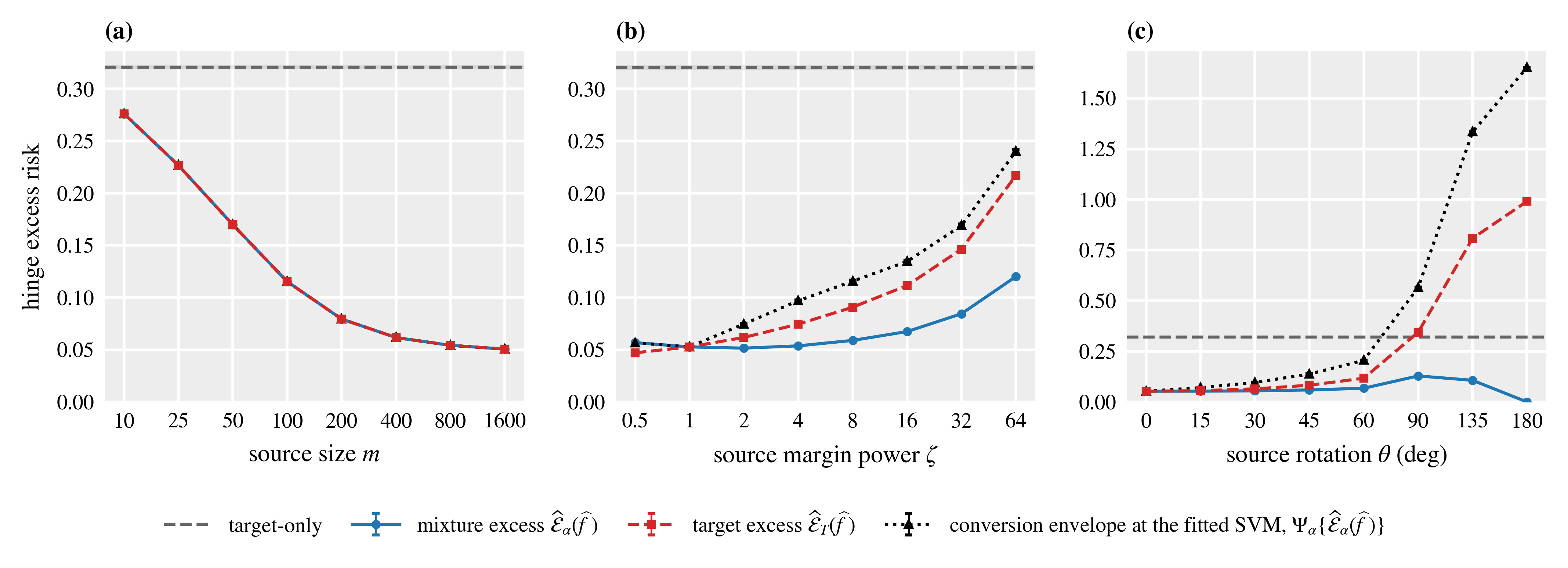}
\caption{Mixture learning and target performance at fixed weight
\(\alpha=1/2\).  For the same fitted SVM, the panels report the
mixture hinge excess \(\widehat{\mathcal E}_\alpha(\widehat f)\)
(solid, circles), the target hinge excess
\(\widehat{\mathcal E}_T(\widehat f)\) (dashed, squares), and the
conversion envelope \(\Psi_\alpha\{\widehat{\mathcal E}_\alpha(\widehat f)\}\)
(dotted, triangles) evaluated at the fitted mixture excess.  The
horizontal dashed line is the target-only reference.  Panels vary
(a) source size under equal laws, (b) source margin power at zero
rotation, and (c) source rotation at unit margin power.  Curves
are means over \(1{,}000\) replications with pointwise Monte Carlo
standard errors; in panel (a) the three curves coincide because
\(P_S=P_T\).  Panel grids and the full data-generating protocol
are given in Section~\ref{cmp:numerical-dgp}.}
\label{fig:sim_bridge}
\end{figH}

Figure~\ref{fig:sim_bridge} distinguishes changes in source
strength from changes in its preferred label.  Equal laws make the
three curves in panel (a) coincide, whereas changing the margin
power (panel (b)) can separate mixture and target excess without
changing the preferred labels.  In panel (c), both target excess
and the conversion envelope cross the target-only reference between
\(60^\circ\) and \(90^\circ\); at \(180^\circ\), equal weighting
gives \(m_\alpha\equiv0\), so every clipped score has zero mixture
excess.

\subsection{Spatial diagnostics for the prediction maps}
\label{cmp:spatial-experiment}

The spatial study of Figure~\ref{fig:sim_votemap} in the main text
uses \(100\) additional training replications at each of the
angles \(0^\circ\), \(90^\circ\), and \(135^\circ\)
(\(n=50\), \(m=1{,}000\)), fitted and evaluated on one reference
sample of \(20{,}000\) covariates shared across the repeated fits.
It compares the same feature locations across
repeated fits.  For method \(M\), let
\(\widehat p_M(x)=B^{-1}\sum_{b=1}^B
\mathbf 1\{\widehat h_{M,b}(x)=1\}\) be the fraction of the
\(B=100\) repeated fits \(\widehat h_{M,1},\ldots,\widehat h_{M,B}\)
that predict label \(+1\) at \(x\).  The target-aligned consensus
\(S_M(x)\) compares these predictions with the target Bayes sign
\(g_T^*(x)=\operatorname{sign}m_T(x)\), and the instability
\(V_M(x)\) records how much the predictions vary across the fitted
samples:
\begin{equation}
  S_M(x)
  =g_T^*(x)\,\{2\widehat p_M(x)-1\},
  \qquad
  V_M(x)
  =4\,\widehat p_M(x)\{1-\widehat p_M(x)\}
  =1-S_M(x)^2.
  \label{eq:sim_consensus}
\end{equation}
In the prediction maps of Figure~\ref{fig:sim_votemap}, blue
corresponds to \(S_M\) near \(1\), red to \(S_M\) near \(-1\), and
near-white to \(S_M\) near \(0\), where \(V_M\) is largest.  The
population row of the same figure carries thin contours of the
cancellation fraction
\(1-|s_\alpha(x)|/\{\alpha|s_T(x)|+u|s_S(x)|\}\), computed from the
known margins on the plotting grid (the figure legend writes it with
the margins \(m_D\), which is the same quantity because the two
domains share the covariate design); it equals \(0\) where the two
domains agree, \(1\) at a tie, and on the disagreement region it is
\(2\min\{\alpha|s_T|,u|s_S|\}/\{\alpha|s_T|+u|s_S|\}\), the
integrand of the cancellation identity
\eqref{eq:cancellation_identity} divided by the total weighted
vote.
Instability discards the direction of the error: consensus values
of \(0.9\) and \(-0.9\) have the same instability \(V_M=0.19\),
but correspond to target-sign error frequencies \(\{1-S_M\}/2\)
of \(0.05\) and \(0.95\).  Multiplication by \(|m_T(x)|\) turns
these error frequencies into their local contributions to target
excess, and the maps connect exactly to the average
classification excess of the repeated fits:
\[
  \frac1B\sum_{b=1}^{B}\mathcal E^{01}_T(\widehat h_{M,b})
  =\mathbb E_X\Bigl[|m_T(X)|\,\frac{1-S_M(X)}{2}\Bigr].
\]
Persistent errors are therefore most costly where the target
strongly favors one label.  The unstable and stable-wrong shares
reported in the figure panels estimate \(P_X(V_M>0.8)\) and
\(P_X(S_M<-0.8)\) on the shared reference sample.

\subsection{Cross-validation rules and implementation details}
\label{supp:numerical_protocol_details}

This section specifies the data-generating model, selection rule,
and implementation choices shared across the comparisons.

\phantomsection\label{cmp:numerical-dgp}%
\emph{Data-generating model and path grids.}  Covariates are
\(X=Z\,\mathbf c(U)+0.1\,\xi\) with \(\mathbf c(t)=(t/\pi)(\cos t,\sin t)^\top\),
\(U\sim\operatorname{Unif}[0,3\pi]\), \(Z\in\{-1,+1\}\) with
equal probability, and \(\xi\) with independent Student-\(t_2\)
coordinates, all independent, as in Section~\ref{sec:numerical}.
Let \(\mathcal C_+\) contain \(2{,}000\) points on
\(\mathbf c([0,3\pi])\), equally spaced in the parameter \(t\) and
including both endpoints, and let \(\mathcal C_-=-\mathcal C_+\).
The signed distance contrast used in \eqref{eq:sim_dgp} is
\[
  d_\theta(x)
  =\min_{z\in\mathcal C_-}\|R_{-\theta}x-z\|_2
  -\min_{z\in\mathcal C_+}\|R_{-\theta}x-z\|_2,
\]
where \(R_\theta\) rotates a point about the origin by \(\theta\);
positive values indicate proximity to the positive arm, and the
zero set is the decision boundary.  The reference sets are finite,
so \(d_\theta\) is the distance contrast to the discretized arms,
not to the continuous curves.  The three paths of
Figures~\ref{fig:sim_bridge} and~\ref{fig:sim_adaptive} use, with
\(n=50\) and \(\beta_T=\beta_S=5\): (a) source size
\(m\in\{10,25,50,100,200,400,800,1600\}\) with \(\theta=0\),
\(\zeta=1\); (b) source margin power
\(\zeta\in\{0.5,1,2,4,8,16,32,64\}\) at \(\theta=0\),
\(m=1{,}000\); (c) source rotation
\(\theta\in\{0,15,30,45,60,90,135,180\}\) degrees at \(\zeta=1\),
\(m=1{,}000\).
The target-only reference in Figure~\ref{fig:sim_bridge} uses one
target-only risk value per replication, and its Monte Carlo
standard error is computed across replications.

The CV-adaptive rule of Section~\ref{subsec:sim_adaptive} compares
the candidate total target weights
\(\mathcal W_{n,m}=\{n/(n+rm):r\in\{0\}\cup\{2^k:k=-6,\ldots,5\}\}
\cup\{1/2\}\) by five-fold cross-validation on the target sample,
repeated five times with fresh fold assignments; within a setting
the same folds are used for every candidate.  Each fold retains
the full source sample and holds the candidate total target weight
\(\alpha\) fixed, so the held-out loss is the target clipped-hinge
loss of an \(\alpha\)-weighted fit on four target folds and all
source observations.  Among source-using candidates, ties in mean
held-out loss are broken in favor of the larger \(\alpha\).  The
best source-using candidate is accepted only when its mean
validation loss is more than \(0.01\) below that of target-only
training; otherwise \(\alpha=1\) is selected.  The reported
classifier is the full-data fit at the selected weight, not an
average of fold-specific classifiers.

Source-size sweeps use nested source prefixes; margin and rotation
sweeps reuse the source covariates and the label uniforms across
grid values, so that differences along a path are not inflated by
sampling variation.  Proportional pooling \(r=1\) gives the source
a total weight share of \(m/(n+m)\), about \(0.952\) for
\(n=50\), \(m=1{,}000\).

For the aggregate experiments in
Figures~\ref{fig:sim_bridge} and~\ref{fig:sim_adaptive}, each
replication draws an independent reference sample of
\(20{,}000\) covariates from \(P_X\), shared across methods and
settings, on which excess risks and the conversion envelope are
evaluated with the known margins; the spatial study of
Section~\ref{cmp:spatial-experiment} instead uses one reference
sample shared across its repeated fits.

The kernel width
\(k(x,x')=\exp\{-2\|x-x'\|_2^2\}\) and regularization
\(\lambda=3\times10^{-4}\) were fixed by an independent pilot
minimizing mean mixture clipped-hinge excess at \(\alpha=1/2\),
\(m=1{,}000\), \(\theta=0\), \(\zeta=1\);
the software cost parameter is \(C=1/(2\lambda)\) with
per-observation weights \(\alpha/n\) and \((1-\alpha)/m\), and the
solver tolerance is \(10^{-4}\) (\texttt{svc\_tol}); with these
per-observation weights the software objective matches
\eqref{eq:transfer_svm_empirical} with an unpenalized intercept.
The plotting grid of Figure~\ref{fig:sim_votemap} is for
visualization only.

\section{A Nonparametric Reinforcement Benchmark}
\label{supp:reinforcement}
We study the deterministically tuned Gaussian SVM on the
reinforcement class \(\mathcal P_\Theta(\nu,q_T,\zeta)\) of
Remark~\ref{rem:reinforcement_benchmark}. The exact kernel-vote
representation yields the upper bound in
Theorem~\ref{thm:gaussian_reinforcement_minimax}; a testing family
in the fixed class yields the matching lower bound under its stated
design and compatibility conditions.
\subsection{The reinforcement class and the benchmark theorem (Theorem~\ref{thm:gaussian_reinforcement_minimax})}
\label{supp:reinforcement_statement}

Under global reinforcement, target and source observations
identify the same decision boundary with different signal
strengths.  The analysis uses a common spatial resolution and
weights the two samples according to the signal available at that
resolution; at scale \(\gamma\), the local volume, the target and
source margin strengths, and the two sample sizes together
determine how well the boundary can be resolved.  The kernel-vote
representation then gives a direct classification analysis of the
Gaussian SVM.  The Gaussian corollaries of
Section~\ref{sec:gaussian} give a general upper guarantee under
arbitrary three-region geometry; this section isolates a separate,
structured benchmark on the fully specified pure-reinforcement
class below, which combines power reinforcement with a fixed
separation between the marginal distributions.
\citet[Theorems~1--2]{cai2021transfer} established the rate form below
for weighted nearest-neighbor classification under posterior drift;
their Theorem~3 gives an adaptive rate up to logarithmic factors.
Here we prove Gaussian-SVM attainment and a matching lower bound on
the specified class under the conditions stated below.

The reinforcement class is defined as follows.  Fix \(d\ge1\), \(0<\nu\le1\), \(q_T\ge0\), and \(0<\zeta<1\).
Fix the class constants
\[
 \Theta=(L_T,C_{q_T},\mu_-,\mu_+,c_{\mathcal X},r_{\mathcal X},c_-,c_+,\delta)
\]
independently of \(n,m\), with
\(L_T,C_{q_T},c_{\mathcal X},r_{\mathcal X},\delta>0\),
\(0<\mu_-\le\mu_+<\infty\), and \(0<c_-\le c_+<\infty\).
Let \(\mathcal P_\Theta(\nu,q_T,\zeta)\) be the class of
pairs \((P_T,P_S)\) satisfying the following conditions, with these
prescribed constants fixed uniformly over the class.  Here \(p_T\) and \(p_S\)
are densities with respect to Lebesgue
measure \(\lambda_d\) on \(\mathbb R^d\), restricted to their common
support.  A support is \emph{regular} for the prescribed constants if
\[
  \lambda_d(\mathcal X\cap B(x,r))
  \ge c_{\mathcal X}\lambda_d(B(x,r))
\]
for every \(x\in\mathcal X\) and every
\(0<r\le r_{\mathcal X}\).  Specifically:
\begin{enumerate}
\item the common support \(\mathcal X\subset[0,1]^d\) is regular in
  this sense and
  \(0<\mu_-\le p_T,p_S\le\mu_+<\infty\);
\item the marginal densities have a fixed separation:
  \(\|p_T-p_S\|_{L_1(\mathcal X)}\ge\delta>0\);
\item the target margin satisfies the prescribed H\"older bound
  \[
    |m_T(x)-m_T(z)|\le L_T\|x-z\|_2^\nu,
    \qquad x,z\in\mathcal X,
  \]
  and
  \[
    P_{X,T}\{0<|m_T|\le t\}\le C_{q_T}t^{q_T},\qquad t>0;
  \]
\item the source reinforces the target globally and has comparable
  power shape:
  \begin{equation}
    m_Tm_S\ge0,\qquad
    c_-|m_T|^\zeta\le |m_S|\le c_+|m_T|^\zeta .
    \label{eq:gaussian_TS}
  \end{equation}
\end{enumerate}
The margin event excludes exact ties, where the classification
excess integrand is zero. The upper power bound in
\eqref{eq:gaussian_TS} controls source contributions from across
the target boundary through the Gaussian tail.  Replacing the upper power inequality by a fixed
\(\zeta\nu\)-H\"older seminorm bound \(L_S>0\) for \(m_S\)
preserves the upper bound. The corresponding lower bound has the
constant qualification given after the lower-bound proof in
Section~\ref{pf:gmm_fixed_class_lower}.

\paragraph*{Sufficient conditions for the matching lower bound}
The upper bound below needs no additional restriction on the fixed constants.
For its matching lower bound, assume
\begin{equation}
\label{eq:t8-design-compatibility}
 0<\delta\le\tfrac12,\quad
 \mu_-\le1-\delta,\quad \mu_+\ge1+\delta,\quad
 c_{\mathcal X}\le2^{-d},\quad r_{\mathcal X}\le\tfrac12.
\end{equation}
These conditions admit the cube and the fixed shifted densities
used in the lower-bound proof. Assume also $\nu q_T\le1$.
If $d=\nu q_T=1$, impose
\begin{equation}
\label{eq:t8-critical-compatibility}
 C_{q_T}L_T^{q_T}\ge1.
\end{equation}
If $d>\nu q_T$, no additional lower bound on $C_{q_T}L_T^{q_T}$ is needed.
The critical condition ensures membership of the one-sided
boundary construction used below.

For a classifier \(g:\mathcal X\to\{-1,1\}\), write
\[
  \mathcal E_T^{01}(g)
  :=
  \int |m_T(x)|
  \mathbf 1\{g(x)\ne\operatorname{sign}m_T(x)\}\,dP_{X,T}(x)
\]
for the target classification excess---the left-hand side of the
calibration identity \eqref{eq:hinge_classification_calibration}---and
define the all-estimator minimax risk
\[
  \mathcal R^{*}_{n,m}(\mathcal P)
  :=
  \inf_{\widehat g}\sup_{(P_T,P_S)\in\mathcal P}
  \mathbb E\,\mathcal E_T^{01}(\widehat g).
\]
The infimum ranges over all measurable randomized classifiers using
the full domain-labeled samples. The benchmark contains no mixing weight: \(\alpha\) enters only as
a tuning parameter of the attaining estimator.

\begin{lemma}[Exact Gaussian vote representation]
\label{lem:exact_gaussian_vote}
Let \(k\) satisfy \(|k(x,z)|\le1\) and \(k(x,x)\le1\), and let
nonnegative sample weights \(w_\ell\) sum to one.  If
\(\lambda_0\ge1/2\), the unique minimizer of
\[
  \lambda_0\|f\|_{H_k}^2+
  \sum_\ell w_\ell(1-Y_\ell f(X_\ell))_+
\]
is
\[
  \widehat f=\frac1{2\lambda_0}
  \sum_\ell w_\ell Y_\ell k(X_\ell,\cdot).
\]
\end{lemma}

\begin{proof}[Proof of Lemma~\ref{lem:exact_gaussian_vote}]
Set
\[
  S_D:=\sum_\ell w_\ell Y_\ell k(X_\ell,\cdot).
\]
The reproducing property gives
\[
  \|S_D\|_{H_k}
  \le\sum_\ell w_\ell k(X_\ell,X_\ell)^{1/2}\le1.
\]
Thus \(f_0:=S_D/(2\lambda_0)\) satisfies
\(\|f_0\|_\infty\le1/(2\lambda_0)\le1\).  For every \(f\in H_k\),
the inequality \((1-a)_+\ge1-a\) yields
\begin{align*}
  \lambda_0\|f\|_{H_k}^2+
  \sum_\ell w_\ell(1-Y_\ell f(X_\ell))_+
  &\ge
  \lambda_0\|f\|_{H_k}^2+1-\langle f,S_D\rangle_{H_k}\\
  &=
  \lambda_0\left\|f-\frac{S_D}{2\lambda_0}\right\|_{H_k}^2
  +1-\frac{\|S_D\|_{H_k}^2}{4\lambda_0}.
\end{align*}
At \(f=f_0\), every hinge term lies on its affine branch, including
the endpoint, so equality holds.  Strict convexity of the
regularizer gives uniqueness.
\end{proof}

For the Gaussian kernel, put
\[
  \chi_\zeta:=\frac{2\nu+d}{2\zeta\nu+d},
  \qquad
  N_\zeta:=n+m^{\chi_\zeta},
  \qquad
  h_{n,m}:=N_\zeta^{-1/(2\nu+d)}.
\]
For known \((\nu,\zeta)\), choose a deterministic bandwidth with
fixed comparison constants
\(0<c_\gamma^-\le c_\gamma^+<\infty\),
\(c_\gamma^-h_{n,m}\le\gamma\le c_\gamma^+h_{n,m}\),
any fixed \(\lambda_0\ge1/2\), and
\begin{equation}
  \alpha_\gamma
  :=
  \frac{n\gamma^\nu}
       {n\gamma^\nu+m\gamma^{\zeta\nu}}.
  \label{eq:minimax_alpha_gamma}
\end{equation}
Lemma~\ref{lem:exact_gaussian_vote} shows that the induced classifier
is the sign of
\[
  \gamma^\nu\sum_{i=1}^nY_i^Tk_\gamma(X_i^T,\cdot)
  +
  \gamma^{\zeta\nu}\sum_{j=1}^mY_j^Sk_\gamma(X_j^S,\cdot).
\]

\begin{theorem}[Gaussian reinforcement minimaxity]
\leavevmode\label{thm:gaussian_reinforcement_minimax}
Fix \((d,\nu,q_T,\zeta,\Theta)\) independently of \(n,m\).
For the tuning above, there is a class-uniform \(C<\infty\),
allowed to depend on the fixed bandwidth comparison constants, such that
for all \(n,m\ge1\),
\begin{equation}
  \sup_{\mathcal P_\Theta(\nu,q_T,\zeta)}
  \mathbb E\,\mathcal E_T^{01}
  \bigl(\operatorname{sign}\widehat f_\gamma\bigr)
  \le
  C\left(
  n+m^{\frac{2\nu+d}{2\zeta\nu+d}}
  \right)^{-\frac{\nu(q_T+1)}{2\nu+d}}.
  \label{eq:gaussian_reinforcement_upper}
\end{equation}
If additionally \(\nu q_T\le1\), the design conditions
\eqref{eq:t8-design-compatibility} hold, and
\eqref{eq:t8-critical-compatibility} holds when \(d=\nu q_T=1\),
then there is \(c>0\), independent of \(n,m\), such that
\begin{align}
  c\left(
  n+m^{\frac{2\nu+d}{2\zeta\nu+d}}
  \right)^{-\frac{\nu(q_T+1)}{2\nu+d}}
  &\le
  \mathcal R^{*}_{n,m}
  \bigl(\mathcal P_\Theta(\nu,q_T,\zeta)\bigr)
  \notag\\
  &\le
  \sup_{\mathcal P_\Theta(\nu,q_T,\zeta)}
  \mathbb E\,\mathcal E_T^{01}
  \bigl(\operatorname{sign}\widehat f_\gamma\bigr),
  \label{eq:gaussian_reinforcement_minimax}
\end{align}
so the same Gaussian transfer SVM family is minimax rate-optimal on
this fixed reinforcement class under the stated sufficient lower-bound
conditions.
\end{theorem}

Under the stated lower-bound conditions, the Gaussian SVM attains
the all-estimator lower bound on the reinforcement class. Theorem~\ref{thm:gaussian_reinforcement_minimax}
uses a direct classification analysis of the exact empirical Gaussian vote.
It shares the separation of target and source variances with
Section~\ref{sec:mixture_to_target}, but is not obtained by specializing
the positive-exponent hinge variance profile of
Section~\ref{pf:two_block_hinge_variance}. In particular, its margin
condition permits positive target-tie mass. On a common tie region,
zero expected hinge excess need not imply zero variance of the
loss difference; classification excess instead gives that region
zero weight. The direct analysis leads to the information equation
\begin{equation}
  n h^{d+2\nu}+m h^{d+2\zeta\nu}\asymp1.
  \label{eq:gaussian_local_information_equation}
\end{equation}
At resolution \(h\), the two variances are separate, the target and
source signals are \(h^\nu\) and \(h^{\zeta\nu}\), and both decide the
same boundary cell.  Equivalently, before optimizing the weight the
local signal-to-noise ratio is
\[
  \mathcal I_\gamma(\alpha)
  \asymp
  \gamma^d
  \frac{(\alpha\gamma^\nu+
  (1-\alpha)\gamma^{\zeta\nu})^2}
  {\alpha^2/n+(1-\alpha)^2/m}.
\]
Cauchy--Schwarz maximizes this expression at
\eqref{eq:minimax_alpha_gamma} and yields
\[
  \max_\alpha\mathcal I_\gamma(\alpha)
  \asymp
  n\gamma^{d+2\nu}+m\gamma^{d+2\zeta\nu}.
\]
This calculation explains the bandwidth and domain-weight choice
within the direct vote analysis. It is a structural comparison with
the two-block oracle, not a derivation from its hinge-variance assumptions.

\begin{remark}[Classification risk and the kernel-vote regime]
\label{rem:reinforcement_scope}
The minimax statement concerns the induced classification risk of
\(g=\operatorname{sign}\widehat f_\gamma\); the hard-thresholded
score has hinge excess exactly \(2\mathcal E_T^{01}(g)\), so it
attains the same rate up to a factor two, while the raw
real-valued score's hinge excess is a different quantity.  At
\(\lambda_0\ge1/2\) the optimizer is the exact weighted kernel
vote: the fixed regularization scales the score magnitude by a
positive factor without changing the classifier's sign, and the
bandwidth and domain weights control the local signal-to-noise
ratio, which is why a fixed \(\lambda_0\) is compatible with the
minimax classification rate. The pointwise vote analysis and
classification-margin integration apply directly on this reinforcement
class, including target ties. Since the class is globally reinforcing,
\(\alpha_\gamma\) may be below \(1/2\) without conflicting with the
no-flip window used when opposition is present.
\end{remark}

\begin{remark}[No acceleration without a relation]
\label{rem:no_relation_no_acceleration}
The source acceleration in
Theorem~\ref{thm:gaussian_reinforcement_minimax} relies on the
relation \eqref{eq:gaussian_TS} coupling \(P_S\) to \(P_T\), and
some such relation is necessary.  Let the target satisfy the
regular-support, density, H\"older, and margin conditions of the
class \(\mathcal P_\Theta(\nu,q_T,\zeta)\), and let \(P_S\) be
otherwise arbitrary. Under \(\nu q_T\le1\), the sufficient
design conditions \eqref{eq:t8-design-compatibility}, and
\eqref{eq:t8-critical-compatibility} when \(d=\nu q_T=1\),
\[
  \inf_{\widehat h}\ \sup_{(P_T,P_S)}
  \mathbb E\,\mathcal E_T^{01}(\widehat h)
  \asymp
  n^{-\frac{\nu(q_T+1)}{2\nu+d}},
  \qquad\text{uniformly in }m.
\]
The upper bound, which needs none of these additional lower-bound
conditions, is attained by ignoring the source; for the lower
bound, one source distribution is held fixed across a target
Assouad family, so the source sample carries no information about
the hard target bits (Supplement,
Section~\ref{pf:unrestricted_source}).  A relation assumption is
therefore necessary for a uniform source-acceleration guarantee:
the region partition of Section~\ref{sec:problem_setup} determines
conversion geometry; it does not determine minimax difficulty.  This does not rule
out acceleration under other structured relations, including known
opposition (Supplement, Section~\ref{supp:adaptive_negative}).
\end{remark}

The complete upper proof, the Assouad lower bound, the primitive
Gaussian-tail compatibility lemma, a counterexample showing why
one-sided dominance alone is insufficient, and the proof of
Remark~\ref{rem:no_relation_no_acceleration} are in
Section~\ref{supp:gaussian_minimax}.
Together, Theorem~\ref{thm:gaussian_reinforcement_minimax} and
Remark~\ref{rem:no_relation_no_acceleration} establish two distinct
benchmarks under the stated fixed-class lower-bound conditions.
On the globally reinforcing class with those conditions, the tuned SVM attains the matching
two-block minimax rate at \(\alpha_\gamma\).  With the target class
fixed under the same sufficient compatibility conditions and the
source unrestricted, the uniform minimax rate remains
target-only, uniformly in \(m\).  The latter statement does not say
that every particular source is useless or describe a fitted rule's
realized risk.

\subsection{Proof of Theorem~\ref{thm:gaussian_reinforcement_minimax}}
\label{supp:gaussian_minimax}

Throughout this subsection all constants may depend on the fixed
parameters \((d,\nu,q_T,\zeta,\Theta)\) and the fixed bandwidth
comparison constants (and on the fixed seminorm bound \(L_S\)
in the source-H\"older alternative), but not on
\(n,m,\gamma\), or the particular distribution pair.  Write
\(s_x:=\operatorname{sign}m_T(x)\) whenever \(m_T(x)\ne0\).  The
letters \(a,b,h,\kappa,r,K\) are local to the proofs of this
subsection and unrelated to their uses in the main text.

\begin{lemma}[Gaussian moments on a regular support]
\label{lem:gaussian_regular_moments}
Under the regular-support and density conditions of
\(\mathcal P_\Theta(\nu,q_T,\zeta)\), for every \(a>0\) there are
positive constants
\(c_0,C_a,C_2\) such that for \(D\in\{T,S\}\), \(x\in\mathcal X\),
and all \(0<\gamma\le\min\{1,r_{\mathcal X}\}\),
\begin{align}
  \int_{\mathcal X}k_\gamma(x,z)p_D(z)\,dz
  &\ge c_0\gamma^d, \label{eq:gmm_mass}\\
  \int_{\mathcal X}k_\gamma(x,z)\|x-z\|_2^a p_D(z)\,dz
  &\le C_a\gamma^{d+a}, \label{eq:gmm_moment}\\
  \int_{\mathcal X}k_\gamma(x,z)^2p_D(z)\,dz
  &\le C_2\gamma^d. \label{eq:gmm_square}
\end{align}
\end{lemma}

\begin{proof}
For \(0<\gamma\le r_{\mathcal X}\), the regular-support definition
and \(k_\gamma(x,z)\ge e^{-1}\) on \(B(x,\gamma)\) give
\[
  \int_{\mathcal X}k_\gamma(x,z)p_D(z)\,dz
  \ge e^{-1}\mu_-c_{\mathcal X}
  \lambda_d(B(0,1))\gamma^d,
\]
which is \eqref{eq:gmm_mass} uniformly at every support point.  For
\eqref{eq:gmm_moment}, use the density
upper bound, enlarge the integration domain to \(\mathbb R^d\), and
make the change of variables \(v=(z-x)/\gamma\).  The resulting
Gaussian moment is finite.  The same argument with
\(e^{-2\|v\|_2^2}\) proves \eqref{eq:gmm_square}.
\end{proof}

\begin{lemma}[Primitive source conditions imply Gaussian
compatibility]
\label{lem:gaussian_compatibility}
Under the regular-support, target H\"older, sign-agreement, and
two-sided power conditions, with \(\mu_-\le p_S\le\mu_+\) on
\(\mathcal X\), there are \(a_G,b_G>0\) and
\(\gamma_0:=\min\{1,r_{\mathcal X}\}\) such
that, for \(x\in\mathcal X\) with \(m_T(x)\ne0\) and
\(0<\gamma\le\gamma_0\),
\begin{equation}
  s_x\int_{\mathcal X}k_\gamma(x,z)m_S(z)p_S(z)\,dz
  \ge
  a_G\gamma^d|m_T(x)|^\zeta
  -b_G\gamma^{d+\zeta\nu}.
  \label{eq:GRC}
\end{equation}
The same conclusion holds if the upper power inequality in
\eqref{eq:gaussian_TS} is replaced by
\[
  |m_S(x)-m_S(z)|\le L_S\|x-z\|_2^{\zeta\nu}.
\]
\end{lemma}

\begin{proof}
Fix \(x\), put \(t=|m_T(x)|\), and \(s=s_x\).
If \(m_T(z)\) has the same sign as \(m_T(x)\), target H\"older
regularity, source lower dominance, and
\((a-b)_+^\zeta\ge a^\zeta-b^\zeta\) give
\[
  s\,m_S(z)
  \ge c_-t^\zeta-C\|x-z\|_2^{\zeta\nu}.
\]
If the signs are opposite, then
\[
  t+|m_T(z)|\le L_T\|x-z\|_2^\nu.
\]
The upper power condition and
\[
  t^\zeta+|m_T(z)|^\zeta
  \le2^{1-\zeta}\{t+|m_T(z)|\}^\zeta
\]
again give, after increasing \(C\),
\[
  s\,m_S(z)=-|m_S(z)|
  \ge c\,t^\zeta-C\|x-z\|_2^{\zeta\nu}.
\]
The same bound holds at target ties.  Integrating it and applying
\eqref{eq:gmm_mass}--\eqref{eq:gmm_moment} proves
\eqref{eq:GRC}.  Under source H\"older regularity, the shorter
pointwise argument is
\[
  s\,m_S(z)
  \ge s\,m_S(x)-|m_S(z)-m_S(x)|
  \ge c_-t^\zeta-L_S\|x-z\|_2^{\zeta\nu},
\]
followed by the same integration.
\end{proof}

\subsubsection{Direct classification upper bound}
\label{pf:gaussian_reinforcement_upper}
\begin{proof}[Proof of the upper bound]
We prove \eqref{eq:gaussian_reinforcement_upper} in three steps,
initially for \(0<\gamma\le\min\{1,r_{\mathcal X}\}\). Write
\(\widehat g_\gamma=\operatorname{sign}\widehat f_\gamma\) and
\(g_T^*=\operatorname{sign}m_T\).

\emph{Step 1: common boundary-resolution scale.}\label{lem:gaussian_fixed_scale}
Let
\[
  A_T(\gamma):=n\gamma^{d+2\nu},\qquad
  A_S(\gamma):=m\gamma^{d+2\zeta\nu},\qquad
  A(\gamma):=A_T(\gamma)+A_S(\gamma).
\]
For the prescribed \(\gamma\asymp h_{n,m}\), we have
\[
  0<c_A\le A(\gamma)\le C_A<\infty.
\]
At \(\gamma=h_{n,m}\),
\[
  A_T=\frac{n}{N_\zeta}\le1,\qquad
  A_S=\frac{m}{N_\zeta^{1/\chi_\zeta}}\le1.
\]
If \(n\ge m^{\chi_\zeta}\), then \(A_T\ge1/2\).  Otherwise
\(N_\zeta<2m^{\chi_\zeta}\), so
\(A_S>2^{-1/\chi_\zeta}\).  Multiplying \(\gamma\) by a fixed
constant changes the two terms only by fixed powers.

\emph{Step 2: pointwise vote error.}\label{lem:gaussian_pointwise_vote}
For \(\gamma\asymp h_{n,m}\) with the fixed comparison constants
and \(0<\gamma\le\min\{1,r_{\mathcal X}\}\), there are constants
\(B\ge1,c,C>0\) such that, whenever
\(t:=|m_T(x)|\ge B\gamma^\nu\) and
\(v:=t/\gamma^\nu\),
\[
  P\{\operatorname{sign}\widehat f_\gamma(x)
  \ne\operatorname{sign}m_T(x)\}
  \le C\exp(-cv^{2\zeta}).
\]
By Lemma~\ref{lem:exact_gaussian_vote} and
\eqref{eq:minimax_alpha_gamma}, multiplying the score by a positive
constant shows that its sign is that of
\[
  V_\gamma(x)
  :=
  \gamma^\nu\sum_{i=1}^nY_i^Tk_\gamma(X_i^T,x)
  +
  \gamma^{\zeta\nu}\sum_{j=1}^mY_j^Sk_\gamma(X_j^S,x).
\]
Put \(Z_x=s_xV_\gamma(x)\).  Target H\"older regularity and
Lemma~\ref{lem:gaussian_regular_moments} give
\[
  s_x\int k_\gamma(x,z)m_T(z)p_T(z)\,dz
  \ge c\gamma^dt-C\gamma^{d+\nu}.
\]
Lemma~\ref{lem:gaussian_compatibility} gives
\[
  s_x\int k_\gamma(x,z)m_S(z)p_S(z)\,dz
  \ge c\gamma^dt^\zeta-C\gamma^{d+\zeta\nu}.
\]
Choose \(B\) large enough.  Since \(t=v\gamma^\nu\), for
\(v\ge B\),
\begin{align}
  \mu_x:=EZ_x
  &\ge
  c\gamma^d\{n\gamma^\nu t+
    m\gamma^{\zeta\nu}t^\zeta\}\notag\\
  &=c\{A_T(\gamma)v+A_S(\gamma)v^\zeta\}
  \ge c'A(\gamma)v^\zeta
  \ge c''v^\zeta .
  \label{eq:gmm_vote_mean}
\end{align}
The independent-summand variance is, by
\eqref{eq:gmm_square},
\[
  \operatorname{Var}(Z_x)
  \le C\{A_T(\gamma)+A_S(\gamma)\}\le C'.
\]
Every centered target summand is bounded by \(2\gamma^\nu\), and
every centered source summand by \(2\gamma^{\zeta\nu}\); hence the
common increment bound is \(2\gamma^{\zeta\nu}\).  Bernstein's
inequality yields
\[
  P(Z_x\le0)
  \le
  \exp\left\{
  -\frac{\mu_x^2}
  {C A(\gamma)+C\gamma^{\zeta\nu}\mu_x}
  \right\}.
\]
The map \(r\mapsto r^2/(a+br)\) is increasing.  We may therefore
insert the lower bound \eqref{eq:gmm_vote_mean}; moreover
\(\gamma^{\zeta\nu}v^\zeta=t^\zeta\le1\).  The bounds from Step 1 then give the asserted exponent.

\emph{Step 3: classification-margin integration.}\label{lem:gaussian_margin_shell}
Fubini gives
\[
  \mathbb E\,\mathcal E_T^{01}
  (\operatorname{sign}\widehat f_\gamma)
  =
  \int |m_T(x)|P\{\widehat g_\gamma(x)\ne g_T^*(x)\}
  \,dP_{X,T}(x).
\]
On \(0<|m_T|\le B\gamma^\nu\), the contribution is at most
\[
  B\gamma^\nu P_{X,T}\{0<|m_T|\le B\gamma^\nu\}
  \le C\gamma^{\nu(q_T+1)}.
\]
For \(k\ge0\), let
\[
  \mathcal S_k
  :=\{2^kB\gamma^\nu<|m_T|
  \le2^{k+1}B\gamma^\nu\}.
\]
The \(k\)-th contribution is at most
\[
  C\gamma^{\nu(q_T+1)}
  2^{(k+1)(q_T+1)}e^{-c2^{2\zeta k}}.
\]
The series is summable for every \(\zeta>0\).  Exact target ties have
zero excess-risk weight.
Consequently, \(\gamma\asymp h_{n,m}\) gives
\[
 \mathbb E\,\mathcal E_T^{01}(\operatorname{sign}\widehat f_\gamma)
 \le C\gamma^{\nu(q_T+1)}
 \le C N_\zeta^{-\nu(q_T+1)/(2\nu+d)}.
\]
For completeness, put
\(\gamma_*:=\min\{1,r_{\mathcal X}\}\) and
\(N_0:=\max\{2,(c_\gamma^+/\gamma_*)^{2\nu+d}\}\).
The preceding argument applies whenever \(N_\zeta\ge N_0\).
For \(N_\zeta<N_0\), classification excess risk is at most one,
so increasing the uniform constant to at least
\(N_0^{\nu(q_T+1)/(2\nu+d)}\) gives the same upper bound.
All these constants depend only on the fixed class and bandwidth
comparison parameters. Alternatively,
\(\gamma=c_\gamma h_{n,m}\) with
\(0<c_\gamma\le\min\{1,r_{\mathcal X}\}\) places every sample
size directly in the range of Steps 1--3.
\end{proof}

\subsubsection{Matching lower bound with fixed class constants}\label{pf:gmm_fixed_class_lower}
Write $q=q_T$ and $C_q=C_{q_T}$ in this proof. We keep the design pair fixed across all scales
and hypothesis signs. On $\mathcal X=[0,1]^d$, let
\[
 p_T=1,\qquad p_S=1+\delta\psi,\qquad
 \psi(x)=\begin{cases}1,&x_1\le1/2,\\-1,&x_1>1/2.\end{cases}
\]
The function $\psi$ has integral zero and $L^1$ norm one. Consequently these
are probability densities in $[\mu_-,\mu_+]$, their $L^1$ distance is exactly
$\delta$, and $p_S\le3/2$. The cube is regular with constants $2^{-d}$ and
$1/2$: at each point, the orthant of a radius-$r\le1/2$ ball pointing toward
the longer side of each coordinate interval lies in the cube. Thus
\eqref{eq:t8-design-compatibility} verifies the prescribed design constants.

\begin{lemma}[Membership of the hard family with prescribed constants]
\label{lem:gmm_hard_membership}
Assume $\nu q\le1$ and \eqref{eq:t8-design-compatibility}.
Let $v_d$ be the volume of the unit ball and set $b=c_-$.

If $r:=d-\nu q>0$, put
\[
 a=L_T/2,\qquad
 \kappa=\frac12\min\left\{16^{-d},\frac{C_q a^q}{d v_d}\right\},
\]
\[
 h_0=\min\left\{\frac1{16},(\kappa/2)^{1/r},
                  (2a)^{-1/\nu},(2ba^\zeta)^{-1/(\zeta\nu)}\right\}.
\]
For every $0<h\le h_0$, let $K=\lfloor\kappa h^{-r}\rfloor$. There are
$K$ radius-$h$ balls inside the cube, with center separation at least $4h$,
such that, for every $\sigma\in\{-1,1\}^K$,
\begin{equation}
\begin{gathered}
 \phi_k(x)=a(h-\|x-x_k\|_2)_+^\nu,\qquad
 m_{T,\sigma}=\sum_{k=1}^K\sigma_k\phi_k,\\
 m_{S,\sigma}=b\sum_{k=1}^K\sigma_k\phi_k^\zeta.
\end{gathered}
\label{eq:gmm_hard_margins}
\end{equation}
define a member of $\mathcal P_\Theta(\nu,q,\zeta)$ with the designs above.
Moreover $K\ge\kappa h^{-r}/2\ge1$.

If $d=\nu q=1$ and \eqref{eq:t8-critical-compatibility} holds, put $a=L_T$,
\[
 h_0=\min\{1/16,(2a)^{-1/\nu},(2ba^\zeta)^{-1/(\zeta\nu)}\}.
\]
For $0<h\le h_0$, use the two hypotheses
\[
 \phi_h(x)=a(h-x)_+^\nu,\qquad
 m_{T,\sigma}=\sigma\phi_h,\qquad
 m_{S,\sigma}=b\sigma\phi_h^\zeta,\qquad \sigma\in\{-1,1\},
\]
on $[0,1]$. They also belong to the same prescribed class.
\end{lemma}

\begin{proof}
Suppose first that $r>0$. Since $h\le(\kappa/2)^{1/r}$,
$K\ge\kappa h^{-r}/2\ge1$. Let $J=\lfloor1/(8h)\rfloor$ and take coordinate
positions $1/4+2h+4hj$, $j=0,\ldots,J-1$. Their Cartesian grid has at least
$16^{-d}h^{-d}$ points, radius-$h$ balls contained in the cube, and the
required separation. Since $K\le\kappa h^{-r}\le16^{-d}h^{-d}$, select
any deterministic $K$ of these centers.

A single bump is $\nu$-H\"older with constant $a$. If two points lie in
different balls with depths $s,t$, their distance is at least $s+t$, and
\[
 |m_{T,\sigma}(x)-m_{T,\sigma}(z)|
 \le a(s^\nu+t^\nu)
 \le2^{1-\nu}a(s+t)^\nu\le L_T\|x-z\|_2^\nu.
\]
The within-ball and outside-ball cases follow from the single-bump bound.
The choice of $h_0$ gives $|m_T|,|m_S|\le1/2$, so the label probabilities
$(1+m_D)/2$ are valid. Disjointness gives sign agreement and
$|m_S|=b|m_T|^\zeta$. Since $b=c_-\le c_+$, this verifies both power
inequalities without assuming that $[c_-,c_+]$ contains one.

For the margin constant, if $0<t<ah^\nu$, put $s=(t/a)^{1/\nu}<h$.
The small-positive-margin shell in one ball has volume
\[
 v_d\{h^d-(h-s)^d\}\le d v_d h^{d-1}s
 \le d v_d a^{-q}h^{d-\nu q}t^q,
\]
where the last inequality follows from
$h^{d-1}s=h^{d-\nu q}s^{\nu q}(s/h)^{1-\nu q}$ and $\nu q\le1$.
For $t\ge ah^\nu$, the full ball volume is at most
$v_da^{-q}h^{d-\nu q}t^q$. Summing over balls and using $p_T=1$ yields
\[
 P_{X,T}\{0<|m_{T,\sigma}|\le t\}
 \le d v_d K a^{-q}h^{d-\nu q}t^q
 \le d v_d\kappa a^{-q}t^q\le C_qt^q.
\]
This also holds at $q=0$. The class constants have not been enlarged.

In the critical case, $q=1/\nu$ and there is one one-sided bump with
H\"older constant $a=L_T$. For $0<t<ah^\nu$, the margin event has length
$(t/a)^{1/\nu}=a^{-q}t^q$; for $t\ge ah^\nu$, its length is
$h\le a^{-q}t^q$. Condition \eqref{eq:t8-critical-compatibility} therefore
gives the required constant. The remaining checks are the same. In
particular the one-sided construction needs $C_q a^q\ge1$, whereas a
two-sided interior tent would need $C_q a^q\ge2$.
\end{proof}

\begin{proof}[Proof of the lower bound]
Use the fixed choices $a,b,h_0$ in Lemma~\ref{lem:gmm_hard_membership}, and put
\[
\begin{gathered}
 D_0=2v_da^2+3v_db^2a^{2\zeta},\qquad h=c_hh_{n,m},\\
 c_h=\min\{1,h_0,(4D_0)^{-1/(d+2\zeta\nu)}\}.
\end{gathered}
\]
As $h_{n,m}\le1$, the hard family is admissible for every $n,m\ge1$.
Let $Q_\sigma=P_{T,\sigma}^{\otimes n}\otimes P_{S,\sigma}^{\otimes m}$.
Use the convention $H^2(P,Q)=\int(\sqrt p-\sqrt q)^2$. For opposite
Bernoulli margins $u,-u$,
\[
 H^2\!\left(\operatorname{Ber}\tfrac{1+u}{2},
             \operatorname{Ber}\tfrac{1-u}{2}\right)
 =2(1-\sqrt{1-u^2})\le2u^2.
\]
For neighbors differing at one bit, the design laws coincide. Integrating
this inequality and using $p_S\le3/2$ gives
\[
 H^2(Q_\sigma,Q_{\sigma^{(k)}})
 \le2v_da^2 n h^{d+2\nu}
       +3v_db^2a^{2\zeta} m h^{d+2\zeta\nu}.
\]
Here product Hellinger squares are bounded by the sum of factor squares:
the affinities $1-H^2/2$ multiply. Also $\operatorname{TV}(P,Q)\le H(P,Q)$
by Cauchy--Schwarz, since
\[
 \frac12\int|\sqrt p-\sqrt q|(\sqrt p+\sqrt q)
 \le\frac12 H(P,Q)\sqrt{4-H^2(P,Q)}\le H(P,Q).
\]
Since
\[
 n h_{n,m}^{d+2\nu}\le1,\qquad
 m h_{n,m}^{d+2\zeta\nu}\le1,
\]
and $c_h\le1$ with $d+2\nu\ge d+2\zeta\nu$,
the last display is at most $D_0c_h^{d+2\zeta\nu}\le1/4$. Thus every
neighboring total variation distance is at most $1/2$. The same bound
holds for the single critical boundary interval.

For an interior ball, write
\[
 I_{d,\nu}=\int_{B(0,1)}(1-\|u\|_2)^\nu\,du>0,\qquad
 \rho_h=\int\phi_kp_T=a I_{d,\nu}h^{d+\nu}.
\]
At each point with positive bump, neighboring hypotheses prefer opposite
labels. The testing inequality $Q(A^c)+Q'(A)=1-[Q(A)-Q'(A)]\ge1-\operatorname{TV}(Q,Q')$,
multiplied by the bump and integrated, implies that the two local expected
excess risks sum to at least $\rho_h/2$. Averaging over all signs and
summing over the $K$ disjoint balls gives, for every classifier,
\[
 \sup_\sigma\mathbb E_\sigma\mathcal E_{T,\sigma}^{01}(\widehat g)
 \ge K\rho_h/4
 \ge\frac{\kappa aI_{d,\nu}}8 h^{\nu(q+1)}.
\]
This remains valid for randomized classifiers by adjoining the same
independent random seed to both experiments.

In the critical case the one-bit loss mass is
\[
 \rho_h=\int_0^h a(h-x)^\nu\,dx
       =\frac{a}{\nu+1}h^{\nu+1}.
\]
The maximal risk is at least $\rho_h/4$, and $\nu+1=\nu(q+1)$.
Taking the infimum over all classifiers and using hard-family inclusion
proves the desired lower bound, with the explicit positive constant
\[
 c=\begin{cases}
 (\kappa aI_{d,\nu}/8)c_h^{\nu(q+1)},&d>\nu q,\\
 \{a/[4(\nu+1)]\}c_h^{\nu+1},&d=\nu q=1.
 \end{cases}
\]
All constants were chosen before $n,m$. There is no finite exceptional
sample-size range left to treat. Together with the upper bound, this
proves Theorem~\ref{thm:gaussian_reinforcement_minimax}.
\end{proof}

\paragraph*{Alternative source-H\"older formulation}
\phantomsection\label{pf:gmm_source_holder}%
Replace only the upper power inequality in \eqref{eq:gaussian_TS}
by the fixed bound
\[
 |m_S(x)-m_S(z)|\le L_S\|x-z\|_2^{\zeta\nu},\qquad L_S>0,
\]
retaining sign agreement and the lower power inequality. The
alternative in Lemma~\ref{lem:gaussian_compatibility} gives the
same uniform upper rate without additional lower-bound conditions.
For the lower rate, retain \eqref{eq:t8-design-compatibility} and
\(\nu q_T\le1\). If \(d>\nu q_T\), repeat the interior construction with
\[
 b=c_-,\qquad a=\min\{L_T/2,(L_S/(2b))^{1/\zeta}\},
\]
and recompute \(\kappa,h_0,D_0,c_h\) using this fixed amplitude.
The target seminorm is at most \(L_T\), the source seminorm is at
most \(2^{1-\zeta\nu}ba^\zeta\le L_S\), and the prescribed margin
bound follows from the same choice of \(\kappa\).
Thus the membership and testing calculations above apply unchanged.
If \(d=\nu q_T=1\), use the boundary bump with
\[
 a=\min\{L_T,(L_S/b)^{1/\zeta}\}
\]
and require instead the sufficient condition
\begin{equation}
 C_{q_T}\bigl[\min\{L_T,(L_S/c_-)^{1/\zeta}\}\bigr]^{q_T}\ge1.
 \label{eq:t8-holder-critical-compatibility}
\end{equation}
Its target and source seminorms are bounded by \(L_T\) and
\(L_S\), respectively, and the one-sided margin coefficient is
\(a^{-q_T}\). Recomputing the fixed scale threshold with this
amplitude completes the same two-point lower bound.

\begin{remark}[Why one-sided reinforcement is insufficient for this
Gaussian proof]
\label{rem:gaussian_one_sided_counterexample}
On \([0,1]\), take \(m_T(x)=x-1/2\) and
\[
  m_S(x)=
  \begin{cases}
  -1,&x<1/2,\\
  (x-1/2)^\zeta,&x\ge1/2.
  \end{cases}
\]
This satisfies sign agreement and the one-sided lower power bound,
but not the upper power condition.  At \(x=1/2+a\gamma\), the
negative and positive parts of the population source vote have
orders
\[
  -\gamma\,\frac{e^{-a^2}}{a}
  \quad\text{and}\quad
  \gamma^{1+\zeta}a^\zeta,
\]
uniformly for \(1\le a=o(1/\gamma)\).  Their balance occurs at
\[
  a^2=\zeta\log(1/\gamma)+O(\log\log(1/\gamma)),
\]
so the smoothed boundary is displaced by
\(\gamma\sqrt{\log(1/\gamma)}\).  A mean-zero density perturbation
away from the boundary makes \(p_T\ne p_S\) without changing this
local calculation.  Thus one-sided dominance alone does not imply
the uniform no-log Gaussian compatibility
\eqref{eq:GRC}.  This counterexample invalidates that argument,
not every possible nonlinear tuning with \(\lambda_0<1/2\).
\end{remark}

\subsubsection{No acceleration without a relation}
\label{pf:unrestricted_source}
We prove Remark~\ref{rem:no_relation_no_acceleration}.
Let the target satisfy the regular-support, density,
\(\nu\)-H\"older, and margin conditions of
\(\mathcal P_\Theta(\nu,q_T,\zeta)\), with \(P_S\) unrestricted, and
assume \(\nu q_T\le1\), \eqref{eq:t8-design-compatibility}, and
\eqref{eq:t8-critical-compatibility} when \(d=\nu q_T=1\).
Only the lower bound below uses these additional compatibility
conditions.

\emph{Upper bound.}  The estimator may discard the source sample:
at \(\alpha=1\) and \(\gamma\asymp n^{-1/(2\nu+d)}\), the
upper-bound argument of Section~\ref{pf:gaussian_reinforcement_upper}
applies verbatim with the source block absent (every source term is
zero and the information scale reduces to \(nh^{d+2\nu}\asymp1\)),
giving
\(\sup\mathbb E\,\mathcal E_T^{01}\lesssim
n^{-\nu(q_T+1)/(2\nu+d)}\) under the target conditions alone.

\emph{Lower bound.} Use the target hard family of
Lemma~\ref{lem:gmm_hard_membership} at
\(h=c_h n^{-1/(2\nu+d)}\), with the fixed \(c_h\) chosen in the
matching lower-bound proof. Since \(n^{-1/(2\nu+d)}\le1\), its
members meet the same prescribed target conditions for every \(n\ge1\).
Pair every sign vector with the single source law
\[
 p_S^0=1+\delta\psi,\qquad m_S^0=0,
\]
where \(\psi\) is the half-cube sign function used above.
This source retains the common support, density bounds, and genuine
covariate shift, even if those design restrictions are retained
after removing the regression relation. Its law does not depend on
the target hard signs. The joint sample laws are
\(P_T^{\sigma,\otimes n}\otimes P_S^{0,\otimes m}\), so the source
factor contributes zero Hellinger distance for every \(m\).
For neighboring target hypotheses the remaining bound is
\[
 H^2(Q_\sigma,Q_{\sigma^{(k)}})
 \le 2v_da^2nh^{d+2\nu}\le\tfrac14
\]
by the same fixed choice of \(c_h\). The testing and risk-mass
calculation of the preceding lower bound therefore gives
\[
 \inf_{\widehat h}\sup\mathbb E\,\mathcal E_T^{01}(\widehat h)
 \ge c n^{-\nu(q_T+1)/(2\nu+d)},
\]
including the critical one-sided boundary construction and
randomized classifiers. The constant is independent of \(m\).
Together with the target-only upper bound, this proves
Remark~\ref{rem:no_relation_no_acceleration} under its stated
fixed-class qualifications.

\end{document}